\documentclass[12pt]{mynewarticle}
\usepackage{array,bbm,amsmath,bxsc,times,mathptm,mathrsfs,t1enc}
\usepackage{hyperref}
\usepackage[0.87]{fsscale}
\usepackage[a4]{format}
\usepackage[latin1]{inputenc}
\author{Volkmar Liebscher\thanks{GSF --- National Research Centre for Environment and 
Health, Institute of Biomathematics and Biometry,
Ingolstädter Landstr.1, D--85758 Neuherberg, Germany,
email:{\tt liebscher\string @gsf.de}}
}
\title{\bf Random Sets and Invariants for (Type II) Continuous Tensor Product Systems of Hilbert Spaces}
\date{\day25\month6\year2003\relax\today}
\numberwithin{equation}{section}
\numberwithin{example}{section}
\numberwithin{lemma}{section}
\numberwithin{definition}{section}
\numberwithin{remark}{section}
\begin{document}
\titlepage
\maketitle
\setcounter{page}{1}
\makeatother
\begin{abstract}
In a series of papers \cite{Tsi03,Tsi00,Tsi00b,Tsi99a,Tsi99b} \textsc{Tsirelson}   constructed from measure types of random sets or (generalised)  random processes a new range of examples for continuous tensor product systems of Hilbert spaces introduced by \textsc{Arveson} \cite{Arv89} for classifying $E_0$-semigroups upto cocycle conjugacy.  This paper starts from establishing the converse. So   we connect  each continuous tensor product systems of Hilbert spaces  with  measure types of distributions of random (closed) sets in $[0,1]$ or $\NRp$. These measure types  are  stationary and  factorise over disjoint intervals. In a special case of this construction,  the corresponding  measure type is an invariant of the product system. This shows, completing in a more systematic way  the \textsc{Tsirelson} examples, that the classification scheme for product systems into types $\mathrm{I}_n$, $\mathrm{II}_n$ and $\mathrm{III}$ is not complete.  Moreover, based on a detailed study of this kind of measure types, we construct for  each stationary factorizing  measure type a continuous tensor product systems of Hilbert spaces such that this measure type arises as the before mentioned invariant. 

 These results are  a further step  in the  classification of  all (separable) continuous tensor product systems of Hilbert spaces  of type $\mathrm{II}$ in completion to the classification of type $\mathrm{I}$ done by \cite{Arv89} and combine well with other invariants like the lattice of product subsystems of a given product system. Although  these invariants  relate to type  $\mathrm{II}$ product systems mainly, they are of general  importance. Namely, the measure types of the above described kind are  connected with representations of the corresponding $L^\infty$-spaces. This leads to direct integral representations of the elements of a given product system which combine well under tensor products.  Using this structure in a constructive way, we can  relate to any (type $\mathrm{III}$) product system a product system of type $\mathrm{II}_0$  preserving isomorphy classes.  Thus,  the classification of type $\mathrm{III}$ product systems   reduces to that of type $\mathrm{II}$ (and even type $\mathrm{II}_0$) ones. 

In this circle of ideas it proves   useful that we reduce the problem of finding a compatible measurable structure for product systems to prove continuity of one periodic unitary group on a single  Hilbert space. As a consequence,  all admissible measurable structures (if there are any) on an algebraic  continuous tensor product systems of Hilbert spaces yield isomorphic product systems. Thus the measurable structure of a continuous tensor product systems of Hilbert spaces is essentially determined by its algebraic one.  
\end{abstract}

\clearpage
\tableofcontents
\clearpage

\section{Introduction}
\label{sec:intro}

In his seminal work \cite{Arv89} \textsc{Arveson} studied $E_0$-semigroups $\sg t\alpha$, which are weakly continuous semigroups of $\ast$-endomorphisms  of some $\B(\H)$, where $\H$ is a separable Hilbert space. One  basic result of \cite{Arv89} was the following: Up to so-called outer (or cocycle) conjugacy, $\sg t\alpha$ is determined by the family $\sg t\E$ of Hilbert spaces
\begin{equation}
\label{eq:def ps from E_0(Arveson)}
  \E_t=\set{u\in\B(\H):\alpha_t(b)u=ub\forall b\in\B(\H)}
\end{equation}
with the inner product $\scpro uv_{\E_t}\unit=u^*v$. Moreover, $\E=\sg t\E$ fulfils $\E_s\otimes\E_t\cong\E_{s+t}$ under $\psi_s\otimes \psi_t\mapsto \psi_s\psi_t$, it is a \emph{continuous tensor product system of Hilbert spaces}, briefly product system. Conversely,  to every  product system there corresponds an $E_0$-semigroup in this way, see  \cite{OP:Arv90} or  Theorem \ref{th:E_0-semigroup} below. \cite{Arv89} classified product systems of type $\mathrm{I}$ (i.e.\ product systems, which are generated by their units, see Definition \ref{def:unit} below) to be isomorphic to one of the exponential product systems $\Gammai(\K)$, $\K$ a separable Hilbert space, where $\Gammai(\K)_t$ is the symmetric Fock space $\Gamma(L^2([0,t],\ell,\K))$ with its  usual tensor product structure.   

Symmetric Fock spaces have natural relations to random processes with independent increments like Brownian motion and  Poisson process, see e.g.\ \cite{Sch93,Arv89,Mey93,Lie98e}. More explicitely, let $M_t=\set{Z\subset[0,t]\times\set{1,\dots,N}:\#Z<\infty}$ be the space of finite point configurations on $[0,t]\times\set{1,\dots,N}$ equipped with the exponential measure $F_t$ \cite{PP:CP72},
\begin{equation}
\label{eq:defexponentialmeasure}
  F_t=\delta_{\emptyset}+\sum_{n\in\NN}\frac 1{n!}\int_{([0,t]\times\set{1,\dots,N})^n}(\ell\otimes\#)^n(\d x_1,\dots,\d x_n)\delta_{\set{ x_1,\dots, x_n}}.
\end{equation}
$F_t$ is (for finite $N$) just a scaled version of the Poisson process on $[0,t]\times\set{1,\dots,N}$ with intensity measure $\ell\otimes\#$. Then for $N\in\set{1,2,\dots}\cup\set\infty$, $L^2(M_t,F_t)$ is isomorphic to $\Gamma(L^2([0,t],\ell,l^2(\set{1,\dots,N})))$ \cite{Gui72,LP89}. For more general product systems respectively $E_0$-semigroups there existed until 1999  only examples of product systems of type $\mathrm{II}$ (with at least one unit) \cite{OP:Pow99} and type $\mathrm{III}$ (without unit)\cite{OP:Pow87,OP:PP90}. For the definition of types we refer to Definition  \ref{def:types of ps} below. 

\textsc{Tsirelson} \cite{Tsi00}  considered random (closed) sets $Z\subset[0,\infty)$ coming from diffusion processes on (subsets) of $\NR^d$.  Setting $Z_t=Z\cap[0,t]$, for all examples the law $\mathscr{L}(Z)$ of $Z$ fulfilled 
\begin{equation}
\label{eq:equZs+t}
  \mathscr{L}(Z_{s+t})\sim\mathscr{L}(Z'_s\cup (Z''_t+s))\dmf{s,t\ge0}
\end{equation}
where $Z',Z''$ are independent copies of $Z$ and $Z+t$ is the set $Z$ shifted by $t$. I.e., the law of the random set $Z\cap[0,s+t]$ is equivalent  (has the same null sets) to the product of the laws of $Z\cap[0,t]$ and the law of $Z\cap[0,s]$ shifted by $t$.  The key point of the construction is the change of  independence properties like  $\mathscr{L}(Z_{s+t})=\mathscr{L}(Z'_t\cup (Z''_s+t))$  to equivalence relations like \ref{eq:equZs+t}, see Proposition \ref{prop:independentmeasure} below. Such equivalence conditions   that one  measure has a positive Radon-Nikodym derivative with respect to  a product measure, are  e.g.\  used in the theory of Markov Random Fields, cf.\ e.g.\ \cite{S:Lau98}. They express that there a no trivial dependencies between the realisations of the field in  different parts, here between $Z\cap[0,t]$ and $Z\cap[t,s+t]$.     

This equivalence relation implies that $\H_t=L^2(\mathscr{L}(Z_t))$  forms a product system.
These  examples  open a new field  in the classification of  product systems. Formerly,  \cite{OP:Pow97} proposed a  classification of product systems similar  to  \textsc{Connes'} classification of  (hyperfinite) $W^*$-factors \cite{OP:Con73} into  types $\mathrm{I}_n$, $\mathrm{II}_n$, and $\mathrm{III}$.   There $n$ is the \emph{numerical index}, an invariant   introduced by \textsc{Arveson} \cite{Arv89}. It is defined as  the number $n\in\NN$ such that the  product subsystem generated by the units is isomorphic to $\Gammai(\NC^n)$, where a \emph{unit} $\sg tu$ is  any (measurable) nonzero section through $\sg t\E$ which factorises: $u_{s+t}=u_s\otimes u_t$.  If that classification had been complete then any  product system of type $\mathrm{II}_n$ would have been  isomorphic to the tensor product of one  product system of type  $\mathrm{II}_0$ and one of type  $\Gammai(\NC^n)$.   By the results of \textsc{Tsirelson} there are even uncountably many nonequivalent examples of  type $\mathrm{II}$ product systems generated by random sets (and random processes), even with  a fixed dimension of the Fock space generated by the units. Thus, the classification of \cite{OP:Pow97} is too coarse to cover all invariants of (type II) product systems and  one should find more invariants of product systems. 

\bigskip
One important  result of the present work is the fact that random sets did not appear by chance in dealing with examples for type $\mathrm{II}$ product systems. Theorem \ref{th:RACS} below shows that  they are intrinsic to them.  The idea is the following. A closed set $Z\subseteq[0,1]$ is characterised by the values $X_{s,t}=\chfc{\set{Z':Z'\cap[s,t]=\emptyset}}(Z)$ for $0\le s<t\le 1$. For random closed sets, $X_{s,t}$ are random variables fulfilling  the relation
\begin{displaymath}
  X_{r,s}X_{s,t}=X_{r,t}.
\end{displaymath}
Under what circumstances  and by which means we can recover from such a family of random variables the set $Z$ will be discussed in section \ref{sec:nonsep}.
\textsc{Tsirelsons} observation was that for every unit $u$, the orthogonal projections $\mathrm{P}^u_{s,t}$ onto $\unit_{\E_s}\otimes \NC u_{t-s}\otimes\unit_{\E_{1-t}}\subset\E_1$  fulfil the similar relation 
\begin{displaymath}
  \mathrm{P}^u_{r,s}\mathrm{P}^u_{s,t}=\mathrm{P}^u_{r,t}.
\end{displaymath}
The result of Theorem \ref{th:RACS} is that we can fix in a sense a distribution of these (quantum stochastic)  random variables and even  a distribution of a random closed set. Namely, for all normal states $\eta$ on $\B(\E_1)$ there is a unique distribution  $\mu_\eta$ of a random closed set  $Z\subseteq[0,1]$ with 
  \begin{displaymath}
    \int \chfc{\set{Z:Z\cap[s_1,t_1]=\emptyset}}\cdots\chfc{\set{Z:Z\cap[s_k,t_k]=\emptyset}}\mathrm{d}\mu_\eta=\eta(\mathrm{P}^u_{s_1,t_1}\cdots \mathrm{P}^u_{s_k,t_k})\dmf{(s_i,t_i)\subseteq[0,1]}.
  \end{displaymath} 
Further,  among  these measures $\mu_\eta$ there is at least one  dominating  all others with respect to   absolute continuity. Thus we fixed the maximal so-called \emph{measure type} of $\mu_\eta$.

This technique works with other  projections too, let $\mathrm{P}^\Us_{s,t}$ project onto $\unit_{\E_s}\otimes \E^\Us_{t-s}\otimes\unit_{\E_{1-t}}\subset\E_1$, where $\E^\Us_t$ is generated by all units of $\E$. Again, they fulfil the same relations and   we can identify  the related  (maximal)   measure type of random sets  as an invariant of the product system (see Theorem \ref{th:ps2randomset}).  Such techniques of associating with a quantum state (here $\eta$) its restriction to an abelian von Neumann algebra proved already useful in the context of locally normal states of a Boson field in vacuum representation, which can be related to point processes, see e.g.\ \cite{FF87}. Here, its proves valuable to extend the structural investigations in product systems from the analysis of the product subsystem generated by the units (which was done by \textsc{Arveson}) to the study in which way this subsystem is embedded into the whole product system, e.g.\ the structure encoded in  the algebraic properties of the projections $\mathrm{P}^\Us_{s,t}$, e.g.  

These results   forces us to go deeper into structure theory of measure types on the compact subsets of compact intervals, which are both stationary and factorizing. Extending  \textsc{Tsirelsons} examples, we are able to show that all such measure types are connected with a product system in Proposition  \ref{prop:statfactmeastype=ps}.  An interesting new structural  aspect is connected with  a construction of a new random closed set by filling holes left by the original random set (cf.\ Corollary \ref{cor:Zalphabeta}). This serves us to show that  all  stationary factorizing measure types of random closed sets in $[0,1]$ appear as that invariant of some product system derived in Theorem \ref{th:ps2randomset}. Further,  this hole-filling procedure can be iterated ad infinitum giving more and more complex measure types and  a tree-like structure of  the set of all stationary factorizing measure types. This tree structure is supported by the fact that such  measure types are a lattice with respect to absolute continuity (see Proposition \ref{prop:measuretypeslattice}). Thus, stationary factorizing measure types of random closed sets in $[0,1]$ build  a rich, fascinating structure in their own right. Although they fulfil seemingly weak equivalence conditions, a lot of results can be obtained by importing ideas from product systems, e.g.\ lattice structures in Proposition \ref{prop:measuretypeslattice} and non-determinism  in the proof of Proposition \ref{prop:statfactmeastype=ps}.  On the other side, the hole-filling procedure still lacks a counterpart for product systems.

This could all together give the impression that we can hope for a complete structure theory at least of type $\mathrm{II}$ product systems. Unfortunately, things turn out to be not so fine. Connected with representations of abelian $W^*$-algebras (which are provided here by Theorem \ref{th:RACS} too) are direct integrals of Hilbert spaces. Using this tool (or its variant, the Hahn-Hellinger Theorem) we succeed in Proposition \ref{prop:typeIIIintotypeII} to relate  any (type $\mathrm{III}$) product system to a type $\mathrm{II}$ one in such a manner that isomorphy classes are preserved. Consequently, the structure theory of  type $\mathrm{II}$ product systems is at least as complex as that of  type $\mathrm{III}$ product systems.  In more positive words, this  parallels the   relation between type $\mathrm{III}$   and type $\mathrm{II}$ $W^*$-factors \cite{OP:Tak73} and shows that for  completing the classification of product systems one can  concentrate on the type $\mathrm{II}$ case  now. Further, the outstanding r\^ole of tensor products in the theory of product systems is reduced since there exist product systems of type $\mathrm{II}_1$ not isomorphic to any nontrivial tensor product, see Proposition \ref{prop:exampleII_2}. The gap should be filled by direct integral decompositions established  in section \ref{sec:direct integrals ps}.    

\textsc{Tsirelson} and coauthors established interesting new examples of product systems coming  from random sets and random processes. These examples fit well into the  general structure of factorizing Standard Borel spaces and stationary factorizing measure types, see Theorem \ref{th:generalproductsystemfrom stationaryfactorizingmeasuretype}.  There we obtain a quite general existence result of corresponding product systems covering all examples of \textsc{Tsirelson} \cite{Tsi03,Tsi00,Tsi99a,C:TV98}.  These results make us  feel  that there are only few  examples from classical probability theory which are readily  applicable  in this context, as e.g.\ (generalised) Gaussian processes in \cite{Tsi03,Tsi00}. The  reason for this is that, typically, classical probability theory studies situation  with a high degree of independence. The examples of \textsc{Tsirelson} tell us  that, at least in the context of product systems, \emph{dependent} structures become important and tractable as well. The notion of a factorizing measure type inherently has not much in common with independence. It concentrates in identifying 0-1-laws characterising a certain type of models. To quantify  the degree of this dependence we develop a further invariant of product systems, a kind of algebra-valued Hausdorff measure. The dependence  degree of a product system   could then be measured by a certain Hausdorff dimension, see Note  \ref{rem:dependence degree}.\label{page:dependence measure} There  $1$ corresponds to the highest degree of  independence, like present in  Levy processes. The less the Hausdorff dimension, the stronger is the dependence structure  in the model. In that direction, we expect  positive impulses from the theory of product systems for classical probability theory too. Appearance of dependent structures emphasises that not product constructions (like  the tensor product of product systems, which is related to independence) but procedures related to conditioning  (like the hole-filling procedure and direct integrals) are of most importance.

There is yet another motivation for these studies from classical probability theory. In order to classify the information structures of random processes it is natural to consider filtrations, i.e.\ increasing families of $\sigma$-fields. Necessarily, one should  complete all  $\sigma$-fields of interest. Equivalently, one could study the corresponding $W^*$-algebras of almost surely bounded functions, giving a more unified view on this problem. In the case of interest, the $W^*$-algebras factorise over disjoint time intervals, the deep reason for tensor product systems of Hilbert spaces to show up. In the present work, we relate  a   tensor product systems of Hilbert spaces closely  to some  tensor product systems of  $W^*$-algebras equipped with a factorizing set of normal states on them, see Proposition  \ref{prop:tensorproductW*=ps}. This approach is related to  another connection between product systems and $E_0$-semigroups. \textsc{Bhat} shows in \cite{OP:Bha96} that  besides \ref{eq:def ps from E_0(Arveson)} $\E$ could equally be described by
\begin{displaymath}
  \E_t=\alpha_t(\Pr\psi) \H
\end{displaymath}
where $\psi\in\H$ is an arbitrary unit vector. In the present paper, we propose and use  even a third variant. There, the main observation is that $\E_t$ is also the GNS Hilbert space \cite[Section 2.3.3]{BR87} of the restriction of a pure normal  state on $\B(\H)$ to the algebra $\alpha_t(\B(\H))'$, if the state factorises with respect to $\B(\H)=\alpha_t(\B(\H))'\otimes\alpha_t(\B(\H))$. To use this fact in terms of product systems, we would like to have states which factorise in this way for all $t\in\NRp$. Unfortunately, such states do not exist if the product system is of type $\mathrm{III}$, see \cite{OP:Arv97} or Proposition \ref{cor:decomposable ps is type I} below. Therefore, we have to connect the different GNS-Hilbert spaces in a natural way. The right structure to achieve  this goal is the following. Suppose  $(H_p)_{p\in P}$ is a family of Hilbert spaces equipped with unitaries $(U_{p,p'})_{p,p'\in P}$, $\map{U_{p,p'}}{H_p}{H_{p'}}$ fulfilling $U_{p,p}=\unit$ and   the cocycle relations 
\begin{displaymath}
  U_{p',p''} U_{p,p'} =U_{p,p''}\dmf{p,p',p''\in P}.
\end{displaymath}
Then 
\begin{displaymath}
  H(P,U)=\set{(\psi_p)_{p\in P}:\psi_p\in H_p\forall p\in P, \psi_{p'}=U_{p,p'}\psi_p\forall p,p'\in P}
\end{displaymath}
is a Hilbert space isomorphic to any $H_p$. This very flexible Hilbert space bundle structure reminds very much of manifolds. In fact, in terms of $L^2$ spaces, such a  structure was first considered  by \textsc{Accardi} \cite{Q:Acc76} to define the intrinsic  Hilbert space  of half-densities associated with a  manifold. In a more general way, it was  used by \textsc{Tsirelson} to describe his  example product systems in an elegant  way. The  relation  of product systems or $E_0$-semigroups with  such bundles of  (GNS-)Hilbert spaces    proves useful. E.g., we can  settle the messy question whether an algebraic product system of Hilbert spaces has a compatible measurable structure. We derive an equivalent condition in terms of continuity of a unitary group of local shifts in Theorem \ref{th:intrinsicmeasurability}. By  this result   we can show that the  measurable structure of product systems is unique upto isomorphy and  can be derived from its algebraic structure provided  it exists. This result should have a strong impact on the axiomatics of product systems since it puts more emphasis  on  the algebraic than the measurable structure.  To close this circle of ideas, a continuous tensor product system of $W^*$-algebras, represented as von Neumann subalgebras of some $\B(\H)$, forms a quantum stochastic filtration.  If each representative is a type $\mathrm{I}$ factor, the canonical shift induces an $E_0$-semigroup on $\B(\H)$.

\bigskip
We expect very fascinating  things coming up from the theory of tensor product systems of  $W^*$-algebras and Hilbert spaces. It is clear that these structures play a crucial r\^ole in theory of quantum Markov processes and quantum Markov random fields \cite{Q:Lie01b,Q:AF01a,Q:AF01b,Q:LS01}. Quantum Markov processes and quantum stochastic differential equations are a main tool in the dilation theory of quantum dynamical semigroups \cite{Par92,Mey93,OP:Bha01}. Unfortunately, the most studied  examples of quantum stochastic calculi are related to symmetric, antisymmetric or full Fock spaces, where the corresponding time shift semigroup is of type $\mathrm{I}$ \cite{OP:Pow99}. Below, there are indications that there should exist  a  quantum stochastic calculus on type $\mathrm{II}$ product systems too. E.g., we can derive a variant of the $\rlap{$\DS\sum$}\DS\int$ lemma of Boson stochastic calculus \cite{Lin93} for stationary factorizing measure types of random closed sets, see Corollary \ref{cor:sumintegrallemma}.  In the context of dilation theory, we should mention that the recent notion of tensor product systems of  Hilbert modules of \textsc{Bhat} and \textsc{Skeide} \cite{Q:BS00} is deeply related to this topic too. We expect that certain techniques developed in the present  paper can prove useful in the context of these more complicated structures too. At least, there develops a theory of type $\mathrm{I}$ continuous tensor product systems of  Hilbert modules \cite{Q:LS01,Q:BBLS01} paralleling that of \textsc{Arveson} \cite{Arv89}. 

\bigskip
To keep this work readable as much as possible  without knowing  other works on product systems   we provide several alternative  proofs  of results of \textsc{Arveson}, \textsc{Bhat} and \textsc{Powers}, especially when we think the proofs  are instructive. From the  point of technique we use essentially the above mentioned GNS-representations and do not touch such powerful tools like the spectral $C^*$-algebra \cite{OP:Arv90a} and the boundary representation of spatial $E_0$-semigroups \cite{OP:Pow91}. Further, we concentrate on invariants of product systems, although most of the results have counterparts in the classification of $E_0$-semigroups upto  conjugacy.

\bigskip 
At the end of this introduction, we want to explain the general structure of this work. After introducing preliminary notions, we show directly, how general product systems are intimitately connected with a stationary factorizing measure type of random closed  sets in $[0,1]$. This leads us to the study of such measure types in section \ref{sec:rs2ps}. Section \ref{sec:hierarchy} and Section \ref{sec:direct integrals} develop a finer structure theory of product systems completing Section \ref{sec:ps2rs}. Especially, Section \ref{sec:hierarchy} studies the implications for the lattice of product subsystems and derives a corresponding lattice structure on the stationary factorizing measure types of random closed sets in $[0,1]$. Section \ref{sec:direct integrals} studies implications of direct integral representations, and provides especially another characterisation of type $\mathrm{I}$ product systems and the reduction of the classification of type  $\mathrm{III}$ product systems to that of type $\mathrm{II}$ ones.    Section \ref{sec:algebraic approach} is devoted to the  study of   measurability questions in product systems by means of GNS representations. This gives us the opportunity to derive in section \ref{sec:randommeasuresandincrements} a   criterion for  existence of product systems built from stationary factorizing measure types on general Standard Borel spaces together with some  analysis of special cases of interest: random closed sets, random measures and random increment processes. Section \ref{sec:nonsep} analyses the analogue of Theorem \ref{th:RACS} if the fibres of the product system are not separable any more and we obtain random bisets rather than random sets. A further, algebraic invariant, derived from ideas in \cite{Tsi00}, is presented in section \ref{sec:algebraic invariant} and computed for product systems with unit as function of the previously introduced invariants. Section \ref{sec:conclusio} contains  some conclusions and presents ideas for further research.

\section{Basics}
\label{sec:basic}
For the natural numbers, positive natural numbers, integers, rational, real, positive real and complex numbers respectively we use the symbols $\NN=\set{0,1,2,\ldots}$, $\NNp=\NN\setminus\set0$, $\NZ$, $\NQ$, $\NR$, $\NRp=[0,\infty)$ and $\NC$ respectively. Let $\NT=\set{z\in\NC:\abs z=1}$ denote the circle. Further, $\ell$ denotes Lebesgue measure, with base space varying between $[0,1],\NRp,\NR$ and $\NT$ according to our needs. 

A \emph{Polish space} is a separable topological space which can be  metrized in a way such that it is complete \cite{Bau91,C:Coh80,OP:Tak79}.  A \emph{Standard Borel space} is a measurable space which is Borel isomorphic to a  Polish space equipped with its Borel field.

All Hilbert spaces are complex.

\bigskip
For a random variable $\xi$ let $\mathscr{L}(\xi)$  denote its distribution. If $\mu$ is a measure on measurable space $(X,\XG)$ and $Y\in\XG$ we use the notation $\mu|Y$ for the measure $\mu(\cdot\cap Y)$. A  \emph{measure type} on   $(X,\XG)$ is an equivalence class $\M$ of (probability) measures on $(X,\XG)$. Equivalence of two measures $\mu,\mu'$ (symbol $\mu\sim\mu'$) means $\mu(Y)=0$ iff $\mu'(Y)=0$ for all $Y\in\XG$. This means both $\mu\ll\mu'$ and $\mu'\ll\mu$, i.e.\ mutual absolute continuity.  We translate all operations on probability measures to operations on measure types in a natural way. E.g., let 
$q$ be   a stochastic kernel. I.e., $\map q{X\times\XG'}{[0,1]}$ such that  $q(x,\cdot)$ is a probability measure and $q(\cdot,Y)$ is a measurable function for all $Y\in\XG'$. Then for a probability measure $\mu$ the construction $\mu\circ q(\cdot)=\int\mu(\d x)q(x,\cdot)$ yields again a probability measure and we can define $\M\circ q=\set{\mu:\mu\sim \mu'\circ q\text{~for one and thus all~}\mu'\in\M}$.

 A measure type  $\M$ determines the algebra $L^\infty(\M)$ of (equivalence classes of) $\M$-essentially bounded functions. Moreover, we  construct also a Hilbert space denoted $L^2(\M)$.  For this goal define  unitaries  $\map{U_{\mu,\mu'}}{L^2(\mu)}{L^2(\mu')}$ for any $\mu,\mu'\in\M$  through
\begin{equation}
  \label{eq:defUmumu'}
U_{\mu,\mu'}\psi(x)=\sqrt{\frac{\d\mu'}{\d\mu\phantom{'}}(x)}\psi(x)  \dmf{\psi\in L^2(\mu),\mu-\text{a.a.~}x\in X}.
\end{equation}
Then 
 \begin{equation}
\label{eq:defL2measuretype}
   L^2(\M)=\set{(\psi_\mu)_{\mu\in\M}: \psi_\mu\in L^2(\mu)\forall \mu\in\M,\psi_{\mu'}=U_{\mu,\mu'}\psi_\mu \forall \mu,\mu'\in\M }
 \end{equation}
is a Hilbert space with the inner product
\begin{displaymath}
  \scpro{\psi}{\psi'}_{L^2(\M)}=\int \ovl{\psi_\mu}\psi'_\mu\d\mu
\end{displaymath}
being independent of the choice of $\mu\in\M$. 
\begin{remark}
  For this construction we find an analogue for general $W^*$-algebras, analogous to independence of the space $L^\infty(\M)=L^\infty(\mu)$ from the choice of  $\mu\in\M$, in section \ref{sec:2measurability}, see \ref{eq:defH(R)} and Example \ref{ex:GNSLinfty}.

Further details about the above construction which originated in  \cite{Q:Acc76} can be found in \cite{Tsi00b,Tsi00}.
\end{remark}
\begin{definition}
  \label{def:meashilbertspacefamily}
  Let $(X,\XG)$ be a measurable space. A family $\H=(\H_x)_{x\in X}$ is a \emph{measurable family of Hilbert spaces}, if each $x\in X$, $\H_x$ is a Hilbert space and there exists a set $\H^0$ of sections $X\ni x\mapsto h_x\in \H_x$ such that $x\mapsto\scpro{h_x}{h'_x}_{\H_x}$ is measurable for  all $h,h'\in\H^0 $ and $\set{h_x:h\in\H^0}$ is total in $\H_x$ for all $x\in X$. Then a \emph{measurable  section} $x\mapsto h_x$ is such that $x\mapsto\scpro{h_x}{h'_x}$ is measurable for  all $h'\in\H^0 $.

For a measurable map  $\map fX{\tilde X}$  and  measurable families $\H=(\H_x)_{x\in X}$, $\tilde \H=(\tilde\H_{\tilde x})_{\tilde x\in \tilde X}$  of Hilbert spaces  we call a family $(A_x)_{x\in X}$ consisting of bounded operators $A_x\in\B(\H_x,\tilde\H_{f(x)})$  measurable  if   $x\mapsto\scpro{A_xh_x}{\tilde h_{f(x)}}$ is measurable for  all $h\in\H^0$, $\tilde h\in{\tilde \H}^0$. Unless stated otherwise, we assume implicitly $f(x)=x$ if $X=\tilde X$.
\end{definition}
\begin{remark}
  From definition,  the structure of $\H$ depends  on the set $\H^0$. But, one can show \cite[4.4.2]{BR87} that if all $\H_x$ are separable  then all  measurable families (corresponding to different  sets $\H^0$) are isomorphic. Since we deal with separable Hilbert spaces only  we need not discuss further subtleties here. 
\end{remark}

\bigskip
In the sequel we use frequently the interval sets $I_{s,t}=\set{(s',t'):s\le s'<t'\le t}$ for $0\le s<t<\infty$. Similarly, we set $I_{0,\infty}=\set{(s',t'):0\le s'<t'< \infty}$.

\bigskip
Let $K$ be a locally compact second countable Hausdorff space. The space $\FG_K$ of all closed subsets  of $K$, topologized  by the myopic hit-and-miss topology  generated by the sets $\set{Z:Z\cap K'=\emptyset}$, $K'$ compact and $\set{Z:Z\cap G\ne\emptyset}$, $G$ open, is a   compact second countable Hausdorff space \cite[section 1-4]{C:Mat75}.  A \emph{random closed  set in $K$}  is just a $\FG_K$ valued random variable, i.e.\ its distribution is a Borel probability measure on $\FG_K$. Observe that $\FG_K$ is a continuous semigroup under $\cup$, leading to a convolution $\ast$ of probability measures by $\mu_1\ast\mu_2(Y)=\mu_1\otimes\mu_2(\set{(Z_1,Z_2):Z_1\cup Z_2\in Y})$. 

\bigskip
Let $\H$ be a Hilbert space. If $\H'\subseteq\H$ is a closed subspace denote the orthogonal projection onto $\H'$ by $\Pr{\H'}$. Especially, for $h\in\H$ we set $\Pr h=\Pr{\NC h}$. For a family of orthogonal projections $(P_i)_{i\in I}\subset\B(\H)$, let $\bigwedge_{i\in I}P_i$ and $\bigvee_{i\in I}P_i$ denote their greatest lower and least upper bound respectively. For limits in the strong and weak topology on $\B(\H)$ we use the symbols $\slim$ or $\limitsto{s}{}$ and $\wlim$ or $\limitsto{w}{}$ respectively. 

Von Neumann algebras are  weakly closed subalgebras of some $\B(\H)$. For any set $S\subset\B(\H)$ its commutant is $S'=\set{a\in\B(\H):as=sa\forall s\in S}$ and $S''$ its bicommutant.  Both are von Neumann algebras if $a\in S\Rightarrow a^*\in S$. All used operator algebras possess a unit denoted $\unit$  and homomorphisms of operator algebras are unital $\ast$-homomorphisms. Further details could be found in \cite[Section 2.4]{BR87}. Again, von Neumann algebras in $\B(\H)$ form a lattice, denote $\bigwedge_{i\in I}\B_i$ and $\bigvee_{i\in I}\B_i$  the greatest lower and least upper bound respectively. A map $\alpha$ between von Neumann algebras  is called normal if it is $\sigma$-strongly continuous, or, equivalently, if $\sup_{d\in D}\alpha(a_d)=\alpha(\sup_{d\in D}a_d)$ for every bounded increasing net $(a_d)_{d\in D}$ of positive operators. 

A functional $\map\eta{\B}\NC$, $\B\subset\B(H)$ a von Neumann algebra, is called  \emph{normal state} if for some positive trace-class operator $\varrho$ of trace $1$ we have $\eta(a)=\Tr\varrho a$ for all $a\in\B$. $\eta$ is \emph{faithful}, if $\eta(a^*a)=0$ implies $a=0$ or, equivalently, if $\varrho$ could be chosen to have  trivial kernel.
 
For any Hilbert space  $\H$, denote $\N(\H)$ the set of all von Neumann subalgebras of $\B(\H)$. We introduce a topology on $\N(\H)$  identifying $\A\in\N(\H)$ with $B_1(\A)=\set{a\in\A:\norm a\le1}$ and using  on $\FG_{B_1(\B(\H))}$ the hit-and-miss topology with respect to the weak topology on  $B_1(\B(\H))$ which turns it into  a compact metric  space. 

A representation of an algebra $\B$ on $\H$ is an homomorphism from $\A$ into $\B(\H)$.

\bigskip 
The symbol \proofend\ indicates the end of a proof or a statement which is either trivial or proven elsewhere.

\section{From Product Systems to Random Sets}
\label{sec:ps2rs}

\subsection{Product Systems}
\label{sec:psbasic}

Loosely speaking, a product system (of Hilbert spaces) is a family $\sg t\E$ with $\E_{s+t}=\E_s\otimes\E_t$ in a consistent manner for all $s,t\ge0$. 
\begin{definition}
  \label{def:ps}
  A \emph{continuous tensor product system of Hilbert spaces} (briefly: product system) is a pair  $\E=(\sg t\E,\sg{s,t}V)$, where $\sg t\E$ is a measurable family of separable Hilbert spaces with $\E_0=\NC$ and  $(V_{s,t})_{s,t\in\NRp}$ is a  family of unitaries  $\map{V_{s,t}}{\E_s\otimes\E_t}{\E_{s+t}}$ with 
  \begin{displaymath}
    V_{0,t}z\otimes \psi_t=z\psi_t=V_{t,0} \psi_t\otimes z \dmf{t\ge0,\psi_t\in\E_t,z\in\E_0=\NC}
  \end{displaymath}
which is  associative in the sense of 
\begin{equation}
  \label{eq:associativity Ust}
 V_{r,s+t}\circ(\unit_{\E_r}\otimes V_{s,t})= V_{r+s,t}\circ( V_{r,s}\otimes\unit_{\E_t})\dmf{r,s,t\in\NRp},
\end{equation}
and  measurable with respect to $\map f{\NRp\times\NRp}\NRp$, $f(s,t)=s+t$, and   the measurable families $(\E_t\otimes\E_s)_{s,t\in\NRp}$, $\sg t\E$ of Hilbert spaces. 

Two product systems $\E$, $\E'$ are \emph{isomorphic} if there is a measurable family $\sg t\theta$ of unitaries $\map{\theta_t}{\E_t}{\E'_t}$ with
\begin{displaymath}
  \theta_{s+t}V_{s,t}=V'_{s,t}\theta_{s}\otimes\theta_{t}\dmf{s,t\in\NRp}.
\end{displaymath}
\end{definition}
\begin{remark}
The above  definition differs  from that of \textsc{Arveson} \cite{Arv89}. There  a Standard Borel structure on $\bigcup_{t>0}\E_t$ is required and the multiplication  $(\bigcup_{t\in\NRp}\E_t)^2\ni(x,y)\mapsto xy$, $x_sy_t:=V_{s,t}x_s\otimes y_t$ is associative, measurable and fibrewise bilinear. In fact, both approaches are essentially equivalent as is shown in  Lemma \ref{lem:both measurability definitions} below. We prefer the above version since we encounter similar, but more complicated structures in section \ref{sec:direct integrals} which are not so easy to capture by the old definition.   Further, differently from \cite{Arv89}, we  incorporate $t=0$ by choosing $\E_0=\NC$ and  allow for  (\emph{trivial}) product systems like $\sgp t\NC$ and $\sgp t{\set0}$. Further, we find the measurable family approach more convenient since we need those for the more complex structures in section \ref{sec:direct integrals}. This change of definition  does not affect the classification task, to characterize all equivalence classes of product systems under isomorphy. The present work amounts to collect more invariants for product systems, i.e.\ maps from product systems to some space which are constant on  isomorphy classes.
\end{remark}
\begin{remark}
Presently, one task in the construction of product systems consists in  proving measurability of the  multiplication, see e.g.\ \cite{Tsi00}.  Below, in Corollary \ref{cor:onlyonemeasurablestructure}, we show that for any algebraic (i.e.\ without any of the measurability requirements above) product system  all  consistent measurable structures (if there is one) lead to isomorphic product systems. Moreover, we can decide easily whether there is a consistent measurable structure.  Clearly, this   is a more complicated version of  the above cited  result for general measurable families of Hilbert spaces.
\end{remark}
\begin{example}
  Consider the trivial one dimensional case  $\E_t\cong\NC$ for all $t\in\NRp$ not covered by the definition in \cite{Arv89}. Then $V_{s,t}$ is encoded by the number $m(s,t)=V_{s,t}1\otimes1\in\NT$. Clearly, \ref{eq:associativity Ust} is equivalent to  
  \begin{displaymath}
      m(r,s+t)m(s,t)=m(r+s,t)m(r,s).
  \end{displaymath}
We will solve  this equation in Lemma \ref{lem:multiplieronR} below, it turns out that there is essentially only the usual  product structure on this  product system. A measurable structure is equivalent to choose an arbitrary (but to be called measurable) curve  $\sg t{z^0}$ in $\NT$. $\sg {s,t}V$ is then measurable if $(s,t)\mapsto m(s,t)z^0_sz^0_t/z^0_{s+t}$ is measurable. Again, there is (upto isomorphism) only one such measurable structure on $\E$ \cite{OP:Arv89}. Further, dimensional considerations show that this is the only  product system \cite[p.19, proof of Proposition 2.2]{Arv89} with $0<\dim \E_t<\infty$. 
\end{example}
\begin{example}
\label{ex:Fockasps}
Let $\H$ be a separable Hilbert space and define its  $n$-fold symmetric tensor product  $\H^{\otimes n}_{\mathrm{sym}}=\ovl{\mathrm{lh}\set{u^{\otimes n}:u\in\H}}\subset\H^{\otimes n}$. The \emph{symmetric Fock space over $\H$} is the Hilbert space 
\begin{displaymath}
  \Gamma(\H)=\NC\oplus \bigoplus_{n=1}^{\infty}\frac1{\sqrt{n!}}\H^{\otimes n}_{\mathrm{sym}}.
\end{displaymath}
Introducing \emph{exponential vectors} $\psi_h\in\Gamma(\H)$,
\begin{displaymath}
\psi_ h=1\oplus 
\bigoplus_{n=1}^{\infty}h^{\otimes 
n}\dmf{h\in\H} 
\end{displaymath}
we have   $\Gamma(\H_1\oplus\H_2)\cong\Gamma(\H_1)\otimes\Gamma(\H_2)$ under the isomorphism extending $\psi_{h_1\oplus h_2}\mapsto\psi_{h_1}\otimes\psi_{h_2}$ \cite[Proposition 19.6]{Par92}. Fix a separable Hilbert space $\K$ and set $\Gammai(\K)_t=\Gamma(L^2([0,t],\K))$. We find that  the multiplication given by $\sg{s,t}V$,
\begin{displaymath}
  V_{s,t}\psi_f\otimes\psi_g=\psi_{f+\rho_sg}\dmf{s,t\in\NRp, f\in L^2([0,s],\K), g\in L^2([0,t],\K)},
\end{displaymath}
where $\sg t\rho$ are the right shift isometries, $\rho_sg(s+t)=g(t)$, is associative. Measurability of the multiplication follows from continuity of $\sg t\rho$ such that  $\Gammai(\K)=(\sg t{\Gammai(\K)},\sg{s,t}V)$  is a product system.  
\end{example}
Further important notions in  product systems are given by 
\begin{definition}
  \label{def:unit}
  A \emph{unit} of a product system $\E$ is a measurable section   $u=\sg tu$ through $\sg t\E$  fulfilling
  \begin{equation}
   \label{eq:defunit}
 u_{s+t}=u_s\otimes u_t=V_{s,t}u_s\otimes u_t\dmf{s,t\ge0}
  \end{equation}
and $u_t\ne0\forall t\ge0$. The set of all units of $\E$ is denoted by $\Usk(\E)$.

  \label{def:subps}  
We say that $\F=\sg t\F$ is a \emph{product subsystem} of $\E$ if $\F_t\subset\E_t$ is a closed subspace for all $t\in\NRp$  and multiplication on $\F$ is given by the operators $(V^\F_{s,t})_{s,t\in\NRp}$,  $V^\F_{s,t}=V^\E_{s,t}\restriction\F_s\otimes\F_t$ which are  unitaries.
\end{definition}
For every product system $\E$ its units generate   a   product subsystem $\E^{\Us}$  of $\E$ by
\begin{equation}
\label{eq:def subsystem generated by units}
 \E^{\Us}_t=\overline{ \mathrm{lh}\set{u^1_{t_1}\otimes u^2_{t_2}\otimes \cdots\otimes u^k_{t_k}:{t_1+\cdots+t_k=t,u^1,\dots,u^k\in\Usk(\E)}}}\subseteq\E_t\dmf{t\in\NRp}.
\end{equation}
In fact, $\E^{\Us}$ is the smallest product subsystem  of $\E$ containing all of the units of $\E$. \textsc{Arveson} proved in \cite[section 6]{Arv89} that $\E^{\Us}$ is isomorphic to a  product system $\Gammai(\K)$ of Fock spaces for  some separable Hilbert space $\K$ (for details see Example \ref{ex:Fockasps}).  A simple proof  based only  on  results of the present paper  is given in Corollary \ref{cor:exponentialhilbertspace}. The different cases concerning the relative position of $\E$ and $\E^\Us$ motivate
\def\citt{\cite{OP:Pow97a}}
\begin{definition}[\citt]
\label{def:types of ps}
A product system $\E$ is said to have \emph{type $\mathrm{I}$}, if $\E=\E^{\Us}$. If neither $\E^\Us=\E$ nor $\E^\Us=\set0$, $\E$ has  \emph{type $\mathrm{II}$}. Product systems with $\E^\Us\ne\set0$ are called \emph{spatial} product systems.

In both cases, if $\E^\Us$ is isomorphic to $\Gammai(\K)$, we speak of \emph{type $\mathrm{I}_{\dim\K}$ or $\mathrm{II}_{\dim\K}$} respectively. $\dim\K$ is called \emph{numerical index} of the product system. 

\emph{Type $\mathrm{III}$} product systems fulfil  $\E^\Us=\set0$. Their numerical index is set to $-\aleph_1$.
\end{definition}
\begin{example}
 It is easy to derive that each unit $u\in\Usk(\Gammai(\K))$ of a  Fock product system $\Gammai(\K)$ has the form $u=u^{z,k}$,
\begin{equation}
\label{eq:defuzk}
  u^{z,k}_t=e^{zt}\psi_{\chfc{[0,t]}k}\dmf{t\in\NRp},
\end{equation}
for some  $k\in\K$ and $z\in\NC$. Denseness of step functions in each $L^2([0,t],\K)$ and totality  of the exponential vectors shows  that $\Gammai(\K)^\Us=\Gammai(\K)$, i.e.\ $\Gammai(\K)$ is of type $\mathrm{I}_{\dim\K}$.  

In \cite{OP:Arv89} it was shown that the numerical index is an invariant of the product system being additive under tensor products of product systems.  Our goal is to find further invariants of product systems being nontrivial if $\E$ is not type $\mathrm{I}$. 
\end{example}

In the following, we want to focus our considerations to one Hilbert space in the family $\sg t\E$, say $\E_1$. A first, useful reduction is achieved by  the following proposition, which we prove for the sake of readability in section \ref{sec:2measurability} on page \ref{page:proof prop:psbyE_1}. 
\begin{proposition}
\label{prop:psbyE_1}
  Product systems $\sg t\E$  correspond via restriction of the parameter space one-to-one to measurable  families $\sgi t\E$  of Hilbert spaces equipped with a measurable family of unitaries $(V_{s,t})_{s,t\in[0,1],s+t\le1}$ fulfilling \ref{eq:associativity Ust} restricted to $r,s,t\ge0$, $r+s+t\le1$. 
\end{proposition}
\begin{remark}
  From \cite{OP:Arv90} we know that we can represent any product system as product system associated with  an $E_0$-semigroup $\sg t\alpha$ on some Hilbert space $\H$. From the results  in \cite{Arv89,OP:PR89,OP:Arv90} it follows that there is a unitary group $\sg tU$ on $\H\otimes\H$ such that $\unit\otimes \alpha_t(b)=U^*_t\unit\otimes b U_t$ for all $t\ge0$, $b\in\B(\H)$. Both results are reproved  in Theorem \ref{th:E_0-semigroup} below. Since this representation is not intrinsic, we find it easier to work with the intrinsic structure on the  Hilbert space $\E_1$. Moreover, the measure types of random sets, we associate with  $\E$ below, would have   depended  on the specific representation if we had  considered the whole half-line. These structures are only  equivalent, if we restrict our considerations to a compact interval, for our convenience  $[0,1]$. This reflects the fact that there are at least two   natural equivalence relations for $E_0$-semigroups: Cocycle conjugacy (related to product systems) and conjugacy. We deal here  with the former only. 
\end{remark}

\subsection{Random Sets in Product Systems}
\label{sec:rsinps}

Next we study structures on $\E_1$ induced by the (restricted) multiplication. 

  There is an important set of unitaries $(\tau_t)_{t\in(0,1)}\subset\B(\E_1)$ acting with regard to  the representations
$\E_{1-t}\otimes\E_t\cong\E_1\cong\E_t\otimes\E_{1-t}$ as \emph{flip}:
\begin{equation}
\label{eq:defflipgroup}
  \tau_t x_{1-t}\otimes x_t=x_t\otimes  x_{1-t}\dmf{x_{1-t}\in\E_{1-t}, x_t\in\E_t}.
\end{equation}
Moreover, setting $\tau_0=\unit$ and   $\tau_{t+k}=\tau_{t}$ for any $t\in[0,1]$, $k\in\NZ$, we get a periodic  1-parameter  group $\gr t\tau$: 
\begin{proposition}
\label{prop:tau_tiscontinuous}
  $\gr t\tau$ is a strongly  continuous 1-parameter group of unitaries.
\end{proposition}
\begin{proof}
First we have to show $\tau_t\tau_s=\tau_{s+t}$. It  is enough to prove  this for $0\le s,t<1$. If $s+t\le1$, we obtain for all $x_s\in\E_s$, $x_t\in\E_t$ and $x_{1-t-s}\in\E_{1-t-s}$
\begin{displaymath}
  \tau_t\tau_s x_{1-t-s}\otimes x_t\otimes x_s=\tau_tx_s\otimes x_{1-t-s}\otimes x_t=x_t \otimes x_s\otimes x_{1-t-s}=\tau_{s+t}x_{1-t-s}\otimes x_t\otimes x_s.
\end{displaymath}
If $s+t>1$, we find for $x_{1-t}\in\E_{1-t}$, $x_{1-s}\in\E_{1-s}$ and $x_{s+t-1}\in\E_{s+t-1}$
\begin{eqnarray*}
  \tau_t\tau_s x_{1-s}\otimes x_{1-t}\otimes x_{s+t-1}&=&\tau_tx_{1-t}\otimes x_{s+t-1}\otimes x_{1-s}\\
&=&x_{s+t-1}\otimes x_{1-s}\otimes x_{1-t}\\
&=&\tau_{s+t} x_{1-s}\otimes x_{1-t}\otimes x_{s+t-1}.
\end{eqnarray*}
This establishes $\tau_{t+s}=\tau_t\tau_s$ for all $s,t\in\NR$. Further, the flip on tensor products of Hilbert spaces is unitary.

  To prove that $\gr t\tau$ is strongly continuous it is sufficient to show  that $(\tau_t)_{t\in(0,1)}$  is (weakly) measurable \cite[Theorem VIII.9]{RS72}.  Without loss of generality assume $\E$ is not trivial and let $\sg t{e^n}$, $n\in\NN$, be measurable sections of $\E$ such that for all $t>0$, $(e^n_t)_{n\in\NN}$ is a complete orthonormal system of $\E_t$. From measurability  of the multiplication it follows that $(s,t)\mapsto\scpro{e^k_{s+t}}{V_{s,t}e^n_s\otimes e^m_t}$ is measurable. Thus  $t\mapsto \scpro{e^k_{1}}{V_{1-t,t}e^n_{1-t}\otimes e^m_t}$ is measurable for all $k,n,m\in\NN$ too. Consequently,  $t\mapsto V_{t,1-t}e^n_t\otimes e^m_{1-t}$ is a measurable function in $\E_1$ for all $n,m\in\NN$. Since $\tau_t V_{1-t,t} e^n_{1-t}\otimes e^m_{t}= V_{t,1-t}e^m_{t}\otimes e^n_{1-t}$, we derive  from complete orthonormality of $( V_{1-t,t}e^n_{1-t}\otimes e^m_t)_{n,m\in\NN}$
  \begin{eqnarray*}
    \scpro{e^k_1}{\tau_te^l_1}&=&\sum_{n,m\in\NN}\scpro{e^k_1}{\tau_t  V_{1-t,t}e^n_{1-t}\otimes e^m_t}\scpro{ V_{1-t,t}e^n_{1-t}\otimes e^m_t}{e^l_1}\\
&=&\sum_{n,m\in\NN}\scpro{e^k_1}{ V_{t,1-t}e^m_{t}\otimes e^n_{1-t}}\scpro{ V_{1-t,t}e^n_{1-t}\otimes e^m_t}{e^l_1}.
  \end{eqnarray*}
This shows that $(\tau_t)_{t\in(0,1)}$  is weakly  measurable, i.e.\ strongly continuous. 
\end{proof}
\begin{remark}
\label{rem:spectral theory of flip as invariant}
It is easy to see that $\gr t\tau$ behaves well under isomorphisms of product systems. Thus it seems reasonable to study the spectrum of (the generator of) $\gr t\tau$ as  a further  invariant for product systems. But, at the moment, there is no indication that there is a nontrivial  product system for which  this invariant is different from $\NZ$ with infinite multiplicity. 
\end{remark}

With  a  product system  $\E$  there also come   two kinds   of  projection families  in $\E_1$. For a fixed unit $u=\sg tu$ set, regarding the representation $\E_1=\E_s\otimes \E_{t-s}\otimes\E_{1-t}$,
\begin{equation}
  \label{eq:defPust}
  \mathrm{P}_{s,t}^u=\Pr{\E_s\otimes \NC u_{t-s}\otimes\E_{1-t}}=\unit_{\E_s}\otimes\Pr{\NC u_{t-s}}\otimes\unit_{\E_{1-t}}\dmf{(s,t)\in I_{0,1}}
\end{equation}
Similarly, with $\E^{\Us}_t$ given by \ref{eq:def subsystem generated by units} define
\begin{equation}
  \label{eq:defPUst}
  \mathrm{P}_{s,t}^{\Us}=\Pr{\E_s\otimes\E^{\Us}_{t-s}\otimes\E_{1-t}}=\unit_{\E_s}\otimes\Pr{\E^{\Us} _{t-s}}\otimes\unit_{\E_{1-t}}\dmf{(s,t)\in I_{0,1}}.
\end{equation}
Clearly, these projections belong to the von Neumann subalgebras
\begin{equation}
\label{eq:defAst}
  \A_{s,t}=\unit_{\E_s}\otimes\B(\E_{t-s})\otimes\unit_{\E_{1-t}}\dmf{(s,t)\in I_{0,1}}.
\end{equation}
Like $\mathrm{P}^u$, $\mathrm{P}^\Us$ we call a set of operators $(b_{s,t})_{(s,t)\in I_{0,1}}\subset\B(\E_1)$ \emph{adapted}, if $b_{s,t}\in\A_{s,t}$. Similar constructions are valid for arbitrary product subsystems $\F$ instead of  $\E^\Us$ or $\sg t{\NC u}$, see \ref{eq:defPFst} in section \ref{sec:hierarchy}.

The flip group acts in a simple manner  on these projections. This action is best written in terms of the automorphism group $\gr t\sigma$,
\begin{equation}
  \label{eq:defshiftautomorphisms}
\sigma_t(a)=\tau_t^*a\tau_t\dmf{a\in\B(\E_1),t\in\NR}.
\end{equation}
With $a=a_t\otimes\unit$, $a_t\in\B(\E_t)$ we obtain 
\begin{eqnarray*}
  \sigma_t(a_t\otimes\unit)x_{1-t}\otimes x_t&=&\tau_t^*(a_t\otimes\unit)\tau_tx_{1-t}\otimes x_t=\tau_t^*(a_t\otimes\unit) x_t\otimes x_{1-t}\\
&=&\tau_t^*a_t x_t\otimes x_{1-t}=x_{1-t}\otimes a_t x_t\\
&=&(\unit\otimes a_t)x_{1-t}\otimes x_t
\end{eqnarray*}
what shows $\sigma_t(a_t\otimes\unit)=\unit\otimes a_t$ and $\sigma_t(\A_{r,s})=\A_{r+t,s+t}$ for $(r,s)\in I_{0,1-t}$. By similar calculations we derive further
\begin{lemma}
  \label{lem:Pstrepr}
  Both families $(\mathrm{P}_{s,t}^u)_{(s,t)\in I_{0,1}}$ and $(\mathrm{P}_{s,t}^{\Us})_{(s,t)\in I_{0,1}}$ consist  of adapted projections $(\mathrm{P}_{s,t})_{(s,t)\in I_{0,1}}\subset\B(\E_1)$ fulfilling  the relations
  \begin{equation}
    \label{eq:Prst}
    \mathrm{P}_{r,t} =\mathrm{P}_{r,s}  \mathrm{P}_{s,t}\dmf{(r,s),(s,t)\in I_{0,1}}
  \end{equation}
and
  \begin{equation}
\label{eq:Pstshift}
    \sigma_r( \mathrm{P}_{s,t})=\left\{
  \begin{array}[c]{cl}
\mathrm{P}_{s+r,t+r}&\text{~if~}t+r\le1\\
\mathrm{P}_{s+r-1,t+r-1}&\text{~if~}s+r\ge1\\
\mathrm{P}_{s+r,1}\mathrm{P}_{0,t+r-1}&\text{~otherwise~}
  \end{array}
\right.\dmf{r\in[0,1],(s,t)\in I_{0,1}}.\proofend
  \end{equation}
\end{lemma}
\begin{remark}
  The shift relation \ref{eq:Pstshift} is  best understood in terms of the unit circle $\NT$,  represented as $[0,1)$ with addition modulo 1. Then the interval $(s,t)$ is  the (open) clockwise arc from $\e^{2\pi\i s}$ to  $\e^{2\pi\i t}$. Thus, if $s>t$, $(s,t)=(s,1)\cup[0,t)$ leading to the convention  $\mathrm{P}_{s,t}=\mathrm{P}_{s,1}\mathrm{P}_{0,t}$ which  simplifies \ref{eq:Pstshift} to 
  \begin{displaymath}
    \sigma_r (\mathrm{P}_{s,t})=\mathrm{P}_{s+r,t+r}\dmf{r,s,t\in \NT,s\ne t}.
  \end{displaymath}
Below, we choose at any occasion that space among $\NT$ or $[0,1]$, which yields  more appropriate notions. 
\end{remark}
\begin{remark}
If $\theta$ is an isomorphism of product systems, it maps units into units. Since the same is true of $\theta^*$, the family $(\mathrm{P}_{s,t}^{\Us})_{(s,t)\in I_{0,1}}$ is an invariant structure.

On the other side,   the following problem is open \cite{OP:Bha01}.  For any pair of  units $u,v\in\Usk(\E)$ with $\norm{u_t}=\norm {v_t}=1$, does  there exist an automorphism of $\E$  mapping $\mathrm{P}_{s,t}^{u}$ into $\mathrm{P}_{s,t}^{v}$ or $u$ into $v$? The importance of this question lead to the notion of \emph{amenable} product system, for which the answer is affirmative \cite[Definition 8.2]{OP:Bha01}. Below we find indications that all  product systems belong to this class, but no proof is available.    Thus, at present, we have to consider the whole collection $\set{(\mathrm{P}_{s,t}^{u})_{(s,t)\in I_{0,1}}:u\in\Usk(\E)}$ to obtain an invariant object. This will be considered  in more detail and in  an appropriate  setting in section \ref{sec:hierarchy}.

All together,  one  should  ask for a structure theory of  projections with \ref{eq:Prst}, which must be necessarily commuting. A good framework to discuss such questions   is representation theory  of  $C^*$-algebras. In the special case of commuting orthogonal projections fulfilling \ref{eq:Prst} we would call the  related objects  (random) bisets, see section \ref{sec:nonsep}.   But, separability of $\E_1$ implies continuity of $\gr t\tau$ (Proposition \ref{prop:tau_tiscontinuous}) and  $(\mathrm{P}_{s,t})_{(s,t)\in I_{0,1}}$ (Proposition \ref{prop:singlepoint}). The latter is essential to get distributions of random \emph{closed} sets, since this requires  some kind of  continuity, see Lemma \ref{lem:biset upper semicontinuous closed}. Due to these facts, the following theorem  can present  a much stronger (von Neumann algebraic) result.
\end{remark}
\begin{theorem}
  \label{th:RACS}
  Let $(\mathrm{P}_{s,t})_{(s,t)\in I_{0,1}}$  be adapted projections in $\B(\E_1)$ which fulfil \ref{eq:Prst} and $\mathrm{P}_{0,1}\ne0$. Then for all normal states $\eta$ on $\B(\E_1)$ there exists a unique probability measure $\mu_\eta$ on $\FG_{[0,1]}$ with 
  \begin{equation}
    \label{eq:mueta}
    \mu_\eta(\set{Z:Z\cap[s_i,t_i]=\emptyset, i=1,\dots,k})=\eta(\mathrm{P}_{s_1,t_1}\cdots \mathrm{P}_{s_k,t_k})\dmf{(s_i,t_i)\in I_{0,1}}
  \end{equation}

If $\eta$ is faithful then  $\mu_{\eta'}\ll\mu_\eta$ for all  normal states $\eta'$ on $\B(\E_1)$. Moreover, the correspondence 
\begin{displaymath}
  \chfc{\set{Z:Z\cap[s,t]=\emptyset}}\mapsto
\mathrm{P}_{s,t}\dmf{(s,t)\in I_{0,1}}
\end{displaymath}
extends to an injective  normal representation $J_{\mathrm{P}}$ of $L^\infty(\mu_\eta)$ on $\E_1$ with  image $\set{\mathrm{P}_{s,t}:{(s,t)\in I_{0,1}}}''$.
\end{theorem}
\begin{remark}
 \ref{eq:Pstshift} is not necessary for the conclusions in the theorem. This will prove useful when we deal with so-called decomposable product systems, see Corollary \ref{cor:decomposable ps is type I}. Similarly, adaptedness is not strictly necessary, it enters the proof only through the next proposition.  

Note that a weaker version of this result was obtained in \cite[Lemma 2.9]{Tsi03} without the non-atomicity results collected in Corollary \ref{cor:open=closed}.   
\end{remark}
For the proof, we need a continuity property of the  family $(\mathrm{P}_{s,t})_{(s,t)\in I_{0,1}}$ (cf.\ Proposition \ref{prop:RACS on Hilbert space}). Define
\begin{equation}
\label{eq:defPcirc}
  \mathrm{P}_s^\circ=\bigvee_{\varepsilon>0}\mathrm{P}_{s-\varepsilon,s+\varepsilon}\dmf{0< s<1}
\end{equation}
The right hand side is  in fact a  $\sigma$-strong limit as $\varepsilon\downarrow0$ since $\varepsilon\mapsto\mathrm{P}_{s-\varepsilon,s+\varepsilon} $ is  increasing by \ref{eq:Prst}. These new projections belong to the von Neumann algebras $\A^\circ_s\subset\B(\E_1)$ defined through 
\begin{equation}
\label{eq:defAcirc}
  \A_s^\circ=\bigwedge_{\varepsilon>0}\A_{s-\varepsilon,s+\varepsilon}\dmf{0< s<1}.
\end{equation}
\begin{proposition}
\label{prop:singlepoint}
For all $s\in(0,1)$ the algebra  $\A_s^\circ$ is $\NC\unit$. More generally, $\lim_{n\to\infty}s_n=s$ and $\lim_{n\to\infty}t_n=t$ force $\lim_{n\to\infty}\A_{s_n,t_n}=\A_{s,t}$ with understanding $\A_{s,s}=\NC\unit$.

Suppose  $(\mathrm{P}_{s,t})_{(s,t)\in I_{0,1}}$ are adapted projections with \ref{eq:Prst}. Then  $\mathrm{P}_{s}^\circ\in\set{0,\unit}$ for all $s\in(0,1)$  and  $\mathrm{P}_{s,t}\ne0$ for some $(s,t)\in I_{0,1}$ implies $\mathrm{P}_{s'}^\circ=\unit$ for all  $s'\in(s,t)$.

If, additionally, $\mathrm{P}_{0,1}\ne0$ then  $s_n\limitsto{}{n\to\infty}s$, $t_n\limitsto{}{n\to\infty}t$  implies 
\begin{displaymath}
  \mathrm{P}_{s_n,t_n}\limitsto{s}{n\to\infty}\left\{
    \begin{array}[c]{c>$l<$}
\mathrm{P}_{s,t}&if $t>s$\\
\unit&if $t=s$
    \end{array}\right.
\end{displaymath}
\end{proposition}
\begin{proof} 
For the proof of the first assertion, fix $a\in\A^\circ_{\set s}$.  For $\delta>0$ we obtain from $\sigma_\delta(\A_{s,t})=\A_{s+\delta,t+\delta}$ that  $a\in\A_{s-\delta,s+\delta}$ implies $\sigma_{\delta}(a)\in\A_{s,s+2\delta}\subseteq\A_{s,1}$. Continuity of $\gr t\tau$ shows $a=\slim_{\delta\downarrow 0}\sigma_\delta(a)\in\A_{s,1}$.  Similarly it follows from considering negative $\delta$ that   $a\in\A_{0,s}$. Thus  $\A^\circ_{\set s}\subseteq\A_{0,s}\cap\A_{s,1}=\NC\unit$, i.e.\ $\A^\circ_{\set s}=\NC\unit$. 

For the proof of the continuity statement observe that by monotony of $(s,t)\mapsto\A_{s,t}$ it is enough to prove the monotone case. \cite[Corollary 3 of Theorem 1-2-2]{C:Mat75} and the Kaplansky density theorem \cite[Theorem 2.4.16]{BR87} show  that this is equivalent to prove 
\begin{displaymath}
  \overline{\bigcup_{\varepsilon>0}\A_{s+\varepsilon,t-\varepsilon}}^w=\A_{s,t}=\bigcap_{\varepsilon>0}\A_{s-\varepsilon,t+\varepsilon}.
\end{displaymath}
Since 
\begin{displaymath}
  a=\slim_{\varepsilon\downarrow0}\sigma_{-\varepsilon}(\sigma_{\varepsilon}(a))\dmf{a\in\A_{s,(s+t)/2},\varepsilon\text{~small}}
\end{displaymath}
we can approximate any element of $\A_{s,(s+t)/2}$ weakly with elements from  $\bigcup_{\varepsilon>0}\A_{s+\varepsilon,t-\varepsilon}$. The same conclusions is valid  for $\A_{(s+t)/2,t}$. Since  $\bigcup_{\varepsilon>0}\A_{s+\varepsilon,t-\varepsilon}$ is an algebra it is weakly dense in $\A_{s,t}$. $\A_{s,t}=\bigcap_{\varepsilon>0}\A_{s-\varepsilon,t+\varepsilon}$ follows from $\overline{\bigcup_{\varepsilon>0}\A_{s+\varepsilon,t-\varepsilon}}^w=\A_{s,t}$  by taking commutants. 

 Adaptedness of $(\mathrm{P}_{s,t})_{(s,t)\in I_{0,1}}$ implies $\mathrm{P}_{s}^\circ\in\A^\circ_{\set s}$ and we derive $\mathrm{P}_{s}^\circ\in\set{0,\unit}$.  By  $\mathrm{P}_{s'}^\circ=0$ for a single $s'\in(0,1)$ we get $\mathrm{P}_{s'-\varepsilon,s'+\varepsilon}=0$ for all $\varepsilon>0$. This shows $\mathrm{P}_{s,t}=0$ for all $0\le s<s'<t\le1$.

Suppose $\mathrm{P}_{0,1}\ne0$ what implies $\mathrm{P}_{s}^\circ=\unit$ for all $s\in(0,1)$. If $s_n\uparrow_{n\to\infty}s$, we get 
\begin{displaymath}
  \unit\ge\mathrm{P}_{s_n,s}\ge\mathrm{P}_{s_n,2s_n-s}\limitsto{s}{n\to\infty}\mathrm{P}_{s}^\circ=\unit
\end{displaymath}
or $\mathrm{P}_{s_n,s}\limitsto{s}{n\to\infty}\unit$. Thus, for $t>s$
\begin{displaymath}
  \slim_{n\to\infty}\mathrm{P}_{s_n,t}=\slim_{n\to\infty}\mathrm{P}_{s_n,s}\mathrm{P}_{s,t}= \mathrm{P}_{s,t}.
\end{displaymath}
Similarly, $t_n\downarrow_{n\to\infty}t$ implies $\mathrm{P}_{s,t_n}\tlimitsto{s}{n\to\infty}\unit$ if $t=s$ and $\mathrm{P}_{s,t_n}\tlimitsto{s}{n\to\infty}\mathrm{P}_{s,t}$ if $t>s$.

From  $s_n\downarrow_{n\to\infty}s$ we obtain from the above relations for $t>s_0$
\begin{displaymath}
  \mathrm{P}_{s,t}=\slim_{n\to\infty}\mathrm{P}_{s,s_n}\mathrm{P}_{s_n,t}=\slim_{n\to\infty}\mathrm{P}_{s,s_n}\slim_{n\to\infty}\mathrm{P}_{s_n,t}=\slim_{n\to\infty}\mathrm{P}_{s_n,t}
\end{displaymath}
The general case follows easily from subsequence arguments.
\end{proof}
\begin{proof}[~of Theorem \ref{th:RACS}]
Observe that \ref{eq:Prst} implies that the projections $(\mathrm{P}_{s,t})_{(s,t)\in I_{0,1}}$ commute with each other. We define for all open sets $G\subseteq[0,1]$   a projection $\mathrm{P}_{G}^\circ $ through the formula
\begin{displaymath}
\mathrm{P}_{G}^\circ =\bigwedge_{
(s,t)\in I_{0,1}, [s,t]\subseteq G
  }\mathrm{P}_{s,t}.
\end{displaymath}
We get  directly that $G_n\uparrow_{n\to\infty} G$ for open sets $G$, $\sequ nG$ leads to $\mathrm{P}_{G_n}^\circ\downarrow\mathrm{P}_{G}^\circ$. Consistently, the meaning of the above formula in the case $G=\emptyset$ is  $\mathrm{P}_\emptyset^\circ=\unit$.

Next we want to prove for all open sets $G,G'\subseteq[0,1]$ that $\mathrm{P}_{G\cup G'}^\circ=\mathrm{P}_{G}^\circ\mathrm{P}_{G'}^\circ$. Since $\mathrm{P}_{G\cup G'}^\circ\le\mathrm{P}_{G}^\circ\mathrm{P}_{G'}^\circ$ is trivially fulfilled,  it suffices to show that any interval $[s,t]\subseteq G\cup G'$ is the finite union of closed intervals which are either contained in $G$ or $G'$. Both $G$ and $G'$ are countable unions of (for each set separately) disjoint open intervals. By compactness, there are sets $G_0\subseteq G$ and $G_0'\subseteq G'$ such that $[s,t]\subseteq G_0\cup G_0'$ and both $G_0$ and $G_0'$ are unions of finitely many open interval. For any open set  $\hat G$ and $\varepsilon>0$ construct the open set 
\begin{displaymath}
  (\hat G)^\varepsilon=\set{t\in \hat G:|t-t'|\le\varepsilon\Rightarrow t'\in\hat G\forall t'\in[0,1]}.
\end{displaymath}
Since $G$ was assumed to be open we derive  $\hat G=\bigcup_{\varepsilon>0}(\hat G)^\varepsilon$. Using again compactness and the fact that $\varepsilon\mapsto(\hat G)^\varepsilon$ is decreasing, we obtain some $\varepsilon>0$ with $[s,t]\subseteq (G_0)^\varepsilon\cup (G_0')^\varepsilon\subseteq\overline{(G_0)^\varepsilon}\cup \overline{(G_0')^\varepsilon}$.  At the end, the definition of $( \hat G)^\varepsilon$ shows that $\overline{(G_0)^\varepsilon}\subseteq G_0\subseteq G$ and  $\overline{(G'_0)^\varepsilon}\subseteq G'_0\subseteq G'$ are unions of finitely many closed intervals and the relation  $\mathrm{P}_{G\cup G'}^\circ=\mathrm{P}_{G}^\circ\mathrm{P}_{G'}^\circ$ is proven.

Now fix a normal state $\eta$ on $\B(\E_1)$. We want to show that the function $T$
\begin{displaymath}
  T(G)=1-\eta(\mathrm{P}_G^\circ)\dmf{G\subseteq[0,1],\text{~open}}
\end{displaymath}
is a Choquet capacity of infinite order \cite{C:Mat75}. This means  that 
\begin{enumerate}
\item\label{cho1} $T(\emptyset)=0$,
\item\label{cho2} $\lim_{n\to\infty}T(G_n)=T(G)$ if $G_n\uparrow G$ and
\item\label{cho3} the functions $\sequ nS$ defined recursively by  $S_0=T$ and
\begin{displaymath}
  S_n(G,G_1,\dots,G_n)=S_{n-1}(G\cup G_1,G_2\dots,G_n)-S_{n-1}(G,G_2,\dots,G_n)
\end{displaymath}
are all positive.
\end{enumerate}
 \ref{cho1} is fulfilled by the normalisation $\eta(\unit)=1$. \ref{cho2} follows from  normality of $\eta$ and the fact that $\mathrm{P}_{G_n}^\circ\downarrow\mathrm{P}_{G}^\circ$ implies $\mathrm{P}_{G}^\circ=\slim_{n\to\infty}\mathrm{P}_{G_n}^\circ$. For the proof of \ref{cho3} observe that $\mathrm{P}_{G\cup G_1}^\circ=\mathrm{P}_{G}^\circ\mathrm{P}_{G_1}^\circ$ implies
\begin{displaymath}
  S_1(G,G_1)=T(G\cup G_1)-T(G)=\eta(\mathrm{P}_{G}^\circ-\mathrm{P}_{G\cup G_1}^\circ)=\eta(\mathrm{P}_{G}^\circ(\unit-\mathrm{P}_{G_1}^\circ))\ge0
\end{displaymath}
and similarly
\begin{displaymath}
  S_n(G,G_1,\dots,G_n)=\eta(\mathrm{P}_{G}^\circ(\unit-\mathrm{P}_{G_1}^\circ)\cdots(\unit-\mathrm{P}_{G_n}^\circ))\ge0.
\end{displaymath}
Now the Choquet theorem \cite[Theorem 2-2-1]{C:Mat75} asserts  that there is a unique Borel probability $\mu_\eta$ on $\FG_{[0,1]}$ such that $\mu_\eta(\set{Z:Z\cap G\ne\emptyset})=T(G)=\eta(\mathrm{P}_{G}^\circ)$. 

Using Proposition \ref{prop:singlepoint} we derive  $\mathrm{P}_{(s,t)}^\circ=\slim_{n\to\infty}\mathrm{P}_{s+1/n,t-1/n}=\mathrm{P}_{s,t}$ and similarly $\mathrm{P}_{[0,t)}^\circ=\mathrm{P}_{0,t}$ and  $\mathrm{P}_{(s,1]}^\circ=\mathrm{P}_{s,1}$. This shows $\mathrm{P}^\circ_{(s_1,t_1)\cup(s_2,t_2)\cup\dots\cup(s_k,t_k)}=\mathrm{P}_{s_1,t_1}\cdots \mathrm{P}_{s_k,t_k}$ what implies
\begin{displaymath}
    \mu_\eta(\set{Z:Z\cap(s_i,t_i)=\emptyset, i=1,\dots,k})=\eta(\mathrm{P}_{s_1,t_1}\cdots \mathrm{P}_{s_k,t_k})\dmf{(s_i,t_i)\in I_{0,1}} 
\end{displaymath}
and the same formula with obvious modifications of the type of intervals if $s_i=0$ or $t_i=1$. Since we are working with closed sets, $Z\cap[s,t]=\emptyset$ iff $Z\cap(s-1/n,t+1/n)=\emptyset$ for some $n\in\NN$. Thus by continuity of (probability) measures, $\sigma$-strong continuity of $\eta$ and $\sigma$-strong continuity of multiplication on $\B(\E_1)$  
\begin{eqnarray*}
\lefteqn{ \mu_\eta(\set{Z:Z\cap[s_i,t_i]=\emptyset, i=1,\dots,k})}&=&\lim_{n\to\infty}\mu_\eta(\set{Z:Z\cap(s_i-1/n,t_i+1/n)=\emptyset, i=1,\dots,k})\\
&=&\lim_{n\to\infty}\eta(\mathrm{P}_{s_1-1/n,t_1+1/n}\cdots \mathrm{P}_{s_k-1/n,t_k+1/n})\\
&=&\eta(\slim_{n\to\infty}\mathrm{P}_{s_1-1/n,t_1+1/n}\cdots\slim_{n\to\infty}\mathrm{P}_{s_k-1/n,t_k+1/n})\\
&=&\eta(\mathrm{P}_{s_1,t_1}\cdots \mathrm{P}_{s_k,t_k}).
\end{eqnarray*}
Further, it is easy to see that the probabilities  $\mu_\eta(\set{Z:Z\cap G\ne\emptyset})$ are already determined by the values $\mu_\eta(\set{Z:Z\cap(s_i,t_i)=\emptyset, i=1,\dots,k})$ and the first part of the theorem is proven.  

\bigskip
For the existence of $J_{\mathrm{P}}$  construct the GNS-representation $(H_{\eta_0},\pi_{\eta_0})$, see section  \ref{GNS-representation}, of the restriction $\eta_0$ of $\eta$ to the abelian von Neumann algebra $\set{\mathrm{P}_{s,t}:{(s,t)\in I_{0,1}}}''\subset\B(\E_1)$. The  map $[\mathrm{P}^\circ_G]_{\eta_0}\mapsto \chfc{\set{Z:Z\cap G=\emptyset}}\in L^2(\mu_\eta)$ extends   to a unitary $\map U{H_{\eta_0}}{L^2(\mu_\eta)}$ since 
\begin{eqnarray*}
  \scpro{[\mathrm{P}_{G'}^\circ]_{\eta_0}}{[\mathrm{P}_{G}^\circ]_{\eta_0}}_{\eta_0}&=&\eta_0(\mathrm{P}_{G'}^\circ\mathrm{P}_{G}^\circ)\\
&=&\eta_0(\mathrm{P}_{G\cup G'}^\circ)\\
&=&\mu_\eta(\set{Z:Z\cap (G\cup G')=\emptyset})\\
&=&\mu_\eta(\set{Z:Z\cap G=\emptyset,Z\cap G'=\emptyset})\\
&=&\int\d\mu_\eta\chfc{\set{Z:Z\cap\smash{G'}=\emptyset}}\chfc{\set{Z:Z\cap G=\emptyset}}\\
&=&\scpro{\chfc{\set{Z:Z\cap\smash{ G'}=\emptyset}}}{\chfc{\set{Z:Z\cap G=\emptyset}}}_{L^2(\mu_\eta)}.
\end{eqnarray*}
Further, the GNS-representation is  unitarily equivalent to the representation of $L^\infty(\mu_\eta)$ on $L^2(\mu_\eta)$, see Example \ref{ex:GNSLinfty}. Since $\eta$ is faithful, $\eta_0$ is too. Therefore, the GNS representation $\pi_{\eta_0}$, $\pi_{\eta_0}(a)[b]_{\eta_0}=[ab]_{\eta_0}$ is faithful (i.e.\ $\pi_{\eta_0}(a)=0$ for positive $a$ implies $a=0$) and its  image is  a von Neumann algebra \cite[Theorem 2.4.24]{BR87}. Thus $\pi_\eta$ has a $\sigma$-weakly continuous faithful inverse and setting  $J_{\mathrm{P}}(\cdot)=\pi_\eta^{-1}(U^*\cdot U)$ the third part is proven. 

Let $\eta'$ be another normal state on $\B(\E_1)$. Then $\eta'\circ J_{\mathrm{P}}$ is a  normal state on $L^\infty(\mu_\eta)$. Therefore, there is a probability measure $\mu'$ on $\FG_{[0,1]}$ with $\mu'\ll\mu_\eta$ for which this normal state is the expectation functional. Now the formula 
\begin{displaymath}
  \mu'(\set{Z:Z\cap G=\emptyset})=\eta'\circ J_{\mathrm{P}}(\chfc{\set{Z:Z\cap G=\emptyset}})=\eta'(\mathrm{P}^\circ_G)=\mu_{\eta'}(\set{Z:Z\cap G=\emptyset})
\end{displaymath}
and the Choquet theorem complete the proof.
\end{proof}
\begin{example}
\label{ex:fullset}
 Suppose $\mathrm{P}_{s,t}=0$ for all $(s,t)\in I_{0,1}$, which fulfils definitely \ref{eq:Prst} but also $\mathrm{P}_{0,1}=0$. We want to show that the conclusions from Theorem \ref{th:RACS} remain valid.   Thus we search for probability measures $\mu_\eta$ with $\mu_\eta(\set{Z:Z\cap[s,t]\ne\emptyset})=1$ for all $(s,t)\in I_{0,1}$. As a consequence, almost all  $Z$ should have points in any interval with rational endpoints. We conclude that $Z$ is dense almost surely and closedness of $Z$ implies $Z=[0,1]$ $\mu_\eta$-a.s.\ or $\mu_\eta=\delta_{[0,1]}$. It is easy to see that this choice fulfils \ref{eq:mueta}.  
\end{example}
\begin{example}
\label{ex:emptyset}
  $\mathrm{P}_{s,t}=\unit$ for all $(s,t)\in I_{0,1}$ fulfils  \ref{eq:Prst} too. Now $\mu_\eta(\set{Z:Z\cap[0,1]=\emptyset})=\eta(\mathrm{P}_{0,1})=\eta(\unit)=1$ shows $\mu_\eta=\delta_{\emptyset}$. 
\end{example}
For later, we add some properties already collected in the proof of Theorem \ref{th:RACS}.
\begin{corollary}
\label{cor:open=closed}
  The same conditions imply the relations  
  \begin{displaymath}
    \mu_\eta(\set{Z:t\in Z})=0\dmf{t\in[0,1]},
  \end{displaymath}
and 
\begin{displaymath}
  \mu_\eta(\set{Z:Z\cap\smash{\bigcup_{i=1}^k}(s_i,t_i)=\emptyset})=\mu_\eta(\set{Z:Z\cap\smash{\bigcup_{i=1}^k}[s_i,t_i]=\emptyset})\dmf{(s_i,t_i)\in I_{0,1}}\end{displaymath}
as well as
\begin{displaymath}
  J_{\mathrm{P}}(\chfc{\set{Z:Z\cap(s,t)=\emptyset}})=  J_{\mathrm{P}}(\chfc{\set{Z:Z\cap[s,t]=\emptyset}})\dmf{(s,t)\in I_{0,1}}.\proofend
\end{displaymath}
\end{corollary}
\subsection{Measure Types as Invariants}
\label{sec:invariant}
Since isomorphisms of product systems send units into units  the system $(\mathrm{P}_{s,t}^{\Us})_{(s,t)\in I_{0,1}}$ is left invariant by them. For that reason we study now the special case of random sets associated with these projections. Although $\mathrm{P}_{s,t}^{\Us}=0$ for type $\mathrm{III}$ product systems such that Theorem \ref{th:RACS} could not be used, Example \ref{ex:fullset} showed that the results extend to this case.  We denote the corresponding  probability  measures  on $\FG_{[0,1]}$ by  $\mu_\eta^\Us$.

Besides their pure existence, the tensor product structure creates additional properties of the measures $\mu_\eta^\Us$. So let  $\mu$ be a probability measure on $\FG_{[0,1]}$.  Define  $\mu_{s,t}$ to be the image of $\mu$ under the map $Z\mapsto Z\cap[s,t]$, and, for $t\in[0,1]$, $\mu+t$ to be the image under the map $Z\mapsto Z+t=(\set{z+t:z\in Z}\cup\set{z+t-1:z\in Z})\cap[0,1]$. These definitions carry over to measure types.
\begin{theorem}
  \label{th:ps2randomset}
  For every product system $\E$ the measures $\mu^\Us_\eta$  belong for all  faithful normal states $\eta$  on $\B(\E_1)$  to the same   measure type. Denoting this measure type by  $\M^{\E,\Us}$ the relations 
\begin{equation}
  \label{eq:MUsfactorizing}
  \M^{\E,\Us}_{r,t}=\M^{\E,\Us}_{r,s}\ast\M^{\E,\Us}_{s,t}\dmf{(r,s),(s,t)\in I_{0,1}}
\end{equation}
and 
\begin{equation}
\label{eq:MUsstationary}
  \M^{\E,\Us}+t=\M^{\E,\Us}\dmf{t\in[0,1]}
\end{equation}
hold.

Moreover, if two product systems  $\E$ and $\E'$ are isomorphic then  $\M^{\E,\Us}=\M^{\E',\Us}$.
\end{theorem}
\begin{remark}
  Remember our conventions about measure types: \ref{eq:MUsfactorizing}  and \ref{eq:MUsstationary} mean for one  and thus all  $\mu\in\M$
  \begin{equation}
  \label{eq:quasifactorizingmeasure}
  \mu_{r,t}\sim\mu_{r,s}\ast\mu_{s,t}\dmf{(r,s),(s,t)\in I_{0,1}}
\end{equation}
and 
\begin{equation}
\label{eq:quasistationarymeasure}
  \mu+t\sim\mu\dmf{t\in\NR}.
\end{equation}
respectively. Such measures we call \emph{quasifactorizing} and \emph{quasistationary} respectively. 
\end{remark}
\begin{remark}
Similar results could be obtained  for $E_0$-semigroups $\sg t\alpha$ on some $\B(\H)$. Then we have to fix    projections $(\mathrm{P}_{s,t})_{(s,t)\in I_{0,\infty}}\subset\B(\H)$ with    \ref{eq:Prst} for all $(r,s),(s,t)\in I_{0,\infty}$ and 
 \begin{displaymath}
   \alpha_r(\mathrm{P}_{s,t})=\mathrm{P}_{s+r,t+r}\dmf{r\in\NRp,(s,t)\in I_{0,\infty}}.
 \end{displaymath}
We do not follow this line here since the analogue of the above result provides us a  measure type  on $\FG_\NRp$ which is \emph{not an invariant} of $\sg t\alpha$ with respect to cocycle conjugacy, the counterpart of isomorphy of product systems for $E_0$-semigroups. It is  only invariant under  conjugacy. 
\end{remark}
\begin{proof}
  From Theorem  \ref{th:RACS} we know for two faithful normal states $\eta,\eta'$ that both $\mu_{\eta'}\ll\mu_\eta$ and $\mu_{\eta}\ll\mu_{\eta'}$ are true, i.e.\ 
 $\mu_{\eta'}\sim\mu_\eta$. This shows independence  of the measure type $\M^{\E,\Us}$ from the choice of $\eta$. 

To prove \ref{eq:quasifactorizingmeasure} we may restrict to the choice  $r=0$ and $t=1$ and $\mu=\mu_\eta$ for some faithful normal state $\eta$ on $\B(\E_1)$. Observe that the representation $\E_1\cong \E_s\otimes\E_{1-s}$ provides us together  with $\eta$ with the restrictions $\eta_{0,s}$ of $\eta$ to $\A_{0,s}$ and $\eta_{s,1}$ to $\A_{s,1}$.  But the normal state  $\eta'=\eta_{0,s}\otimes\eta_{s,1}$, defined via
\begin{displaymath}
  \eta_{0,s}\otimes\eta_{s,1}(b_s\otimes b_{1-s})=\eta(b_s\otimes\unit_{\E_{1-s}})\eta(\unit_{\E_s}\otimes  b_{1-s})\dmf{b_s\in\B(\E_s), b_{1-s}\in\B(\E_{1-s})}
\end{displaymath}
is again faithful on $\B(\E_1)$. Fix intervals  $[s_i,t_i]$, $i=1,\dots,k$, which are contained in $[0,s]$ for $i\le l\le k$ and which are contained in $[s,1]$ for $i>l$. We find 
\begin{eqnarray*}
\lefteqn{\mu_{\eta'}(\set{Z:Z\cap[s_i,t_i]=\emptyset, i=1,\dots,k})}
&=&\eta'(\mathrm{P}_{s_1,t_1}\cdots \mathrm{P}_{s_k,t_k})\\
&=&\eta_{0,s}(\mathrm{P}_{s_1,t_1}\cdots \mathrm{P}_{s_l,t_l})\eta_{s,1}(\mathrm{P}_{s_{l+1},t_{l+1}}\cdots \mathrm{P}_{s_k,t_k})\\
&=&  (\mu_{\eta})_{0,s}(\set{Z:Z\cap[s_i,t_i]=\emptyset, i=1,\dots,k})  (\mu_{\eta})_{s,1}(\set{Z:Z\cap[s_i,t_i]=\emptyset, i=1,\dots,k})\\
&=&  (\mu_{\eta'})_{0,s}(\set{Z:Z\cap[s_i,t_i]=\emptyset, i=1,\dots,k})  (\mu_{\eta'})_{s,1}(\set{Z:Z\cap[s_i,t_i]=\emptyset, i=1,\dots,k})\\
&=&(\mu_{\eta'})_{0,s}\ast(\mu_{\eta'})_{s,1}(\set{Z:Z\cap[s_i,t_i]=\emptyset, i=1,\dots,k}).
\end{eqnarray*}
 Theorem \ref{th:RACS}  proves \ref{eq:quasifactorizingmeasure}. 

To prove equation \ref{eq:quasistationarymeasure}, define for $t\in[0,1]$ and  a faithful normal state  $\eta$ on $\B(\E_1)$  the state  $\eta^t=\eta\circ\sigma_t=\eta(\tau_t\cdot\tau_t^*)$ which is again faithful. Then we see from \ref{eq:Pstshift} for $(s_1,t_1),\dots,(s_l,t_l)\subset[0,1-t]$ and $(s_{l+1},t_{l+1}),\dots,(s_k,t_k) \subset[1-t,1]$
\begin{eqnarray*}
\lefteqn{\mu_{\eta^t}(\set{Z:Z\cap\smash{\bigcup_{i=1}^k}(s_i,t_i)=\emptyset})}
&=&\eta^t(\mathrm{P}_{s_1,t_1}\cdots \mathrm{P}_{s_k,t_k})\\
&=&\eta(\sigma_t(\mathrm{P}_{s_1,t_1}\cdots \mathrm{P}_{s_l,t_l}\mathrm{P}_{s_{l+1},t_{l+1}}\cdots \mathrm{P}_{s_k,t_k}))\\
&=&\eta(\sigma_t(\mathrm{P}_{s_1,t_1})\cdots\sigma_t( \mathrm{P}_{s_l,t_l})\sigma_t(\mathrm{P}_{s_{l+1},t_{l+1}})\cdots \sigma_t(\mathrm{P}_{s_k,t_k}))\\
&=&\eta(\mathrm{P}_{s_1+t,t_1+t}\cdots \mathrm{P}_{s_l+t,t_l+t}\mathrm{P}_{s_{l+1}+t-1,t_{l+1}+t-1}\cdots \mathrm{P}_{s_k+t-1,t_k+t-1})\\
&=& \mu_\eta(\set{Z:Z\cap\smash{\bigcup_{i=1}^l}(s_i+t,t_i+t)=\emptyset,Z\cap\smash{\bigcup_{i=l+1}^k}(s_i+t-1,t_i+t-1)=\emptyset}) \\
&=& \mu_\eta (\set{Z:(Z-t)\cap(s_i,t_i)=\emptyset, i=1,\dots,k}).
\end{eqnarray*}
In these calculations we can use open instead of  closed intervals by Corollary \ref{cor:open=closed}. So we can conclude from  Theorem \ref{th:RACS} that $\mu_\eta-t=\mu_{\eta_t}\sim\mu_\eta$ for all $t\in[0,1]$.

If two product systems $\E$ and $\E'$ are isomorphic under $\sg t\theta$, $\theta$ and $\theta^*$  map units into units. Thus we obtain $\theta_t\E^\Us_t=(\E')^\Us_t$ and  the families $(\mathrm{P}_{s,t}^{\E,\Us})_{(s,t)\in I_{0,1}}$ and $(\mathrm{P}_{s,t}^{\E',\Us})_{(s,t)\in I_{0,1}}$ are unitarily equivalent:
\begin{displaymath}
  \theta_1^*\mathrm{P}_{s,t}^{\E,\Us}\theta_1=\mathrm{P}_{s,t}^{\E',\Us}\dmf{(s,t)\in I_{0,1}}.
\end{displaymath}
Setting  $\eta'=\eta(\theta_1^*\cdot\theta_1)$ on $\B(\E'_1)$, it is straight forward to see $\mu^{\E,\Us}_\eta=\mu^{\E',\Us}_{\eta'}$ and the proof is complete.  
\end{proof}
\begin{remark}
  We prove in Corollary  \ref{cor:muetafull} below that $\set{\mu_\eta:\eta\text{~faithful~}}$ is already a measure type, i.e.\ it is equal to  $\M^{\E,\Us}$.  
\end{remark}
We  present now  two simple examples.
\begin{example}
\label{ex:measure type III}
  If $\E$ is of type $\mathrm{III}$, i.e.\ it has no unit, we know $\mathrm{P}^\Us_{s,t}=0$ for all $(s,t)\in I_{0,1}$. Thus Example \ref{ex:fullset} shows  $\mu^\Us_\eta=\delta_{[0,1]}$ and  $\M^{\E,\Us}=\set{\delta_{\set{[0,1]}}}$.  
\end{example}
\begin{example}
Consider $\E=\Gammai(\K)$ which  is of type $\mathrm{I}$.  Thus  $\mathrm{P}^{\Gammai(\K),\Us}_{0,1}=\unit$ and  Example \ref{ex:emptyset} yields  $\M^{\E,\Us}=\set{\delta_{\emptyset}}$.
\end{example}
\subsection{Measure Types Related to Units}
\label{sec:measure types units}
It is clear, that for the  first part of Theorem \ref{th:ps2randomset}   the family $(\mathrm{P}_{s,t})_{(s,t)\in I_{0,1}}$ need not be associated with $\E^\Us$. So we obtain  in a similar, but   due to $\mathrm{P}^u_{0,1}=\Pr{u_1}\ne 0$ more simple, fashion
\begin{proposition}
  \label{prop:psunit2randomset}
  For every product system $\E$ and every unit $u\in\Usk(\E)$ the measures $\mu^u_\eta$ defined via Theorem \ref{th:RACS} from the system $(\mathrm{P}_{s,t}^{u})_{(s,t)\in I_{0,1}}$ with varying  faithful normal state $\eta$  on $\B(\E_1)$ belong to a unique measure type $\M^{\E,u}$    fulfilling \ref{eq:MUsfactorizing} and \ref{eq:MUsstationary}. \proofend
\end{proposition}
\begin{remark}
\label{rem:ME,u as invariant}
We would like to have, like for the product subsystem  $\E^\Us$, that $\M^{\E,u}$ is an invariant of the product system $\E$.  So the main problem is whether $\M^{\E,u}$ is a well-defined invariant, i.e.\ whether for two normalized units $u,v$ the measure types $\M^{\E,u}$ and $\M^{\E,v}$ are equal. Copying the relevant part of  the proof of Theorem \ref{th:ps2randomset} gives us only  $\M^{\E,u}=\M^{\E',\theta_1 u}$ if $\theta$ is an isomorphism between $\E$ and $\E'$. Thus  $\M^{\E,u}=\M^{\E,v}$ would result for all $u,v\in\Usk(\E)$ if there were  an automorphism of $\E$ such that $v_t=\theta_t u_t$. Whether the latter is true is an open problem until now \cite{OP:Bha01}, although we have strong indications (see Remark \ref{rem:zerocapacitiesofM^u=M^v} and Proposition \ref{prop:unitandunitalprojections}) that the answer is affirmative.
\end{remark}
To derive the connection between the measure types $\M^{\E,u}$ and $\M^{\E,\Us}$ we need to collect more properties of   units. It is easy to prove (cf.\ \cite[Theorem 4.1]{Arv89}) that for two units $u,v\in\Usk(\E)$  there is some $\gamma(u,v)\in\NC$  such that 
\begin{equation}
\label{eq:defcovariance}
  \scpro{u_t}{v_t}_{\E_t}=\e^{-\gamma(u,v)t}\dmf{t\in\NRp}.
\end{equation}
Thereby, $\map\gamma{\Usk(\E)\times\Usk(\E)}\NC$ is a conditionally positive kernel and called \emph{covariance function} \cite{Arv89}.

\emph{Poisson processes} $\Pi_\nu$  are distributions of (finite) random closed sets  defined by
\begin{equation}
\label{eq:defPoisson}
  \Pi_\nu=\e^{-\nu([0,1])}\left(\delta_\emptyset+\sum_{n\in\NN}\frac1{n!}\int\nu^n(\d t_1,\dots,\d t_n)\delta_{\set{t_1,\dots,t_n}}\right).
\end{equation}
where $\nu$ is  a finite measure on $[0,1]$ with $\nu(\set t)=0\forall t$. Observe that $  \Pi_\nu(\set{Z:Z\cap G=\emptyset})=\e^{-\nu(G)}$ for all (open) sets $G\subset[0,1]$ which is  characteristic for $\Pi_\nu$ by the Choquet theorem \cite[Theorem 2-2-1]{C:Mat75}. Observe that this theorem shows too that $\Pi_\nu$ is factorizing and stationary if $\nu=\ell$. It is not hard to prove from \ref{eq:defPoisson} that it is quasifactorizing for all $\nu\sim \ell$.
\begin{lemma}
  \label{lem:unitalmeasurePoisson}
Let $u,v$ be two units with $\norm{v_1}=1$. Then  the measure $\mu^u_\eta$ associated with the projection  family  $(\mathrm{P}_{s,t}^{u})_{(s,t)\in I_{0,1}}$ and the pure normal state  $\eta(a)=\Tr\Pr{v_1}a$ $\forall a\in\B(\E_1)$ by Theorem \ref{th:RACS} is the Poisson process  $\Pi_{2\Re\gamma(u,v)\ell}$.
\end{lemma}
\begin{proof}
If we  normalize $u$ and $v$, we obtain from \ref{eq:defcovariance} for disjoint intervals ${(s_i,t_i)\in I_{0,1}}$ 
\begin{eqnarray*}
  \mu_\eta(\set{Z:Z\cap[s_i,t_i]=\emptyset, i=1,\dots,k})&=& \eta(\mathrm{P}^u_{s_1,t_1}\cdots \mathrm{P}^u_{s_k,t_k})\\
&=&\prod_{i=1}^k\e^{2\Re\gamma(u,v)(t_i-s_i)}\\
&=&\Pi_{2\Re\gamma(u,v)\ell}(\set{Z:Z\cap[s_i,t_i]=\emptyset, i=1,\dots,k}).
\end{eqnarray*}
Theorem \ref{th:RACS} completes the proof.
\end{proof}
For mutually commuting projection families, it is straight forward to prove the following analogue of Theorem \ref{th:RACS}.
\begin{proposition}
  \label{prop:RACSmult}
  If $N\in\NN^*\cup\set\infty$ and $(\mathrm{P}^n_{s,t})_{(s,t)\in I_{0,1}}$, $n=1,\dots,N$   are commuting families of  adapted projections with \ref{eq:Prst}, \ref{eq:Pstshift}, there exists for all normal states $\eta$ on $\B(\E_1)$ a unique probability measure $\mu^{1,\dots,N}_\eta$ on $\FG_{[0,1]\times\set{1,\dots,N}}$ with 
  \begin{equation}
    \label{eq:muetamult}
    \mu_\eta^{1,\dots,N}(\set{Z:Z\cap\smash{\bigcup_{i=1}^k}([s_i,t_i]\times \set{n_i})=\emptyset})=\eta(\mathrm{P}^{n_1}_{s_1,t_1}\cdots \mathrm{P}^{n_k}_{s_k,t_k})\dmf{(s_i,t_i)\in I_{0,1},1\le n_i\le N}
  \end{equation}

If $\eta$ is faithful $\mu^{1,\dots,N}_{\eta'}\ll\mu^{1,\dots,N}_\eta$ for all  normal states $\eta'$. Moreover, the correspondence 
\begin{displaymath}
  \chfc{\set{Z:Z\cap([s,t]\times \set{n})=\emptyset}}\mapsto\mathrm{P}^{n}_{s,t}
\end{displaymath}
extends to a normal representation of $L^\infty(\mu^{1,\dots,N}_\eta)$ on $\B(\E_1)$ with  image $\set{\mathrm{P}^n_{s,t}:(s,t)\in I_{0,1},n=1,\dots,N}''$.\proofend
\end{proposition}
For any $Z\in\FG_K$, denote $\check{Z}$ the set of its limit points:
\begin{displaymath}
  \check{Z}=\set{t\in Z:t\in\overline{Z\setminus\set t}}.
\end{displaymath}
With these preparations we can  derive an interesting relation between the measure types $\M^{\E,u}$ and $\M^{\E,\Us}$, we prove for simplicity reasons in section \ref{sec:direct integrals ps}, on page \pageref{page:proof unitunitalprojections}. 
\begin{proposition}
  \label{prop:unitandunitalprojections}
  Let $u$ be a unit of a product system $\E$ and define the  measure type  $\M^{u,\Us}$ of random sets on $\FG_{[0,1]\times\set{1,2}}$ by
  \begin{displaymath}
   \M^{u,\Us}=\set{\mu:\mu\sim\mu_\eta^{1,2}\text{~for all faithful normal states $\eta$  on $\B(\E_1)$}}, 
  \end{displaymath}
where the measures  $\mu_\eta^{1,2}$ are  associated with  the families  $(\mathrm{P}^1_{s,t})=(\mathrm{P}^u_{s,t})$,  $(\mathrm{P}^2_{s,t})=(\mathrm{P}^\Us_{s,t})$ by the preceding proposition.

 Then $\mathrm{P}^\Us_{s,t}=J_{\mathrm{P}^u}(\chfc{\set{Z:\#(Z\cap[s,t])<\infty}})$ and therefore  $\M^{u,\Us}$ is concentrated on the set $\set{(Z\times\set1)\cup(\check{Z}\times\set2):Z\in\FG_{[0,1]}}$.  Consequently, $\M^{\E,\Us}=\M^{\E,u}\circ l^{-1}$, with  $l$  mapping $Z$ into $\check{Z}$.
\end{proposition}
\begin{remark}
Surely, both the set of measure types on $\FG_{[0,1]}$ as well as its subset of stationary factorizing ones are monoids under $\ast$.  If we equip the (isomorphy classes of) product systems  with the (commutative) operation of taking the tensor product, the map $\E\mapsto\M^{\E,\Us}$ becomes an homomorphism by this result. This parallels the similar result from \cite{OP:Arv89} on additivity of the numerical index.
\end{remark}
\begin{remark}
  This result is essential since the measure type $\M^{\E,u}$ is much easier to compute than $\M^{\E,\Us}$, what  is used several times below. 
\end{remark}

\subsection{Tensor Products (I)}
\label{sec:tensor products(I)}
One main operation on  product systems is the tensor product. I.e., if  $\E,\E'$ are product systems then  $\E\otimes\E'=(\sgp t{\E_t\otimes\E'_t},\sgp{s,t}{V_{s,t}\otimes V'_{s,t}}$ is again a product system.  It is interesting to ask, what implications this operation has in terms of the measure types  associated with $\E$ and $\E'$.

 \textsc{Arveson} proved in \cite[Corollary 3.9]{OP:Arv89}, see Corollary \ref{cor:unitstensoraretensorofunits} below, that all units of $\E\otimes\E'$ have the form $u\otimes u'$. Then 
  \begin{displaymath}
    \mathrm{P}^{\E\otimes\E',u\otimes u'}_{s,t}= \mathrm{P}^{\E,u}_{s,t}\otimes\mathrm{P}^{\E',u'}_{s,t}\dmf{(s,t)\in I_{0,1}}.
  \end{displaymath}
Further, this result implies also
  \begin{displaymath}
    \mathrm{P}^{\E\otimes\E',\Us}_{s,t}= \mathrm{P}^{\E,\Us}_{s,t}\otimes\mathrm{P}^{\E',\Us}_{s,t}\dmf{(s,t)\in I_{0,1}}.
  \end{displaymath}
\begin{proposition}
  \label{prop:tensorproduct}
Let $\E$ and $\E'$ be product systems. Then $\M^{\E\otimes\E',\Us}=\M^{\E,\Us}\ast\M^{\E',\Us}$. If $u$ and $u'$ are units of $\E$ and $\E'$ respectively then $\M^{\E\otimes\E',u\otimes u'}=\M^{\E,u}\ast\M^{\E',u'}$. 
\end{proposition}
\begin{proof}
For proving the first relation we use Theorem \ref{th:RACS} for a faithful state $\eta\otimes\eta'$ on $\B(\E_1\otimes\E_1')\cong\B(\E_1)\otimes\B(\E_1')$ with $\eta$ being faithful on $\B(\E_1)$ and $\eta'$ being faithful on $\B(\E_1')$.  Then  we obtain for all $(s_i,t_i)\in I_{0,1}$, $i=1,\dots,k$
\begin{eqnarray*}
 \lefteqn{ \mu_{\eta\otimes\eta'}^\Usk(\set{Z:Z\cap[s_i,t_i]=\emptyset, i=1,\dots,k})}
&=&\eta\otimes\eta'(\mathrm{P}^{\E\otimes\E',\Us}_{s_1,t_1}\cdots \mathrm{P}^{\E\otimes\E',\Us}_{s_k,t_k})\\
&=&\eta\otimes\eta'(\mathrm{P}^{\E,\Us}_{s_1,t_1}\otimes\mathrm{P}^{\E',\Us}_{s_1,t_1}\cdots\mathrm{P}^{\E,\Us}_{s_k,t_k}\otimes\mathrm{P}^{\E',\Us}_{s_k,t_k} )\\
&=&\eta(\mathrm{P}^{\E,\Us}_{s_1,t_1}\cdots\mathrm{P}^{\E,\Us}_{s_k,t_k} )\eta'(\mathrm{P}^{\E',\Us}_{s_1,t_1}\cdots\mathrm{P}^{\E',\Us}_{s_k,t_k} )\\
&=&\eta(\mathrm{P}^{\E,\Us}_{s_1,t_1}\cdots \mathrm{P}^{\E,\Us}_{s_k,t_k})\eta'(\mathrm{P}^{\E',\Us}_{s_1,t_1}\cdots \mathrm{P}^{\E',\Us}_{s_k,t_k})\\
&=& \mu_\eta^{\E,\Us}(\set{Z:Z\cap[s_i,t_i]=\emptyset, i=1,\dots,k}) \mu_{\eta'}^{\E',\Us}(\set{Z':Z'\cap[s_i,t_i]=\emptyset, i=1,\dots,k})\\
&=& \mu_\eta^{\E,\Us}\otimes \mu_{\eta'}^{\E',\Us}(\set{(Z,Z'):(Z\cup Z')\cap[s_i,t_i]=\emptyset, i=1,\dots,k})\\
&=& \mu_\eta^{\E,\Us}\ast \mu_{\eta'}^{\E',\Us}(\set{Z:Z\cap[s_i,t_i]=\emptyset, i=1,\dots,k})
\end{eqnarray*}
what shows the assertion on $\M^{\E\otimes\E',\Us}$ by Theorem \ref{th:RACS}. The proof for $\M^{\E\otimes\E',u\otimes u'}$ is the same.
\end{proof}
A further analysis of these measure types   follows in section \ref{sec:tensor products(II)}.

\section{From Random Sets to Product Systems}
\label{sec:rs2ps}

In the last section we established that any  product system corresponds to at least one measure type of random sets. To get a deeper understanding of this relation we study now the reverse question: How can we  derive product systems from measure types.  As a side-effect, we show that both $\E\mapsto \M^{\E,\Us}$ and $\E\mapsto \M^{\E,u}$  are onto. Surprisingly, we can use some already established facts about product systems to derive properties of the measure types of interest.

\subsection{General Theory}
\label{sec:theoryrs2ps}

Let $L$ be any locally compact Hausdorff space.  On $\FG_{[0,1]\times L}$ we have the natural operations $Z\mapsto Z_{s,t}=Z\cap([s,t]\times L)$ and $Z\mapsto Z+t=(\set{(z+t,l):(z,l)\in Z}\cup\set{(z+t-1,l):(z,l)\in Z})\cap[0,1]\times L$. The convolution on the probability measures of $\FG_{[0,1]\times L}$ associated with $\cup$ is denoted $\ast$. In generalisation of the results in Theorem  \ref{th:RACS} we make the following   
\begin{definition}
  \label{def:statfactmeastype}
  We call a measure type  $\M$ on $\FG_{[0,1]\times L}$ \emph{stationary factorizing measure type} if  
\begin{equation}
  \label{eq:factorizingmeasuretype}
  \M_{r,t}=\M_{r,s}\ast\M_{s,t}\dmf{(r,s),(s,t)\in I_{0,1}}
\end{equation}
and 
\begin{equation}
\label{eq:stationarymeasuretype}
  \M+t=\M\dmf{t\in\NR}.
\end{equation}  
\end{definition}
\begin{remark}
  We remind the reader that \ref{eq:factorizingmeasuretype} and \ref{eq:stationarymeasuretype} mean \ref{eq:quasifactorizingmeasure} and \ref{eq:quasistationarymeasure} for all $\mu\in\M$.  
\end{remark}
\begin{proposition}
  \label{prop:statfactmeastype=ps}
  Every stationary factorizing measure type $\M$  on $\FG_{[0,1]\times L}$ determines  a product system $\E^\M$ associated with $(\sgi t{\E^\M},(V_{s,t})_{s,t\in[0,1],s+t\le1})$ by Proposition \ref{prop:psbyE_1}, where $\E^\M_t=L^2(\M_{0,t})$, $t\in[0,1]$, and
  \begin{displaymath}
    (V_{s,t}\psi_s\otimes\psi_t)_{\mu_{0,s}\otimes\mu'_{s,s+t}}(Z)=(\psi_s)_{\mu_{0,s}}(Z_{0,t})(\psi_t)_{\mu_{s,s+t}-t}(Z_{s,s+t}-t)\dmf{s,t\in [0,1],s+t\le1}
  \end{displaymath}
A corresponding measurable structure is given by sections $\sgi t\psi$ for which there is a measurable function $\map{\hat\psi}{[0,1]\times\FG_{[0,1]}}\NC$ such that for some $\mu\in\M$
\begin{displaymath}
  (\psi_t)_{\mu_{0,t}}(Z)=\hat\psi(t,Z_{0,t})\dmf{t\in[0,1],\mu_{0,t}-\text{a.a.~}Z\in\FG_{[0,1]}}.
\end{displaymath}
\end{proposition}
The proof of this result requires some results from section \ref{sec:Polish Product System}. Therefore, it is  postponed to section \ref{sec:psfromRACS}, page \pageref{page:proof prop:statfactmeastype=ps}. Above we  introduced a  general $L$ since we  wanted to have an analogue of   Proposition \ref{prop:RACSmult} too. In the following, we concentrate on the case where $L$ is a singleton, i.e.\ $Z\in\FG_{[0,1]}$. First  we summarize some properties of stationary factorizing measure types on $\FG_{[0,1]}$, being partly derived in the proof of the previous proposition.
\begin{corollary}
  \label{cor:propmu}
Suppose $\M$ is a stationary factorizing measure type on $[0,1]$ different from $\set{\delta_{[0,1]}}$. Then for all $\mu\in\M$
  \begin{enumerate}
  \item\label{mu1} $\mu(\set{Z:t\in Z})=0$ for all $t\in[0,1]$. More generally, $\mu(\set{Z:Z_0\cap Z\ne\emptyset})=0$ for all countable   $Z_0\subset[0,1]$.
  \item\label{mu2} $\mu(\set\emptyset)>0$.
  \item\label{mu3} $\ell(Z)=0$ for $\mu$--a.a.\ $Z\in\FG_{[0,1]}$. Thus for $\mu$-a.a.\ $Z$ the complement $\cmpl{Z}$ is a dense open set with full Lebesgue measure.
  \end{enumerate}
\end{corollary}
\begin{proof}
The first part of  \ref{mu1} is equivalently stated and proved in  Lemma \ref{lem:threedifferentpossibilities}. The second part of \ref{mu1}  follows easily.   

For $\varepsilon>0$ consider  $f_\varepsilon(Z)=\chfc{\set{Z':Z'\cap[t-\varepsilon,t+\varepsilon]=\emptyset}}(Z)$. Now $\lim_{\varepsilon\downarrow0}f_\varepsilon(Z)=1-\chfc Z(t)$ together with \ref{mu1} shows $\chfc Z(t)=1$ and  $\lim_{\varepsilon\downarrow0}\chfc{\set{Z:Z\cap[t-\varepsilon,t+\varepsilon]=\emptyset}}=1$ $\mu$-a.s., i.e.\ there is some $\varepsilon>0$ such that $\mu(\set{Z:Z\cap[t-\varepsilon,t+\varepsilon]=\emptyset})=\mu_{t-\varepsilon,t+\varepsilon}(\set\emptyset)>0$. Covering $[0,1]$ by finitely many intervals of length less than $2\varepsilon$ and using \ref{eq:quasifactorizingmeasure} yields \ref{mu2}. 

To prove  \ref{mu3} observe for $\mu$-a.a.\ $Z\in\FG_{[0,1]}$
\begin{displaymath}
  0=\int^1_0 \ell(\d t)\chfc Z(t)=\ell(Z)
\end{displaymath}
what proves the third  assertion.
\end{proof}
\begin{remark}
It is not complicated to get   similar results  for general  $L$. Then \ref{mu1} reads $\mu(\set{Z:(\set t\times L)\cap Z\ne\emptyset})=0$ after excluding cases similar to $\mu=\delta_{[0,1]}$. 

\ref{mu2} has to be converted into $\mu(\set{Z:Z\cap([0,1]\times L')=\emptyset})>0$ for all compact $L'\subseteq L$. Note that  this may be false for $L'=L$, see the next example. This fact is the reason we explore the structure of the product system $\E^\M$ in the case of singleton $L$ only. 

Using the canonical projection $\map\pi{{[0,1]}\times L}{[0,1]}$,  \ref{mu3} shall mean $\ell(\pi(Z))=0$.
\end{remark}
\begin{example}
\label{ex:infinite Poisson}
   Consider the Poisson process $\Pi_{\ell\otimes\#}$ on $\FG_{[0,1]\times\NN}$ with intensity measure $\ell\otimes\#$. Since $\FG_{[0,1]\times\NN}\cong(\FG_{[0,1]})^\NN$ under $Z\mapsto (\pi(Z\cap [0,1]\times\set n))_{n\in\NN}$ we can define $\Pi_{\ell\otimes\#}=\Pi_\ell^{\otimes \NN}$. This product measure  is stationary and factorizing since $\Pi_ell$ is so. On the other side, by the law of large numbers, for any Borel set $Y\subset[0,1]$ with $\ell(Y)>0$  we obtain $\mu(\set{Z:Z\cap(Y\times \NN)=\emptyset})=0$ since $\Pi_\ell(\set{Z:Z\cap Y=\emptyset})=\mathrm{e}^{-\ell(Y)}<1$.
\end{example}
From now on, we will convert notations from the interval  $[0,1]$ to the circle $\NT$ viewed as $[0,1)$ with the topology derived from the map $t\mapsto\mathrm{e}^{2\pi\mathrm{i}t}\in\NC$. This conversion is possible as the above corollary, \ref{mu1} shows that for any quasistationary quasifactorizing measure $\mu$ its  image under $Z\mapsto Z\cap[0,1)\subseteq\NT$ contains all information about $Z$ and yields almost surely a closed set.  The canonical structure of  $\NT$ as a topological group (with respect to addition modulo 1) yields for $\mu$-a.a.\ $Z$ that $Z+t=\set{z+t:z\in Z}$ which is  the canonical notation. This change is done  for convenience only,  all our results stated for $[0,1]$ have a counterpart on $\NT$ and vice versa. 
\begin{corollary}
  \label{cor:psrstypeII}
Let $\M$ be a stationary factorizing measure type on $\FG_\NT$ different from $\set{\delta_\NT}$.
Then the  product system   $\E=\E^\M$ has at least one unit, corresponding to 
  \begin{displaymath}
    (u_t)_\mu=\mu(\set\emptyset)^{-1/2}\chfc{\set\emptyset}\dmf{t\in[0,1],\mu\in\M_{0,t}}.
  \end{displaymath}
The product system $\E^\Us$ generated by all units of $\E^\M$ is given by 
\begin{displaymath}
\E^\Us_t=\set{\psi\in L^2(\M_{0,t}):\chfc{\set{Z:\#Z=\infty}}\psi_\mu=0\text{~for all $\mu\in\mu_{0,t}$}}\dmf{t\in\NRp}.  
\end{displaymath}

Moreover, $\M^{\E,u}=\M$ and $\M^{\E,\Us}=\M\circ l^{-1}$, denoting $l$ the map  $Z\mapsto\check{Z}$. Consequently, for the type of $\E$ there are the following possibilities:
\begin{displaymath}
  \begin{array}[c]{|c|l|}\hline
\mathrm{I}_0&\text{~if~}\M=\set{\delta_\emptyset}\\\hline
\mathrm{I}_1&\text{~if~}\M=\set{\mu:\mu\sim\Pi_\ell}\\\hline
\mathrm{II}_0&\text{~if~$\M$-a.s. $\check Z=Z$  and~}\M\ne\set{\delta_\emptyset}\\\hline
\mathrm{II}_1&\text{~otherwise~}\\\hline
  \end{array}
\end{displaymath}
\end{corollary}
The symbol $\FG^f_K$ denotes  the set of  finite sets:
\begin{equation}
\label{eq:defFfinite}
  \FG^f_K=\set{Z\subseteq K:\# Z<\infty}.
\end{equation}
\begin{proof}
   Clearly, 
  \begin{eqnarray*}
   \lefteqn{\sqrt{\mu_{0,s+t}(\set{\emptyset})}\chfc{\set\emptyset}(Z_{0,s+t})}&=&\sqrt{\mu_{0,s+t}(\set{\emptyset})}\chfc{\set\emptyset}(Z_{0,s})\chfc{\set\emptyset}(Z_{s,s+t}) \\
&=&\sqrt{\frac{\mu_{0,s+t}(\set{\emptyset})}{\mu_{0,s}(\set{\emptyset})\mu_{s,s+t}(\set{\emptyset})}}\sqrt{\mu_{0,s}(\set{\emptyset})}\chfc{\set\emptyset}(Z_{0,s})\sqrt{\mu_{s,s+t}(\set{\emptyset})}\chfc{\set\emptyset}(Z_{s,s+t}) \\
&=&V_{s,t}\sqrt{\mu_{0,s}(\set{\emptyset})}\chfc{\set\emptyset}\otimes\sqrt{\mu_{s,s+t}(\set{\emptyset})}\chfc{\set\emptyset}(Z)
  \end{eqnarray*}
This shows that $u$ factorizes. Clearly, the map $t\mapsto\mu_{0,t}(\set\emptyset)$ is decreasing and therefore measurable. Since $\chfc{\set\emptyset}$ is measurable too, $u$ is a unit of $\E^\M$.

For the second part, we see that $\scpro{u_t}{\psi}=\mu^{1/2}(\emptyset)\psi_\mu(\emptyset)$ such that $(\Pr{u_t}\psi)_\mu=\psi_\mu(\emptyset)\chfc{\set\emptyset}$. Fix $0\le s<t\le1$, a measure $\mu\in\M$ with $\mu=\mu_{0,s}\ast\mu_{s,t}\ast\mu_{t,1}$ and $\psi\in L^2(\mu_{0,s})$, $\psi'\in L^2(\mu_{s,t}-s)$, $\psi''\in L^2(\mu_{t,1}-t)$. We define $\tilde\psi\in L^2(\M)$ by 
\begin{displaymath}
  \tilde\psi_\mu(Z)=\psi(Z_{0,s})\psi'(Z_{s,t}-s)\psi''(Z_{t,1}-t).
\end{displaymath}
Then 
\begin{eqnarray*}
  (\mathrm{P}^u_{s,t}\tilde\psi)_\mu(Z)&=&\psi(Z_{0,s})\Pr{(u_t)_{\mu_{s,t}-s}}\psi'(Z_{s,t}-s)\psi''(Z_{t,1}-t)\\
&=&\psi(Z_{0,s})\psi'(\emptyset)\chfc{\set\emptyset}(Z_{s,t}-s)\psi''(Z_{t,1}-t)\\
&=&\psi(Z_{0,s})\psi'(Z_{s,t}-s)\chfc{\set\emptyset}(Z_{s,t}-s)\psi''(Z_{t,1}-t)\\
&=&\chfc{\set{Z':Z'\cap[s,t]=\emptyset}}(Z)\psi(Z_{0,s})\psi'(Z_{s,t}-s)\psi''(Z_{t,1}-t).
\end{eqnarray*}
This shows that $\mathrm{P}^u_{s,t}$ is multiplication by $\chfc{\set{Z:Z\cap[s,t]=\emptyset}}$. By normality, $J_{\mathrm{P}^u}(f)$ is multiplication by $f$ for all $f\in L^\infty(\M)$ what implies $\M^{\E,u}=\M$. The relation $\M^{\E,\Us}=\M\circ l^{-1}$ follows from  Proposition \ref{prop:unitandunitalprojections}. This  implies also  $\mathrm{P}^\Us_{s,t}=J_{\mathrm{P}^u}(\chfc{\set{Z:\#(Z\cap[s,t])<\infty}})$. 

For  a stationary factorizing measure type  $\M$ on $\FG_\NT$ Proposition \ref{prop:finitepoisson} shows that the restriction of $\M$ to $\FG^f_\NT$ is equivalent either to $\Pi_{0}$  or $\Pi_\ell$. This  implies either  $\E^\Us=\E^{\set{\delta_\emptyset}}=\Gammai(\set0)$ or $\E^\Us=\E^{\set{\mu:\mu\sim\Pi_\ell}}=\Gammai(\NC)$ respectively. If $\E=\E^\Us$, this amounts to the types $\mathrm{I}_0$ and $\mathrm{I}_1$ respectively. Otherwise, $\E$ has type $\mathrm{II}_0$ or $\mathrm{II}_1$. From   Proposition \ref{prop:unitandunitalprojections} we derive that $\M|\FG^f_\NT\sim \delta_{\emptyset}$ is equivalent to $\check Z=Z$ $\M$-a.s.  This completes the proof.
\end{proof}
\begin{remark}
  Similarly to $\M=\set{\delta_\emptyset}$, $\M=\set{\delta_\NT}$ yields a type $\mathrm{I}_0$ product system.
\end{remark}
\begin{remark}
  Similarly, one  gets a type $\mathrm{II}$ product system  for compact $L$, but then there are more choices for the numerical index of $\E^\M$, see proposition \ref{prop:exampleII_2} and example \ref{ex:three Poisson in two}.

If $L$ is not compact things are much more complicated since \ref{mu2} from Corollary \ref{cor:propmu} is not valid, see Example \ref{ex:infinite Poisson}. Therefore even type $\mathrm{III}$ product system could arise, although we have no example for that behaviour. The measure type of the Poisson process in  example \ref{ex:infinite Poisson} yields a type $\mathrm{I}_\infty$ product system and alteration of finitely many factors yields type $\mathrm{II}$. Whether or not type $\mathrm{III}$ is possible should be first discussed in the context of infinite tensor products of product systems \textsc{Bhat} (personal communication).  
\end{remark}
We want to underline the complex structure of type II product systems by several examples, following some ideas of \textsc{Tsirelson} \cite{Tsi00}.
\subsection{Example 1: Finite Random Sets}
\label{sec:finite}
First we  look at product systems $\F$ for which   the measure type $\M^\F$ is concentrated on finite sets $\FG^f_\NT$ (defined like in  \ref{eq:defFfinite}). The main examples are Poisson processes $\Pi_\nu$ defined in \ref{eq:defPoisson}.
\begin{proposition}
  \label{prop:finitepoisson}
  Suppose $\mu$ is  a quasistationary quasifactorizing probability measure on $\NT$ with $\mu(\FG^f_\NT)=1$. Then either $\mu\sim \Pi_\ell$ or $\mu=\delta_\emptyset$.
\end{proposition}
\begin{proof}
Consider the sets $Y^{s,t}_{\set{s_1,\dots,s_k}}\subset\FG^f_{\NT}$ for $0\le s=s_0<s_1<\dots<s_k<s_{k+1}= t\le1$ given by
\begin{displaymath}
  Y^{s,t}_{\set{s_1,\dots,s_k}}=\set{Z:\#(Z\cap[s_i,s_{i+1}])\le1\forall i=0,\dots,k}.
\end{displaymath}
Especially, $Y^{s,t}_{\emptyset}=\set{Z: \#(Z\cap[s,t])\le1}$. Then it is easy to see that 
\begin{displaymath}
  \FG^f_{\NT}=\bigcap_{S\in\FG^f_{(0,1)},~S\subset\NQ}Y^{0,1}_S
\end{displaymath}
since for any finite set with at least two elements the is a positive minimal distance between each pair of such elements. Therefore it is sufficient to fix the equivalence class of all $\mu|Y^{0,1}_S$. Further, by \ref{eq:quasifactorizingmeasure}, we derive
\begin{displaymath}
  \mu|Y^{0,1}_{\set{s_1,\dots,s_k}}\sim \mu_{0,s_1}\ast\mu_{s_1,s_2}\ast\dots\ast\mu_{s_k,1}|Y^{0,1}_{\set{s_1,\dots,s_k}}=\mu_{0,s_1}|Y^{0,s_1}_\emptyset\ast\mu_{s_1,s_2}|Y^{s_1,s_2}_\emptyset\ast\dots\ast\mu_{s_k,1}|Y^{s_k,1}_\emptyset
\end{displaymath}
  Thus we need to determine   the equivalence classes of $\mu_{s_1,s_2}|Y^{s_1,s_2}_\emptyset$ only. By  \ref{eq:quasifactorizingmeasure}, it is enough to  do this for $\mu|Y^{0,1}_\emptyset$. From assumption we know $\mu(\set\NT)=0$ and Corollary  \ref{cor:propmu} shows $\mu(\set\emptyset)>0$. Thus we need to determine   $\mu|\set{Z:\#Z=1}$ only. Especially, we have  to show that $\mu|\set{Z:\#Z=1}\sim\Pi_\ell|\set{Z:\#Z=1}$ or $\mu|\set{Z:\#Z=1}=0$.  

We can identify $\set{Z:\#Z=1}$ as Borel subset of $\FG_\NT$ with $\NT$ via $t\mapsto\set{t}$. Correspondingly,  there exists a Borel measure $\mu'$ on $\NT$ such that $\mu'(Z')=\mu(\set{Z:Z\subseteq Z',\# Z=1})$.   But \ref{eq:quasistationarymeasure} shows that $\mu$ is quasiinvariant under the shift on $\NT$. 
Since $\set{Z:\#Z=1}$ is invariant under this shift, $\mu'$ is equivalent to the  shift invariant measure $\int_\NT\ell(\d t)(\mu'+t)$.  The latter   must be  a multiple of  Lebesgue measure on $\NT$. Thus  either $\mu'=0$ or $\mu'\sim\ell$ and the proof is  complete.
\end{proof}
The proof of the following proposition  is formulated and proved in a more general form  in  section \ref{sec:direct integrals ps}, Corollary \ref{cor:exponentialhilbertspace}.
\begin{proposition}
  \label{prop:exponentialhilbertspace}
  Suppose $\E$ is a product system with a unit $\sg tu$ for which  the measure type  $\M^{\E,u}$ from Proposition \ref{prop:psunit2randomset} is concentrated on $\FG^f_\NT$. Then $\E$ is an exponential product system $\Gammai(\K)$. 
\end{proposition}
We want to show here, that the key point of the new type of examples from \cite{Tsi00} is the switch from independence to quasifactorisation  in equation \ref{eq:quasifactorizingmeasure}. A similar result is given by \cite[Lemma 5.6]{Tsi03}.
\begin{proposition}
  \label{prop:independentmeasure}
  Suppose $\mu$ is a quasistationary  probability measure on $\FG_\NT$ with
\begin{equation}
  \label{eq:strictlyfactorizingmeasure}
  \mu_{r,t}=\mu_{r,s}\ast\mu_{s,t}\dmf{(r,s),(s,t)\in I_{0,1}}.
\end{equation}
Then $\mu$ is a Poisson process or $\delta_\NT$.
\end{proposition}
\begin{proof}
  According to Corollary \ref{cor:propmu}, we need to consider only the case $\mu(\set\emptyset)>0$. Define  the functional $Q(s,t)=\mu(\set{Z:Z\cap[s,t]=\emptyset})=\mu(\set{Z:Z\cap(s,t)=\emptyset})$. Then 
$Q$ is continuous and fulfils  
\begin{displaymath}
Q(r,t)=Q(r,s)Q(s,t)\dmf{(r,s),(s,t)\in I_{0,1}}.
\end{displaymath}
Further, $\mu(\set\emptyset)>0$ and \ref{eq:strictlyfactorizingmeasure} imply $Q(s,t)\ne0$ for all $s\ne t$. Thus, $-\ln Q$ is nonnegative, additive and continuous. Consequently,   there is a (diffuse) measure $\nu$ on $[0,1]$ such that $\nu((0,t))=-\ln Q(0,t)$, or $Q(s,t)=\e^{-\nu((s,t))}$. By \ref{eq:strictlyfactorizingmeasure}, we see that $\mu(\set{Z:Z\cap G=\emptyset})= \e^{-\nu(G)}=\Pi_\nu(\set{Z:Z\cap G=\emptyset}) $ for all open $G\subset\NT$. The  Choquet theorem \cite[Theorem 2-2-1]{C:Mat75} implies that $\mu=\Pi_\nu$.
\end{proof}
\subsection{Example 2: Countable Random Sets}
\label{sec:countable}
Now we want to look at stationary factorizing measure types on $\FG_\NT$ being concentrated on special countable closed  sets. This study advances from the construction in   \cite[section 5]{Tsi00} which  used the jump times of a continuous time  Markov process on a discrete set. Since the process was regular, the jump times formed a countable (closed) set. Recall that $\check{Z}$ is the set of  limit points of $Z$. The following result can be derived from standard set theory.
\begin{lemma}
  \label{lem:countableRACSnotperfect}
  Let $Z$ be a countable closed set. Then the sequence $\sequ nZ$ defined through  $Z_0=Z$, $Z_{n+1}=\check{Z_n}$, $n\ge1$ is strictly  decreasing  unless  $Z_n=\emptyset$ for some $n\in\NN$.
\end{lemma}
\begin{proof}
Suppose  $\check{Z'}=Z'$, i.e.\ $Z'$  is a perfect set. It follows from   \cite[Theorem 6.65]{A:HS75} that $Z'$  is uncountable or empty. Since $Z$ is countable all $Z_n$ are and the assertion is proven. 
\end{proof}
Thus   $N(Z)=\inf\set{n\in\NN:Z_n=\emptyset}$ may serve as a degree of complexity of $Z$. But, observe  firstly that $N(Z)=\infty$ is possible. Further,  if $Z$ is random and almost surely countable, $N(Z)$ need not be uniformly bounded over all realisations. Example \ref{ex:complexcountable} provides some illustrations of such situations.

We restrict our attention here to the simplest case  to show that  there is a huge set of tractable examples. Assume that $\mu$  is concentrated on  the set of all $Z$ which fulfil
\begin{condition}{C}
\item \label{count1} $\check{Z}\in\FG^f_\NT$  (may be, it is empty) such that $\check{\check{Z}}=\emptyset$ and
\item \label{count2} for all $t\in[0,1)$ there is some $\varepsilon(t)>0$ such that $(t,t+\varepsilon(t)]\cap Z=\emptyset$, i.e.\ $Z$ has no limit point from the right.
\end{condition}
If $Z\in\FG_\NT$, $\check{Z}\ne\emptyset$,  obeys these two conditions, we can construct a sequence  $\sequ n\lambda=\sequp n{\lambda_n(Z)}\subset[0,1]$ in the following way. Let $Z\cap [0,\min\check{Z}]=\set{\xi_n:n\in\NN}$ where $\sequ n\xi$ is strictly increasing.  Then  set
\begin{displaymath}
  \lambda_0=\frac{\min\check{Z}-\xi_0}{\min\check{Z}}
\end{displaymath}
and 
\begin{displaymath}
  \lambda_{n+1}\xi_n+(1-\lambda_{n+1})\min\check{Z}=\xi_{n+1}\dmf{n\in\NN},
\end{displaymath}
i.e.\ $\lambda_{n+1}=\frac{\min\check{Z}-\xi_{n+1}}{\min\check{Z}-\xi_n}$.

Now consider for a  probability measure $\nu$ on  $[0,1]^\NN$ the following three conditions.
\begin{condition}{F}
\item \label{fact1} For all $k\in \NN$ the distribution of $(\lambda_0,\dots,\lambda_k)$ under $\nu$ is equivalent to Lebesgue measure on $[0,1]^{k+1}$.
\item \label{fact2} The law of the sequence $\sequp n{\lambda_{n+1}}$ under $\nu$ is equivalent to $\nu$.
\item \label{fact3} Almost surely under $\nu$, the relation  $\prod_{n\in\NN}\lambda_n=0$ is valid.  
\end{condition}
\begin{proposition}
  \label{prop:simplecountable}
Let $\mu$ and $\mu'$ be two quasistationary quasifactorizing measures on $\FG_\NT$ for which almost all $Z\in\FG_\NT$ fulfil \ref{count1} and \ref{count2} and $\check{Z}\ne\emptyset$ occurs with positive probability. Then $\mu$ and $\mu'$ are equivalent iff the laws of $\lambda=\sequ n\lambda$, conditional on $\check{Z}\ne\emptyset$, under $\mu$ and $\mu'$ are equivalent. 

Moreover, for a measure type $\N$ on $[0,1]^\NN$ there is a quasistationary quasifactorizing measure $\mu$ on $\FG_\NT$ with realisations obeying \ref{count1} and \ref{count2} almost surely and $\mu(\check{Z}\ne\emptyset)>0$ such that the distribution of $\sequ n\lambda$, conditional on $\check{Z}\ne\emptyset$,  belongs to $\N$ iff $\N$ fulfils \ref{fact1}--\ref{fact3}.
\end{proposition}
\begin{proof}
Equation  \ref{eq:quasifactorizingmeasure} implies like in the proof of Proposition \ref{prop:finitepoisson} that it is enough to fix the measure types of the law of $Z$ conditional on $\# \check Z=0$ and conditional on  $\# \check Z=1$. From the same  proposition it follows that the former is the measure type of a Poisson process, since $\# \check Z=0$ is the same as $Z\in \FG_\NT^f$. Similar arguments show for  conditioning on $\# \check Z=1$ that it is sufficient to fix (the equivalence class of) the law of $Z$ conditional on $\# \check Z=1$ and $\max Z=\max\check Z$. By quasistationarity the  equivalence class of this law is fully determined by the distribution of $\lambda=\sequ n\lambda$, conditional on $\check Z\ne\emptyset$ and vice versa.

Moreover, there is almost surely  some rational $t>0$ such that $Z_{0,t}=Z\cap[0,t]$  is finite. Again \ref{eq:quasifactorizingmeasure}  and Proposition \ref{prop:finitepoisson} imply that the law of $Z_{0,t}$, conditional on a fixed  $t$, belongs to the measure type of  a Poisson process with intensity measure $\ell(\cdot\cap[0,t])$.  Conditioning on the numbers of points in $[0,k]$ we find that the law of $(\xi_0,\dots,\xi_k)$ is equivalent to Lebesgue measure on the set $0<\xi_0<\dots<\xi_k<t$. Since $\min \check Z$ has also a law equivalent to Lebesgue measure, conditional on $\# Z\cap[0,t]<\infty$, we obtain  that the law of $(\lambda_0,\dots,\lambda_k)$ is equivalent to Lebesgue measure on $[0,1]^{k+1}$, what is condition \ref{fact1}. Observe that for any $t\in\NT$ the laws of $Z_{t,1}$ conditional on $Z_{0,t}=\emptyset$ and $\# Z_{0,t}=1$ are equivalent. But, the latter impose the shift by one in the sequence $\sequ n\lambda$ with regard to the former. This shows \ref{fact2}. \ref{fact3} follows from $\xi_n\limitsto{}{n\to\infty}\min\check Z$. 

 Sufficiency  follows from  Corollary \ref{cor:Zalphabeta}, similar to the proof of Theorem \ref{th:surj}. We may  choose $Q_{s,t}$ as the law of $\bigcup_{n\in\NN}\set{s+\lambda_0\cdots\lambda_n(t-s)}\cup\set{t}$ under some law from $\N$ for $\sequ n\lambda$ and $Q_\emptyset$ as Poisson process. Then \ref{fact1} establishes \ref{q3}, \ref{fact2} forces \ref{q2a}, both together yield \ref{q2} whereas \ref{q1} is fulfilled by construction. \ref{fact3} proves that the realisations have really limits in the predefined  points of $\check Z$. This completes the proof. 
\end{proof}
\begin{corollary}
  \label{cor:uncountablymanycountable}
  The set of   different stationary factorizing measure types on $\FG_\NT$ concentrated on countable sets has the cardinality of the continuum.
\end{corollary}
\begin{proof}
  For any $\alpha>0$ set $\lambda^\alpha_n=U_n^{\alpha}$, where $\sequ nU$ is i.i.d.\ uniformly distributed in $[0,1]$. Clearly, the distribution of $\lambda^\alpha$ fulfils \ref{fact1} and \ref{fact2}. From  the law of large numbers for $\sequ n{\ln \lambda^\alpha}$ we derive \ref{fact3}.  Since the distributions of $\lambda^\alpha_n$ and $\lambda^\beta_n$ are different for $\alpha\ne\beta$, the laws of $\lambda^\alpha$ and $\lambda^\beta$ are singular, what can again be seen from the law of large numbers. 

On the other side, $\FG_\NT$ is a Polish space, so the set of all probability measures on it has at most the cardinality of the continuum. This completes the proof. 
\end{proof}
\subsection{Example 3: Random Cantor Sets}
\label{sec:perfect}
Here we want to deal with stationary factorizing measure types $\M$ which are such that $\E^\M$ has type $\mathrm{II}_0$. Corollary \ref{cor:psrstypeII} and Proposition \ref{prop:unitandunitalprojections} show that this happens if and only if   $\check{Z}=Z$ almost surely and $Z\ne\emptyset$ with positive probability.   We construct for such nonempty $Z$ random numbers $(\lambda^l_{i,j})_{0\le j\le 2^i-1}\subset [0,1]$, $(\lambda^r_{i,j})_{0\le j\le 2^i-1}\subset [0,1]$ in the following way:

Set $x_{0,0}=\min Z$, $y_{0,0}=\max Z$ and choose $z_{0,0}$ uniformly in $[x_0,y_0]$. Then we set $x_{1,0}=x_{0,0}$, $y_{1,0}=\max (Z\cap[x_0,z_0])$ and  $x_{1,1}=\min (Z\cap[z_0,y_0])$ and $y_{1,1}=y_{0,0}$. We proceed by choosing all $z_{i,j}$ independent and uniformly in $[x_{i,j},y_{i,j}]$ and setting
\begin{displaymath}
  \begin{array}{lll}
  x_{i+1,2j}&=&x_{i,j}\\
y_{i+1,2j}&=&\max Z\cap[x_{i,j},z_{i,j}]\\
 x_{i+1,2j+1}&=& \min Z\cap[z_{i,j},y_{i,j}]\\
 y_{i+1,2j+1}&=&y_{i,j}
  \end{array}
\end{displaymath}
Now, set $\lambda_{0,0}^l=x_{0,0}$,  $\lambda_{0,0}^r=(1-y_{0,0})/(1-x_{0,0})$, 
\begin{displaymath}
  \begin{array}{lll}
\lambda^l_{i,j}&=&\frac{y_{i+1,2j}-x_{i+1,2j}}{z_{i,j}-x_{i,j}}\\
\lambda^r_{i,j}&=&\frac{y_{i+1,2j+1}-x_{i+1,2j+1}}{y_{i,j}-z_{i,j}}
  \end{array}
\end{displaymath}
Observe 
\begin{equation}
\label{eq:random Cantor set}
  Z=\bigcap_{i\in\NN}\bigcup_{j=0}^{2^i-1}[x_{i,j},y_{i,j}].
\end{equation}
\begin{remark}
  Equation \ref{eq:random Cantor set} shows that every nonempty realisation is the limit of iterating the procedure of dividing every existing interval into two new ones by discarding an interval inside the old one. This procedure is the same as for the construction of the standard Cantor set in $[0,1]$, except of using random numbers $\lambda^l_{i,j},\lambda^r_{i,j}$ instead of the fixed $1/3$ and $2/3$ for the standard Cantor dust. This analogy motivated the title of this section.
\end{remark}
\begin{example}
\label{ex:zerosBt}
 Let $\M$ be constructed from the set of zeros of Brownian motion. Consider the set of zeros of a Brownian motion $\sg tB$ in $\NR$ with the  law  $\mathscr{L}(B_0)$ of the starting point being  equivalent to Lebesgue measure, i.e.\ $Z=\set{t\in[0,1]:B_t=0}$ and $\mu=\mathscr{L}(Z)$. It is well-known \cite{C:Per98} that $Z$ is almost surely an uncountable closed  set with no isolated point if $B_t$ starts in $0$.  I.e., conditional on $Z\ne\emptyset$, $Z$ is a perfect set.   In  \cite[section 2]{Tsi00} there was proven that the law of $Z_a=\set{t\in[0,1]:B_t=a}$ is quasistationary and quasifactorizing for any Brownian motion starting in $0$ and $a\ne0$. But it is clear that the laws $\mathscr{L}(Z)$ and $\mathscr{L}(Z_a)$ are equivalent.  We prefer the former way to describe this measure type.
 
Since the distribution of the hitting time $T_a=\min\set{s:B^x_s=a}$ is for all $x\ne 0$  equivalent to Lebesgue measure  on $\NRp$ \cite{C:RY99}, we see from the independence of increments that the distribution of  $(x_{0,0},y_{0,0})$ is equivalent to Lebesgue measure $\ell\otimes\ell$ restricted to $I_{0,1}$. Conditionally on $Z\ne\emptyset$  and $(x_{0,0},y_{0,0})$, $Z$ is the set of zeros of a Brownian bridge starting at 0 at time $x_{0,0}$ and ending in 0 at time $y_{0,0}$. Similar  arguments apply to the following steps, we don't go into the messy details here.
\end{example}
\begin{proposition}
  \label{prop:perfectLebesgueequivalence}
 The joint distribution of the  sequences $\lambda^r,\lambda^l$, characterizes $\M$. This joint  law is locally absolutely continuous with respect to  Lebesgue measure $\ell$ on $[0,1]$ and $\mu$ and $\mu'$ are equivalent iff these laws are equivalent.  
\end{proposition}
\begin{proof}
  The proof is similar of the case of countable sets, but we use another  construction of a finite set.  

First observe that if we had chosen instead of $0$ another point in $\NT$ a similar  construction would have produced other random variables $(\lambda^l,\lambda^r)$ with,   by quasistationarity, equivalent  law. Thus we could choose this point uniformly on $\NT$. 

We enumerate the points $x_{i,j},y_{i,j}$ in a distinct manner, say $t_{2^k+2l}=x_{k,2l+1}$ and $t_{2^k+2l+1}=y_{k,2l}$ for $l=0,\dots,2^{k-1}-1$. Choose independently a random variable $\xi$ geometrically distributed with mean 1 and define $Z'=\bigcup_{i=1}^\xi\set{t_i}$. Clearly, the correspondence $Z\mapsto Z'$ defines a stochastic kernel $q$ on $\FG_\NT$ which, extended by $q(\emptyset,\cdot)=\delta_\emptyset$, fulfils \ref{qq1}   and \ref{qq2}. Proposition \ref{prop:extensionRACS} shows that the measure type $\M\circ q$ is again stationary and factorizing. By Proposition \ref{prop:finitepoisson}, it is the measure type of the Poisson process. Thus the law of $\sequ it$ is locally absolutely continuous with respect to Lebesgue measure restricted to a certain simplex. Now the construction implies that  $(\lambda^l_{i,j},\lambda^r_{i,j})_{0\le i\le k,0\le j\le 2^i-1}$ has a strictly positive density with respect to Lebesgue measure on $[0,1]^{2^{k+2}-2}$.
\end{proof}

\begin{remark}
Similar to the case of countable closed sets one can characterize the possible laws of $(\lambda^l,\lambda^r)$ completely. By these means one can show that the set of   different stationary factorizing measure types concentrated on perfect sets has the cardinality of the continuum. We do not pursue these ideas  here since \cite[section 4]{Tsi00} and \cite{Tsi03} already established the latter result.
\end{remark}

\subsection{Tensor Products (II)}
\label{sec:tensor products(II)}

In this section, we want to analyse the question  whether  tensor products of product systems can  yield new measure types $\M^{\E^\M,\Us}$.
\begin{proposition}
  \label{prop:convolutionmeasuretypes}
  Let $\M$ be a  stationary factorizing measure type on $\FG_\NT$ with 
  \begin{equation}
 \label{eq:empty intersection}
  Z_1\cap Z_2=\emptyset \dmf{\M\otimes\M-\text{a.a.~}(Z_1,Z_2)}.
  \end{equation}
 Then $\M\ast\M=\M$. 
\end{proposition}
\begin{proof}
Fix any $\mu\in\M$ and suppose $\mu(Y)>0$.  Then  
\begin{displaymath}
  \mu\ast\mu(Y)=\mu\otimes\mu(\set{(Z_1,Z_2):Z_1\cup Z_2\in Y})\ge\mu\otimes\mu(\set{(Z_1,\emptyset):Z_1\in Y})=\mu(Y)\mu(\set\emptyset)>0.
\end{displaymath}
Thus $\mu\ll \mu\ast\mu$.

Now assume $\mu(Y)=0$ for some Borel set $Y\subset\FG_\NT$. We introduce the Borel sets $ Y_n\subset\FG_{[0,1]}\times\FG_{[0,1]}$,
\begin{displaymath}
  Y_n=\set{(Z_1,Z_2):Z_1\cap[\frac{k-1}{2^n},\frac{k}{2^n}]=\emptyset\text{~or~}Z_2\cap[\frac{k-1}{2^n},\frac{k}{2^n}]=\emptyset\text{~for all~}1\le k\le 2^n}.
\end{displaymath}
Clearly, $\bigcup_{n\in\NN}Y_n=\set{(Z_1,Z_2):Z_1\cap Z_2=\emptyset}$ and this set has, due to our assumption, full measure with respect to  $\mu\otimes\mu$. For two disjoint sets $S_1,S_2\subset\set{1,\dots, 2^n}$ of indeces define 
\begin{displaymath}
  Y^n_{S_1,S_2}=\set{(Z_1,Z_2):\text{for all~$i=1,2$~ and~}k\in\set{1,\dots, 2^n}\setminus S_i\text{~}Z_i\cap[\frac{k-1}{2^n},\frac{k}{2^n}]=\emptyset}.
\end{displaymath}
Obviously, 
\begin{displaymath}
  Y_n=\bigcup_{
S_1,S_2\subset\set{1,\dots, 2^n},
S_1\cap S_2=\emptyset,
S_1\cup S_2=\set{1,\dots, 2^n}
}Y^n_{S_1,S_2}\dmf{n\in\NN}.
\end{displaymath}
 Fix $n\in\NN$, the measures $\mu_k=\mu_{\frac{k-1}{2^n},\frac{k}{2^n}}$ and  appropriate $S_1,S_2\subset\set{1,\dots, 2^n}$. We find
\begin{eqnarray*}
\lefteqn{  \bigast_{1\le k\le 2^n}\mu_k\otimes\bigast_{1\le l\le 2^n}\mu_l(\set{(Z_1,Z_2):Z_1\cup Z_2\in Y,(Z_1,Z_2)\in Y^n_{S_1,S_2}})}
&=&\bigast_{k\in S_1}\mu_k\otimes\bigast_{l\in S_2}\mu_l(\set{(Z_1,Z_2):Z_1\cup Z_2\in Y})\prod_{k\in \cmpl{S_1}}\mu_k(\set\emptyset)\prod_{l\in \cmpl{S_2}}\mu_l(\set\emptyset)\\&=& \bigast_{1\le k\le 2^n}\mu_k(Y) \bigast_{1\le l\le 2^n}\mu_l(\set\emptyset).
\end{eqnarray*}
From $\mu\sim\bigast_{1\le k\le 2^n}\mu_k$ we derive $\mu\otimes\mu(\set{(Z_1,Z_2):Z_1\cup Z_2\in Y,(Z_1,Z_2)\in Y^n_{S_1,S_2}})=0$ and therefore $\mu\ast\mu(Y)=0$. The proof is complete.
\end{proof}
\begin{remark}
  We conjecture that the above result is true regardless of \ref{eq:empty intersection}.   
\end{remark}
\begin{remark}
In contrast to product systems of Hilbert spaces a general tensor product construction for product systems of Hilbert modules still it not known to exist.   An analogue of the procedure to  alternate local products of the measures to get their convolution, was used by \textsc{Skeide} \cite{Q:Ske01} to construct such a tensor product at least for  product systems of Hilbert modules with a so-called central unit. There the empty set is replaced by a projection onto a so-called central unit. Unfortunately,  his construction does not necessarily yield  the tensor product of   $\E^\M$ and $\E^{\M}$ if \ref{eq:empty intersection} is not satisfies, see the forthcoming \cite{Q:BLS03}. This is another reason for looking for techniques to decide whether \ref{eq:empty intersection} happens or not.    
\end{remark}
\begin{remark}
\label{rem:Kahane}
Concerning \ref{eq:empty intersection}, there is a general theorem by \textsc{Kahane} \cite{C:Kah86}, which states that 
\begin{equation}
\label{eq:Kahane}
  \dim_{\mathscr{H}}(Z_1\cap (Z_2+r))\ge \dim_{\mathscr{H}} Z_1+\dim_{\mathscr{H}} Z_2 -1-\varepsilon
\end{equation}
for all $\varepsilon>0$ for a set of $r\in\NR$ of positive Lebesgue measure, denoting $\dim_{\mathscr{H}} Z$ the Hausdorff dimension of $Z$ defined as follows.  The \emph{Hausdorff measure}  of a Borel set $F$ with respect to a function $\map h\NRp\NRp$ is defined as 
\begin{equation}
  \label{eq:defHausdorff}
  \mathscr{H}^h(F)=\sup_{\varepsilon>0}\mathscr{H}^h_\varepsilon(F),
\end{equation}
where
  \begin{equation}
  \label{eq:defHausdorffa}
 \mathscr{H}^h_\varepsilon(F)= \inf\set{\sum_{i\in\NN}h(d(B_i)):\text{$\sequ iB$ are balls with  $d(B_i)\le\varepsilon$ and $\bigcup_{i\in\NN}B_i\supseteq F$}},
  \end{equation}
denoting $d(B)$ the diameter of $B$. Then  the \emph{Hausdorff dimension} $\dim_{\mathscr{H}}F$ of a Borel set $F$ is  defined by 
\begin{displaymath}
  \dim_{\mathscr{H}} F=\inf\set{\alpha>0:\mathscr{H}^{h^\alpha}(F)>0 \text{~for~}h^\alpha(\varepsilon)=\varepsilon^\alpha  }.
\end{displaymath}

By quasistationarity we derive  $\mu\otimes\mu(\set{(Z_1,Z_2):Z_1\cap Z_2\ne\emptyset})>0$ if $\dim_{\mathscr{H}} Z>1/2$ for a set of $Z$ with positive $\mu$-measure.  Consequently, random closed sets with high Hausdorff dimension  will not obey \ref{eq:empty intersection}.
\end{remark}
\begin{example}
\label{ex:zerosBt1}
 We continue Example \ref{ex:zerosBt}, denote $Z=\set{t\in[0,1]:B_t=0}$ the set of zeros of a Brownian motion $\sg tB$ in $\NR$ with   $\mathscr{L}(B_0)\sim\ell$. Then we see  for an independent copy  $Z'$  of $Z$ that $Z\cap Z'$  is the distribution of $\set{t\in[0,1]:B^{(2)}_t=0}$, where $\sg t{B^{(2)}}$ is a Brownian motion in $\NR^2$ with the initial law  $\mathscr{L}(B^{(2)}_0)$ being  equivalent to Lebesgue measure $\ell^2$. But, $\set{t\in[0,1]:B^{(2)}_t=0}$ is almost surely void conditional on the fact   that the  Brownian motion does not start in $(0,0)$ \cite{C:Per98}. Since the latter happens almost surely, we find $\mu\ast\mu\sim\mu$ in this case.  The above formula \ref{eq:Kahane} does not allow the same conclusion  since $\dim_{\mathscr{H}} Z=1/2$ for $\mu$-a.a.\ nonempty $Z$ in this case \cite{C:Per98}.
\end{example}
\begin{example}
  Following \cite[Section 4]{Tsi99a},  suppose $\mu_\delta$ is the distribution of the set of zeros of a  Bessel process with parameter $\delta\in(0,2)$. This random set is  constructed for this Bessel process in the same scheme which was used for Brownian motion in Example \ref{ex:zerosBt}. It is well-known that  $\dim_{\mathscr{H}}(Z)=1-\delta/2$ almost surely conditional on $Z\ne\emptyset$ \cite{Tsi99a,C:BG60}.    From Note \ref{rem:Kahane} we derive that  for $\delta<1$ the relation  \ref{eq:empty intersection} is not almost surely valid. Nevertheless, it is not clear whether $\mu\sim\mu\ast\mu$ is true or not in this case. Since almost surely $\check{Z}=Z$, this would be implied by  $\E^\M\cong\E^\M\otimes \E^\M$.
\end{example}
\begin{example}
  We can also consider  the measure $F_1$ introduced in \ref{eq:defexponentialmeasure}. Since $F_1$ is a multiple of  the Poisson process $\Pi_\ell$, Corollary \ref{cor:propmu}, \ref{mu1} implies $F_1\otimes F_1(\set{(Z_1,Z_2):Z_1\cap Z_2\ne\emptyset})=0$. Consequently, we derive    $F_1\ast F_1\sim F_1$. On the other side,  one can compute directly 
  \begin{displaymath}
    \int F_1(\d Z)\int F_1(\d Z')f(Z\cup Z')=\int F_1(\d Z'')2^{\# Z''}f(Z'').
  \end{displaymath}
This is a direct conclusion of the  $\rlap{$\DS\sum$}\DS\int$ lemma, see \cite{FF91b,LP89}, 
  \begin{equation}
\label{eq:sigmaintegrallemma}
    \int F_1(\d Z)\int F_1(\d Z')f(Z,Z')=\int F_1(\d Z'')\sum_{Z\cup Z'=Z'',Z\cap Z'=\emptyset}f(Z,Z').
  \end{equation}
That lemma can serve as one basis of quantum stochastic calculus on the symmetric Fock space, i.e.\ exponential product systems, cf.\ \cite{Lin93}.  Since we obtain  a similar   $\rlap{$\DS\sum$}\DS\int$-lemma for quasistationary quasifactorizing measures $\mu$ (at least if $\mu\ast\mu\sim\mu$) below, we expect  a calculus on the product systems $\E^\M$ too.
\end{example}
\begin{corollary}
  \label{cor:sumintegrallemma}
  For all quasistationary quasifactorizing measures $\mu$ there is a unique (upto $\mu$-equivalence) stochastic kernel $q^\mu$ from $\FG_\NT$ to $\FG_\NT\times\FG_\NT$ with $q^\mu(Z'',\set{(Z,Z'):Z\cup Z'=Z''})=1$ for all $Z''\in\FG_\NT$ such that for all positive measurable functions $\map{f}{\FG_\NT\times\FG_\NT}{\NR}$, $\map{g}{\FG_\NT}{\NR}$
  \begin{equation}
\label{eq:sumintegralformula}
    \int \mu(\d Z)\int \mu(\d Z')f(Z,Z')g(Z\cup Z')=\int \mu\ast\mu(\d Z'')g(Z'')\int q^\mu(Z'',\d(Z,Z'))f(Z,Z').
  \end{equation}
\end{corollary}
\begin{proof}
  Since $\FG_\NT$ is a Polish space, we can disintegrate the measure $\mu\otimes\mu$  conditional on  $Z\cup Z'$ \cite{Par67}. Since $Z\cup Z'$ has distribution $\mu\ast\mu$, the disintegrating kernel fulfils \ref{eq:sumintegralformula}.  By construction it satisfies the support condition $\mu\ast\mu$-a.s. It is a  standard procedure now  to construct  a version of this disintegrating kernel which fulfils  the support condition everywhere. 
\end{proof}
\begin{example}
  We consider the Poisson process $\Pi_\ell=\e^{-1} F_1$.  From the $\rlap{$\DS\sum$}\DS\int$-lemma \ref{eq:sigmaintegrallemma} we find
  \begin{eqnarray*}
\lefteqn{\int\Pi_\ell(\d Z_1)\int\Pi_\ell(\d Z_2)g(Z_1\cup Z_2)f(Z_1,Z_2)}&=&\e^{-2}\int F_1(\d Z_1)\int F_1(\d Z_2)g(Z_1\cup Z_2)f(Z_1,Z_2)\\
&=&\e^{-2}\int F_1(\d Z)g(Z)\sum_{Z_1\cup Z_2=Z,Z_1\cap Z_2=\emptyset}f(Z_1,Z_2)\\
&=&\e^{-2}\int F_1\ast F_1(\d Z)g(Z) 2^{-\# Z}\sum_{Z_1\cup Z_2=Z,Z_1\cap Z_2=\emptyset}f(Z_1,Z_2)\\    
&=&\int\Pi_\ell\ast \Pi_\ell(\d Z)g(Z) 2^{-\# Z}\sum_{Z_1\cup Z_2=Z,Z_1\cap Z_2=\emptyset}f(Z_1,Z_2)\\    
  \end{eqnarray*}
This shows that $q^{\Pi_\ell}(Z,\p)$ is the uniform distribution on the finite set 
\begin{displaymath}
  \set{{(Z_1,Z_2):Z_1\cup Z_2=Z,Z_1\cap Z_2=\emptyset}}.
\end{displaymath}
It is easy to see that this distribution corresponds exactly to the random  mechanism where  each point in $Z$  chooses independently  from all other points whether to belong to the set $Z_1$ or $Z_2$ with probability $1/2$. It was proven  in \cite{Fic75a} that this kernel $q^{\Pi_\ell}$ characterizes $\Pi_\ell$ among all probability measures on $\FG^f_\NT$.   We want to remark that $\set{\mu:\mu\sim\Pi_\ell}$ is the only stationary factorizing measure type with $q^\mu$ being concentrated on a discrete set. This shows that the kernel  $q^{\Pi_\ell}$ characterizes $\Pi_\ell$ even among all probability measures on $\FG_\NT$.
\end{example}
\begin{corollary}
  \label{cor:tensorproduct no new intrinsic measure type}
 \ref{eq:empty intersection} implies  $\M^{\E^\M\otimes\E^\M,\Us}=\M^{\E^\M,\Us}$ and $\M^{\E^\M\otimes\E^\M,u\otimes u'}=\M^{\E^\M,u}$.\proofend 
\end{corollary}
\begin{remark}
  This result is a simple proof of the fact   that the measure types $\M^{\E,\Us}$ and $\M^{\E,u}$ do not  characterize a type $\mathrm{II}$ product system completely. Another indication is given in Proposition  \ref{prop:typeIIIintotypeII} below. That result shows that it is even not enough to consider the measure types of all product subsystems, how described in Theorem \ref{th:spect}.   
\end{remark}

We want to mention the following simple  description for tensor products of the special product systems $\E^\M$.
\begin{proposition}
  \label{prop:tensorproductto measure types}
Suppose  $\M$ and $\M'$ are stationary factorizing measure types on $\FG_{\NT\times L}$ and $\FG_{\NT\times L'}$ respectively. Define the measure type $\M\otimes\M'$ on $\FG_{\NT\times (L\cup L')}=\FG_{\NT\times L}\times\FG_{\NT\times  L'}$ by
\begin{displaymath}
  \M\otimes\M'=\set{\tilde\mu:\tilde\mu\sim\mu\otimes\mu'\text{~for one and thus all~}\mu\in\M,\mu'\in\M'}.
\end{displaymath}
Then $\E^\M\otimes\E^{\M'}\cong\E^{\M\otimes\M'}$.\proofend
\end{proposition}
\label{sec:II_2neII_0otimesI_2}
At the end of this section, we want show by a simple   example  that the r\^ole of tensor product is, besides   additivity relations for  the numerical index and the measure type $\M^{\E,\Us}$,   a restricted one if one is concerned with the structure of the category of  product systems.   Recall that a product system $\E$ is of type $\mathrm{II}_k$ or $\mathrm{I}_k$ if $\E^\Us\cong\Gammai(\NC^k)$ (set $\NC^\infty=l^2$) and $\E\ne\E^\Us$ or $\E=\E^\Us$ respectively. Before the work of \textsc{Tsirelson} there  was hope that any type $\mathrm{II}_k$ product system is isomorphic to some tensor product $\F\otimes\Gammai(\NC^k)$ where $\F$ is type $\mathrm{II}_0$ \cite{OP:Pow97}.  We show  now  for the  type   $\mathrm{II}_1$ product systems studied  in section \ref{sec:countable} that this is not  true.
\begin{proposition}
  \label{prop:exampleII_2}
Let $\M$ be a  stationary factorizing measure type on $\FG_{\NT}$ with realisations which are countable sets but $Z\notin\FG^f_\NT$  with positive probability. Then   $\E=\E^\M$ is a type $\mathrm{II}_1$ product system which is not isomorphic to some $\F\otimes \Gammai(\NC)$, where $\F$ is any (type $\mathrm{II}_0$) product system. It is even not isomorphic to any nontrivial tensor product  $\F\otimes\tilde\F$, $\F\not\cong\Gammai(\set0)\not\cong\tilde\F$.
\end{proposition}
\begin{proof}
The statement about the type was contained in Corollary \ref{cor:psrstypeII}.

 Assume $\E\cong\F\otimes\tilde\F$ with nontrivial $\F,\tilde\F$ and  let $u$ be the unit in $\E$ corresponding to $\chfc{\set\emptyset}$. From Corollary \ref{cor:unitstensoraretensorofunits} we obtain $u\cong v\otimes\tilde v$ for units $v\in\F$, $\tilde v\in\tilde\F$. Since $\M=\M^{\E,u}=\M^{\F,v}\ast\M^{\tilde\F,\tilde v}$ is concentrated on at most countable sets, both  $\M^{\F,v}$ and $\M^{\tilde F,v}$ are so. At least one of them must have nonfinite realisations, say $\M^{\F,v}$. Then we see analogous to Corollary \ref{cor:psrstypeII} that $\F$ is of type $\mathrm{II}_1$. By additivity of the numerical index, $\tilde \F$ must be either type $\mathrm{II}_0$ or $\mathrm{I}_0$. The former is not possible since by  Corollary \ref{cor:psrstypeII} $\M^{\tilde F,v}$ would have uncountable realisations.   This shows  that $\tilde\F$ is trivial and that $\E$ can be only  a trivial tensor product of product systems.  
\end{proof}
\begin{remark}
  There exist similar  examples for product systems of type $\mathrm{II}_k$, $k\in\NN$, $k\ge2$, derived from measure types on $\FG_{\NT\times\set{1,\dots,k}}$.
\end{remark}
Summarisingly, tensor product do not impose enough structure on the set of product systems (apart from the type $\mathrm{I}$ ones).  As a consequence, we consider other procedures in the next section and in section \ref{sec:direct integrals ps}.

\subsection{The map $\E\mapsto\M^{\E,\Us}$ is surjective}
\label{sec:surjective}

Which  stationary factorizing measure types may arise in the two ways considered in Section \ref{sec:ps2rs}? For $\M^{\E,u}$, Corollary \ref{cor:psrstypeII} shows that  the answer is \emph{all but $\set{\delta_\NT}$}. It is a bit surprising that   the same is true (without the exception) for the measure type  $\M^{\E,\Us}$.
\begin{theorem}
  \label{th:surj}
  The maps $\E\mapsto\M^{\E,\Us}$ and  $(\E,u)\mapsto\M^{\E,u}$ are surjective onto the set of all stationary factorizing measure types and all stationary factorizing measure types different from $\set{\delta_\NT}$ respectively.
\end{theorem}
 For  this result we need to establish that  any quasistationary quasifactorizing random set on $\NT$ is the law of the set of limit points of another quasistationary quasifactorizing random set, so the proof is given on page \pageref{page:proof th:surj} below. This   follows from the more general results of this section, which were already used in section \ref{sec:countable}. Recall  for a stochastic kernel  $q$ on $\FG_\NT$ that $q(Z)$ is a probability measure on $\FG_\NT$. On $[0,1)\cong\NT$, we use the intervals $(s,t)=(s,1)\cup[0,t)$ and $[s,t]=[s,1)\cup[0,t]$, if $s\ge t$, and  extend the notations $\mu_{s,t}$ and $Z_{s,t}$ in the obvious way.
\begin{proposition}
\label{prop:extensionRACS}
Let $q$ be a stochastic kernel from $\FG_\NT$ to $\FG_\NT$ such that
\begin{enumerate}
\item\label{qq1} For all  $s,t\in\NT$  there are two stochastic kernels $q'_{s,t},q''_{t,s}$ on $\FG_{[s,t]}$ and  $\FG_{[t,s]}$ respectively fulfilling   $q(Z)\sim q'_{s,t}(Z_{s,t})\ast q''_{t,s}(Z_{t,s})$  for all  $Z\in\FG_\NT$, $s,t\not\in Z$.
\item\label{qq2} $q(Z+t)\sim q(Z)+t$ for all $t\in\NT$.
\end{enumerate}
Then, for any  stationary factorizing measure type $\M\ne\set{\delta_\NT}$ on $\FG_\NT$ also $\M\circ q=\set{\mu':\mu'\sim\mu\circ q,\mu\in\M}$ is a  stationary factorizing measure type.
\end{proposition}
\begin{proof}
  Clearly, $\tilde\mu\in\M\circ q$ is quasistationary by \ref{qq2}. Further, it is quasifactorizing   iff 
  \begin{equation}
   \label{eq:quasifactinterval}
 \tilde\mu\sim\tilde\mu_{s,t}\ast\tilde\mu_{t,s}\dmf{s,t\in\NT, s\ne t}.
  \end{equation}
To show this, fix $\mu\in\M$ and $s\ne t\in\NT$. If $\mu=\delta_\NT$ equation \ref{eq:quasifactinterval} follows immediately from \ref{qq1}. Otherwise, 
  \begin{displaymath}
    \mu\circ q\sim (\mu_{s,t}\ast\mu_{t,s})\circ q\sim (\mu_{s,t}\ast\mu_{t,s})\circ(q'_{s,t}\ast q''_{t,s})=\mu'_{s,t}\ast\mu''_{t,s}
  \end{displaymath}
for two measures $\mu'$ and $\mu''$ since $s,t\notin Z$ $\mu$-a.s.  We derive from this that $(\mu\circ q)_{s,t}\sim\mu'_{s,t}$ and $(\mu\circ q)_{t,s}\sim\mu''_{t,s}$ what implies \ref{eq:quasifactinterval}.
\end{proof}
Now we develop  a procedure to fill the (many) holes  of a typical quasistationary quasifactorizing random set (different from $\NT$).  
\begin{corollary}
  \label{cor:Zalphabeta}
  Suppose that there are given  a quasistationary quasifactorizing measure $Q_\emptyset$ and a family $(Q_{s,t})_{s,t\in\NT}$ where  $Q_{s,t}$ is a measure probability on $\FG_{[s,t]}$ which fulfil 
  \begin{condition}{Q}
  \item\label{q1}    $Q_{s+r,t+r}=Q_{s,t}+r$ for all $s,t\in\NT$, $r\in\NT$,
  \item\label{q2a} $(Q_{r,t})_{r,s}\sim (Q_{r,t'})_{r,s}$ and $(Q_{r,t})_{s,t}\sim (Q_{r',t})_{s,t}$ for all $r,r',s,t,t'\in\NT$, $s\in(r,t)\cap(r,t')\cap(r',t)\cap(r',t')$,
  \item\label{q2} $Q_{r,t}\sim (Q_{r,t})_{r,s}\ast(Q_{r,t})_{s,t}$ for all $r,t\in\NT$, $s\in (r,t)$,  and
  \item\label{q3} $(Q_{r,t})_{s,s'}\sim(Q_\emptyset)_{s,s'}$ for all $r,s,s',t\in\NT$, $[s,s']\subset (r,t)$.
  \end{condition}
Define  the stochastic kernel $q$ by 
\begin{equation}
\label{eq:infinite convolution holes}
  q(Z)=\left\lbrace
  \begin{array}[c]{cl}
\delta_Z\ast\bigast_{(\alpha,\beta)\subseteq \cmpl Z,\text{maximal}}Q_{\alpha,\beta} &\text{~if $Z\ne\emptyset$}\\
Q_\emptyset&\text{~otherwise}
  \end{array}\right.,
\end{equation}
where the convolution is taken over the family of maximal open subintervals of $\cmpl Z$. 

Then for all stationary factorizing measure types $\M$ on $\FG_\NT$ the measure type  $\M\circ q$ is again stationary and factorizing.
\end{corollary}
\begin{proof}
If $\M=\set{\delta_\NT}$ the statement is clearly fulfilled. So let us assume  $\M\ne\set{\delta_\NT}$.

First we show that the infinite convolution leads to a closed set. Take $t\in\cmplk{Z'}$ if $Z'$ is the sample from $q(Z)$for fixed  $Z$. Then $t\in\cmpl{Z}$, i.e.\ $t\in(\alpha,\beta)$ for some maximal $\alpha,\beta$. Thus $t$ is in the complement of the  random closed set $Z_{\alpha,\beta}=(Z'\cap(\alpha,\beta))\cup\set{\alpha,\beta}$, generated according to  $Q_{\alpha,\beta}$. Since all other sets  $Z_{\alpha',\beta'}$ do not meet the open interval $(\alpha,\beta)$,  there is an open neighbourhood of $t$ not only in $\cmplk{Z_{\alpha,\beta}}$ but in $\cmplk{Z'}$ too. This shows that $Z'$ is closed.

  Now  quasifactorisation  of $Q_\emptyset$ implies \ref{qq1} for $Z=\emptyset$, whereas \ref{q2a} and  \ref{q2} provide that condition in the case where neither $Z_{s,t}$ nor $Z_{t,s}$ is empty. \ref{q3} covers the situation where exactly one of these two sets is empty.  Clearly, \ref{q1} implies \ref{qq2} and the proof is over. 
\end{proof}
\begin{proof}[~of Theorem \ref{th:surj}]
\label{page:proof th:surj}
Clearly, for a unit we have never $P^u_{s,t}=0$ for some $0\le s<t\le1$. This shows in comparison with example \ref{ex:measure type III} that $\set{\delta_\NT}$ is not possible for $\M^{\E,u}$. For each stationary factorizing measure type $\M\ne\set{\delta_\NT}$ Corollary \ref{cor:psrstypeII} shows that there is the product system $\E^\M$ which has at least one unit with $\M^{\E,u}=\M$.

As far as $\M^{\E,\Us}$ is concerned, this measure type is  $\set{\delta_\NT}$ iff $\E$ is of type $\mathrm{III}$ (see Example \ref{ex:measure type III}) and there are examples for product systems of type $\mathrm{III}$, e.g.\ in  \cite{Tsi00}. Thus we may assume again that $\M\ne\set{\delta_\NT}$. In the following we construct a stochastic kernel $q$ on $\FG_\NT$ according to Corollary  \ref{cor:Zalphabeta} such that the measure type $\M^{\E^{\M\circ q},\Us}$, which can be easily calculated  using Corollary \ref{cor:psrstypeII}, is $\M$. This $q$ is universal for all $\M\ne\set{\delta_\NT}$.

 We construct a specific $Q_{s,t}$ by $Q_{s,t}(Y)=Q(\set{Z:\set{s+(t-s)z:z\in Z}\in Y})$, where $Q$ is the law of 
$\set{\xi^l_n:n\in\NN}\cup\set{\xi^r_n:n\in\NN}\cup\set{0,1}$ for the following random variables $\sequ n{\lambda^l}$, $\sequ n{\lambda^r}$.
Choose one point $r_0$ uniformly in $(0,1)$ and two sequences $\sequ n{\lambda^l}$, $\sequ n{\lambda^r}$ uniformly distributed in $[0,1]$. 
These sequences define  via
\begin{eqnarray*}
  \xi^l_0&=&r_0\\
  \xi^l_{n+1}&=&  \lambda^l_n\xi^l_n\dmf{n\in\NN}\\
  \xi^r_0&=&r_0\\
  \xi^r_{n+1}&=&  \lambda^r_n\xi^r_n+(1-\lambda^r_n)\dmf{n\in\NN}
\end{eqnarray*}
two other sequences $\sequ n{\xi^l}$, $\sequ n{\xi^r}$ and consequently $Q$.  Clearly, \ref{q1} is fulfilled  and \ref{q2} is  checked like follows. Define a new enumeration $(x_n)_{n\in\NZ}$ of $\xi$, namely $x_n=\left\{
  \begin{array}[c]{cl}
\xi^r_n&\text{~if~} n\ge0\\
\xi^l_{-n}&\text{~if~} n<0
  \end{array}\right.
$. Then each finite subsequence  of $(x_n)_{n\in\NZ}$ has a law equivalent to Lebesgue measure restricted to a simplex of ordered numbers.  Further, the law of $x$ is invariant under the shift $\sequ nx\mapsto \sequp n{x_{n-1}}$. From this  it is easy to see that $(Q_{r,t})_{r,s}\sim\mathscr{L}(\set r\cup \set{r+(s-r)\xi^r_n:n\in\NN})$. Conditioning on the value of $\max\set{n\in\NZ:x_n\le s}$ yields  \ref{q2}. Choosing $Q_\emptyset=\Pi_\ell$ we get \ref{q3}.
\label{page:countableqdef}
 
Now take any stationary factorizing measure type $\M$ and consider $\E^{\M\circ q}$. We know that  $\M^{\E^{\M\circ q},\Us}$ is the law of $\check {Z'}$, if $Z'$ is distributed according to $\mu\circ q$. But any limit point of $Z'$ is a point of $Z$, since any $Q_{s,t}$ is concentrated on sets with  limit points exactly in $s$ and $t$. Further, any point of $Z$ is a limit point of $Z$ or it is some $\alpha$ for a maximal interval  $(\alpha,\beta)\subset\cmpl Z$. In the latter case, it is a limit point of $Z'$ and we conclude that $\check{Z'}=Z$ almost surely. Corollary \ref{cor:psrstypeII} shows  $\M^{\E^{\M\circ q},\Us}=\M$ and the proof is complete.  
\end{proof}
\begin{remark}
\label{rem:measuretypeasitsowninvariant}
We established the product system $\E^\M$ for all stationary factorizing measure types $\M$ on $\FG_{[0,1]\times L}$.  But, upto now, we could only prove that $\M\circ l^{-1}$,  $l(Z)=\check{Z}$, is an invariant of this product system if $\M$ is a measure type on $\FG_{[0,1]}$ (being $\M^{\E^\M,\Usk}$). One might  conjecture that,  actually, $\E^\M$ is isomorphic to $\E^{\M'}$ iff the lattices $\FG_L$ and $\FG_{L'}$ are isomorphic  and $\M=\M'$ with this identification. For singletons $L,L'$ this is related to  the question whether  automorphisms of product systems act transitively on the set of normalized units. Unfortunately, even if the latter  would be true in general, things are not so simple, as the following  example shows.     
\end{remark}
\begin{example}\label{ex:three Poisson in two}
  Let $Z_1,Z_2,Z_3\in\FG^f_\NT$ be independent random sets distributed according to a Poisson process and define
  \begin{displaymath}
    Z=((Z_1\cup Z_2)\times\set1)\cup((Z_1\cup Z_3)\times\set2)\in\FG^f_{\NT\times\set{1,2}}.
  \end{displaymath}
Then the law  $\mathscr{L}(Z)$ is  stationary and  factorizing since the distribution $\mathscr{L}((Z_1\times\set1)\cup (Z_2\times\set2)\cup (Z_3\times\set3))$ is (see Proposition \ref{prop:extensionRACS}). On the other side, we can recover $Z_1,Z_2,Z_3$ from $Z$ by
\begin{eqnarray*}
  Z_1&=&\pi(Z\cap(\NT\times\set1))\cap\pi(Z\cap(\NT\times\set2))\\
Z_2&=&\pi(Z\cap(\NT\times\set1))\setminus Z_1\\
Z_3&=&\pi(Z\cap(\NT\times\set2))\setminus Z_1,
\end{eqnarray*}
where $\map\pi{\NT\times\set{1,2}}\NT$ is the canonical projection. Thus the product systems corresponding to   $\mathscr{L}(Z)$ and  $\mathscr{L}((Z_1\times\set1)\cup (Z_2\times\set2)\cup (Z_3\times\set3))$ are isomorphic (to $\Gammai(\NC^3)$), although $L$ is $\set{1,2}$ for the first and $\set{1,2,3}$ for the second distribution. 
\end{example}
We are now in a good position to  establish an example of an arbitrarily complex countable quasistationary quasifactorizing random set.
\begin{example}
\label{ex:complexcountable}
  Define a sequence $\sequp n{ \mu^{(n)}}$ of quasistationary quasifactorizing measures on $\FG_\NT$ by
  \begin{displaymath}
   \mu^{(0)}=\delta_\emptyset ,\qquad \mu^{(n+1)}=\mu^{(n)}\circ q\dmf{n\in\NN}
  \end{displaymath}
where $q$ is the stochastic kernel derived in the proof of Theorem \ref{th:surj}. We will see below, in the proof of  Proposition \ref{prop:measuretypeslattice}, that 
\begin{displaymath}
  \mu=\sum_{n\in\NN}2^{-n-1}\mu^{(n)}
\end{displaymath}
is a quasistationary quasifactorizing  measure. It is  concentrated on countable closed sets, since each $\mu^{(n)}$ is so. Moreover, none of the random sets $Z_n$ defined in Lemma \ref{lem:countableRACSnotperfect} is  $\mu$-a.s.\ empty. 

From $\mu$ we derive a family $(Q_{s,t})_{s,t\in\NT}$ fitting into Corollary \ref{cor:Zalphabeta}  as images under the maps $Z\mapsto \set{(t-s)\frac{z}{\max Z}+s:z\in Z}$ (with appropriate meaning for $t\le s$). $Q_\emptyset$ is determined by this choice upto equivalence. Let $Z'$ be a quasistationary quasifactorizing random set and $\tilde Z$ be constructed from $Z'$ by a realisation of the random transition behind the stochastic kernel $q$ from \ref{eq:infinite convolution holes}. Suppose $\#Z'=\infty$ with positive probability and fix such a realisation. Then, by the law of large numbers, for all $z\in Z'$, $n\in\NN$ the set $\tilde Z_n$ in Lemma \ref{lem:countableRACSnotperfect} contains a strictly increasing sequence converging towards $z$. This shows $Z'\subseteq\tilde Z_n$ for all $n\in\NN$. 

Clearly, we can iterate this procedures to advance further and further in the hierarchy of (random) countable closed sets in $\NT$ ordered by inclusion.    
\end{example}
\section{An Hierarchy of Random Sets}
\label{sec:hierarchy}

\subsection{Factorizing Projections and Product Subsystems}
\label{sec:factorizing projections}

Almost all results presented in Section \ref{sec:ps2rs} depended only on the existence of adapted projections $(\mathrm{P}_{s,t})_{(s,t)\in I_{0,1}}$ fulfilling \ref{eq:Prst} and \ref{eq:Pstshift}.  Existence of such projections is  stable  under tensor products:  E.g., if  $\E$ has a unit  but $\F$ has none (i.e., it is of type III) then  $(\mathrm{P}_{s,t}^{\E,\Us}\otimes\unit_{\F_1})_{(s,t)\in I_{0,1}}$  are projections on $\E_1\otimes\F_1$ fitting into Theorem \ref{th:RACS}. But, $\E\otimes\F$ is of type III, see Corollary \ref{cor:XotimesIIIisIII} below. Thus the above developed technique has applications even to product systems of type III. From an abstract point of view,  the crucial notion is that of a \emph{product subsystem}.  Therefore, we will analyse in this section  the whole set $\mathscr{S}(\E)=\set{\F:\text{~$\F$ is a subsystem of $\E$~}}$ which is invariant under isomorphisms.
\begin{remark}
  \label{rem:spectrum}
$\mathscr{S}(\E)$ is never void, since it contains both  $\F=\E$ and $\F$, $\F_t=\set0\subset\E_t$.
\end{remark}
For a product subsystem $\F$ we define projections $(\mathrm{P}_{s,t}^{\F})_{(s,t)\in I_{0,1}}\subset\B(\E_1)$,
\begin{equation}
  \label{eq:defPFst}
  \mathrm{P}_{s,t}^{\F}=\Pr{\E_s\otimes\F_{t-s}\otimes\E_{1-t}}=\unit_{\E_s}\otimes\Pr{\F _{t-s}}\otimes\unit_{\E_{1-t}}\dmf{(s,t)\in I_{0,1}}.
\end{equation}
\vskip-20pt
\begin{proposition}
  \label{prop:Pst=subsystem}
 $\mathscr{S}(\E)$ is  under the map $\F\mapsto(\mathrm{P}^\F_{s,t})_{(s,t)\in I_{0,1}}$ in one-to-one correspondence with the set of families $(\mathrm{P}_{s,t})_{(s,t)\in I_{0,1}}\subset\B(\E_1)$ of adapted  projections fulfilling both \ref{eq:Prst} and \ref{eq:Pstshift}. 

With this identification, the product subsystems $\sg t{\NC u}$ generated by units $\sg tu$ correspond one-to-one to families $(\mathrm{P}_{s,t})_{(s,t)\in I_{0,1}}$ for which  $\mathrm{P}_{0,1}$ is a one-dimensional projection.
\end{proposition}
\begin{proof}
  Clearly, to a projection family $(\mathrm{P}_{s,t})_{(s,t)\in I_{0,1}}$ there correspond subspaces $\F'_{s,t}\subseteq\E_1$. Since 
$\mathrm{P}_{s,t}\in\A_{s,t}$, there is some $\F''_{s,t}\subseteq\E_{t-s}$ such that $\F'_{s,t}=\E_s\otimes \F''_{s,t}\otimes\E_{1-t}$. \ref{eq:Prst} implies $\F''_{r,s}\otimes\F''_{s,t}=\F''_{r,t}$ and \ref{eq:Pstshift} shows that $\F''_{s,t}=\F''_{0,t-s}$. The reverse conclusion is obvious.

For the proof of the second assertion fix vectors $\sgi tv$ with $\mathrm{P}_{0,t}=\Pr{v_t}\otimes \unit_{\E_{1-t}}$. Since $t\mapsto \mathrm{P}_{0,t}$ is continuous by Proposition \ref{prop:singlepoint}, we could choose $\sgi tv$ as measurable section.  Equation \ref{eq:Pstshift} shows $\mathrm{P}_{s,t}=\unit_{\E_s}\otimes\Pr{v_{t-s}}\otimes \unit_{\E_{1-t}}$ for all $s<t$.  \ref{eq:Prst} implies that for $s,t\ge0$, $s+t\le1$ there are complex numbers $c(s,t)\in\NT$ such that $v_{s+t}=c(s,t)v_s\otimes v_t$. An elementary calculation gives 
\begin{displaymath}
    c(s+t,r)c(s,t)=c(s,t+r)c(t,r)\dmf{s,t,r\ge0, s+t+r\le1}.
\end{displaymath}
By Corollary \ref{cor:multipliersimplex}, there are numbers $\sgi tz$ such that $c(s,t)=z_sz_t/z_{s+t}$ for all $s,t\ge0$, $s+t\le1$. Thus $u_t=z_tv_t$ fulfils $u_{s+t}=u_s\otimes u_t$ for all $s,t\ge0$, $s+t\le1$ and this relation   extends trivially to all $s,t\in\NRp$. Applying \cite{OP:Arv89}, we can choose $\sgi tz$ and thus $\sg tu$ measurably. This means that the  latter is a unit and the proof is  complete.
\end{proof}
\begin{corollary}
  \label{cor:nonmeasurable units}
For any section $\sg tu$ through a product system $\E=\sg t\E$ with $u_{t_0}\ne 0$ for some $t_0>0$  and \ref{eq:defunit} there are a unit $v\in\Usk(\E)$ and complex numbers $\sg tc\subset\NC$ with $c_{s+t}=c_sc_t$, $s,t\in\NRp$, such that
\begin{displaymath}
  u_t=c_tv_t\dmf{t\in\NRp}.
\end{displaymath}
\end{corollary}
\begin{proof}
Define the family  $(\mathrm{P}_{s,t}^u)_{(s,t)\in I_{0,1}}$ exactly like for a unit. Clearly, this family  fulfils   \ref{eq:Prst} and \ref{eq:Pstshift} too.   Now the above proposition determines a unit  $v\in\Usk(\E)$ such that $\mathrm{P}_{s,t}^u=\mathrm{P}_{s,t}^v$. This shows $\NC u_t=\NC v_t$ and the proof is complete.
\end{proof}
In this circle of problems, we want to mention the following result proved in   \cite[Corollary 3.9]{OP:Arv89}.
\begin{corollary}
  \label{cor:unitstensoraretensorofunits}
  If $\E,\E'$ are product systems,  a unit $w=\sg tw$ of  $\E\otimes\E'$ has necessarily the form $w_t=u_t\otimes v_t$ for two units $u\in\Usk(\E)$, $v\in\Usk(\E')$.
\end{corollary}
\begin{proof}
  Consider the operators $(\mathrm{P}^w_{s,t})_{(s,t)\in I_{0,1}}\subset\B(\E_1\otimes\E'_1)$. We can build operators $a_{s,t}\in\B(\E_1)$, $(s,t)\in I_{0,1}$,  being the partial trace of $\mathrm{P}^w_{s,t}$ over $\B(\E'_{t-s})$. I.e., $a_{s,t}\in\A^\E_{s,t}$ is uniquely determined by $\Tr b\otimes\unit \mathrm{P}^w_{s,t}=\Tr ba_{s,t}$ for all $b\in\B(\E_{t-s})$. It follows that
  \begin{equation}
\label{eq:arst}
    a_{r,t}=a_{r,s}a_{s,t}\dmf{(r,s),(s,t)\in I_{0,1}}.
  \end{equation}
Now we use  the method of \cite[Corollary 3.9]{OP:Arv89} to show that all  $a_{s,t}$ are  projections and $a_{0,1}$ is one-dimensional. Observe that  the operators $a_{s,t}$ are positive, trace-class and nonzero. Thus, denoting $P'_{s,t}$ the eigenprojection corresponding to  the largest eigenvalue of $a_{s,t}$, we find \ref{eq:Prst} for the family   $(P'_{s,t})_{(s,t)\in I_{0,1}}$. Further, since $a_{0,1}$ is trace-class, $P'_{0,1}$ is finite-dimensional. This means that the product subsystem related to $(P'_{s,t})_{(s,t)\in I_{0,1}}$ consists of finite-dimensional Hilbert spaces. But $\dim P'_{s,t}\dim P'_{t,r}=\dim P'_{s,r}$ shows that  $P'_{0,1}$ is 1-dimensional.  Now denote $P''_{s,t}$ the eigenprojection of the second largest eigenvalue of $a_{s,t}$. These projections fulfil
\begin{displaymath}
  P''_{0,2t}= P'_{0,t}P''_{t,2t} +P''_{0,t}P'_{t,2t}\dmf{t\in[0,1/2]}. 
\end{displaymath}
Thus $\dim P''_{0,2t}=2\dim P''_{0,t}$ or $\dim P''_{0,t}=2^n\dim P''_{0,2^{-n}t}\ge 2^n$ for all $n,t$. Compactness and positivity  of $a_{0,t}$ show that $P''_{0,t}$ has only one  eigenvalue different from zero. Since $P'_{0,1}$ is onedimensional, in fact $a_{s,t}=P'_{s,t}$. Thus $a$ corresponds to a unit and $\mathrm{P}^w_{0,1}$ factorizes to $P'_{0,1}\otimes \tilde P_{0,1}$. Similar arguments show that $\tilde P_{0,1}$ is onedimensional and $w$ is the tensor product of the units corresponding to $P'_{0,1}$ and $\tilde P_{0,1}$.       
\end{proof} 
\begin{corollary}
  \label{cor:XotimesIIIisIII}
  If $\E,\F$ are product systems, where  $\F$ is of type $\mathrm{III}$, $\E\otimes\F$ is of type $\mathrm{III}$ too.\proofend
\end{corollary}
\begin{remark}
  Clearly, \ref{eq:arst} is an extension of \ref{eq:Prst}. Extending Theorem \ref{th:RACS} one could ask, which (random) structures correspond to  this relation. A sufficient answer is given in section \ref{sec:random increments} below. At the moment we just mention that there is again a probability measure behind such families being the distribution of a pair of a random closed set and an increment process.  Unfortunately, the set of positive operator families $(a_{s,t})_{(s,t)\in I_{0,1}}$ obeying  \ref{eq:arst} is not a lattice with respect to the natural order such that this larger (invariant) set  would be much more difficult to use. 
\end{remark}
\begin{theorem}
\label{th:S(E)islattice}
$\mathscr{S}(\E)$ is a complete lattice with respect to the inclusion relation, i.e.\ for any family $(\F^i)_{i\in I}$ of product subsystems there are a least upper bound  product subsystem $\bigvee_{i\in I}\F^i$  and a greatest lower bound $\bigwedge_{i\in I}\F^i$.

This lattice is an invariant of the product system.
\end{theorem}
\begin{proof}
 The fibrewise  intersection of a family of product subsystems is again a product subsystem. Consequently, the greatest lower bound of a  family of product subsystems corresponds to the intersection of the corresponding subspaces in $\E_1$. Further,  the set of all closed subspaces in $\E_1$ is a complete lattice.  Thus the assertion  follows from \cite[Theorem I.6]{S:Bir84} which states that any subset of a complete lattice for which any subsubset has a greatest lower bound is again a complete lattice (with, may be, a different least upper bound functional).

Since any isomorphism of product systems maps product subsystems one-to-one into  product subsystems, $\mathscr{S}(\E)$ is clearly an invariant of $\E$.
\end{proof}
\begin{remark}
  Set $\mathrm{P}^i_{s,t}=\mathrm{P}^{\F^i}_{s,t}$. Following the above proof, we find
  \begin{equation}
\label{eq:maximumP^i}    \bigwedge_{i\in I}\mathrm{P}^i_{r,s}\bigwedge_{i\in I}\mathrm{P}^i_{s,t}=\bigwedge_{i\in I}\mathrm{P}^i_{r,t}\dmf{(r,s),(s,t)\in I_{0,1}},
  \end{equation}
i.e.\   $\bigwedge_{i\in I}\F^i$ corresponds to $(\bigwedge_{i\in I}\mathrm{P}^i_{s,t})_{(s,t\in I_{0,1})}$. 
On the other hand, for  the operators 
\begin{equation}
\label{eq:minimumP^i}
\mathrm{P}^{s_0,\dots,s_N}_{s,t}= \prod_{k=1}^N\bigvee_{i\in I}\mathrm{P}^i_{s_{k-1},s_k}
\end{equation}
the strong  limit $\mathrm{s}-\lim_{s=s_0<s_1<\cdots<s_N=t} \mathrm{P}^{s_0,\dots,s_N}_{s,t}=\mathrm{P}^\infty_{s,t}$ for finer and finer partitions, which exists due to monotony, defines the  least upper bound of $\set{\mathrm{P}^i_{s,t}{(s,t\in I_{0,1})}i\in I}$. 
\end{remark}
\begin{proposition}
  \label{prop:unitinF}
  Let $\E$ be a  product system with at least one unit and $\F\ne \sgp t{\set0}$  be a subsystem of it. Then $\F$ contains a unit too. Especially, if $\sg tu$ is a  unit the section $\sgi tv$ determined  by   $v_t\otimes u_{1-t}=\mathrm{P}^\F_{0,t}u_1$ extends to  a unit.
\end{proposition}
\begin{proof}
  Continuity of $(\mathrm{P}^\F_{s,t})_{(s,t)\in I_{0,1}}$ shows that $(\mathrm{P}^\F_{0,t})_{0\le t\le1}$ is a measurable family of projections. This implies that $\sg tv$ is a measurable section. From \ref{eq:Prst} and \ref{eq:Pstshift} it follows that $v_{s+t}=v_s\otimes v_t$ for all $s,t\ge 0$, $s+t\le1$.

 Assume $v_t=0$ for some $t>0$. Then $v_{t/2^n}=0$ for all $n\in\NN$. This implies on $\E_1$ that $\mathrm{P}^\F_{0,t/2^n}u_1=v_{t/2^n}\otimes u_{1-t/2^n}=0$. On the other side, by Proposition \ref{prop:singlepoint}, $\slim_{n\to\infty}\mathrm{P}^\F_{0,t/2^n}=\unit$. We derive   $\lim_{n\to\infty}\mathrm{P}^\F_{0,t/2^n}u_1=u_1$ which is  a contradiction.  Thus the proof is complete.
\end{proof}
\begin{remark}
 The study of families of projection with \ref{eq:Prst} and \ref{eq:Pstshift}  and their lattice structure originates in work of \textsc{Powers}, cf.\ \cite{OP:Pow99,OP:Bha01} for similar    results  on $\mathscr{S}(\E)$ stated in terms of $E_0$-semigroups and dominated $e_0$-semigroups respectively cocycles of these semigroups.  

 A complete analysis of the lattice $\mathscr{S}(\E)$ for type $\mathrm{I}$ product systems can be found in \cite[section 7, page 65]{OP:Bha01}, see also the next section. Below, we will consider several other examples, built as product systems $\E^\M$ for stationary factorizing measure types $\M$ on $\FG_\NT$. 
\end{remark}

\subsection{Subsystems of $\E^\M$}
\label{sec:subsystems EM}

In the simplest case the  product systems $\E^\M$ are  type $\mathrm{I}$. The  discussion of  $\mathscr{S}(\Gammai(\K))$ is simplified by the following result on units in an arbitrary product system. 
\label{sec:units}
We want to  show that $\Usk(\E)$ carries  a certain affine  structure reflecting that one of the  Hilbert space $\K$ if $\E^\Us\cong\Gammai(\K)$. So consider for two units $u,u'\in\Usk(\E)$, $\lambda\in\NC$ and $t,s_1,\dots,s_k\in\NRp$ with  $s_1+\dots+s_k=t$ (denote this fact by  $(s_1,\dots,s_k)\in\Delta_t$)   the vectors $v^\lambda_{s_1,\dots,s_k}$ defined through
\begin{displaymath}
v^\lambda_{s_1,\dots,s_k}=(\lambda u_{s_1}+(1-\lambda)u'_{s_1} )\otimes\cdots\otimes(\lambda u_{s_k}+(1-\lambda)u'_{s_k})\in\E_t.
\end{displaymath}
Further, we call $(s'_1,\dots,s'_l)\in\Delta_t$  a \emph{refinement} of  $(s_1,\dots,s_k)\in\Delta_t$, if there are natural numbers $l_0=0,l_1,\dots,l_k=l$ such that $s'_{l_i+1}+\dots+s'_{l_{i+1}}=s_{i+1}$ for $i=0,\dots,k-1$. This way, $\Delta_t$ becomes a directed set.  Limit procedures  over this directed set we denote by $(s_1,\dots,s_k)\to\Delta_t$.
\begin{proposition}
\label{prop:affineunit}
  If $u,u'\in\Usk(\E)$  and $\lambda\in\NC$ then there is another unit $v\in\Usk(\E)$ such that 
  \begin{displaymath}
    \lim_{\set{s_1,\dots,s_k}\to\Delta_t}v^\lambda_{s_1,\dots,s_k}=v_t\dmf{t\in\NRp}.
  \end{displaymath}
 Denoting the  covariance function (see \ref{eq:defcovariance}) of $\E$ by $\gamma$ this  unit $v$ fulfils 
\begin{displaymath}
  \gamma(v,w)=\lambda\gamma(u,w)+(1-\lambda)\gamma(u',w)\dmf{w\in\Usk(\E)}.
\end{displaymath}
\end{proposition}
\begin{proof}
  First we consider $\norm{v^\lambda_{s_1,\dots,s_k}-v^\lambda_t}$ for  $(s_1,\dots,s_k)\in\Delta_t$.  Then
\begin{displaymath}
  \norm{v^\lambda_{s_1,\dots,s_k}}^2=\prod_{i=1}^k\norm{\lambda u_{s_i}+(1-\lambda)u'_{s_i}}^2.
\end{displaymath}
From equation   \ref{eq:defcovariance} we get  that $\norm{u'_{s}}^2=1+\gamma(u',u')s+\OO(s^2)$ as well as $\scpro{u_{s}}{u'_{s}}=1+\gamma(u',u')s+\OO(s^2)$. Throughout this proof, $\OO(x)$ denotes terms $T(x)$ with  $\abs{T(x)}\le Cx$ for a universal constant $C$, which depends on $t,u,u'$ and $\lambda$ only. The elementary relation $\prod_{i=1}^k(1+cs_i+\OO(s_i^2))=1+c\sum_{i=1}^ks_i+\OO((\sum_{i=1}^ks_i)^2)$ shows for  $(s_1,\dots,s_k)\in\Delta_s$, 
\begin{eqnarray*}
  \norm{v^\lambda_{s_1,\dots,s_k}}^2&=&\prod_{i=1}^k\big(\absq\lambda(1+\gamma(u,u)s_i+\OO(s^2_i))+2\Re \ovl\lambda(1-\lambda)(1+\gamma(u,u')s_i+\OO(s^2_i))\\&&\qquad+\absq{1-\lambda}(1+\gamma(u',u')s_i+\OO(s^2_i))\big)\\
&=&\prod_{i=1}^k(1+(\absq\lambda\gamma(u,u)+2\Re \ovl\lambda(1-\lambda)\gamma(u',u')\\&&\qquad+\absq{1-\lambda}\gamma(u',u'))s_i+\OO(s^2_i))\\&=&1+(\absq\lambda\gamma(u,u)+2\Re \ovl\lambda(1-\lambda)\gamma(u',u')+\absq{1-\lambda}\gamma(u',u'))t+\OO(s^2).
\end{eqnarray*}
This implies that $\norm{v^\lambda_{s_1,\dots,s_k}}$ is bounded as long as  $(s_1,\dots,s_k)\in\Delta_s$, $s\le t$, and 
\begin{displaymath}
  \norm{v^\lambda_s}^2=1+(\absq\lambda \gamma(u,u)+2\Re\ovl\lambda(1-\lambda)\gamma(u,u')+\absq{1-\lambda}\gamma(u',u'))s+\OO(s^2).
\end{displaymath}
Similarly, we obtain 
\begin{eqnarray*}
  \scpro{v^\lambda_s}{v^\lambda_{s_1,\dots,s_k}}
&=&1+(\absq\lambda \gamma(u,u)+2\Re\ovl\lambda(1-\lambda)\gamma(u,u')+\absq{1-\lambda}\gamma(u',u'))s+\OO(s^2).
\end{eqnarray*}
As a consequence we get  $\norm{v^\lambda_{s_1,\dots,s_k}-v^\lambda_s}=\OO(s^2)$.  Boundedness of $ v^\lambda_{s_1,\dots,s_k}$ implies for any two partitions $(s_1,\dots,s_k),(s'_1,\dots,s'_l)\in\Delta_t$, the latter being finer than the former, that 
\begin{displaymath}
  \norm{v^\lambda_{s_1,\dots,s_k}-v^\lambda_{s'_1,\dots,s'_l}}=\OO(\sum_{i=1}^ks^2_i).
\end{displaymath}
This  proves the asserted convergence. 

In a similar fashion, we derive for any unit $w$ and $(s_1,\dots,s_k)\in\Delta_t$, $t\in[0,T]$,
\begin{displaymath}
  \scpro{v^\lambda_{s_1,\dots,s_k}}{w_t}=1+(\lambda\gamma(u,w)+(1-\lambda)\gamma(u',w))t+\OO(t^2),
\end{displaymath}
what implies the covariance formula and completes the proof.
\end{proof}
\begin{remark}
\textsc{Arveson } \cite[Proposition 6.3]{Arv89} proves  a similar result, showing that 
\begin{displaymath}
  v^n=u_{t/2n}u'_{t/2n}u_{t/2n}u'_{t/2n}\dots u_{t/2n}u'_{t/2n}\in\E_t
\end{displaymath}
  converges weakly to another unit $v$ with the covariance function fulfilling $\gamma(v,w)=\frac12\gamma(u,w)+\frac12\gamma(u',w)$ for all units $w$. This defines  a certain convex structure on $\Usk(\E)$ too.  We have   shown above quite briefly (with establishing  strong convergence) that there is an affine  structure on $\Usk(\E)$. This extended result  could avoid much of the  work in  \cite[Section 6]{Arv89} to  prove that any type $\mathrm{I}$ product system is isomorphic to some $\Gammai(\K)$. We do not repeat  the remaining  arguments from there  here since the same  fact is derived  in  Corollary \ref{cor:exponentialhilbertspace} by computing direct integrals.
\end{remark}
\begin{example}
  Let $\E=\Gammai(\K)$ be a type $\mathrm{I}$ product system. We know that all units $u\in\Usk(\Gammai(\K))$ have the form $u=u^{z,k}$ given in \ref{eq:defuzk}. Then it is easy to derive  $v=u^{z'',k''}$ with $z=\lambda z+(1-\lambda)z'$, $k''=\lambda k+(1-\lambda)k'$.  Thus, if  $\F$ a subsystem of $\Gammai(\K)$, the units in $\F$ correspond to a closed  affine subspace of $\NC\oplus\K$.  Consequently, $\F$ is either $\set 0$ or again type $\mathrm{I}$. Moreover, any closed affine subspace of $\NC\oplus\K$ can arise in this way. Thus  $\mathscr{S}(\E)$  is in one-to-one correspondence with the closed affine subspaces of $\NC\oplus\K$.   
\end{example}
Now  we want to analyse  $\mathscr{S}(\E^\M)$ for a general stationary factorizing measure type $\M$ on $\FG_{[0,1]}$. As a first step, the following lemma yields a vivid description.
\begin{lemma}
  \label{lem:PstinEM}
  Suppose that $(\mathrm{P}_{s,t})_{(s,t)\in I_{0,1}}\subset L^\infty(\M)$ are adapted projections fulfilling 
\ref{eq:Prst} and \ref{eq:Pstshift}. Then there exists a measurable function $\map f{\FG_{[0,1]}}{\FG_{[0,1]}}$ such that $\mathrm{P}_{s,t}$ is the operator of  multiplication by  
\begin{math}
  \chfc{\set{Z:f(Z)\cap[s,t]=\emptyset}}(\cdot)
\end{math}. This function can be chosen to fulfil either $f(Z)\equiv[0,1]$ or $f(Z)\subseteq Z$. Further, $f(Z+t)=f(Z)+t$ will be true  for all $t\in\NR$ for almost all $Z$ and $f(Z)_{s,t}=f(Z_{s,t})$ for all $(s,t)\in I_{0,1}$ for almost all $Z\in\FG_{[0,1]}$. 

Conversely, any such $f$ corresponds to an adapted projection family  with \ref{eq:Prst} and \ref{eq:Pstshift}.
\end{lemma}
\begin{proof}
If $\M=\set{\delta_{[0,1]}}$ there is nothing to prove. Thus, by  Corollary \ref{cor:propmu},  we may assume that $\set\emptyset$ is an atom for $\M$.

   Considering pairs   $(s,t)\in I_{0,1}$ of dyadic numbers we derive random variables $\map{\xi_{k,n}}{\FG_{[0,1]}}{\set{0,1}}$, $n\in\NN$, $k\in\set{0,\dots,2^n-1}$, such that $\mathrm{P}_{k2^{-n},(k+1)2^{-n}}$ is multiplication by $\xi_{k,n}$ and $\xi_{k,n}(Z)=\xi_{k,n}(Z\cap[k2^{-n},(k+1)2^{-n}])$.  We derive that almost surely
  \begin{equation}
\label{eq:dyadic factorization}
    \xi_{2k,n+1}\xi_{2k+1,n+1}=\xi_{k,n}\dmf{n\in\NN,k\in\set{0,\dots,2^n-1}}
  \end{equation}
as well as $\xi_{k,n}(Z+2^{-n})=\xi_{k-1,n}(Z)$. Since $\set\emptyset$ is an atom of $\M$, we derive that 
\begin{eqnarray*}
  \xi_{k,n}(\emptyset)&=&\xi_{2k,n+1}(\emptyset)\xi_{2k+1,n+1}(\emptyset)=\xi_{2k,n+1}(\emptyset)\xi_{2k,n+1}(\emptyset-2^{-n-1})\\
&=&\xi_{2k,n+1}(\emptyset)^2=\xi_{2k+1,n+1}(\emptyset)^2\dmf{n\in\NN,k\in\set{0,\dots,2^n-1}}.
\end{eqnarray*}
This shows that there are only two cases: $\xi_{k,n}(\emptyset)=0$ for all $k,n$ or $\xi_{k,n}(\emptyset)=1$. 

In the former case, we obtain from \ref{eq:dyadic factorization} for almost any $Z\in \FG_{[0,1]}$ such that $[l2^{-m},(l+1)2^{-m}]\cap Z=\emptyset$ for some $l,m$ that $\xi_{k,n}(Z)=0$ if $[l2^{-m},(l+1)2^{-m}]\subseteq[k2^{-n},(k+1)2^{-n}]$. From    Corollary \ref{cor:propmu}, \ref{mu3} we know that for \Maa\ $Z$ for any $k,n$ there exists such $l,m$. We conclude that  $\xi_{k,n}(\emptyset)\equiv0$ implies $\xi_{k,n}\equiv0$ \Mas\ and by continuity, $\mathrm{P}_{s,t}=0$, $(s,t)\in I_{0,1}$. Setting   $f(Z)\equiv[0,1]$ we derive the assertion.

For the rest of the proof, we assume  $\xi_{k,n}(\emptyset)=1$ for all $k,n$.  We change $\xi$ on a null set to obtain the above derived  relations everywhere. Define $\map f{\FG_{[0,1]}}{\FG_{[0,1]}}$ by 
\begin{displaymath}
  f(Z)=\bigcap_{n\in\NN}\bigcup_{\def\arraystretch{0.7}
    \begin{array}[c]{>\scriptstyle c}
k=0,\dots,2^n-1\\\xi_{k,n}(Z)=0
    \end{array}
}[k2^{-n},(k+1)2^{-n}]
\end{displaymath}
Clearly, $n\mapsto \bigcup_{\def\arraystretch{0.5}
    \begin{array}[c]{>\scriptscriptstyle c}
k=0,\dots,2^n-1\\\xi_{k,n}(Z)=0
    \end{array}
}[k2^{-n},(k+1)2^{-n}]$ is decreasing due to \ref{eq:dyadic factorization} such that $f$ maps into $\FG_{[0,1]}$. $f$ is  measurable since finite unions and countable intersections are so. Further, $t\in f(Z)$  iff $\xi_{k,n}(Z)=0$ for all $k,n$ with $t\in[k2^{-n},(k+1)2^{-n}]$. By  definition and the above  compatibility relation for $\xi$ we get  for dyadic $s,t$, say $s=k2^{-n}$, $t=l2^{-n}$, that 
\begin{displaymath}
  f(Z)\cap(s,t)\ne\emptyset\iff\prod_{p=k}^{l-1}\xi_{p,n}(Z)=0.
\end{displaymath}
By definition of $\xi$, this shows that multiplication by \begin{math}
  \chfc{\set{Z:f(Z\cap[s,t]=\emptyset)}}(\cdot)
\end{math} is the projection $\prod_{p=k}^{l-1}\mathrm{P}_{p2^{-n},(p+1)2^{-n}}=\mathrm{P}_{s,t}$. Upper continuity and monotony  of  $(s,t)\mapsto \chfc{\set{Z:f(Z)\cap[s,t]=\emptyset)}} $ and continuity of $(\mathrm{P}_{s,t})_{(s,t)\in I_{0,1}}$ (see Proposition \ref{prop:singlepoint}) imply this fact  for all other pairs $(s,t)$. Since $f$ is almost surely determined by this property  $f(Z)_{s,t}=f(Z_{s,t})$ for all $(s,t)\in I_{0,1}$ for \Maa\ $Z\in\FG_{[0,1]}$ follows from  the fact that $\set\emptyset$ is an $\M$-atom. $f(Z+t)=f(Z)+t$ for  $t\in \NR$ \Mas\ follows from \ref{eq:Pstshift}. Lastly, $f(Z)\cap(k2^{-n},(k+1)2^{-n})\ne\emptyset$ implies $f(Z\cap[k2^{-n},(k+1)2^{-n}])\ne\emptyset$, i.e. $Z\cap[k2^{-n},(k+1)2^{-n}]\ne\emptyset$. Therefore, $f(Z)\subseteq Z$. 

The converse statement is obvious  and the proof is complete.
\end{proof}
\begin{example}
  The simplest example of such local maps is the  map $f(Z)=\check{Z}$, mapping $Z$ into the set of its limit points.   This map appeared  in the analysis of  the projections $(J_{\mathrm{P}^u}^{-1}(\mathrm{P}^\Us_{s,t}))_{(s,t)\in I_{0,1}}\subset L^\infty(\M^{\E,u})$ for some unit $u$ of a  (spatial) product system $\E$, see Proposition \ref{prop:unitandunitalprojections}. It should be easy to prove  that this $f$ is the largest possible  in the sense that  any  such $f$ with $\check{Z}\subseteq f(Z)\subseteq Z $ \Mas\ is either $f(Z)=Z$ or $f(Z)=\check{Z}$. 
\end{example}
\begin{example}
\label{ex:product subsystems for Tsirelson}
If the structure of $Z$ is more complicated more evolved selection criteria for the points of $f(Z)$ arise, are in fact necessary. So    
  consider the special product system  $\E=\E^\M$, where $\M$ is constructed from the zeros of a Brownian motion, see example \ref{ex:zerosBt}. For this product system we know $\E^\Us\cong\Gammai(\set0)$. Thus, Proposition \ref{prop:unitinF} shows for any  product subsystem $\F$ that $\mathrm{P}^\F_{s,t}\ge\mathrm{P}^\Us_{s,t}$, $(s,t)\in I_{0,1}$. I.e.,  $\mathrm{P}^\F_{s,t}$ and $\mathrm{P}^\Us_{s,t}$ commute. Since the latter operators generate a maximal abelian subalgebra of $\B(\E_1)$, there is a function $f$ described in the previous lemma.  Possible candidates for such local maps are the sets of so-called \emph{fast points} like 
  \begin{displaymath}
    f_c(Z)=\set{z\in Z:\limsup_{\varepsilon\downarrow0}\frac{\mathscr{H}^h(Z\cap(z-\varepsilon,z+\varepsilon))}{\sqrt{\varepsilon|\log \varepsilon|}}\ge c}
  \end{displaymath}
for suitable $c>0$, where  $\mathscr{H}^h$ is  the Hausdorff measure corresponding to the function $h(\varepsilon)=\sqrt{\varepsilon\log|\log\varepsilon|}$, see equation \ref{eq:defHausdorff}.  Observe that almost all points of $Z$ fulfil a law of iterated logarithm, these are not contained in $f_c(Z)$. It is known from \cite{C:ST98}, see also \cite{C:OT74,C:Mar99}, that under suitable choices of $c$, $f_c(Z)$ is not almost surely either $Z$ or $\emptyset$. But, this result is not sufficient to establish that $\mathscr{S}(\E^\M)$ consists of more than two elements for we need that $f(Z)$ is closed, too. This is typically not the case for $f_c(Z)$, for most   of the $c$ they are known to be  dense in $Z$ almost surely. So we can only  conjecture that $\mathscr{S}(\E^\M)$ is not trivial, but a lattice with at least one  chain of cardinality of the continuum. Another description of $\mathscr{S}(\E^\M)$ can be  found in Proposition \ref{prop:measuretypeslattice} below, but  that result  does not help  to answer the above  question either.
\end{example}
Another description proves more adequate.
\begin{lemma}
\label{lem:measure types and projections}
   The  families $(\mathrm{P}_{s,t})_{(s,t)\in I_{0,1}}\subset L^\infty(\M)$  of  adapted projections with  \ref{eq:Prst}, \ref{eq:Pstshift}, and $\mathrm{P}_{0,1}\ne0$ stay in one-to-one correspondence with stationary factorizing measure types $\M'$, $\M'\ll\M$, under the map $ \M'\mapsto (\mathrm{P}^{\M'}_{s,t})_{(s,t)\in I_{0,1}}$, 
\begin{displaymath}
 \mathrm{P}^{\M'}_{s,t} = Z\mapsto\chfc{(0,\infty)}(\frac{\d\M'_{s,t}}{\d\M_{s,t}}(Z_{s,t}))\dmf{(s,t)\in I_{0,1}}.
\end{displaymath}
\end{lemma}
\begin{remark}
  Observe, that by the chain rule for  Radon-Nikodym derivatives, $\chfc{(0,\infty)}(\frac{\d\mu'_{s,t}}{\d\mu_{s,t}})$ does  not (essentially) depend on the choice $\mu\in\M$, $\mu'\in\M'$.  
\end{remark}
\begin{proof}
  Assume that  $\M'\ll\M$ are stationary factorizing measure types and  fix $\mu\in\M$, $\mu'\in\M'$. By assumption, $\mu_{s,t}'\ll\mu_{s,t}$  and we can  define a family of random variables
\begin{displaymath}
  \xi_{s,t}(Z)=\chfc{(0,\infty)}(\frac{\d\mu'_{s,t}}{\d\mu_{s,t}}(Z_{s,t}))\dmf{(s,t)\in I_{0,1}}.
\end{displaymath}
Obviously,  $\xi_{s,t}$ does not depend upto $\M$-equivalence on the specific choice of  $\mu,\mu'$ and therefore we obtain for $\M$-a.a.\ $Z\in\FG_\NT$
\begin{equation}
\label{eq:factorization xi}
   \xi_{r,t}(Z)= \xi_{r,s}(Z) \xi_{s,t}(Z)\dmf{(r,s),(s,t)\in I_{0,1}}
\end{equation}
and
\begin{equation}
\label{eq:shift xi}
   \xi_{s+r,t+r}(Z+r)=  \xi_{s,t}(Z)\dmf{(s,t)\in I_{0,1}, r\in\NR}
\end{equation}
This shows that we can associate  $\M'$ with  the system of projections $(\mathrm{P}^{\M'}_{s,t})_{(s,t)\in I_{0,1}}$, corresponding to the $\set{0,1}$ valued random variables $(\xi_{s,t})_{(s,t)\in I_{0,1}}$. Observe that $\xi_{0,1}\mu\sim\mu'$ such that we can recover $\M'$ from $\xi$.

Conversely, let $(\mathrm{P}_{s,t})_{(s,t)\in I_{0,1}}\subset L^\infty(\M)$ be a family of adapted projections with  \ref{eq:Prst}, \ref{eq:Pstshift}, and $\mathrm{P}_{0,1}\ne0$.  Therefore,  we find measurable ($\set{0,1}$ valued) functions $(\xi_{s,t})_{(s,t)\in I_{0,1}}$ such that $\mathrm{P}_{s,t}$ is multiplication by $\xi_{s,t}$. If $\mu\in\M$, we form the measure $\mu'$,
\begin{displaymath}
  \mu'=\frac1{\int \xi_{0,1}\d\mu}\xi_{0,1}\mu.
\end{displaymath}
Observe that $\int \xi_{0,1}\d\mu>0$ since $\mathrm{P}_{0,1}\ne0$. From $\xi_{r,t}=\xi_{r,s}\xi_{s,t}$ $\M$-a.s.\ we get that  $\M'=\set{\mu'':\mu''\sim\mu'}$ is again a stationary factorizing measure type. By definition, $\M'\ll\M$ and $\mathrm{P}^{\M'}_{s,t}=\mathrm{P}_{s,t}$.  Consequently, families $(\mathrm{P}_{s,t})_{(s,t)\in I_{0,1}}\subset L^\infty(\M)$ with  \ref{eq:Prst}, \ref{eq:Pstshift},  and $\mathrm{P}_{0,1}\ne0$, are in one-to-one correspondence with  stationary factorizing measure types $\M'$, $\M'\ll\M$. 
\end{proof}
As a corollary, we can compute $\mathscr{S}(\E^\M)$ for certain stationary factorizing measure types. For the sake of simplicity, let $\mathscr{S}^*(\E^\M)$ denote all of $\mathscr{S}(\E^\M)$ except the zero subsystem.
\begin{corollary}
  \label{cor:S(E)forE^M}
  Let $\M$ be a stationary factorizing measure type on $\FG_\NT$ with $Z=\emptyset$ or $\# Z=\infty$ $\M$-a.s. Then 
  \begin{displaymath}
\mathscr{S}^*(\E^\M)\cong \set{\M':\M'\text{~is a stationary factorizing measure type with~}\M'\ll\M}
  \end{displaymath}
\end{corollary}
\begin{proof}
 Since  $ \M'\mapsto (\mathrm{P}^{\M'}_{s,t})_{(s,t)\in I_{0,1}}\subset L^\infty(\M)\hookrightarrow \B(\E^\M_1)$ is injective, it remains to prove it is onto. Fix any family $(\mathrm{P}_{s,t})_{(s,t)\in I_{0,1}}\subset\B(\E^\M_1)$ of non-zero projections with  \ref{eq:Prst} and \ref{eq:Pstshift}. From Proposition \ref{prop:unitinF} we know that $\mathrm{P}_{s,t}\ge\mathrm{P}^u_{s,t}$ for some unit of $\E^\M$. In Corollary \ref{cor:psrstypeII} we proved that $\E^\M$ has only one family of unital projections $(\mathrm{P}^u_{s,t})_{(s,t)\in I_{0,1}}$, coming from the unit $\sg tu$ associated with $(u_1)_\mu=\mu(\set\emptyset)^{-1/2}\chfc{\set\emptyset}$. Consequently, $(\mathrm{P}_{s,t})_{(s,t)\in I_{0,1}}\subseteq L^\infty(\M)$ and the proof of the above lemma shows that $\mathrm{P}_{s,t}=\mathrm{P}^{\M'}_{s,t}$ for some stationary factorizing measure type $\M'$.  This completes the proof. 
\end{proof}
\begin{remark}
  One should ask for the connection of the two characterisations. This is easily given, setting the random variable $\xi_{0,1}$ in the proof of Lemma \ref{lem:measure types and projections} to
  \begin{displaymath}
    \xi_{0,1}(Z)=\chfc{\set{Z':f(Z')=\emptyset}}(Z),
  \end{displaymath}
i.e.\ $\M'_f=\xi_{0,1}\M$. Again, this is in line with the approach  in Theorem \ref{th:RACS} and Lemma \ref{lem:PstinEM} that a closed set in $[0,1]$ is characterised by the intervals where it is empty. Of course, it is more easy to work instead of the whole $f(Z)$ only with the information where it is empty.   
\end{remark}
\subsection{The Lattice of Stationary Factorizing Measure Types}
\label{sec:MT lattice}
We want to  use Theorem \ref{th:RACS} to get more insight into  the lattice $\mathscr{S}(\E)$ by considering the  map $\F\mapsto\M^\F$ mapping  product subsystems $\F$  into  the  stationary factorizing measure type $\M^\F$ on $\NT$ associated with the family $(\mathrm{P}^\F_{s,t})_{(s,t)\in I_{0,1}}$ by Theorem \ref{th:RACS}. Before  we have to   study a corresponding   concept of order on the stationary factorizing measure types. The previous sections indicates  that the (almost) right one  is induced by  absolute continuity $\ll$. 

Recall that a lattice $S$ is \emph{relatively complete} if any order interval $[e,f]=\set{s\in S:e\le s\le f}$ is a complete lattice, i.e.\ any subset has a greatest lower and a least upper bound. $S$ is  \emph{$\sigma$-complete} if any countable family has a greatest lower and least upper bound. 
\begin{lemma}
  Let $(X,\XG)$ be a measurable space and $\M_1\ll\M_2$ be measure types on it. Then $\set{\M':\text{measure type,~}\M_1\ll\M'\ll\M_2}$ is a complete lattice with respect to $\ll$.
\end{lemma}
\begin{proof}
First  observe that $\ll$ induces a partial order: $\M'\ll\M''$ and $\M''\ll\M'$ imply $\M'\sim\M''$ what is the same as  $\M'=\M''$. 

 For a measure type $\M'$, $\M_1\ll\M'\ll\M_2$ define the measurable function $p_{\M'}$ like in the proof of Lemma \ref{lem:measure types and projections} by   $p_{\M'}=\chfc{(0,\infty)}(\frac{\d\M'}{\d\M_2})$. Then $p_{\M'}$ is $\M_2$-a.s.\ determined. 

Conversely, a $\set {0,1}$ valued measurable function $p\ge p_{\M_1}$  $\M_2$-a.s.\ determines a measure type $\M_1\ll\M'\ll\M_2$ for which $p_{\M'}=p$ $\M_2$-a.s. 

Clearly,  $\M'\ll\M''$ iff  $p_{\M'}\le p_{\M''}$  $\M_2$-a.s. Thus greatest lower bound and least upper bound for measure types $(\M_\alpha)_{\alpha\in A}$ translate into essential greatest lower bound and essential least upper bound for the functions $(p_{\M_\alpha})_{\alpha\in A}$, which exist and preserve the ordering $p_{\M_1}\le p$.
\end{proof}
\begin{proposition}
\label{prop:measuretypeslattice}
The set of all stationary factorizing measure types on $\FG_\NT$ different from $\set{\delta_\NT}$ ordered by $\ll$ is a lattice. It is relatively and $\sigma$-complete and has a least element, $\set{\delta_{\emptyset}}$.  
\end{proposition}
\begin{proof}
  Firstly,  like in the proof of Proposition \ref{prop:convolutionmeasuretypes} we  derive for stationary factorizing measure types $\M_1,\M_2$ that  $\M_1,\M_2\ll\M_1\ast\M_2$. Thus we can dominate any finite  family of stationary factorizing measure types by another one.   

Next we prove that each dominated set of stationary factorizing measure types is a complete lattice. From Corollary \ref{cor:propmu}, \ref{mu2} we know that $\delta_{\emptyset}\ll\M$ for any stationary factorizing measure type $\M\ne\set{\delta_\NT}$. Due to the above lemma and \cite[Theorem I.6]{S:Bir84} it is enough  to prove that for any family $(\M_\alpha)_{\alpha\in A}$ of stationary factorizing measure types  the  greatest lower bound is again   stationary factorizing. This is easily seen from the random variables $(\xi_{s,t})_{(s,t)\in I_{0,1}}$ used in the proof of Lemma \ref{lem:measure types and projections}. The greatest lower bound for the measures corresponds to the essential greatest lower bound of the corresponding  random variables. This new random variable   fulfils again \ref{eq:factorization xi} and \ref{eq:shift xi}. Consequently, each set of projection families associated with  stationary factorizing measure types dominated by a fixed  stationary factorizing measure type $\M$ has a greatest lower and least upper bound. Further, the chain rule for Radon-Nikodym derivatives shows that the measure type associated with these two projection families is independent from the choice of the dominating measure type $\M$. Thus the set of  stationary factorizing measure types on $\FG_\NT$ is a lattice, which is even relatively complete. 

 Clearly, $\set{\delta_{\emptyset}}$ is the least element of this lattice.

Lastly, $\sigma$-completeness of the lattice follows from the fact that  a countable family $\sequp n{\mu^n}$ of quasistationary quasifactorizing probability measures is dominated by another quasistationary quasifactorizing probability measure, e.g.
\begin{displaymath}
\mu=\sum_{n\in\NN}2^{-n}\mu^1\ast\cdots\ast\mu^n.
\end{displaymath}
Clearly, $\mu_i\ll\mu$ and $\mu$ is quasistationary. To prove that $\mu$ is quasifactorizing observe that
\begin{displaymath}
  \mu_{r,s}\ast\mu_{s,t}=\sum_{n,m\in\NN}2^{-n-m}\mu^1_{r,s}\ast\cdots\ast\mu^n_{r,s}\ast\mu^1_{s,t}\ast\cdots\ast\mu^m_{s,t}
\end{displaymath}
Now $\mu^1_{r,s}\ast\cdots\ast\mu^n_{r,s}\ll\mu^1_{r,s}\ast\cdots\ast\mu^{n+1}_{r,s}$ shows that  
\begin{displaymath}
  \sum_{n,m\in\NN}2^{-n-m}\mu^1_{r,s}\ast\cdots\ast\mu^n_{r,s}\ast\mu^1_{s,t}\ast\cdots\ast\mu^m_{s,t}\sim\sum_{n\in\NN}2^{-n}\mu^1_{r,s}\ast\cdots\ast\mu^n_{r,s}\ast\mu^1_{r,s}\ast\cdots\ast\mu^n_{r,s}
\end{displaymath}
and we find 
\begin{eqnarray*}
  \mu_{r,s}\ast\mu_{s,t}&\sim&\sum_{n\in\NN}2^{-n}\mu^1_{r,s}\ast\cdots\ast\mu^n_{r,s}\ast\mu^1_{s,t}\ast\cdots\ast\mu^n_{s,t}\\
&=&\sum_{n\in\NN}2^{-n}\mu^1_{r,s}\ast\mu^1_{s,t}\ast\cdots\ast\mu^n_{r,s}\ast\mu^n_{s,t}\\
&\sim&\sum_{n\in\NN}2^{-n}\mu^1_{r,t}\ast\cdots\ast\mu^n_{r,t}=\mu_{r,t}.
\end{eqnarray*}
This completes the proof.
\end{proof}

As a summary, we derive  the following result, enriching the structure of $\mathscr{S}(\E)$.  We need a new  partial order $\preceq$ on all stationary factorizing measure types on $\FG_\NT$, defined  by $\M'\preceq\M$ iff $\M'\ll\M$ or $\M=\set{\delta_\NT}$.   
\begin{theorem}
  \label{th:spect}
The set of  all stationary factorizing measure types on $\FG_\NT$, ordered by $\preceq$ is a complete lattice with least element $\set{\delta_\emptyset}$ and greatest element $\set{\delta_\NT}$. For any product system $\E$, the map  $m_\E$ on $\mathscr{S}(\E)$, $m_\E(\F)=\M^\F$ provides an  order  antihomomorphism, which is an invariant of  $\E$.  
\end{theorem}
\begin{proof}
The lattice property follows essentially from Proposition \ref{prop:measuretypeslattice} and the observation that any undominated set of  stationary factorizing measure types has $\set{\delta_\NT}$ as least upper bound under $\preceq$.
 For completing the proof, it is enough to show that $\F\subseteq\F'$ implies $\M^\F\gg\M^{\F'}$. This follows readily from Proposition \ref{prop:measuretypeslattice} since $\mathrm{P}^{\F'}_{s,t}\in\set{\mathrm{P}^\F_{s',t'}:(s',t')\in I_{0,1}}''$.
\end{proof}
\begin{remark}
\label{rem:refinementS(E)}
Clearly, the marks adjusted by the map  $m_\E$ to  the elements in  $\mathscr{S}(\E)$ provide more (invariant) information about a product system $\E$.  The only  problem is that $\mathrm{P}^u_{s,t}\wedge\mathrm{P}^v_{s,t}=0 $ for two units $u,v$ which are not multiples of each other. Thus we had to introduce the new order $\preceq$.

On the other hand, we know very little about the image of $\mathscr{S}(\E)$ under $m_\E$. A first step to its characterisation would be the solution of  the question  whether automorphisms of the product system $\E$ act transitively on the normalized units.  An affirmative answer to the latter would imply that all measure types  $\M^{\E,u}$, $u\in\Usk(\E)$, derived in Proposition \ref{prop:psunit2randomset} are equal. I.e., we could not distinguish these product subsystems, which are the minimal ones in   $(\mathscr{S}^*(\E),\subseteq)$. 

Further, we could mark instead of  the \emph{vertices} of the graph of the relation $\subseteq$ in $\mathscr{S}(\E)$ (by the map $m_\E$ from above) the \emph{edges} of that graph, since  by  Theorem \ref{th:RACS} every pair $\F,\F'\in\mathscr{S}(\E)$, $\F\subseteq\F'$, gives rise to a measure type $\M^{\F,\F'}$ too. Again, the relation  between different measure types like  $\M^{\F,\F'},\M^{\F,\F''},\M^{\F',\F''}$ is not simple. This can be  seen in  the examples covered by Proposition \ref{prop:exampleII_2}. There both $\M^{\sg t{\NC u},\E^\Us}$ and $\M^{\E^\Us,\E}$ are equivalent  to a Poisson process, but  $\M^{\sg t{\NC u},\E^\Us}$ may vary among a continuum  of different measure types by Proposition \ref{prop:simplecountable}.
\end{remark}
\section{Direct Integral Representations}
\label{sec:direct integrals}

Since  $\M^{\E,\Us}$ is not a fully characterising  invariant of $\E$ we shall look for further  invariant structures encoded in the projections $(\mathrm{P}^{\E,\Us}_{s,t})_{(s,t)\in I_{0,1}}$. Clearly,  the double commutant algebra $\set{\mathrm{P}^{\E,\Us}_{s,t}:{(s,t)\in I_{0,1}}}''\subset\B(\E_1)$ is an abelian von Neumann algebra. So  we can apply the theory of direct integrals of Hilbert spaces \cite{Sak71,BR87} or the Hahn-Hellinger theorem \cite{Par92}. This makes notations considerably  heavy, but we think the results derived by this technique are worth this effort.

\subsection{Random Sets and Direct Integrals}
\label{sec:direct integrals theory}

 Recall that if $\H=(\H_x)_{x\in X}$ is a measurable family of Hilbert spaces and $\nu$ is a measure on $X$, the direct integral $\int^\oplus\nu(\d x)\H_x$ is defined as the space of square integrable measurable sections. A linear  operator $\map A{\int^\oplus \mu(\d x)\H_x}{\int^\oplus \mu'(\d x)\H'_x}$ on  direct integrals is \emph{diagonal} if it has a representation as $A(\psi_x)_{x\in X}=(A_x\psi_x)_{x\in X}$ for a measurable family of operators $(A_x)_{x\in X}$.   Then we write $A=\int^\oplus\mu(\d x)A_x$.
\begin{proposition}
  \label{prop:HZ}
 Let $\F$ be a product subsystem of $\E$. Then  for all $\mu\in\M^{\F}$ there exists a measurable family $(H_Z)_{Z\in\FG_{[0,1]}}$ of separable Hilbert spaces and a unitary  $\map{U_\mu}{\E_1}{\int^\oplus\mu(\d Z)H_Z}$ such that 
  \begin{displaymath}
   (U_\mu \mathrm{P}^{\F}_{s,t}U_\mu^*\psi)_Z=\chfc{\set{Z':Z'\cap[s,t]=\emptyset}}(Z)\psi_Z.
  \end{displaymath}
If $\mu'\in\M^{\F}$, $(H'_Z)_{Z\in\FG_{[0,1]}}$  and   $\map{U'_{\mu'}}{\E_1}{\int^\oplus\mu'(\d Z)H'_Z}$ have the same property then  there is a diagonal unitary $\map U{\int^\oplus\mu(\d Z)H_Z}{\int^\oplus\mu'(\d Z)H'_Z}$ with  $U U_\mu=U'_{\mu'}$.

Moreover, for all $r\in\NR$ the Hilbert spaces $H_Z$ and $H_{Z+r}$ are isomorphic for $\M^\F$-a.a.\ $Z\in\FG_{[0,1]}$.
\end{proposition}
\begin{proof}
Theorem \ref{th:RACS} yields the normal $\ast$-homomorphism $\map{J_{\mathrm{P}^\F}}{L^\infty(\M^\F)}{\B(\E_1)}$. Applying the Hahn-Hellinger theorem \cite[Theorem 7.6]{Par92} to the abelian von Neumann algebra $J_{\mathrm{P^\F}}(L^\infty(\M^\F))$ we find a unitary $\map{U_0}{\E_1}{L^2(\FG_{[0,1]}\times\NN,\nu)}$ with a probability measure $\nu$ on $\FG_{[0,1]}\times\NN$ such that
\begin{displaymath}
  (U_0 \mathrm{P}^{\F}_{s,t}U_0^*\psi)(Z,n)=\chfc{\set{Z':Z'\cap[s,t]=\emptyset}}(Z)\psi(Z,n).
\end{displaymath}
We want to  show that  under the projection  $\map\pi{\FG_{[0,1]}\times\NN}{\FG_{[0,1]}}$ the measure  $\nu\circ \pi^{-1}$ on $\FG_{[0,1]}$ belongs to  $\M^\F$. Choose a complete orthonormal system $\sequ l\psi$ in $L²(\nu)$. Then $\eta'(\cdot)=\sum_{l\in\NN}2^{-l-1}\scpro{\psi_l}{\cdot\psi_l}$ is a faithful normal state on $\B(L²(\nu))$ and $\eta=\eta'(U_0\cdot U_0^*)$ one on $\B(\E_1)$. We find 
\begin{eqnarray*}
  \eta(\mathrm{P}^{\F}_{s_1,t_1}\cdots \mathrm{P}^{\F}_{s_k,t_k})&=&\eta'(U_0\mathrm{P}^{\F}_{s_1,t_1}\cdots \mathrm{P}^{\F}_{s_k,t_k}U_0^*)\\
&=&\sum_{l\in\NN}2^{-l-1}\scpro{\psi_l}{U_0\mathrm{P}^{\F}_{s_1,t_1}U_0^*\cdots U_0\mathrm{P}^{\F}_{s_k,t_k}U_0^*\psi_l}\\
&=&\sum_{l\in\NN}2^{-l-1}\int\nu(\d Z,\d n)\overline{\psi_l(Z,n)}\prod_{i=1}^k\chfc{\set{Z':Z'\cap[s_i,t_i]=\emptyset}}(Z)\psi_l(Z,n)\\
&=&\int\nu(\d Z,\d n)\sum_{l\in\NN}2^{-l-1}\absq{\psi_l(Z,n)}\chfc{\set{Z':Z'\cap[s_i,t_i]=\emptyset, i=1,\dots,k}}(Z)\\
&=&\nu'\circ \pi^{-1}(\set{Z:Z\cap[s_i,t_i]=\emptyset, i=1,\dots,k})
\end{eqnarray*}
for the measure $\nu'=\sum_{l\in\NN}2^{-l-1}\absq{\psi_l}\nu$. Theorem \ref{th:RACS} yields $\M^\F\ni \mu_\eta=\nu'\circ\pi^{-1}$. Since $\sequ l\psi$ is complete, $\sum_{l\in\NN}2^{-l-1}\absq{\psi_l}>0$ \nuae\ and $\nu'\sim\nu$. This shows $\nu\circ\pi^{-1}\in\M^\F$.

 Disintegration of $\nu$ with respect to $\nu\circ\pi^{-1}$ gives a stochastic kernel $p$ from $\FG_{[0,1]}$ to $\NN$ such that for all \nui\ $f$
 \begin{displaymath}
\int \nu(\d Z,\d n)f(Z,n)=\int  \nu\circ\pi^{-1}(\d Z)\int p(Z,\d n)  f(Z,n).
 \end{displaymath}
We set $H_Z=L^2(\NN,p(Z,\p))$, which is measurable by measurability of $p$. Derive $\map{U_1}{L²(\nu)}{\int^\oplus\nu\circ\pi^{-1}(\d Z)H_Z}$ from 
\begin{displaymath}
  U_1 \psi (Z)=\psi(Z,\cdot).
\end{displaymath}
Since 
\begin{eqnarray*}
  \scpro{U_1\psi'}{U_1\psi}&=&\int  \nu\circ\pi^{-1}(\d Z)\scpro{U_1\psi'(Z)}{U_1\psi(Z)}_{H_Z}\\
&=&\int  \nu\circ\pi^{-1}(\d Z)\int p(Z,\d n)\overline{U_1\psi'(Z)(n)}U_1\psi(Z)(n)\\
&=&\int  \nu\circ\pi^{-1}(\d Z)\int p(Z,\d n)\overline{\psi'(Z,n)}\psi(Z,n)\\
&=&\int  \nu(\d Z,\d n)\overline{\psi'(Z,n)}\psi(Z,n)\\
&=&\scpro{\psi'}{\psi}_{L²(\nu)},
\end{eqnarray*}
$U_1$ is isometric.  Since for any $(\psi_Z)_{Z\in\FG_{[0,1]}}\in\int^\oplus\nu\circ\pi^{-1}(\d Z)H_Z$ the function $\tilde\psi$, 
\begin{math}
  \tilde\psi(Z,n)=\psi_Z(n)\dmf{Z\in\FG_{[0,1]},n\in\NN}
\end{math}
is measurable, in  $L²(\nu)$  and $U_1\tilde\psi=(\psi_Z)_{Z\in\FG_{[0,1]}}$,  $U_1$ is surjective too. Thus $U_{\nu\circ\pi^{-1}}=U_1U_0$ is unitary and the direct integral structure is established for the special $\nu\circ\pi^{-1}$. Observe for any $\mu,\mu'\in\M^\F$ that $U_{\mu,\mu'}$ defined like in \ref{eq:defUmumu'} is a diagonal unitary. Consequently, for all $\mu\in\M^\F$ the operator $U_\mu=U_{\nu\circ\pi^{-1},\mu}U_{\nu\circ\pi^{-1}}$ yields the desired  objects.

Uniqueness is a general result for direct integral decompositions\cite{Par92}.

The isomorphism of  $H_Z$ and $H_{Z+r}$ for $\M^\F$-a.a.\ $Z\in\FG_\NT$ is established by the fact that $\sigma_r$ induces an automorphism on  $J_{\mathrm{P^\F}}(L^\infty(\M^\F))$, implementing $\map{\beta_r}{L^\infty(\M^\F)}{L^\infty(\M^\F)}$,  $\beta_r(f)(Z)=f(Z-r)$. Uniqueness of the direct integral representation completes the proof.
\end{proof}
\begin{corollary}
\label{cor:muetafull}
  If $\E$ is a product system then   $\M^{\E,\Us}=\set{\mu_\eta:\eta\text{~faithful on~}\B(\E_1)}$.
\end{corollary}
\begin{proof}
Fix some $\mu_\eta$ and $\mu\in \M^{\E,\Us}$ and a complete orthonormal system $\sequ lh$ in $L^2(\mu_\eta)$.   From measurability of  $(H_Z)_{Z\in\FG_\NT}$ we get measurable sections $(\psi^n_Z)_{Z\in\FG_\NT}$, $n\in\NN$,  such that $\scpro{\psi^n_Z}{\psi^m_Z}=0$ for all $Z$, if $n\ne m$, and $\sequp n{\psi^n_Z}$ is total in any $H_Z$. Further, we may achieve that
\begin{displaymath}
  \sum_{n\in\NN}\norm{\psi^n_Z}^2_{H_Z}=\frac{\frac{\d \mu}{\d \mu_\eta}(Z)}{\sum_{l\in\NN}2^{-l-1}\absq{h_l(Z)}}\dmf{Z\in\FG_\NT}.
\end{displaymath}
 Now  the sections $(\phi^{l,n})_{{l,n}\in\NN}$, $\phi^{l,n}_Z=h_l(Z)\psi^n_Z$, are total in $\int^\oplus\mu_\eta(\d Z)H_Z$. Namely, let a section $\psi$ fulfil $\scpro\psi{\phi^{l,n}}=0$ for all $l,n\in\NN$. Then $\int \mu_\eta(\d Z)h_l(Z)\scpro{\psi_Z}{\psi^n_Z}_{H_Z}=0$ for all $l,n\in\NN$. Completeness of $\sequ lh$ shows that for all $n\in\NN$ and $\mu_\eta$-a.a.\  $Z\in\FG$, $\scpro{\psi_Z}{\psi^n_Z}_{H_Z}=0$. Totality of $\sequp n{\psi^n_Z}$ in $H_Z$ implies $\psi_Z=0$ $\mu_\eta$-a.s. Consequently, $(\phi_{l,n})_{l,n\in\NN}$ is total in $\int^\oplus\mu_\eta(\d Z)H_Z$ and the functional $\eta'(\cdot)=\sum_{l,n\in\NN}2^{-l-1}\scpro{\phi^{l,n}}{\cdot\phi^{l,n}}$ is a faithful normal state on $\B(\int^\oplus\mu_\eta(\d Z)H_Z)$. 

Further, applying Theorem \ref{th:RACS} to  the faithful normal state  $\eta''(\cdot)=\eta'(U_\mu \cdot U_\mu^*)$ on $\B(\E_1)$ we find
\begin{eqnarray*}
  \lefteqn{ \mu^{\E,\Us}_{\eta''}(\set{Z:Z\cap[s_i,t_i]=\emptyset, i=1,\dots,k})}
&=&\eta''(\mathrm{P}^\Us_{s_1,t_1}\cdots \mathrm{P}^\Us_{s_k,t_k})
=\eta'(U_\mu\mathrm{P}^\Us_{s_1,t_1}\cdots \mathrm{P}^\Us_{s_k,t_k} U_\mu^*)
=\eta'(\chfc{\set{Z:Z\cap[s_i,t_i]=\emptyset, i=1,\dots,k}})\\
&=&\int\mu_\eta(\d Z)\sum_{l,n\in\NN}2^{-l-1}\absq{h_l(Z)}\norm{\psi^n_Z}^2_{H_Z}\chfc{\set{Z':Z'\cap[s_i,t_i]=\emptyset, i=1,\dots,k}}(Z)\\
&=&\int\mu(\d Z)\chfc{\set{Z':Z'\cap[s_i,t_i]=\emptyset, i=1,\dots,k}}(Z)\\
&=&\mu(\set{Z:Z\cap[s_i,t_i]=\emptyset, i=1,\dots,k}).
\end{eqnarray*}
The uniqueness statement of the same  theorem  implies that $\mu^{\E,\Us}_{\eta''}=\mu$.
\end{proof}
\begin{example}
  We consider $\Gammai(\NC^d)$, $d\in\NN^*$.  The relations  $\Gammai(\NC^d)\cong\Gammai(\NC)\otimes\cdots\otimes\Gammai(\NC)$ and $\Gammai(\NC)=\E^{\set{\mu:\mu\sim\Pi_\ell}}$ imply  that $\Gammai(\NC^d)=\E^{\M'}$, where $\M'$ is the measure type on $\FG_{\NT\times\set{1,\dots,d}}$ determined by $\M'=\set{\mu:\mu\sim \Pi_\ell\otimes\cdots\otimes\Pi_\ell}$, see Proposition \ref{prop:tensorproductto measure types}. We choose the unital projections corresponding to $\mathrm{P}^u_{0,1}=\chfc{\set{\emptyset}}$.  Looking on the representation of $\Gamma(L^2([0,1],\NC^d))$ as $L^2(\FG^f_\NT,F_1)$, where $F_1$ was defined in \ref{eq:defexponentialmeasure}, we get $\mathrm{P}^u_{s,t}=\chfc{\set{Z\in\FG_{\NT\times\set{1,\dots,d}}:\pi(Z)\cap[s,t]=\emptyset}}$ where the projection $\pi(Z)$ is give as $\pi(Z)=\set{z\in \NT:\exists k\in\set{1,\dots,d}: (z,k)\in Z}$. Formula  \ref{eq:defexponentialmeasure} shows that the disintegrating kernel $p$ from the proof of Proposition \ref{prop:HZ} is given by
  \begin{displaymath}
    p(\set{t_1,\dots,t_n},\d Z)=d^{-n}\sum_{k_1,\dots,k_n=1}^d\delta_{\set{(t_1,k_1),\dots,(t_n,k_n)}}(\d Z).
  \end{displaymath}
Thus $H_{\set{t_1,\dots,t_n}}=L^2((\frac1d\sum_{k=1}^d\delta_k)^{\otimes n})=(\NC^d)^{\otimes n}$. 
\end{example}
\begin{example}
   Construct $Z$  from any set of probability measures $(Q_{s,t})_{s,t\in \NT,s\ne t}$ used in the proof of Theorem \ref{th:surj} and a random closed set  $Z_0$ distributed  according to the Poisson process $\Pi_\ell$ such  that $\check Z=Z_0$. Let $\M$ be the measure type of $\mathscr{L}(Z)$ and $\E=\E^\M$. Proposition \ref{prop:unitandunitalprojections} shows analogously to Corollary \ref{cor:psrstypeII} that    $\mathrm{P}^\Us_{s,t}=J_{\mathrm{P}^\Us}(\chfc{\set{Z\in\FG_{\NT\times\set{1,\dots,d}}:\check Z\cap[s,t]=\emptyset}})$.  Similar to the above example we derive $H_{\set{t_1,\dots,t_n}}=L^2(Q_{t_1,t_2}\otimes\cdots\otimes Q_{t_n,t_1})$, so all $H_Z$ are infinite dimensional. Thus $\M^\F$ together with the random variables $Z\mapsto\dim H_Z$ does not characterize the product system $\E$ upto equivalence. The additional   associative product structure on $(H_Z)_{Z\in\FG_\NT}$, we derive in the next section, is an essential ingredient. This is the same situation as in product systems. There each fibre is isomorphic to $L^2(\NN)$ but there are nonisomorphic product systems.       
\end{example}

\subsection{Direct Integrals in  Product Systems }
\label{sec:direct integrals ps}

Of course, Proposition \ref{prop:HZ} has an analogue for arbitrary $t$ replacing $1$. Besides  the  direct integral representations for different $t\in\NRp$, there is also the tensor product  structure on $\sg t\E$. In this section we want to analyse the implications of this product structure on the direct integral representations.  

First, we want to approach this structure from  a   constructive point of view. Fix  a stationary factorizing measure type  $\M$ on $\FG_\NRp$ and a measurable  family $H=(H^t_Z)_{t\in\NRp,Z\in\FG_{[0,t]}}$  of Hilbert spaces  with a more complicated product structure. This product structure is defined by unitaries $(V^{s,t}_{Z_s,Z_t})_{s,t\in\NRp,Z_s\in\FG_{[0,s]},Z_t\in\FG_{[0,t]}}$, $\map{V^{s,t}_{Z_s,Z_t}}{H^s_{Z_s}\otimes H^t_{Z_t}}{H^{s+t}_{Z_s\cup (Z_t+s)}}$ which 
\begin{enumerate}
\item are  measurable in $t,Z_t$ with regard to $f(s,Z_s,t,Z_t)=(s+t,Z_s\cup Z_t+s)$ and
\item fulfil the associativity condition 
\begin{equation}
\label{eq:associativity VstZsZt}
 V_{Z_r,Z_s\cup (Z_t+s)}^{r,s+t}\circ(\unit_{H_{Z_r}^r}\otimes V^{s,t}_{Z_s,Z_t})= V_{Z_r\cup( Z_s+r),Z_t}^{r+s,t}\circ(V_{Z_r\cup (Z_s+r)}^{r,s}\otimes\unit_{H_{Z_t}^t})
\end{equation}
for all $r,s,t\in\NRp$ for $\M_p$-a.a.\ $Z_p$, where $p$ varies among $r,s,t$.
\end{enumerate}
Observe that $\sg t{\E^0}$, $\E^0_t=H^t_\emptyset$, is a product system then if $\M $ differs from $\set{\delta_\NT}$. Isomorphisms of this structure are given by a measurable family $(\theta^t_{Z})_{t\in\NRp,Z\in\FG_{[0,t]}}$ of unitaries $\map{\theta^t_{Z}}{H^t_{Z}}{{\tilde H}^t_{Z}}$ with
\begin{displaymath}
  \theta^{s+t}_{Z_s\cup Z_t+s}=\theta^s_{Z_s}\otimes \theta^t_{Z_t}\dmf{s,t\in\NRp,\M^{\F}_{0,p}-\text{a.a.~}Z_p, p=s,t}.
\end{displaymath}

Define the Hilbert spaces 
  \begin{displaymath}
    \E^{\M,H}_t=\set{\psi=(\psi_\mu)_{\mu\in\M_{0,t}}:\psi_\mu\in \int^\oplus\mu(\d Z)H^t_{Z},\psi_{\mu'}=U_{\mu,\mu'}\psi_\mu\forall \mu,\mu'\in\M_{0,t} },
  \end{displaymath}
where the unitaries   $U_{\mu,\mu'}$ are again given by \ref{eq:defUmumu'}, i.e.\
\begin{displaymath}
U_{\mu,\mu'}\psi(Z')=\sqrt{\frac{\d\mu'}{\d\mu\phantom{'}}(Z')}\psi(Z')  \dmf{\psi\in \int^\oplus\mu(\d Z)H^t_{Z},\mu-\text{a.a.~}Z'\in \FG_{[0,1]}}.
\end{displaymath}
Then, the  inner products
\begin{displaymath}
  \scpro{\psi}{\psi'}_{\E^{\M,H}_t}=\int \scpro{\psi_\mu(Z)}{\psi'_\mu(Z)}_{H^t_Z}\mu(\d Z)
\end{displaymath}
do not depend on the choice of $\mu\in\M_{0,t}$ and define an inner product on $\E_t^{\M,H}$.  Multiplication on $\E^{\M,H}$ is given by the unitaries $V_{s,t}$, 
  \begin{equation}
\label{eq:multiplication on EMH}
    (V_{s,t}\psi_s\otimes \psi_t)_{\mu_{0,s}\otimes(\mu'_{0,t}+s)}(Z)=V^{s,t}_{Z_{0,s},Z_{s,s+t}-s}\underbrace{(\psi_s)_{\mu_{0,s}}(Z_{0,s})}_{\in H^s_{Z_{0,s}}}\otimes\underbrace{(\psi_t)_{\mu'_{0,t}}(Z_{s,s+t}-s)}_{\in H^t_{Z_{s,s+t}-s}}\in H^{s+t}_{Z}.
  \end{equation}
The measurable structure is determined by such sections $\psi_s(Z_s)$ for which  for any $\mu\in\M$ the map
\begin{displaymath}
  (s,Z)\mapsto(\psi_s)_{\mu_{0,s}}(Z_{0,s}) 
\end{displaymath}
is almost surely measurable.
\begin{lemma}
  \label{lem:unitalizationistypeII}
$\E=\sg t{\E^{\M,H}}$ is a product system. 

 Assume additionally that  $\M $ is different from $\set{\delta_\NT}$. Then $\E$ has a unit iff $\E^0$ has one.

If $\E^0$ is  type $\mathrm{I}_0$ then 
\begin{displaymath}
  \E^\Us_t=\set{\psi:\psi_\mu(Z)=0 \text{~unless~$\# Z<\infty$}}\dmf{t\in\NRp}.
\end{displaymath}
\end{lemma}
\begin{proof} 
Measurability of $V_{s,t}$ follows from that  of  $(V^{s,t}_{Z_s,Z_t})_{s,t\in\NRp,Z_s\in\FG_{[0,s]},Z_t\in\FG_{[0,t]}}$. In the same manner, we find associativity of the former family by the corresponding property of the latter family. 

If $\M\ne\set{\delta_\NT}$ the family $\sg t{\E^1}$, $\E^1_t=\set{\psi:\psi_\mu(Z)=0 \text{~unless~$Z=\emptyset$}}$ is isomorphic to $\sg t{\E^0}$ and a product subsystem of $\E^{\M,H}$. Thus Proposition \ref{prop:unitinF} shows that  $\E^0$ has a unit if and only if $\E^{\M,H}$ has one.

If $\E^0\cong\Gammai(\set0)$  we  set $(u_t)_\mu(Z)=\mu(\set\emptyset)^{-1/2}\chfc{\set\emptyset}(Z)\in\NC\cong\E^0_t$ and find that $u$ is a unit of $\E^{\M,H}$. It is an easy computation (see  the proof of Corollary \ref{cor:psrstypeII}) that 
\begin{displaymath}
  J_{\mathrm{P}^u}(f)\psi_\mu(Z)=f(Z)\psi_\mu(Z).
\end{displaymath}
Proposition \ref{prop:unitandunitalprojections} completes the proof.
\end{proof}
\begin{theorem}
  \label{th:directintegralps}
  If $\E$ is a product system and $\F$ is a subsystem of it. Then there exists a measurable family of Hilbert spaces $H=(H^t_Z)_{t\in\NRp, Z\in\FG_{[0,t]}}$ such that $\E\cong\E^{\M^\F,H}$ under an isomorphism respecting the natural actions of $J_{\mathrm{P}^\F}(L^\infty(\M^\F))$. 

Further, if $\E$ and $\tilde\E$ are isomorphic under $\sg t\theta$ such that $\theta_t\F_t=\tilde\F_t$ for all $t\in\NRp$ then all $\theta_t$ are isomorphic to  diagonal operators $\theta_t=\int^\oplus\mu_{0,t}(\d Z_t)\theta^t_{Z_t}$ such that $\theta$ is an isomorphism between $H$ and $\tilde H$.

Consequently, $\E$ and $\tilde\E$ are isomorphic iff $\M^{\E,\Us}=\M^{\tilde\E,\Us}$ and the  Hilbert space families  $(H^t_{Z})_{t\in\NRp,Z\in\FG_{[0,t]}}$ and  $(\tilde H^t_{Z})_{t\in\NRp,Z\in\FG_{[0,t]}}$  derived from the subsystems $\E^\Us$ and $\tilde\E^\Us$ are isomorphic.
\end{theorem}
\begin{proof}
Application of Proposition \ref{prop:HZ} yields a family $(H_Z)_{Z\in\FG_{[0,1]}}$ associated with the family $\mathrm{P}^\F$. In a similar way, we can define projections $\mathrm{P}^{\F,t}_{s',t'}\in\B(\E_t)$ through
\begin{displaymath}
 \mathrm{P}_{s',t'}^{\F,t}=\Pr{\E_{s'}\otimes\F_{t'-s'}\otimes\E_{t-t'}}\dmf{(s',t')\in I_{0,t}}. 
\end{displaymath}
The analogue of Proposition \ref{prop:HZ} yields a family $(H^t_Z)_{Z\in\FG_{[0,1]}}$ associated with the family $\mathrm{P}^{\F,t}$ for all $t\in\NRp$. For selecting $H^t_Z$   in a measurable way, with both $t$ and $Z$ free, we will apply  Lemma \ref{lem:measurabilitydirectintegral} below and use further notions introduced for convenience in section \ref{sec:measurability}.  Proposition \ref{prop:singlepoint} shows that  $t\mapsto\mathrm{P}^\F_{ts',tt'}$ is continuous for all $(s',t')\in I_{0,t}''$. Consequently, the family $\sgi t{\pi'}$, $\pi'_t=J_{\mathrm{P}^\F}\circ d_t$ with the dilations $d_tf(Z)=f(\set{z:tz\in Z_{0,t}})$ is a measurable family of representations of $L^\infty(\FG_{[0,1]})$ on $\E_1$. Normality of the representations follows from normality of $J_{\mathrm{P}^\F}$. Fix a measurable family $\sgi t{\hat\eta}$ of normal states on $\B(\E_t)$ each. Then $\pi_t=\unit\otimes\eta_{1-t}\circ \pi'_t$ is again a representation of $L^\infty(\FG_{[0,1]})$ on $\E_t$ and $\sgi t{\pi}$ is measurable.  Now   existence  of a family $H$ with the required measurability is a consequence of Lemma \ref{lem:measurabilitydirectintegral} below.

Denoting $U^t_\mu$ the unitary from Proposition \ref{prop:HZ}, define   unitaries  $\map{V^{s,t}_{Z_{0,s},{Z_{s,s+t}}}}{H^s_{Z_{0,s}}\otimes H^t_{Z_{s,s+t}}}{H^{s+t}_{Z_{0,s}\cup Z_{s,s+t}}}$ by 
  \begin{displaymath}
   U^{s+t}_{\mu_{s+t}}V^{s,t}_{Z_1,Z_2+s}\psi(Z_1\cup Z_2+s)=\sqrt{\frac{\d\mu_s\otimes(\mu_t+s)}{\d \mu_{s+t}}}(Z_1\cup Z_2+s)U^s_{\mu_s}\otimes U^t_{\mu_t}V_{s,t}^{-1}\psi(Z_1, Z_2+s)
  \end{displaymath}
for $\M_{0,s}$-a.a.\ $Z_1\in\FG_{[0,s]}$ and $\M_{0,t}$-a.a.\ $Z_2\in\FG_{[0,t]}$. Since $\E_{s+t}=\E_s\otimes\E_t$ and $\mathrm{P}^\F$ factorizes according to \ref{eq:Prst} the uniqueness part of Proposition \ref{prop:HZ} implies  $H^{s+t}_Z\cong H^s_{Z_{0,s}}\otimes H^t_{Z_{s,s+t}}$ for $\mu_{s+t}$--a.a.\ $Z$ with respect to this unitaries. Naturally, this implies the same factorisation for any $\mu_s,\mu_t,\mu_{s+t}$ from the same measure types which is \ref{eq:multiplication on EMH}. 

By assumption, $\theta_1$ intertwines the families $(\mathrm{P}_{s,t}^\F)_{(s,t)\in I_{0,1}}$ and $(\mathrm{P}_{s,t}^{\F'})_{(s,t)\in I_{0,1}}$. Thus the assertion  on isomorphy follows from  Proposition \ref{prop:HZ}. 

In the special case $\F=\E^\Us$ we know that $\M^{\E,\Us}=\M^{\tilde\E,\Us}$ is necessary for $\E\cong\tilde \E$ and any isomorphism of $\E$ and $\tilde E$ maps $\E^\Us$ onto $\tilde \E^\Us$. The proof is  completed by the preceding results.
\end{proof}
\begin{remark}
  \label{rem:multidimensionaldirectintegral}
  Similarly, there is a multidimensional analogue of this result complementing Proposition \ref{prop:RACSmult}.
\end{remark}
\begin{remark}
On a first sight, the last statement gives us a mean to determine  the whole structure of $\E$. But there are some    obstacles. We cannot provide the relations  $H_Z\cong H_{Z_{0,t}}\otimes H_{Z_{t,1}}$  simultaneously  for all $t\in\NRp$, $Z\in\FG_{[0,t]}$ and we have no direct control over the unitaries encountered in this equivalence. So we replaced the problem of classifying product systems by the more complicated problem of classifying families $(H^t_{Z_t})_{t\in\NRp, Z_t\in\FG_{[0,t]}}$ with much more complicated product structure.  

Nevertheless, the above  result is useful in two directions. Firstly, it gives us a hint to reduce the structure theory of type $\mathrm{III}$ product systems to that of type $\mathrm{II}$. Secondly,  we  can prove now   Corollary \ref{cor:exponentialhilbertspace} which classifies all type $\mathrm{I}$ product systems solely by   results of the present paper, quite differently from the methods used in \cite[section 6]{Arv89}. 
\end{remark}
\subsection{Characterisations of Type I Product Systems}
\label{sec:decomposable}
We want to use the direct integral technique to characterize type $\mathrm{I}$ product systems. 
\begin{corollary}
  \label{cor:exponentialhilbertspace}
  Suppose, $\E$ is a product system with a unit $\sg tu$ for which  the measure type  $\M^{\E,u}$ from Proposition \ref{prop:psunit2randomset} is concentrated on $\FG^f_\NT$. Then $\E$ is an exponential product system $\Gammai(\K)$. Especially, $\ell$-a.a.\ Hilbert spaces $H^1_{\set{t}}$ are isomorphic to $\K$. 
\end{corollary}
\begin{proof}
The proof follows much the lines of that of Proposition \ref{prop:finitepoisson}.
 
  By Corollary  \ref{cor:muetafull} we could assume that $\mu^u_\eta$ is just $\Pi_\ell$. Next we apply Theorem  \ref{th:directintegralps} to the family  $(\mathrm{P}^u_{s,t})_{(s,t\in I_{0,1})}$. In any $\E_T$, there acts a unitary  flip group $\gr t{\tau^T}$,
\begin{equation}
\label{eq:deftauT}
  \tau^T_t x_{T-t}\otimes x_t=x_t\otimes  x_{T-t}\dmf{x_{T-t}\in\E_{T-t}, x_t\in\E_t}.
\end{equation}
Like in Theorem \ref{th:directintegralps} we find that $\tau^T_r$ induces an isomorphism of almost all $H^T_Z$ with $H^T_{Z+r}$. For fixed $k\in\NN\cup\set\infty$ consider the measurable function $\map{f^T_k}{[0,T]}\NRp$, $ f^T_k(t)= \chfc{\set{Z:\dim H^T_Z=k}}(\set t)$. Using Proposition \ref{prop:HZ} we find $f^T_k(t+r)=f^T_k(t)$ for $\ell$-a.a.\ $t$. This shows that $f^T_k\ell$ is a translation invariant measure, i.e.\ $f^T_k$ is constant $\ell$-a.e. Since $f^T_k\cdot f^T_{k'}=0$ $\ell$-a.e.\ for  $k\ne k'\in\NN$ all but one function out of $f^T_k$, $k\in\NN$ vanish  almost surely.  I.e., almost all Hilbert spaces $H^T_{\set t}$ have the same dimension, they are isomorphic to some Hilbert space $\K^T$.  Moreover,  Theorem \ref{th:directintegralps} together with the complete factorisation of $\Pi_\ell$ implies
\begin{displaymath}
  H^{T+S}_{\set t}\cong\left\{
    \begin{array}[c]{cl}
H^T_{\set t}\otimes H^S_\emptyset&\text{~if~$t<T$}\\
 H^T_\emptyset\otimes H^S_{\set t}&\text{~otherwise}
    \end{array}\right.\dmf{\ell-\text{a.a.} t\in[0,T+S]}.
\end{displaymath}
Thus $\K^T\cong\K^S\cong\K$ for all $T,S>0$. Using  \ref{eq:multiplication on EMH} and $H^T_\emptyset=\NC$ we find like in the proof of Proposition \ref{prop:finitepoisson} that  $H^T_{\set{t_1,\dots,t_n}}\cong \K^{\otimes n}$. This  shows $\E_t\cong \Gammai(\K)_t$, since $\mu_\eta=\Pi_\ell$ is concentrated on $\FG^f_{[0,1]}$.
\end{proof}
This result allows us to complete the \label{page:proof unitunitalprojections} 
\begin{proof}[ of Proposition \ref{prop:unitandunitalprojections}]
The measure type  $\M^{u,\Us}$ is well-defined since $\mathrm{P}^u_{s,t}\le\mathrm{P}^\Us_{s,t}$ implies that $\mathrm{P}^u$ and $\mathrm{P}^\Us$ commute. In the following we use that $\FG_{[0,1]\times\set{1,2}}\ni Z\mapsto(\set{t:(t,1)\in Z},\set{t:(t,2)\in Z})\in\FG_{[0,1]}\times\FG_{[0,1]}$ is a Borel isomorphism, and consider instead of $\M^{u,\Us}$ its image. 

We define a new family $(\mathrm{P}_{s,t})_{(s,t)\in I_{0,1}}$ of projections by $\mathrm{P}_{s,t}=J_{\mathrm{P}^u}(\chfc{\set{Z:\#(Z\cap[s,t])<\infty}})$. Observe that $\set{Z\in\FG_{[0,1]}:\#(Z\cap[s,t])<\infty}$ is contained in the $\sigma$-subfield of $\FG_{[0,1]}$ generated by all sets $\set{Z\in\FG_{[0,1]}:Z\cap[s',t']=\emptyset}$, $s\le s'<t'\le t$. From normality of $J_{\mathrm{P}^u}$ we derive that $\mathrm{P}_{s,t}\in\set{\smash{\mathrm{P}^u_{s',t'}}:s\le s'<t'\le t}''\subset\A_{s,t}$, i.e.\ this family is adapted. It fulfils \ref{eq:Prst} since 
\begin{displaymath}
  \# (Z\cap[r,t])<\infty\iff \# (Z\cap[r,s])<\infty\text{~and~}\# (Z\cap[s,t])<\infty\dmf{(r,s),(s,t)\in I_{0,1}}.
\end{displaymath}
By \ref{prop:Pst=subsystem}, there is a product subsystem $\E^0$ of $\E$ such that $\mathrm{P}_{s,t}$ projects onto $\E_s\otimes\E^0_{t-s}\otimes\E_{1-t}$.

Normalize $u$ and consider another normalized unit $v$ and the corresponding normal state $\eta$. We know from Lemma   \ref{lem:unitalmeasurePoisson} that $\mu_\eta$ is the stationary Poisson process on $[0,1]$ with intensity $-2\Re\gamma(u,v)$, which is concentrated on $\FG^f_{[0,1]}$. This implies
\begin{displaymath}
  1=\mu_\eta(\FG^f_{[0,1]})=\eta(\mathrm{P}_{0,1})=\scpro{v_1}{\mathrm{P}_{0,1}v_1}
\end{displaymath}
and $\mathrm{P}_{0,1}v_1=v_1$. Consequently, $\E^\Us$ is contained in $\E^0$.  Further, Corollary \ref{cor:exponentialhilbertspace} shows  that $\E^0$ is isomorphic to some $\Gammai(\K)$.  Since any product system $\Gammai(\K)$ is generated by its units, $\E^0$ is  it too. But  $\E^0$ is a subsystem of $\E$. Thus   all units of $\E^0$ are units of $\E$ and we derive $\E^0=\E^\Us$. This shows $\mathrm{P}^\Us_{s,t}=J_{\mathrm{P}^u}(\chfc{\set{Z\in\FG_{[0,1]}:\# (Z\cap[s,t])<\infty}})$ for all $(s,t)\in I_{0,1}$.

For any faithful normal state $\eta$ we get from $\mathrm{P}^u_{s,t}\le\mathrm{P}^\Us_{s,t}$ and Corollary \ref{cor:open=closed}  that  $\check{Z}\cap(s,t)=\emptyset$ iff $\#(Z\cap[s,t])<\infty$ for $\mu_\eta$-a.a.\ $Z$. Thus, $\mathrm{P}^\Us_{s,t}=J_{\mathrm{P}^u}(\chfc{\set{Z:\#(Z\cap[s,t])<\infty}})$ and  consequently,
\begin{eqnarray*}
 \lefteqn{\mu_\eta(\set{(Z,Z'):Z\cap[s_i,t_i]=\emptyset,Z'\cap[s'_i,t'_i]=\emptyset, i=1,\dots,k})}
 &=& \eta(\mathrm{P}^u_{s_1,t_1}\cdots \mathrm{P}^u_{s_k,t_k}\mathrm{P}^\Us_{s'_1,t'_1}\cdots \mathrm{P}^\Us_{s'_k,t'_k})\\
 &=&\eta(\mathrm{P}^u_{s_1,t_1}\cdots \mathrm{P}^u_{s_k,t_k}J_{\mathrm{P}^u}(\set{Z:\check{Z}\cap[s'_i,t'_i]=\emptyset, i=1,\dots,k}))\\
&=&\mu_\eta(\set{(Z:Z\cap[s_i,t_i]=\emptyset,\check{Z}\cap[s'_i,t'_i]=\emptyset, i=1,\dots,k}).
\end{eqnarray*}
This completes the proof.
\end{proof}
\begin{corollary}
  For  any spatial product system $\E$ there is some Hilbert space $\K$ such that  $\E^\Us\cong\Gammai(\K)$. 
\end{corollary}
\begin{proof}
  This follows from  Proposition \ref{prop:unitandunitalprojections} since for any unit $u\in\Usk(\E)$ the measure type $\M^{\E^\Us,u}$ is concentrated on $\FG^f_{[0,1]}$.  
\end{proof}
Before we collect some concluding remarks, we want to show that another result of \textsc{Arveson} \cite{OP:Arv97,OP:BEF99} can be proved  with the  techniques presented so far.

Call a vector $v\in\E_t$ \emph{decomposable} if for all $0\le s\le t$ there are some vectors $v_s\in\E_s$, $v'_{t-s}\in\E_{t-s}$ with $v=v_s\otimes v'_s$. By Theorem \ref{th:S(E)islattice}, the decomposable vectors generate a product subsystem of $\E$ (may be, the empty one). If this    product subsystem is all of $\E$ it is called \emph{decomposable}.
\begin{corollary}
  \label{cor:decomposable ps is type I}
  Every decomposable product system is of type $\mathrm{I}$.

Consequently, every product system with at least one  decomposable nonzero vector has at least one unit.
\end{corollary}
\begin{proof}
  We use the ideas for the proofs of Theorem \ref{th:ps2randomset}, Proposition \ref{prop:independentmeasure}, Corollary \ref{cor:exponentialhilbertspace} and Theorem  \ref{th:directintegralps}.

Suppose $\E$ is decomposable. Then there is at least one nonzero decomposable vector in $\E_1$, say $v$. If all decomposable vectors are just multiples of $v$, $\E_1$ is one-dimensional and the proof is finished. So let us assume that the set of  decomposable vectors in $\E_1$ is larger. 

To such a $v$ and $s,t\in I_{0,1}$ there correspond vectors $v'_s\in\E_s$, $v_{t-s}\in\E_{t-s}$ and $v''_{1-t}\in\E_{1-t}$ with $v=v'_s\otimes v_{t-s}\otimes v''_{1-t}$. This allows us to define 
\begin{displaymath}
  \mathrm{P}_{s,t}^v=\unit_{\E_s}\otimes\Pr{\NC v_{t-s}}\otimes\unit_{\E_{1-t}}\dmf{(s,t)\in I_{0,1}}.
\end{displaymath}
From the  decomposability of $v$ we obtain  that $(\mathrm{P}^v_{s,t})_{(s,t\in I_{0,1})}$  are adapted projections which fulfil \ref{eq:Prst}. Now  Theorem \ref{th:RACS} provides us for all normal  states $\eta$ on $\B(\E_1)$ with a probability measure $\mu^v_\eta$ on $\FG_\NT$ such that \ref{eq:mueta} is true for $\mathrm{P}_{s,t}=\mathrm{P}^v_{s,t}$. For faithful states $\eta'$, we get even  a faithful representation $J^{\eta'}_{\mathrm{P}^v}$ of $L^\infty(\mu^v_{\eta'})$.

If  $w$ is another normalized   decomposable vector in $\E_1$ set $\eta=\scpro w{\p w}$. Since both $v$ and $w$ factorize, $\mu^v_\eta$  factorizes, i.e.\ it fulfils \ref{eq:strictlyfactorizingmeasure}. Since $\mu^v_\eta$ is not stationary, we cannot derive $\mu^v_\eta(\set{Z:t\in Z})=0$ for all $t\in\NT$ from Corollary \ref{cor:propmu} directly. So suppose there is some $t\in\NT$ such that $\mu^v_\eta(\set{Z:t\in Z})>0$. If $\eta'$ is a faithful normal state on $\B(\E_1)$, Theorem \ref{th:RACS} implies $\mu^v_\eta\ll\mu^v_{\eta'} $. This shows  that $\mu^v_{\eta'}(\set{Z:t\in Z})>0$.  Faithfulness of $\eta'$  implies that $P^{v, \circ}_t< \unit$, where  $P^{v, \circ}_t$ is defined according to \ref{eq:defPcirc}. Proposition \ref{prop:singlepoint} forces  $P^{v}_{0,1}=0$ and equivalently $v=0$, what is a contradiction. Thus  $\mu^v_\eta(\set{Z:t\in Z})=0$ for all $t\in\NT$.  Proceeding like in the proof of Corollary \ref{cor:propmu}, we obtain $\mu^v_\eta(\set\emptyset)>0$. Continuing like in the proof of Proposition \ref{prop:independentmeasure}, we obtain a finite diffuse measure $\lambda^{v,w}$ on $[0,1]$ such that $\lambda^{v,w}([0,t))=-\ln \mu^v_\eta(\set{Z:Z\cap[0,t]=\emptyset})$ and the  Choquet theorem \cite[Theorem 2-2-1]{C:Mat75} implies that $\mu^v_\eta=\Pi_{\lambda^{v,w}}$. Since $\lambda^{v,w}$ is a finite measure, $\Pi_{\lambda^{v,w}}(\FG_\NT^f)=1$. Thus we obtain for a faithful normal states $\eta'$  on $\B(\E_1)$ that 
\begin{displaymath}
  \norm{J^{\eta'}_{\mathrm{P}^v}(\chfc{\FG_\NT^f})w}^2=\eta(J^{\eta'}_{\mathrm{P}^v}(\chfc{\FG_\NT^f}))=\mu^v_\eta(\FG_\NT^f)=\Pi_{\lambda^{v,w}}(\FG_\NT^f)=1,
\end{displaymath}
what implies $J^{\eta'}_{\mathrm{P}^v}(\chfc{\FG_\NT^f})w=w$. Since  $\E$ is decomposable we derive that  $J^{\eta'}_{\mathrm{P}^v}(\chfc{\FG_\NT^f})=\unit$, i.e.\ $\mu^v_{\eta'}(\FG_\NT^f)=1$.   

Using equation \ref{eq:quasifactorizingmeasure} like in the proof of Proposition \ref{prop:finitepoisson} we obtain that $\mu^v_{\eta'}\sim\Pi_{\nu}$ for the finite diffuse  measure $\nu$ on $\NT$ given by $\nu(Y)=\mu^v_{\eta'}(T_1(Y))$, where $T_1(Y)=\set{\set t:t\in Y})$. 

If $(\mathrm{P}^v_{s,t})_{(s,t\in I_{0,1})}$ would fulfil \ref{eq:Pstshift}, we could proceed like in the proof of Corollary \ref{cor:exponentialhilbertspace}.   Thus we will construct in the following a nonzero factorizing vector which is $\gr t\tau$ invariant too (this vector will correspond to a unit then).  

Using the analogue of Theorem \ref{th:directintegralps} we find that $\E_1\cong\int^\oplus_{\FG_\NT^f}\Pi_{\nu}(\d Z)H_Z$ with the identification $\mathrm{P}^v_{s,t}\cong\chfc{\set{Z:Z\cap[s,t]=\emptyset}}$. Proceeding like in the proof of Corollary \ref{cor:exponentialhilbertspace}, we find a measurable family $(H_s)_{s\in\NT}$ of Hilbert spaces such that  $H_{s_1,\dots,s_n}=H_{s_1}\otimes\cdots\otimes H_{s_n}$ for $\nu^n$-a.a.\ $s_1,\dots,s_n$ and $H_\emptyset=\NC$. Moreover, $w$ factorizes iff there is a measurable family   $h=(h_s)_{s\in\NT}$ with $h\in\int^\oplus\nu(\d s)H_s$ and $\alpha\in\NC$ such that $w_{s_1,\dots,s_n}=\alpha h_{s_1}\otimes\cdots\otimes h_{s_n}$  for $\nu^n$-a.a.\ $s_1,\dots,s_n$ and $h_\emptyset=\alpha$. A short calculation shows 
\begin{displaymath}
  \scpro{w}{w'}=\overline{\alpha}\alpha'\e^{\scpro h{h'}}
\end{displaymath}
We assume that $w\ne0$, i.e.\ $\alpha\ne0$.  From the curve $t\mapsto\tau_tw$ we derive  continuous  families $(h^t)_{t\in\NT}\subset\int^\oplus\nu(\d s)H_s$ and $(\alpha^t)_{t\in\NT}\subset\NC\setminus\set0$. In the same manner  let $(k^t)_{t\in\NT}$ and  $(\beta^t)_{t\in\NT}\subset\NC\setminus\set0$ correspond to $t\mapsto\tau_tv$. Using that $\scpro{\tau_tv}{\tau_sw}\ne0$ for all $s,t$ and  $\gr t\tau$ is continuous, we find  continuous versions of $t\mapsto \ln\alpha^t,\ln\beta^t$. The group property of $\gr t\tau$ implies 
\begin{displaymath}
  \overline{\alpha^s}\beta^t\mathrm{e}^{\scpro{h^s}{k^t}}=\mathrm{e}^{\overline{\ln\alpha^s}+\ln\beta^t-\scpro{h^s}{k^t}}=\scpro{\tau_sw}{\tau_tv}=\scpro{w}{\tau_{t-s}v}=\mathrm{e}^{\overline{\ln\alpha^0}+\ln\beta^{t-s}+\scpro{k^0}{h^{t-s}}}.
\end{displaymath}
Now we get from  continuity of $(s,t)\mapsto\scpro{\tau_tv}{\tau_sw}$ and $t\mapsto \ln\alpha^t,\ln\beta^t,(h^t)_{t\in\NT},(k^t)_{t\in\NT}$ that
\begin{displaymath}
{\overline{\ln\alpha^s}+\ln\beta^t-\scpro{h^s}{k^t}}={\overline{\ln\alpha^0}+\ln\beta^{t-s}+\scpro{k^0}{h^{t-s}}}\dmf{s,t\in\NR}.
\end{displaymath}
Let $u$ be the factorizing vector corresponding to  $\gamma=\mathrm{e}^{\int\ell(\d t)\ln\beta^t}$ and $k=\int\ell(\d t)k^t$. We find 
\begin{eqnarray*}
  \scpro{w}{\tau_ru}&=&\scpro{\tau^*_rw} u=\mathrm{e}^{\int\ell(\d t)(\overline{\ln\alpha^{-r}}+\ln\beta^t+\scpro{h^{-r}}{k^t})}\\
&=&\mathrm{e}^{\int\ell(\d t)(\overline{\ln\alpha^0}+\ln\beta^{t+r}+\scpro{h^0}{k^{t+r} })}\\
&=&\mathrm{e}^{\overline{\ln\alpha^0}+\int\ell(\d t)(\ln\beta^t+\scpro{h^0}{k^t})}=\overline{\alpha^0}\gamma\mathrm{e}^{\scpro{h^0}{k}}=\scpro{w}{u}.
\end{eqnarray*}
Since $\E$ is decomposable and $w$ can be chosen freely, this shows $\tau_ru=u$ for all $r\in\NR$. Thus $u$ corresponds to a unit and Corollary \ref{cor:exponentialhilbertspace} completes  the proof.
\end{proof}
\begin{remark}
\label{rem:structuretheorybydirectintegral}
 We see that in principle we can  use the direct integral  structures coming from any set of projections fulfilling \ref{eq:Prst}. E.g., the projection can be related to an arbitrary  product subsystem $\F$ of $\E$.  Then $(\mathrm{P}^u_{s,t})_{(s,t\in I_{0,1})}$ gives a finer structure than $(\mathrm{P}^\Us_{s,t})_{(s,t\in I_{0,1})}$  since $\set{\mathrm{P}^u_{s,t}:s,t\in I_{0,1}}''\subset\set{\mathrm{P}^\Us_{s,t}:s,t\in I_{0,1}}''$. A further advantage of the former is that $H^T_\emptyset\cong\NC$, which allows an effective use of the product structure given by $V^{s,t}_{Z_s,Z_t}$, see the above proofs. On the other side, we do not know whether the direct integral structure  induced by  the former family  is an invariant of the product system (this corresponds to the already mentioned problem whether automorphisms of a product system act transitively on the normalized units). 

 If, fortunately,  the product system is of type $\mathrm{II}_0$, i.e.\ both subsystems coincide, one can say a little bit more about its automorphisms, like the next example shows.
\end{remark}
\begin{example}
  \label{ex:automorphism for Bt0}
  We consider the measure type $\M$ of the zero sets of Brownian motion, introduced in Example \ref{ex:zerosBt1}.  Suppose $\theta=\sg t\theta$ is an automorphism of $\E^\M$. Since  we know that the product system is of type $\mathrm{II}_0$, $\theta$ must leave the projections onto the unit, which are just $\mathrm{P}^\Us_{s,t}$, invariant. By Theorem \ref{th:directintegralps}, $\theta_1$ is a diagonal operator in the direct integral representation associated with $\E^\Us$. But the algebra $\set{\mathrm{P}^\Us_{s,t}:(s,t\in I_{0,1})}''$ is maximal abelian in $\B(\E_1)$ such that almost all fibre spaces are one dimensional in this direct integral space. Thus, $\theta_Z$ is just multiplication by a complex number $z(Z)\in\NT$. Further, the tensor product structure forces 
  \begin{displaymath}
    z(Z)=z(Z_{0,t})z(Z_{t,1})\dmf{0<t<1,\M- \text{a.a.~}Z}
  \end{displaymath}
Examples of such operators can be obtained from  $z(Z)=\mathrm{e}^{\mathrm{i c\mathscr{H}^h(Z)}}$ $\M$-a.s.\ for some $c\in\NR$, where $\mathscr{H}^h$ is the Hausdorff measure with respect to the function $h(\varepsilon)=\sqrt{\varepsilon|\log|\log\varepsilon||}$, see \ref{eq:defHausdorff}. By results of \cite{C:Per81}, the function is well-defined. We do not know whether this set of unitaries exhausts all automorphisms of $\E^\M$.  The unsolved problem to find all product subsystems of $\E^\M$ is seemingly related, see Example \ref{ex:product subsystems for Tsirelson}.
\end{example}
\begin{remark}
By the isomorphy result of Theorem \ref{th:directintegralps},  one can also use the direct integral representations to enrich the lattice $\mathscr{S}(\E)$.   But,  we see no use of this. Conceptionally, it is clearer to consider the lattice of tensor product subsystems of von Neumann algebras of $(\A_{s,t})_{(s,t)\in I_{0,\infty}}$, see Note \ref{rem:productsystemW*-algebras} below.
\end{remark}
\begin{remark}
  In \cite{OP:Pow99} there were constructed $E_0$-semigroups associated with type $\mathrm{II}$ product systems by  domination  of fixed non-unital (so-called $e_0$-)semigroups on $\B(\H)$ which are type $\mathrm{I}$. But the  $e_0$-semigroups dominated by a $E_0$-semigroup correspond to product subsystems, see \cite[Theorem 4.2]{OP:Pow99} and Proposition \ref{prop:Pst=subsystem} above. Thus  we can think of $\mathrm{P}^\F_{0,t}$ as projections onto $H^t_\emptyset$ (which is a Fock space) and the extension amounts in choosing $\M$ and suitable $H^t_Z$. This is similar to the procedure established in Theorem \ref{th:surj}, applied to the measure type of $\Pi_\ell$.
\end{remark}
\begin{remark}
\label{rem:HahnHellingerinvariant}
  Together with the Hahn-Hellinger  theorem there is the notion of spectral type of a representation of an abelian $W^*$-algebra. Namely, we can order the measures $\nu_k(\p)=\nu(\p\cap (\FG_{[0,1]}\times \set k))$ (where $\nu$ is the measure derived from this theorem, see proof of Proposition \ref{prop:HZ}) in such a way  that $\nu_{k+1}\ll\nu_k$, $k\in\NN$. For any $k$, let $\mathrm{mult}(\nu_k)=\#\set{k':\nu_{k'}\sim\nu_k}$ denote the spectral multiplicity of $\nu_k$. Then the set $\set{(\set{\nu:\nu\sim\nu_k},\mathrm{mult}(\nu_k)):k\in\NN}$ is a characteristic  unitary invariant of the representation of $L^2(\M^\F)$ \cite[Theorem 7.6]{Par92}.   Analogous to the fact that $\dim\E_t$ is either $1$ or $\infty$ one can derive that the spectral multiplicity of $\nu_0$ is  either $1$ or $\infty$ too, giving an additional invariant to $\M^{\E,\Us}$, say.  Moreover, this result gives each $\E_t$ the structure of  $L^2(\M'_{0,t})$, $\M'_{0,t}$ being a measure type on $\FG_{[0,t]}\times\NN$. Unfortunately, the relation between  the other spectral multiplicities and measure types is not clear. Thus, we get  no canonical (associative) transformation $\map{c_{s,t}}{(\FG_{[0,s]}\times\NN)\times(\FG_{[0,t]}\times\NN)}{(\FG_{[0,s+t]}\times\NN)}$, such that any type $\mathrm{II}$ product system would fit into Theorem  \ref{th:generalproductsystemfrom stationaryfactorizingmeasuretype} below.
\end{remark}

\subsection{Unitalizing Type III Product Systems}
\label{sec:typeIIItotypeII}
Here we  construct examples of type II product systems parametrised by, besides a stationary factorizing measure type, all possible type III product systems. The main idea is to construct from a product system $\E$ the new product system as $\E^{\M,H}$.  We will consider the special case $H^t_Z=\E_{T(t,Z)}$, where $T(t,Z)\in\NRp$ is a ``new'' (random) time for the old product system.  

For a stationary  factorizing measure type $\M$ on $\FG_\NT$ and  a function $\map h\NRp\NRp$  consider the following condition
\begin{condition*}{\textrm{H}-$h$}
\label{H-h}
  For $\M$-a.a.\ $Z\in\FG_\NT$  either $\mathscr{H}^h(Z)\in(0,\infty)$ or $Z=\emptyset$.
\end{condition*}
Thereby, $\mathscr{H}^h$ is the Hausdorff measure related to $h$ by  \ref{eq:defHausdorff}.
\begin{example}
We may  consider the measure type $\M$ of the set  of   zeros of Brownian motion, see  \cite[section 2]{Tsi00} or the above Example \ref{ex:zerosBt1}. Then, according to  \cite{C:Per81}, $\M$ fulfils \ref{H-h}  for the function  $\map h\varepsilon{\sqrt{\varepsilon\log|\log\varepsilon|}}$.
\end{example}
\begin{proposition}
  \label{prop:typeIIIintotypeII}
  Let $\E$ and $\tilde\E$ be two  product systems and $\M$  a stationary  factorizing measure type $\M$ on $\FG_\NT$ which fulfils \ref{H-h}.  We build two   families $(H^t_Z)_{t\in\NRp, Z\in\FG_{[0,t]}}$, $(\tilde H^t_Z)_{t\in\NRp, Z\in\FG_{[0,t]}}$ of Hilbert spaces by $H^t_Z=\E_{\mathscr{H}^h(Z)}$ and  $\tilde H^t_Z=\tilde\E_{\mathscr{H}^h(Z)}$. Further, we equip $H$ with the multiplication encoded in the unitaries $(V^{s,t}_{Z_s,Z_t})_{s,t\in\NRp,Z_s\in\FG_{[0,s]},Z_t\in\FG_{[0,t]}}$ given by $V^{s,t}_{Z_s,Z_t}=V_{\mathscr{H}^h(Z_s),\mathscr{H}^h(Z_t)}$ and  $\tilde H$ with similarly defined unitaries $(\tilde V^{s,t}_{Z_s,Z_t})_{s,t\in\NRp,Z_s\in\FG_{[0,s]},Z_t\in\FG_{[0,t]}}$.

Then   $\M$ together with both $H$ and $\tilde H $ fulfils the hypothesis of Lemma \ref{lem:unitalizationistypeII}. Moreover,  $\E$ and $\tilde\E$ are isomorphic iff $\E^{\M,H}$ and  $\E^{\M,\tilde H}$ are so.
\end{proposition}
\begin{proof}
 Since $(Z,t)\mapsto\mathscr{H}^h(Z\cap[0,t])$ is measurable, $(H^t_Z)_{Z\in\FG_{[0,t]},t\in\NRp}$ is a measurable family of Hilbert spaces. If we define 
 \begin{math}
   V^{s,t}_{Z_s,Z_t}=V_{\mathscr{H}^h(Z_s),\mathscr{H}^h(Z_t)}
 \end{math},
the associativity relation \ref{eq:associativity VstZsZt} is obvious. Thus both $\E^{\M,H}$ and $\E^{\M, \tilde H}$ exist by  Lemma \ref{lem:unitalizationistypeII}.
 
  The \emph{only if} direction of the isomorphy statement is obvious. 

  For the proof of the  \emph{if} part  observe that the Hausdorff measure of any countable set is $0$. Thus \ref{H-h} implies $\M_{0,t}(\set{Z:0<\#Z<\infty})=0$ for all $t>0$. From $H^t_\emptyset=\E_0=\NC=\tilde\E_0=\tilde H^t_\emptyset$ and   Lemma \ref{lem:unitalizationistypeII} we derive that both  $\E^{\M,H}$ and  $\E^{\M,\tilde H}$ are type $\mathrm{II}_0$, related to the units   
\begin{displaymath}
    (u_t)_{\mu_t}=\mu_t(\set\emptyset)^{-1/2}\chfc{\set\emptyset}=(\tilde u_t)_{\mu_t}\dmf{t\in\NRp,\mu_t\in\M_{0,t}}.
  \end{displaymath}
 Since an isomorphism $\theta$ of $\E^{\M,H}$ and  $\E^{\M,\tilde H}$  fulfils
  \begin{displaymath}
    \theta_1\mathrm{P}^\Us_{s,t}\theta_1^*=\mathrm{P}^{\tilde\Us}_{s,t}\dmf{(s,t)\in I_{0,1}}
  \end{displaymath}
we derive 
  \begin{displaymath}
    \theta_1\mathrm{P}^u_{s,t}\theta_1^*=\mathrm{P}^{\tilde u}_{s,t}\dmf{(s,t)\in I_{0,1}}
  \end{displaymath}
Consequently, by Theorem \ref{th:directintegralps} this   isomorphism is a diagonal operator, mapping fibrewise   $H^t_Z$ into $\tilde H^t_Z$ for almost all $Z$ and tensoring with  respect to the product structure of $H$ and $\tilde H$.  Fix $\mu\in\M_{0,1}$ and regard $\E^{\M,H}_1\cong\int^\oplus\mu(\d Z)H^1_Z$ and similarly $\E^{\M,\tilde H}_1\cong\int^\oplus\mu(\d Z)\tilde H^1_Z$. Then $\theta_1=\int^\oplus\mu(\d Z)\theta^1_{Z}$ where $\map{\theta^t_{Z_t}}{\E^t_{Z_t}}{\tilde\E^t_{Z_t}}$ is $\mu_{0,t}$-a.s.\ unitary. Since $\theta$ is an isomorphism of product systems we find  for all $t\in (0,1)$ and  $\mu$-a.a.\ $Z$  that $\theta^1_Z\cong\theta^t_{Z_{0,t}}\otimes\theta^{1-t}_{Z_{t,1}}$. Consequently,  for $\mu$-a.a.\ $Z$ this factorisation takes place for a countable dense set of $t\in[0,1]$ simultaneously. Conditioning $\mu$ on a fixed value of $\mathscr{H}^h(Z)$  we see  that there is at least one  $T\in\NRp$, such that $T=\mathscr{H}^h(Z)$ for some $Z$ where this factorisation takes place for a countable dense set of $t$ simultaneously.  Since  the Hausdorff measure $\mathscr{H}^h$ is continuous (all singletons have zero measure) the map  $t\mapsto\mathscr{H}^h(Z\cap[0,t])$ is continuous. Thus  there is a unitary $\map{\tilde \theta}{\E_T}{\tilde\E_T}$ such that 
\begin{equation}
  \label{eq:facttildetheta}
  \tilde\theta=\tilde\theta_t\otimes\tilde\theta'_{T-t}
\end{equation}
for a countable dense set of $t\in(0,T)$ and unitaries $\tilde\theta_t,\tilde\theta'_{T-t}$.  Without loss of generality, we may assume that $T=1$ and Corollary \ref{cor:onlyonemeasurablestructure} shows then that  $\E$ and $\E'$ are isomorphic. This completes the proof. 
\end{proof}
\begin{remark}
Clearly,  this is the analogue of the relation between type $\mathrm{III}$   and type $\mathrm{II}$ $W^*$-factors \cite{OP:Tak73}. Somewhat surprisingly, both procedure use direct integral techniques. But observe that different from the classification of factors all constructed product systems of the above kind (there are different choices of $\M$ possible) are type $\mathrm{II}_0$. Consequently, this relation gives us no direct information which is useful  for the classification of type $\mathrm{III}$ product systems. We want also remark that our construction is solely related to product systems, we see no direct construction in terms of related $E_0$-semigroups. 

Summarily, for  completing the classification of product systems one could  concentrate on type $\mathrm{II}$ product systems now.
\end{remark}
\begin{remark}
  This relation to  random (Hausdorff) measures suggests to study, besides stationary factorizing measure types of random closed sets,  stationary factorizing measure types of random measures. That such an approach is useful we will show in the following example.  Moreover, we will show in section  \ref{sec:random measure} that  the latter measure types  lead again to product systems, but less is known about their structure at the moment. 
\end{remark}
\begin{example}
  We follow the above construction for general $\M$ with \ref{H-h} and $\E$ to be the type $\mathrm{I}_1$ product system, i.e.\ $\E_t=L^2(\FG_{[0,t]},\Pi_{\ell})$. Then $H^t_Z=L^2(\FG_{[0,\mathscr{H}^h(Z\cap[0,t])]},\Pi_{\ell})$, but additivity and continuity of the Hausdorff measure show that we may equivalently  choose $H^t_Z=L^2(\FG_{[0,t]},\Pi_{\nu^Z})$, where the measure $\nu^Z$ is defined through $\nu^Z(Y)=\mathscr{H}^h(Z\cap Y)$ and the Poisson process $\Pi_\nu{}$  with  intensity measure $\nu$ is defined in  \ref{eq:defPoisson}. This shows immediately that $\E^{\M,H}$ is isomorphic to $\E^{\M^2}$, where $\M^2$ is the measure type of a random closed set in $[0,1]\times\set{1,2}$, i.e.\ a random pair $(Z_1,Z_2)$, $Z_{1,2}\in\FG_{[0,1]}$, where $\mathscr{L}(Z_1)\in\M$ and $Z_2$, conditionally on $Z_1$,  is distributed according to $\Pi_{\nu^Z}$. Thus, almost surely, $Z_2$ is a finite set contained in $Z_1$, the latter being either empty or uncountable.
\end{example}
\section{Measurability in Product Systems: An  Algebraic Approach}
\label{sec:algebraic approach}
In this section we want to study the r\^ole of the measurability condition in Definition \ref{def:ps}. It will  result that it is more or less  an intrinsic property of the multiplication encoded in the unitaries $(V_{s,t})_{s,t\in\NRp}$. Although the derivations  are quite  technical, we see it as a big  achievement to clarify the secondary r\^ole  of  the measurability condition in the axioms of product systems.

\subsection{GNS-representations}

\label{sec:B(H)representation}
\label{GNS-representation}
The key idea used in this section is the change from  the Hilbert space point of view (looking at $\E_t$ etc.\ directly) to an algebraic one. So  we will focus on the algebras $\A_{s,t}\subseteq\B(\E_1)$ introduced in \ref{eq:defAst}. The relation between these  two aspects is given by the Gelfand-Neimark-Segal (GNS-) representation. We will explain the main points of this representation now, for more details see, e.g., \cite[section 2.3.3]{BR87}.

Let $\eta$ be a  state on a $C^*$-algebra $\A$. I.e.\ $\A$ is a Banach algebra with norm fulfilling $\norm{a^*a}=\norm{a}^2$ and $\eta$ a  norm continuous  functional with $\eta(a^*a)\ge0$ and $\eta(\unit)=1$. On $\A$ we can define an semi-inner product $\scpro ab_\eta=\eta(a^*b)$. The kernel $\mathcal{K}$ of $\scpro\cdot\cdot_\eta$ is a linear subspace (in fact a right ideal in $\A$), dividing it out and completing $\A/\mathcal{K}$ with respect to the respective  inner product determines a Hilbert space $H_\eta$. $\A$ has a representation $\pi_\eta$ on $H_\eta$ determined through $\pi_\eta(a)[b]_\eta=[ab]_\eta$ where we use the convention $[b]_\eta=b+\mathcal{K}$. 

As a special case, this scheme applies to  $W^*$-algebras   equipped with a normal state. $W^*$-algebras, the abstract version of von Neumann-algebras, are  $C^*$-algebras which are the dual of a Banach space, its predual. A normal state on a  $W^*$-algebra is a state which is an element of the predual of the algebra too. In the present paper, we are mainly interested in the following  two examples.
\begin{example}
\label{ex:GNSnormalstateB(H)}
  Let $\eta$ be a normal state on $\B(\H)$. Then 
  \begin{eqnarray*}
   \scpro{[a]_\eta}{[b]_\eta}_{H_\eta}&=& \eta(a^*b)=\Tr a^*b\varrho=\Tr \varrho^{1/2}a^*b\varrho^{1/2}\\
&=&\Tr (a\varrho^{1/2})^*b\varrho^{1/2}=\scpro{a\varrho^{1/2}}{b\varrho^{1/2}}_{\B^2(\H)}
  \end{eqnarray*}
where $\B^2(\H)$ is the Hilbert space of Hilbert-Schmidt operators on $\H$ (i.e.\ for $b\in\B^2(\H)$ the operator $b^*b$ is trace-class) equipped with the scalar product $\scpro ab_{\B^2(\H)}=\Tr a^*b$. This shows that $[a]_\eta\mapsto a\varrho^{1/2}$  defines an isometric embedding of $H_\eta$ into $\B^2(\H)$. 

If $\eta$ is faithful this map is onto and $H_\eta\cong\B^2(\H)$. Denote $\H^*$ the adjoint Hilbert space of $\H$ which is $\H$ as a set with the same addition but scalar multiplication $\lambda\cdot\psi=\overline{\lambda}\psi$ and inner product $\scpro{\psi}{\psi'}_{\H^*}=\scpro{\psi'}{\psi}$. Then $\H\otimes\H^*\cong\B^2(\H)$ under the unitary map $\psi\otimes\psi'\mapsto\lvert\psi\rangle\langle\psi'\rvert=\scpro{\psi'}\cdot\psi$. Further, by the Riesz representation theorem $\H^*\cong\H$ such that   $H_\eta\cong\H\otimes\H$. Observe that this unitary relation is not any more intrinsic but depends on the choice of a conjugation realizing  $\H^*\cong\H$.  

Denote $\P(\H)$ the set of pure normal  states on $\B(\H)$, i.e.\ $\eta\in\P$ is given by $\eta(a)=\scpro\psi{a\psi}$ for a unit vector $\psi\in\H$ or $\varrho=\Pr\psi$.  We see $[a]_\eta=0$ if and only if $a\psi=0$ and  $\scpro{[a]_\eta}{[b]_\eta}_{H_\eta}=\scpro{a\psi}{b\psi}_\H$. Thus $[a]_\eta\mapsto a\psi$ is an isometry from $H_\eta$ into $\H$. The choice $a=\Pr{\psi',\psi}=\scpro\psi\cdot\psi'$ shows that this map is onto and therefore a unitary. We conclude that $H_\eta\cong\H$, a key observation for the following construction. Since there is never a canonical choice of $\eta$ nor $\psi$, it is nice that we can  connect the different spaces $(H_\eta)_{\eta\in\P(\H)}$ isomorphic to $\H$ in a natural manner.
\end{example}
\begin{lemma}
  \label{lem:well-defined Uetaeta'}
  The map $[a]_\eta\mapsto[a\Pr{\psi,\psi'}]_{\eta'}$ is well-defined and extends to a unitary $\map{U_{\eta,\eta'}}{H_\eta}{H_{\eta'}}$.
\end{lemma}
\begin{proof}
  Clearly, the map is onto. Further,
  \begin{eqnarray*}
    \scpro{U_{\eta,\eta'}[a]_\eta}{U_{\eta,\eta'}[b]_\eta}_{H_{\eta'}}&=& \scpro{[a\Pr{\psi,\psi'}]_{\eta'}}{[b\Pr{\psi,\psi'}]_{\eta'}}_{H_{\eta'}}\\
&=&\scpro{a\Pr{\psi,\psi'}\psi'}{b\Pr{\psi,\psi'}\psi'}_\H=\scpro{a\psi}{b\psi}_\H=\scpro{[a]_\eta}{[b]_\eta}_{H_\eta}.
  \end{eqnarray*}
This shows both that $U_{\eta,\eta'}$ respects the kernel of $\scpro\cdot\cdot_\eta$ and that it is isometric. This completes the proof.
\end{proof}
\begin{example}
\label{ex:GNSLinfty}  
Let $\M$ be a measure type on $(X,\XG)$ and $\A=L^\infty(\M)$. Any normal state $\eta$ on $\A$ is given by $\eta(f)=\int f\d\tilde\mu$, where $\tilde \mu$ is a probability measure  with  $\tilde\mu\ll\M$. Since $\scpro{[f]_\eta}{[g]_\eta}_{H_\eta}=
\int\ovl f g\d\tilde\mu=\scpro fg_{L^2(\tilde\mu)}$, we see $H_\eta\cong L^2(\tilde\mu)$. 
\end{example}
Next we want to extend  the construction of $L^2(\M)$ (see \ref{eq:defL2measuretype}), which connects the Hilbert spaces $L^2(\mu)$, $\mu\in\M$, to  families  $(H_\eta)_{\eta\in\S}$ of GNS-representation  spaces. Thereby, $U_{\eta,\eta'}$ replaces $U_{\mu,\mu'}$ defined in \ref{eq:defUmumu'}. The following result is easy to prove.
\begin{lemma}
\label{lem:defH(S)}
  Suppose $\S$ is a set of  states on a $C^*$-algebra $\A$ and $\U=(U_{\eta,\eta'})_{\eta,\eta'\in\S}$ a family of unitary operators $\map{U_{\eta,\eta'}}{H_\eta}{H_{\eta'}}$,  $\eta,\eta'\in\S$  such that $U_{\eta,\eta}=\unit_{H_\eta}$, $\eta\in\S$ and 
  \begin{equation}
\label{eq:productUetaeta'}
    U_{\eta',\eta''}U_{\eta,\eta'}=U_{\eta,\eta''}\dmf{\eta,\eta',\eta''\in\S}.
  \end{equation}
Then the linear space
\begin{equation}
\label{eq:defH(R)}
  H(\A,\S,\U)=\set{(\psi_\eta)_{\eta\in\S}:\psi_\eta\in H_\eta\forall\eta\in\S, \psi_{\eta'}=U_{\eta,\eta'}\psi_\eta\forall\eta,\eta'\in\S}
\end{equation}
equipped with the inner product
\begin{displaymath}
  \scpro\psi{\psi'}_{H(\A,\S,\U)}=\scpro{\psi_\eta}{\psi'_\eta}_{H_\eta}\dmf{\psi,\psi'\in H(\A,\S,\U)},
\end{displaymath}
which is independent of  the choice of $\eta$, is a Hilbert space isomorphic to each $H_\eta$, $\eta\in\S$.

If for all $\eta,\eta'\in\S$ the map $\map{U_{\eta,\eta'}}{H_\eta}{H_{\eta'}}$ intertwines the GNS-representations $\pi_\eta$ and $\pi_{\eta'}$:
  \begin{displaymath}
   U_{\eta,\eta'}\pi_\eta(a) =\pi_{\eta'}(a)U_{\eta,\eta'}\dmf{a\in\A}.
  \end{displaymath} 
then $H(\A,\S,\U)$ carries the representation $\pi$ of $\A$, 
\begin{displaymath}
  \pi(a)(\psi_\eta)_{\eta\in\S}=(\pi_\eta(a)\psi_\eta)_{\eta\in\S}\dmf{a\in\A}.\proofend
\end{displaymath}
\end{lemma}
\begin{example}
  It is easy to see that for any measure type $\M$ on a measurable space $(X,\XG)$, viewed as sets of states on $L^\infty(\M)$, the relation  $H(L^\infty(X,\XG),\M,\U)=L^2(\M)$, where the unitaries $U_{\mu,\mu'}$ are defined in \ref{eq:defUmumu'}, is valid.
\end{example}
\begin{example}
  Example \ref{ex:GNSnormalstateB(H)} shows  that $H(\B(\H),\P(\H),\U)\cong\H$, where $\P(\H)$ is the set of pure normal states on $\H$ and $U_{\eta,\eta'}$ was defined in Lemma \ref{lem:well-defined Uetaeta'}. For getting \ref{eq:productUetaeta'}, we need to fix the vectors $\psi$ corresponding to each $\eta\in\P(\H)$ in advance. That there is no canonical choice to achieve this  will present some problems later. 
\end{example}

\subsection{Algebraic Product Systems and Intrinsic Measurable Structures}
\label{sec:2measurability}

Now we want to examine  \emph{algebraic product systems} $\E=(\sg t\E,\sg{s,t}V)$. I.e., $\E_t$ is separable and  $\sg{s,t}V$ fulfils  \ref{eq:associativity Ust}, but we do not impose any measurability requirement in advance. We want to show in this section that the measurable structure of a (standard, measurable) product system is (upto isomorphy) already determined by its algebraic one.

First we collect some more notions.  Suppose $(X,d)$ is a complete separable metric  space. Then $\FG_X$ equipped with the Fell topology is Polish iff $X$ is locally compact. But $\FG_X$ is  Polish under the Wijsman topology which is  the coarsest topology making   
\begin{displaymath}
  \FG_X\ni Z\mapsto \mathrm{dist}(x,Z)=\inf\set{d(x,z):z\in Z}
\end{displaymath}
continuous for all $x$ \cite{A:Bee93}. Observe that the definition of the topology depends on $d$ but the Borel sets on $\FG_X$ coincide for all choices making $X$ complete. Further, they coincide even with the Borel sets generated by the Fell topology.  
On $\N(\H)$ one introduces the Effros Borel structure \cite[Section IV.8]{OP:Tak79} which is the $\sigma$-field generated by the maps $\A\mapsto \norm{\eta_\A}$, where $\eta$ ranges over the $\sigma$-weakly continuous functionals on $\B(\H)$ and $\eta_\A$ is the restriction of $\eta$ to $\A$. Then, by \cite[Corollary IV.8.5]{OP:Tak79}  $\N(\H)$ is a Standard Borel space. Further, there exist countably many operator valued measurable  sections $\N(\H)\ni\A\mapsto a^n_\A$ such that $\set{a^n_\A:n\in\NN}$ is ($\sigma$-) weakly dense in  $B_1(\A)$ for all $\A\in\N(\H)$. This means, that the Effros Borel structure is the Borel structure for the topology on $\N(\H)$ we introduced in the beginning. 

Let  $\P(\H)$ denote  the set of pure normal states on $\B(\H)$ which is a closed subset of the separable Banach space $\B^1(\H)=\set{\varrho\in\B(\H):\Tr\abs\varrho<\infty}$. 

The central result of  this section is the following
\begin{theorem}
\label{th:intrinsicmeasurability}
  Let $\E$ be an algebraic product system. Then the following statements are equivalent.
  \begin{enumerate}
  \item \label{intrinsic1}There exists a  measurable structure on $\sg t\E$ such that $\E$ is a product system.
  \item \label{intrinsic2}The  unitary group $\gr t{\tau}$, defined by \ref{eq:defflipgroup} is a  continuous one parameter group of unitaries on $\E_1$.
  \item \label{intrinsic3}The  automorphism group $\gr t{\sigma}$ on $\B(\E_1)$, defined by \ref{eq:defshiftautomorphisms}, is weakly  continuous.
  \end{enumerate}
\end{theorem}
\begin{remark}
  It is well-known that weak measurability and continuity coincide  both for  $\gr t\sigma$ and $\gr t\tau$. 

Further, continuity of $\gr t\sigma$ implies that there is continuous 1-parameter unitary group $\gr t{u}$ with $\sigma_t( \cdot)=u_t^*\cdot u_t$. This means that there is a character $\gamma$ on $\NR$ with $\tau_t=\gamma(t)u_t$. The surprising point is that $\tau$ (or $\gamma$) has to be continuous if $\sigma$ is, we have no other arguments in favour of this than the quite complicated proof below. 
\end{remark}
The proof of the essential part of this theorem (which can be found on page \pageref{page:proof th:intrinsicmeasurability}) relies on  Proposition \ref{prop:psbyE_1} which we prove first. 
\begin{proof}[~of Proposition \ref{prop:psbyE_1}]
\label{page:proof prop:psbyE_1}
Suppose we have a measurable family $\sgi t\E$ of Hilbert spaces equipped with a measurable family of unitaries $\map{V_{s,t}}{\E_s\otimes\E_t}{\E_{s+t}}$, $s,t\ge0$, $s+t\le1$ fulfilling \ref{eq:associativity Ust}. Then we fix the family $\sg t\F$ of Hilbert spaces
\begin{displaymath}
  \F_{n+t}=(\E_1)^{\otimes n}\otimes \E_t\dmf{n\in\NN, t\in[0,1)}.
\end{displaymath}
For defining a multiplication by  operators $\sg{s,t}W$, observe that $\F_{n+t}$ is isomorphic to each of the Hilbert spaces $H_\eta$, where $\eta$ is a pure normal state on $\B(\E_1^{\otimes m})$, $m>n$, if $\eta$ factorizes with respect to $\A_{0,n+t}\otimes\A_{n+t,m}=\B(\E_1^{\otimes m})$, defining $\A_{r,s}$ like in \ref{eq:defAst}. For $n+t$, $n'+t'$, $n''+t''$ consider on $H_\eta$, where $\eta$ is a pure normal state on $\B(\E_1^{\otimes m})$, $m>n+n'+n''+2$, the operators
\begin{displaymath}
  W_{n+t,n'+t'}[a_t]_\eta\otimes [a_{t'}]_\eta=[a_t\sigma^m_{n+t}(a_{t'})]_\eta
\end{displaymath}
where $\gr t{\sigma^m}$ is defined on $\B(\E_1^{\otimes m})$ like $\gr t\sigma$ in \ref{eq:defshiftautomorphisms}. The groups property of $\sigma^m$ shows that  the operators $\sg{s,t}W$ are associative in the sense of \ref{eq:associativity Ust}. Their  definition is  even independent from the choice of $\eta$ and $m$, as we derive for $n,n'\in\NN$, $\varepsilon,\varepsilon'\ge0$, $\varepsilon+\varepsilon'\le1$ that
\begin{eqnarray*}
 \lefteqn{ W_{n+\varepsilon,n'+\varepsilon'}(x^1\otimes x^2\otimes\dots\otimes x^n\otimes x_\varepsilon)\otimes((y^1_{1-\varepsilon}\otimes y^1_{\varepsilon})\otimes(y^2_{1-\varepsilon}\otimes y^2_{\varepsilon})\otimes\dots\otimes (y^{n'}_{1-\varepsilon}\otimes y^{n'}_{\varepsilon})\otimes y_{\varepsilon'})}
&=&x^1\otimes x^2\otimes\dots\otimes x^n\otimes (x_\varepsilon\otimes y^1_{1-\varepsilon})\otimes (y^1_{\varepsilon}\otimes y^2_{1-\varepsilon})\otimes\dots\otimes (y^{n'-1}_{\varepsilon}\otimes y^{n'}_{1-\varepsilon})\otimes (y^{n'}_{\varepsilon}\otimes y_{\varepsilon'})\\&&\dmf{x^i\in\E_1,x_\varepsilon,y^i_\varepsilon\in\E_\varepsilon,y^i_{1-\varepsilon}\in\E_{1-\varepsilon},y_{\varepsilon'}\in\E_{\varepsilon'}}
\end{eqnarray*}
and for $\varepsilon+\varepsilon'>1$ that
\begin{eqnarray*}
 \lefteqn{ W_{n+\varepsilon,n'+\varepsilon'}(x^1\otimes\dots\otimes x^n\otimes x_\varepsilon)\otimes((y^1_{1-\varepsilon}\otimes y^1_{\varepsilon})\otimes\dots\otimes (y^{n'}_{1-\varepsilon}\otimes y^{n'}_{\varepsilon})\otimes (y_{1-\varepsilon}\otimes y'_{\varepsilon'+\varepsilon-1}))}
&=&x^1\otimes \dots\otimes x^n\otimes (x_\varepsilon\otimes y^1_{1-\varepsilon})\otimes\dots\otimes (y^{n'-1}_{\varepsilon}\otimes y^{n'}_{1-\varepsilon})\otimes (y^{n'}_{\varepsilon}\otimes y_{1-\varepsilon})\otimes   y'_{\varepsilon'+\varepsilon-1}\\&&\dmf{x^i\in\E_1,x_\varepsilon,y^i_\varepsilon\in\E_\varepsilon,y^i_{1-\varepsilon},y_{1-\varepsilon}\in\E_{1-\varepsilon},y'_{\varepsilon'+\varepsilon-1}\in\E_{\varepsilon'+\varepsilon-1}}.
\end{eqnarray*}
Plugging in  measurable sections for the various $x$'s and $y$'s in these formulae it  is an easy consequence of measurability  of $(V_{s,t})_{s,t\ge0,s+t\le1}$ that $\sg{s,t}W$ is a measurable  operator family too.

Suppose further that  $\sgi t\E$ arises from  a product system $\E$ by restricting $\sg t\E$ to the index set $[0,1]$.  Then it is clear from these formulae (and associativity of $\sg {s,t}V$) that  $\F=(\sg t\F,\sg {s,t}W)$ is isomorphic to $\E$ as a product system. This completes the proof. 
\end{proof}
For the  following lemma remember the conventions of intervals on $\NT$, being closed clockwise arcs. Further, set
\begin{displaymath}
  \A_{[s,t]}=
  \left\{\begin{array}[c]{c>$l<$}
\A_{s,t}& if $s<t$\\
\A_{s,1}\otimes\A_{0,t}& if $s\ge t$
  \end{array}\right.\qquad\qquad\A_{(s,t)}=
  \left\{\begin{array}[c]{c>$l<$}
\A_{s,t}& if $s\le t$\\
\A_{s,1}\otimes\A_{0,t}& if $s> t$
  \end{array}\right..
\end{displaymath}
Note that both definitions coincide except for $s=t$. For any $s,t\in[0,1]$, let  $\P_{[s,t]}$ be the set of all pure normal states $\eta$ on $\B(\E_1)$  which factorize like 
\begin{displaymath}
  \eta(ab)=\eta(a)\eta(b)\dmf{a\in\A_{[s,t]},b\in\A_{(t,s)}}.
\end{displaymath}
We abbreviate $\P_{[0,t]}$ to $\P_t$ and $\P_{t_1}\cap\P_{t_2}\cap\dots\cap \P_{t_n}$ to $\P_{t_1,\dots,t_n}$.
\begin{lemma}
\label{lem:measurable P_t}
Let $\sg t\E$ be an algebraic product system for which $\gr t\sigma$ is $\sigma$-strongly continuous. 

Then  $\set{(s,t,\eta):s,t\in[0,1],\eta\in\P_{[s,t]}}\subset[0,1]\times\P(\E_1)$ is closed  as well as $\P_{[s,t]}$. 

Further,  the  set $\P^1=\bigcup_{t\in(0,1)}\P_t$ is measurable inside $\P(\E_1)$. More generally, for all $n\in\NN$ the set $\P^n=\bigcup_{0<t_1<t_2<\dots< t_n<1}\P_{t_1,\dots,t_n}$ is measurable.

Additionally, there exists a measurable section $\sgi t\eta$  through $\sgi t\P$.
\end{lemma}
\begin{proof}
For simplicity, we prove the results for $\sgi t\P$ only, the general proof goes along the same lines but it is notationally more evolved.  

Consider for $0\le s<s'\le1$ the sets $\Q_{s,s'}=\set{\eta\in\P(\E_1):\eta(ab)=\eta(a)\eta(b)\forall a\in\A_{0,s},b\in\A_{s',1}}$ which are closed and hence measurable.  Since $\bigcup_{s<t}\A_{0,s}$ is $\sigma$-strongly dense in $\A_{0,t}$ by Proposition \ref{prop:singlepoint}, we obtain $\P_t=\bigcap_{
s,s'\in\ND:t\in(s,s')}\Q_{s,s'}$ where $\ND=\set{k2^{-n}:n\in\NN,k=0,\dots,2^n}$ are the dyadic numbers. This shows closedness of $\P_t$. Similarly, $t_n\to t$ and $\P_{t_n}\ni\eta_n\to\eta$ implies $\eta\in \Q_{s,s'}$ for all $s,s'$ with $s<t<s'$ since $\Q_{s,s'}$ is closed. Thus $\eta\in \P_t$ and  $\set{(t,\eta):t\in[0,1],\eta\in\P_t}$ is closed in $[0,1]\times\P(\E_1)$.

Now we consider for  $r,r'\in\ND$ the set  $\R_{r,r'}=\bigcap_{n\in\NN}\bigcup_{s\in 2^{-n}\NN\cap[r,r']}\Q_{s,s+2^{-n}}$. Then the set  $\R_{r,r'}$ is closed and $\eta\in \R_{r,r'}$ iff for all $n\in\NN$ there is an $s_n\in 2^{-n}\NN\cap [r,r']$ such that $\eta\in\Q_{s_n,s_n+2^{-n}}$. Clearly, there is a subsequence $\sequ kn$ such that $k\to [s_{n_k},s_{n_k}+2^{-n_k}] $ is decreasing and $\bigcap_{k\in\NN}[s_{n_k},s_{n_k}+2^{-n_k}]=\set t$ for some $t\in[r,r']$. This shows $\eta\in\P_t$. Conversely, if   $\eta\in\P_t$ for some $t\in[r,r']$ then $\eta\in\R_{r,r'}$. Thus we conclude that $\bigcup_{t\in(0,1)}\P_t=\bigcup_{r,r'\in\ND}\R_{r,r'}$ is measurable. The generalisation follows from formulae similar to 
\begin{displaymath}
  \bigcup_{0<t<t'<1}\P_{t,t'}=\bigcup_{0<r<r'<r''<r'''\in\ND}\R_{r,r'}\cap\R_{r'',r'''}.
\end{displaymath}

Since $\set{(s,t,\eta):s,t\in[0,1],\eta\in\P_{[s,t]}}\subset[0,1]\times\P(\E_1)$ is closed, the map $(s,t)\mapsto \P_{[s,t]}\in\FG_{\P_{\E_1}}$ is upper semicontinuous in the sense that for any open set $G\subset\P(\E_1)$ and  $s_n\tlimitsto{}{n\to\infty}s$, $t_n\tlimitsto{}{n\to\infty}t$ with $G\cap\P_{[s,t]}\ne\emptyset$ fulfils  $G\cap\P_{[s_n,t_n]}\ne\emptyset$ eventually. Thus this map  is  also measurable.   \cite[Theorem 6.6.7]{A:Bee93} completes the proof.
\end{proof}
For $\eta\in\P(\E_1)$ we define  the set $F(\eta)=\set{t\in(0,1):\eta\in\P_t}$ and $\tilde F(\eta)=\set{(s,t):\eta\in\P_{[s,t]}}$.
\begin{lemma}
\label{lem:measurability F(eta)}
   $F(\eta)\in\FG_{(0,1)}$ and  $\tilde F(\eta)\in\FG_{\NT²}$ for all $\varepsilon\in \P(\E_1)$. Further,  the maps  $\eta\in\P(\E_1)\mapsto F(\eta),\tilde F(\eta)\in\FG_{(0,1)}$ are  measurable.
\end{lemma}
\begin{proof}
   The first assertions follows from the fact that the set $\set{(s,t,\eta):s,t\in[0,1],\eta\in\P_{[s,t]}}$ is closed in $[0,1]\times\P(\E_1)$. 

For the second one, we observe for all $0\le s<t\le 1$  that $\set{\eta:F(\eta)\cap(s,t)\ne\emptyset}=\bigcup_{r,r'\in\ND,s<r<r'<t}\R_{r,r'}$ is measurable for all $0<s<t<1$. Since the sets  $\set{Z:Z\cap(s,t)=\emptyset}$, $0\le s<t\le1$, generate the Borel $\sigma$-field on $\FG_{(0,1)}$ the set  $F$ is measurable. The proof for $\tilde F$ is similar.
\end{proof}
We also need the following result on (not necessarily measurable) multipliers on $\NR$ resp.\ $[0,1]$. This is necessary to extend some  results derived in \cite{OP:Arv89} for usual product systems  for  algebraic ones.
\begin{lemma}
\label{lem:multiplieronR}
  Suppose $\map m{\NR^2}\NT$ fulfils for all $s,t,r\in\NR$
  \begin{equation}
  \label{eq:multiplier}
  m(s+t,r)m(s,t)=m(s,t+r)m(t,r).
  \end{equation}
Then there is a function $\map f\NR\NT$ such that 
\begin{equation}
\label{eq:solution multiplier}
  m(s,t)=\frac{f(s)f(t)}{f(s+t)}\dmf{s,t\in\NR}.
\end{equation}
\end{lemma}
\begin{remark}
Formula   \ref{eq:multiplier} restricted to  $\NRpp$   determines all algebraic product structures on $\sgp t\NC$, identifying $\NT$ with the set of unitaries from $\NC\otimes\NC$ to $\NC$.  
\end{remark}
\begin{proof}
 By Tykhonovs theorem,  $\NT^\NR$ is compact in the product topology. For any finite dimensional $\NQ$-vector subspace $H\subset \NR$ we consider 
  \begin{displaymath}
    S_H=\set{\map f\NR\NT:f\text{~fulfils \ref{eq:solution multiplier} for all~}s,t\in H}.
  \end{displaymath}
Clearly, $S_H$ is closed. In the following we show that it is nonvoid.

Let $H$ be fixed, $\H=L^2(H,\sum_{x\in H}\delta_x)$ and define operators $(u_t)_{t\in H} $ on $\H$ by 
\begin{displaymath}
  u_th(x)=m(t,x-t)h(x-t)\dmf{x,t\in H}
\end{displaymath}
It is easy to see that each $u_t$ is unitary and
\begin{displaymath}
  u_su_t=m(s,t)u_{s+t}\dmf{s,t\in H}.
\end{displaymath}
Thus, defining automorphisms  $(\alpha_t)_{t\in H} $ of $\B(\H)$ by
\begin{displaymath}
  \alpha_t(a)=u_t^* au_t\dmf{a\in\B(\H)}
\end{displaymath}
we arrive at
\begin{displaymath}
  \alpha_{s+t}=\alpha_s\circ\alpha_t\dmf{s,t\in H}.
\end{displaymath}
We want to show that there is a collection $(v_t)_{t\in H}$ of unitaries in $\H$ such that 
\begin{equation}
  \label{eq:representation automorphism}
  \alpha_t(a)=v_t^* av_t\dmf{a\in\B(\H)}
\end{equation}
and $v_sv_t=v_{s+t}$. 

Choosing a basis of $H$ (over $\NQ$) it is easy to see that we can restrict ourselves to the case $H=\NQ$ for a while. Now it is well-known that $u_t^* au_t=v_t^* av_t\forall a\in\B(\H)$ implies that $v_t$ is a multiple of $u_t$. Therefore, 
\begin{displaymath}
  T=\set{(v_t)_{t\in \NQ}:\text{\ref{eq:representation automorphism} is valid for all~}t\in \NQ}
\end{displaymath}
is compact under the product topology. Define for $k\in \NN^*$
\begin{displaymath}
  T_k=\set{(v_t)_{t\in \NQ}\in T: v_{p/k}v_{q/k}=v_{(p+q)/k}\text{~ for all~}p,q\in\NZ}
\end{displaymath}
Then $T_k$ is closed.   The choice 
\begin{displaymath}
  v_t=\left\{
    \begin{array}[c]{c>{$~}l<$}
u_{1/k}^p&if $t=p/k$ for some $p\in\NZ$\\
u_t&otherwise
    \end{array}
\right.\dmf{t\in\NQ}
\end{displaymath}
fulfils $v_sv_t=v_{s+t}$ for all $s,t\in\frac1k\NZ$. Further, for $t=p/k$ and $a\in\B(\H)$ we know 
\begin{eqnarray*}
 v_t^* av_t= u_{1/k}^{-p} au_{1/k}^p&=& \underbrace{u_{1/k}^*u_{1/k}^*\cdots u_{1/k}^*}_{p\text{~times~}} a \underbrace{u_{1/k}u_{1/k}\cdots u_{1/k}}_{p\text{~times~}}\\
&=& \underbrace{\alpha_{1/k}\circ\alpha_{1/k}\circ\cdots \circ\alpha_{1/k}}_{p\text{~times~}}(a)=\alpha_{p/k}(a)=\alpha_t(a).
\end{eqnarray*}
Thus  $T_k$ is not void. Since $T_{k'}\subset T_k$ if $k| k'$, the family $\sequ kT$ enjoys the finite intersection property and by compactness, $\bigcap_{k\in\NN^*}T_k$ is nonvoid, let $(v_t)_{t\in\NQ}$ be a member of this set. Then it is easy to see that $v_sv_t=v_{s+t}$ for all $k\in\NN$ for all $s,t\in\frac1k\NZ$, i.e.\ for all $s,t\in\NQ$. This  is our desired property.

Consequently, in the general case there exists a similar representation  $(v_t)_{t\in H}$ of $H$. By the preceding remarks, there are numbers $(f(t))_{t\in H}$ such that $v_t=f(t)u_t$ for all $t\in H$ what implies 
\begin{displaymath}
 f(s+t) u_{s+t}=v_{s+t}=v_sv_t=f(s)u_sf(t)u_t=f(s)f(t)m(s,t)^{-1}u_{s+t}\dmf{s,t\in H}
\end{displaymath}
and henceforth \ref{eq:solution multiplier}. Extending $f$ arbitrarily to a function on $\NR$ shows that $S_H$ is nonvoid. 

Further, $S_{H'}\subset S_H$ for $H'\supset H$ shows the finite intersection property for the sets $S_H$ and we find that $\bigcap_H S_H\ne\emptyset$. This shows the assertion. 
\end{proof}
\begin{corollary}
  \label{cor:multipliersimplex}
Suppose that $\map m{\set{(s,t):s,t\ge0,s+t\le1}}\NT$ is such that \ref{eq:multiplier} is fulfilled for all $s,t,r\ge0,s+t+r\le1$. Then there is a map $\map f{[0,1]}\NT$ with \ref{eq:solution multiplier} for all $s,t\ge0,s+t\le1$.
\end{corollary}
\begin{proof}
  By the preceding lemma, it is enough  to extend $m$ to a map on the whole $\NR^2$ in compliance with \ref{eq:multiplier}.
This is done in a similar way  like in the proof of Proposition \ref{prop:psbyE_1} (compare Note \ref{rem:onedimensional product systems}).  Define the map $\map\kappa{[0,1]}\NT$ by
\begin{displaymath}
  \kappa(\varepsilon)=\frac{m(\varepsilon,1-\varepsilon)}{m(1-\varepsilon,\varepsilon)}
\end{displaymath}
and $\map M{\NR²}\NT$,
\begin{displaymath}
M (n+\varepsilon,n'+\varepsilon')=\left\{
    \begin{array}[c]{c>{$~}l<$}
\kappa(\varepsilon)^{n'}m(\varepsilon,\varepsilon')&if $\varepsilon+\varepsilon'\le 1$\\
\kappa(\varepsilon)^{n'}m(\varepsilon,1-\varepsilon)m(1-\varepsilon,\varepsilon+\varepsilon'-1)^{-1}&otherwise
    \end{array}
\right.
\end{displaymath}
for arbitrary $n,n'\in\NZ$, $\varepsilon,\varepsilon'\in[0,1)$. From \ref{eq:multiplier} we know $m(1,0)=m(0,0)=m(0,1)$ such that $M$ extends $m$ to $\NR^2$. It remains to show \ref{eq:multiplier}, i.e.\ that
\begin{displaymath}
 M (n+\varepsilon,n'+\varepsilon'+n''+\varepsilon'')M (n'+\varepsilon',n''+\varepsilon'') = M (n+\varepsilon+n'+\varepsilon',n''+\varepsilon'')M (n+\varepsilon,n'+\varepsilon')
\end{displaymath}
for all $n,n',n''\in\NZ$ and $\varepsilon,\varepsilon',\varepsilon''\in[0,1)$. Clearly, it is enough to show this for $n=n'=n''=0$, but concerning the relations between $\varepsilon,\varepsilon',\varepsilon''$ we have to distinguish several cases.
\begin{enumerate}\def\labelenumi{\arabic{enumi}.~}
\item $\varepsilon+\varepsilon'+\varepsilon''\le1$ implies $\varepsilon+\varepsilon'\le1$, $\varepsilon'+\varepsilon''\le1$
  \begin{eqnarray*}
    \mathrm{LHS}&=&m (\varepsilon,\varepsilon'+\varepsilon'')m (\varepsilon',\varepsilon'') \\
    \mathrm{RHS}&=&m (\varepsilon+\varepsilon',\varepsilon'')m (\varepsilon,\varepsilon')
  \end{eqnarray*}
Equality follows directly from \ref{eq:multiplier}.
\item $\varepsilon+\varepsilon'+\varepsilon''>1$, $\varepsilon+\varepsilon'\le1$, $\varepsilon'+\varepsilon''\le1$ show $2\ge\varepsilon+\varepsilon'+\varepsilon''$.
  \begin{eqnarray*}
    \mathrm{LHS}&=&m(\varepsilon,1-\varepsilon)m(1-\varepsilon,\varepsilon+\varepsilon'+\varepsilon''-1)^{-1}m (\varepsilon',\varepsilon'') \\
    \mathrm{RHS}&=&m (\varepsilon+\varepsilon',1-\varepsilon-\varepsilon')m(1-\varepsilon-\varepsilon',\varepsilon+\varepsilon'+\varepsilon''-1)^{-1}m (\varepsilon,\varepsilon')
  \end{eqnarray*}
$\mathrm{LHS}=\mathrm{RHS}$ follows from 
\begin{eqnarray*}
 m (\varepsilon+\varepsilon',1-\varepsilon-\varepsilon') m (\varepsilon,\varepsilon')&=&m(\varepsilon,1-\varepsilon)m(\varepsilon',1-\varepsilon-\varepsilon')\\
m(1-\varepsilon-\varepsilon',\varepsilon+\varepsilon'+\varepsilon''-1)m (\varepsilon',\varepsilon'')&=&m(\varepsilon',1-\varepsilon-\varepsilon')m(1-\varepsilon,\varepsilon+\varepsilon'+\varepsilon''-1).
\end{eqnarray*}
\item $\varepsilon+\varepsilon'>1$, $\varepsilon'+\varepsilon''\le1$ imply $2>\varepsilon+\varepsilon'+\varepsilon''>1$
  \begin{eqnarray*}
     \mathrm{LHS}&=&m(\varepsilon,1-\varepsilon)m(1-\varepsilon,\varepsilon+\varepsilon'+\varepsilon''-1)^{-1}m (\varepsilon',\varepsilon'') \\
    \mathrm{RHS}&=&m (\varepsilon+\varepsilon'-1,\varepsilon'')m(\varepsilon,1-\varepsilon)m(1-\varepsilon,\varepsilon+\varepsilon'-1)^{-1}
  \end{eqnarray*}
$\mathrm{LHS}=\mathrm{RHS}$ follows from 
\begin{displaymath}
  m(1-\varepsilon,\varepsilon+\varepsilon'+\varepsilon''-1)m (\varepsilon+\varepsilon'-1,\varepsilon'')=m (\varepsilon',\varepsilon'')m(1-\varepsilon,\varepsilon+\varepsilon'-1).
\end{displaymath}
\item $\varepsilon+\varepsilon'\le1$, $\varepsilon'+\varepsilon'' >1$ imply $2>\varepsilon+\varepsilon'+\varepsilon''>1$
  \begin{eqnarray*}
     \mathrm{LHS}&=&\kappa(\varepsilon)m(\varepsilon,\varepsilon'+\varepsilon''-1)m(\varepsilon',1-\varepsilon')m(1-\varepsilon',\varepsilon'+\varepsilon''-1)^{-1}\\
    \mathrm{RHS}&=&m (\varepsilon+\varepsilon',1-\varepsilon-\varepsilon')m(1-\varepsilon-\varepsilon',\varepsilon+\varepsilon'+\varepsilon''-1)^{-1}m (\varepsilon,\varepsilon')
  \end{eqnarray*}
$\mathrm{LHS}=\mathrm{RHS}$ follows from 
\begin{eqnarray*}
  m(1-\varepsilon-\varepsilon',\varepsilon+\varepsilon'+\varepsilon''-1)m(\varepsilon,\varepsilon'+\varepsilon''-1)&=&m(1-\varepsilon',\varepsilon'+\varepsilon''-1)m(1-\varepsilon-\varepsilon',\varepsilon)\\
  m(\varepsilon',1-\varepsilon')m(1-\varepsilon-\varepsilon',\varepsilon)&=&m(1-\varepsilon,\varepsilon)m(\varepsilon',1-\varepsilon-\varepsilon')\\
  m(\varepsilon,1-\varepsilon)m(\varepsilon',1-\varepsilon-\varepsilon')&=&m(\varepsilon+\varepsilon',1-\varepsilon-\varepsilon')m(\varepsilon,\varepsilon').
\end{eqnarray*}
\item $2\ge\varepsilon+\varepsilon'+\varepsilon''>1$, $\varepsilon+\varepsilon'>1$, $\varepsilon'+\varepsilon''>1$ 
\begin{eqnarray*}
     \mathrm{LHS}&=&\kappa(\varepsilon)m(\varepsilon,\varepsilon'+\varepsilon''-1)m(\varepsilon',1-\varepsilon')m(1-\varepsilon',\varepsilon'+\varepsilon''-1)^{-1}\\
    \mathrm{RHS}&=&m (\varepsilon+\varepsilon'-1,\varepsilon'')m(\varepsilon,1-\varepsilon)m(1-\varepsilon,\varepsilon+\varepsilon'-1)^{-1}
  \end{eqnarray*}
$\mathrm{LHS}=\mathrm{RHS}$ follows from 
\begin{eqnarray*}
m(1-\varepsilon,\varepsilon)m (\varepsilon+\varepsilon'-1,1-\varepsilon') &=& m(\varepsilon',1-\varepsilon')m(1-\varepsilon,\varepsilon+\varepsilon'-1)\\
m (\varepsilon+\varepsilon'-1,\varepsilon'')m(1-\varepsilon',\varepsilon'+\varepsilon''-1)&=& m(\varepsilon,\varepsilon'+\varepsilon''-1) m (\varepsilon+\varepsilon'-1,1-\varepsilon')
\end{eqnarray*}
\item  $\varepsilon+\varepsilon'+\varepsilon''>2$ implies  $\varepsilon+\varepsilon'>1$, $\varepsilon'+\varepsilon''>1$
\begin{eqnarray*}
     \mathrm{LHS}&=&\kappa(\varepsilon)m(\varepsilon,1-\varepsilon)m(1-\varepsilon,\varepsilon+\varepsilon'+\varepsilon''-2)^{-1}m(\varepsilon',1-\varepsilon')m(1-\varepsilon',\varepsilon'+\varepsilon''-1)^{-1}\\
    \mathrm{RHS}&=&m (\varepsilon+\varepsilon'-1,2-\varepsilon-\varepsilon')m(2-\varepsilon-\varepsilon',\varepsilon+\varepsilon'+\varepsilon''-2)^{-1}m(\varepsilon,1-\varepsilon)m(1-\varepsilon,\varepsilon+\varepsilon'-1)^{-1}
  \end{eqnarray*}
$\mathrm{LHS}=\mathrm{RHS}$ follows from 
\begin{eqnarray*}
m(1-\varepsilon',\varepsilon'+\varepsilon''-1)m(1-\varepsilon,\varepsilon+\varepsilon'+\varepsilon''-2)&=&m(2-\varepsilon-\varepsilon',\varepsilon+\varepsilon'+\varepsilon''-2)m(1-\varepsilon',1-\varepsilon)\\
m(1-\varepsilon,\varepsilon)m (\varepsilon+\varepsilon'-1,1-\varepsilon') &=& m(\varepsilon',1-\varepsilon')m(1-\varepsilon,\varepsilon+\varepsilon'-1)\\
m (\varepsilon+\varepsilon'-1,  2-\varepsilon-\varepsilon')m(1-\varepsilon',1-\varepsilon)&=&m(\varepsilon,1-\varepsilon)m (\varepsilon+\varepsilon'-1,1-\varepsilon')
\end{eqnarray*}
\end{enumerate}
Since this covers all possible relations between $\varepsilon,\varepsilon',\varepsilon''$ the proof is complete.  
\end{proof}
\begin{proof}[ of Theorem \ref{th:intrinsicmeasurability}]
\label{page:proof th:intrinsicmeasurability}
The  statement \ref{intrinsic1}$\Rightarrow$\ref{intrinsic2} was proved in   Proposition \ref{prop:tau_tiscontinuous}. The implication \ref{intrinsic2}$\Rightarrow$\ref{intrinsic3} is trivial. So it remains to  show \ref{intrinsic3}$\Rightarrow$\ref{intrinsic1},  we use Proposition \ref{prop:psbyE_1} for this purpose. Thus  it suffices to show that we can equip $\sgi t\E$ with a measurable structure making $(V_{s,t})_{s,t\ge0,s+t\le1}$ measurable. This is done by constructing an algebraically isomorphic family $\sgi t{\E'}$ as $\sgip t{H(\A_{0,t},\P_t,\U^t)}$ where we can exploit the measurable structure coming from $\N(\E_1)$ and $\P(\E_1)$ to define an intrinsic  measurable structure. A crucial point thereby is that measurability of $\sigma$ implies its $\sigma$-weak continuity \cite[Proposition 2.5]{Arv89}.

For the  construction of  the spaces  $H(\A_{0,t},\P_t,\U^t)$ according to Lemma \ref{lem:defH(S)} the problem is the choice of  unitaries $U^t_{\eta,\eta'}$ intertwining the GNS Hilbert spaces. Following Example \ref{ex:GNSnormalstateB(H)}, we assume that there are operators $Q^t_{\eta,\eta'}\in\A_{0,t}$ such that
\begin{displaymath}
  U^t_{\eta,\eta'}[a]_\eta=[aQ^t_{\eta,\eta'}]_{\eta'}\dmf{a\in\A_{0,t}}.
\end{displaymath}
Since the unitaries $(U^t_{\eta,\eta'})_{\eta,\eta'\in\P_t}$ should fulfil \ref{eq:productUetaeta'}, the operators $(Q^t_{\eta,\eta'})_{\eta,\eta'\in\P_t}$ cannot be chosen freely. For $t\in[0,1]$ and $\eta\in\P_t$ denote by $P^t_\eta$  the supporting projection of $\eta$ in $\A_{0,t}$. I.e., $P^t_\eta=\Pr\psi\otimes\unit_{\E_{1-t}}$ for some $\psi\in\E_t$ and  $\eta(P^t_\eta)=1$. Now we assume that the operators $(Q^t_{\eta,\eta'})_{\eta,\eta'\in\P_t}$ are chosen in such a way that
\begin{enumerate}
\item \label{mess0}$P^t_\eta =Q^t_{\eta,\eta'}(Q^t_{\eta,\eta'})^*$ and $P^t_{\eta'}=(Q^t_{\eta,\eta'})^*Q^t_{\eta,\eta'}$ for all $t\in[0,1]$, $\eta,\eta'\in\P_t$.
\item \label{mess3}$Q^t_{\eta,\eta'}Q^t_{\eta',\eta''}=Q^t_{\eta,\eta''}$ for all $t\in[0,1]$, $\eta,\eta',\eta''\in\P_t$.
\item \label{mess1}$(t,\eta,\eta')\mapsto Q^t_{\eta,\eta'}$ is weakly measurable on the closed set  $\set{(t,\eta,\eta'):t\in[0,1],\eta,\eta'\in\P_t}$.
\item \label{mess4}\begin{math}
 Q^s_{\eta,\eta'}\sigma_s(Q^t_{\eta\circ\sigma_{-s},\eta'\circ\sigma_{-s}})=Q^{s+t}_{\eta,\eta'}
\end{math} for all $s,t\in[0,1]$, $s+t\le1$ and $\eta,\eta'\in\P_{s,s+t}$
\end{enumerate}
Below we show how to construct $(Q^t_{\eta,\eta'})_{\eta,\eta'\in\P_t}$ with \ref{mess0}--\ref{mess4}, but first we proceed along the general lines of the proof. Observe that \ref{mess0} and \ref{mess3}  imply $Q^t_{\eta,\eta}=P^t_\eta$, $Q^t_{\eta',\eta}=(Q^t_{\eta,\eta'})^*$ and  fix $Q^t_{\eta,\eta'}$ upto a complex factor. Clearly, these two conditions assure that for all $t\in[0,1]$  the family $(U^t_{\eta,\eta'})_{\eta,\eta'\in\P_t}$  is well-defined and fulfils \ref{eq:productUetaeta'}. Thus  we can set $\E'_t=H(\A_{0,t},\P_t,\U^t)$, $t\in[0,1]$. 

Now we define a product in $\sgi t{\E'}$ by  operators $\sgii st{V'}$, $\map{V'_{s,t}}{\E'_s\otimes\E'_t}{\E'_{s+t}}$ like follows. For   $s,t\ge0$, $s+t\le1$ and   $\eta\in\P_{s,s+t}$ we  set (note that $\eta\circ\sigma_{-s}\in\P_t$)
\begin{equation}
\label{eq:def multiplication V'st}
  \big(V'_{s,t} ([a_{\eta'}]_{\eta'})_{{\eta'}\in\P_s}\otimes ([b_{\eta''}]_{\eta''})_{\eta''\in\P_t}\big)_{\eta}=[a_{\eta}\sigma_s(b_{\eta\circ\sigma_{-s}})]_{\eta}.
\end{equation}
We derive from 
\begin{eqnarray*}
  \lefteqn{\scpro{[a_{\eta}\sigma_s(b_{\eta\circ\sigma_{-s}})]_{\eta}}{[a'_{\eta}\sigma_s(b'_{\eta\circ\sigma_{-s}})]_{\eta}}}&=&\eta(\sigma_s(b_{\eta\circ\sigma_{-s}})^*a_{\eta}^*a'_{\eta}\sigma_s(b'_{\eta\circ\sigma_{-s}}))\\
&=&\eta(\underbrace{a_{\eta}^*a'_{\eta}}_{\in\A_{0,s}}\underbrace{\sigma_s(b_{\eta\circ\sigma_{-s}}^*b'_{\eta\circ\sigma_{-s}})}_{\in\A_{s,s+t}})\\
&=&\eta(a_{\eta}^*a'_{\eta})\eta(\sigma_s(b_{\eta\circ\sigma_{-s}}^*b'_{\eta\circ\sigma_{-s}}))\\
&=&\scpro{([a_{\eta'}]_{\eta'})_{{\eta'}\in\P_s}}{([a'_{\eta'}]_{\eta'})_{{\eta'}\in\P_s}}_{\eta}\scpro{\smash{([b_{\eta''}]_{\eta''})_{\eta''\in\P_t}}}{\smash{([b'_{\eta''}]_{\eta''})_{\eta''\in\P_t}}}_{\eta\circ\sigma_{-s}}
\end{eqnarray*}
that $V'_{s,t}$ is isometric and thus   well-defined on $H_{\eta}\otimes H_{\eta\circ\sigma_{s}}$. Further, $\A_{0,s}\A_{s,s+t}$ is $\sigma$-weakly total in $\A_{0,s+t}$ and $\eta$ is $\sigma$-weakly continuous such that the image of  $V'_{s,t}$ is dense. Thus $V'_{s,t}$ extends to a unitary on $H_{\eta}\otimes H_{\eta\circ\sigma_{s}}$. Now take another  $\eta'\in\P_{s,s+t}$. We find from \ref{mess4} and \ref{mess0}
\begin{eqnarray*}
U_{\eta,\eta'} [a_{\eta}\sigma_s(b_{\eta\circ\sigma_{s}})]_{\eta} &=&[a_{\eta}\sigma_s(b_{\eta\circ\sigma_{s}})Q^{s+t}_{\eta,\eta'}]_{\eta'}\\
 &=&[a_{\eta'}Q^s_{\eta',\eta}\sigma_s(b_{\eta'\circ\sigma_{s}}Q^t_{\eta'\circ\sigma_{s},\eta\circ\sigma_{s}})Q^{s+t}_{\eta,\eta'}]_{\eta'}\\
 &=&[a_{\eta'}\sigma_s(b_{\eta'\circ\sigma_{s}})Q^s_{\eta',\eta}\sigma_s(Q^t_{\eta'\circ\sigma_{s},\eta\circ\sigma_{s}})Q^{s+t}_{\eta,\eta'}]_{\eta'}\\ 
 &=&[a_{\eta'}\sigma_s(b_{\eta'\circ\sigma_{s}})Q^{s+t}_{\eta',\eta}Q^{s+t}_{\eta,\eta'}]_{\eta'}\\ 
&=&[a_{\eta'}\sigma_s(b_{\eta'\circ\sigma_{s}})]_{\eta'}.
\end{eqnarray*}
A density argument shows that $U_{\eta,\eta'}(V'_{s,t}\psi\otimes\psi')_{\eta}=(V'_{s,t}\psi\otimes\psi')_{\eta'}$, i.e.\ the family $V'$ is well-defined on $\E'$. Further, we get  for any  $\eta^*\in\P_{s,s+t,s+t+r}$
\begin{eqnarray*}
\lefteqn{  \big(V'_{s,t+r} ([a_\eta]_\eta)_{\eta\in\P_s}\otimes (V'_{t,r} ([b_{\eta'}]_{\eta'})_{\eta'\in\P_t}\otimes ([c_{\eta''}]_{\eta''})_{\eta''\in\P_r})\big)_{\eta^*}}&=&[a_{\eta^*}\sigma_s(b_{\eta^*\circ\sigma_{s}}\sigma_t(c_{\eta^*\circ\sigma_{s}\circ\sigma_{t}}))]_{\eta^*}=[a_{\eta^*}\sigma_s(b_{\eta^*\circ\sigma_{s}})\sigma_{s+t}(c_{\eta^*\circ\sigma_{s+t}})]_{\eta^*}.
\end{eqnarray*}
This shows that the operators $\sgii st{V'}$ are associative in the sense of \ref{eq:associativity Ust}.

For establishing a consistent measurable structure on $(\sgi t{\E'},\sgii st{V'})$, observe that Proposition \ref{prop:singlepoint} is valid since its  proof relied on  the  continuity of $\gr t\sigma$ only. Thus for any $\sigma$-weakly continuous functional $\eta$ on $\B(\E_1)$ the map $I_{0,1}\ni(s,t)\mapsto\norm{\eta_{\A_{s,t}}}$ is measurable and  consequently  $I_{0,1}\ni(s,t)\mapsto\A_{s,t}$ is  measurable too. This gives us countably many measurable sections $\sgi t{a^n}$ through $\B(\E_1)$   such that $\overline{\mathrm{lh}\set{a^n_t:n\in\NN}}=\A_{0,t}$. Moreover, using Lemma \ref{lem:measurable P_t}  we find  a measurable curve $\sgi t\eta$ through $\P(\E_1)$   such that $\eta_t\in\P_t$.  Define the measurable structure $(\E')_0$ as the set of all sections $\sgip t{\psi^t}$ where $\psi^t\in H(\A_{0,t},\P_t,\U^t)$ and $t\mapsto\scpro{[a^n_t]_{\eta_t}}{\psi^t_{\eta_t}}_{\eta_t}$ is measurable for all $n\in\NN$. We get easily that for another measurable section $\sgi t{b}$ through $\sgip t{\A_{0,t}}$  that the map
\begin{displaymath}
  t\mapsto\scpro{[a^n_t]_{\eta_t}}{[b_t]_{\eta_t}}_{\eta_t}=\eta_t((a^n_t)^*b_t)
\end{displaymath}
is measurable since multiplication and taking adjoints are  weakly measurable operations in $\B(\E_1)$. Similarly, we could use another measurable curve $\sgi t{\eta'}$ with the prescribed properties. Then
\begin{displaymath}
t\mapsto  \scpro{[a^n_t]_{\eta'_t}}{U^t_{\eta'_t,\eta_t}[b_t]_{\eta_t}}_{\eta'_t}
=\eta_t((a^n_t)^*b_tQ^t_{\eta'_t,\eta_t})
\end{displaymath}
is  measurable what shows that our measurable structure does not depend on $\sgi t\eta$. Determine  for $s,t\ge0,s+t\le1$ the normal state $\eta_{s,t}$ on $\B(\E_1)$ by
\begin{displaymath}
  \eta_{s,t}(a_s\otimes a'_t\otimes a''_{1-s-t})=\eta_s(a_s)\eta_t(a'_t)\eta_{1-s-t}(a''_{1-s-t}).
\end{displaymath}
Clearly, the sheet $\sgii st\eta$ is  measurable  with $\eta_{s,t}\in\P_{s,s+t}$. Now it follows for $k,l,m\in\NN$ and $s,t\ge0, s+t\le1$ that
\begin{displaymath}
  \scpro{[a^k_{s+t}]_{\eta_{s,t}}}{[b_{\eta_{s,t}}\sigma_s(c_{\eta_{s,t}\circ\sigma_{s}})]_{\eta_{s,t}}}_{\eta_{s,t}}
=\eta_{s,t}((a^k_{s+t})^*b_{\eta_{s,t}}\sigma_s(c_{\eta_{s,t}\circ\sigma_{s}})),
\end{displaymath}
i.e.,  $\sgii st{V'}$  is measurable.

Next we construct  a (noncanonical) isomorphism between  $\sgi t\E$ and $\sgi t{\E'}$. Choose a  section $\sgi t\psi$ through $\sgi t\E$ such that  
\begin{displaymath}
  \eta_t(a_t\otimes\unit_{\E_{1-t}})=\scpro{\psi_t}{a_t\psi_t}\dmf{t\in[0,1],a_t\in\B(\E_t)}.
\end{displaymath}
We define operators $\map{\theta'_t}{\E_t'}{\E_t}$ through
\begin{displaymath}
  \theta'_t([a_\eta\otimes\unit]_\eta)_{\eta\in\P_t}=a_{\eta_t}\psi_t
\end{displaymath}
It is clear that  $\theta'_t$ is well-defined. Since it is isometric and has full range, it is unitary. 

Further, fix $s,t\ge0$, $s+t\le1$. Then
\begin{displaymath}
  V_{s,t}\theta'_s([a_\eta\otimes\unit]_\eta)_{\eta\in\P_s}\otimes \theta'_t([b_{\overline\eta}\otimes\unit]_{\overline\eta})_{{\overline\eta}\in\P_t}=  V_{s,t}a_{\eta_s}\psi_s\otimes b_{\eta_t}\psi_t
\end{displaymath}
and 
\begin{displaymath}
 \theta'_{s+t}V'_{s,t} ([a_\eta\otimes\unit]_\eta)_{\eta\in\P_s}\otimes ([b_\eta\otimes\unit]_\eta)_{\eta\in\P_t}
=c_{\eta_{s+t}}\psi_{s+t}
\end{displaymath}
where
\begin{displaymath}
  [c_{\eta_{s,t}}]_{\eta_{s,t}}=[a_{\eta_{s,t}}\otimes\unit \sigma_s(b_{\eta_{s,t}\circ\sigma_{s}})]_{\eta_{s,t}}
=[a_{\eta_s}\otimes b_{\eta_t}\otimes\unit_{\E_{1-s-t}}Q^s_{\eta_{s,t},\eta_s}\sigma_s(Q^t_{\eta_{s,t}\circ\sigma_{s},\eta_t})]_{\eta_{s,t}}.
\end{displaymath}
Thus, we may choose $c_{\eta_{s+t}}=a_{\eta_s}\otimes b_{\eta_t}\otimes\unit_{\E_{1-s-t}}Q^s_{\eta_{s,t},\eta_s}\sigma_s(Q^t_{\eta_{s,t}\circ\sigma_{s},\eta_t})Q^{s+t}_{\eta_{s+t},\eta_{s,t}}$. From  
\begin{displaymath}
  Q^s_{\eta_{s,t},\eta_s}\sigma_s(Q^t_{\eta_{s,t}\circ\sigma_{s},\eta_t})Q^{s+t}_{\eta_{s+t},\eta_{s,t}}=m(s,t)P^{s+t}_{\eta_{s,t}}
\end{displaymath}
for some $m(s,t)\in\NT$ we find
\begin{displaymath}
   V_{s,t}\theta'_s\otimes \theta'_t=m(s,t) \theta'_{s+t}V'_{s,t}.
\end{displaymath}
Define a new system $\sgii st{V''}$ by  $V''_{s,t}=(\theta'_{s+t})^{-1} V_{s,t}\theta'_s\otimes \theta'_t$, $s,t\in[0,1]$, $s+t\le1$. Then 
\begin{displaymath}
  V''_{r,s+t}\circ(\unit_{\E_r}\otimes V''_{s,t})=(\theta'_{r+s+t})^{-1} V_{r,s+t}\circ(\unit_{\E_r}\otimes V_{s,t})\theta'_r\otimes \theta'_s\otimes \theta'_t
\end{displaymath}
shows by  associativity (in the sense of \ref{eq:associativity Ust}) of  $V$ that $V''$ is  associative. Since  $V'$ is associative too we derive    that $m$ fulfils \ref{eq:multiplier}. By Corollary \ref{cor:multipliersimplex}, there is some $\map f{[0,1]}\NT$ such that $m(s,t)=f(s)f(t)f(s+t)^{-1}$. Defining new unitaries  $\sgi t\theta$, $\theta_t=f(t)^{-1}\theta'_t$, we obtain  an isomorphism between $\sgi t{\E'}$  and $\sgi t{\E}$. Therefore, the  image of the measurable structure of $\sgi t{\E'}$ under $\sgi t\theta$ is a  measurable structure on  $\sgi t{\E}$ compatible with its product structure.  

At the end, it remains  to show that there are operators $Q^t_{\eta,\eta'}$ which fulfil the  conditions \ref{mess0}--\ref{mess4}. Clearly, these conditions are fulfilled in the sense of sets, i.e.\ if we consider the set of all possible operators $Q^t_{\eta,\eta'}$  instead of a single ones. Thus we have just to show that there is a section through these sets which fulfils \ref{mess0}--\ref{mess4} pointwise.

For this goal, we extend the structure a bit and let $P^{s,t}_\eta$ denote the supporting projection of a normal state $\eta$ in  $\A_{[s,t]}$. Note that the set $Z^{s,t}_{\eta,\eta'}=\set{q\in\A_{[s,t]}:qq^*=P^{s,t}_\eta, q^* q=P^{s,t}_{\eta'}}$ is (weakly) compact for all  $s,t\in\NT$ and  $\eta,\eta'\in\P_{[s,t]}$ (in fact it is an image of $\NT$ in $\A_{[s,t]}$). Assume $s_n\tlimitsto{}{n\to\infty}  s,t_n\tlimitsto{}{n\to\infty}  t$, $\P_{[s_n,t_n]}\ni\eta_n\tlimitsto{}{n\to\infty} \eta\in\P_{[s,t]}$, $\P_{[s_n,t_n]}\ni\eta'_n\tlimitsto{}{n\to\infty} \eta'\in\P_{[s,t]}$. Without loss of generality, we may assume that $\eta_n,\eta,\eta_n',\eta'$ correspond to unit vectors $\varphi_n\otimes\psi_n,\varphi\otimes\psi,\varphi_n\otimes\psi'_n,\varphi\otimes\psi'$ factorizing appropriately and fulfilling $\varphi_n\otimes\psi_n\tlimitsto{}{n\to\infty} \varphi\otimes\psi$ as well as $\varphi_n\otimes\psi'_n\tlimitsto{}{n\to\infty} \varphi\otimes\psi'$. Fix   a weak limit point $p$ of $\sequp n{P^{s_n,t_n}_{\eta_n}}$. Clearly, $\unit\ge p\ge0$ and $p\in\A_{s,t}$ by continuity of $(s,t)\mapsto\A_{s,t}$. From uniform boundedness of $\sequp n{P^{s_n,t_n}_{\eta'_n}}$ we find
\begin{displaymath}
  \eta(p)=\lim_{k\to\infty}\eta_{n_k}(P^{s_{n_k},t_{n_k}}_{\eta_{n_k}})=1
\end{displaymath}
such that $p\ge  P^{s,t}_\eta $. On the other side, we find by similar arguments $(\unit-\Pr{\varphi\otimes\psi})p(\unit-\Pr{\varphi\otimes\psi})=0$. From $p\in\A_{s,t}$ it follows that $p=P^{s,t}_\eta $  and   $P^{s_n,t_n}_{\eta_n}\tlimitsto{w}{n\to\infty}  P^{s,t}_\eta $. since $P^{s_n,t_n}_{\eta_n}$ are projections, it follows  $P^{s_n,t_n}_{\eta_n}\tlimitsto{s}{n\to\infty}  P^{s,t}_{\eta}$ and by analogy $P^{s_n,t_n}_{\eta'_n}\tlimitsto{s}{n\to\infty}  P^{s,t}_{\eta'}$.  Without loss of generality, we may assume that $\eta_n,\eta,\eta_n',\eta'$ correspond to unit vectors $\varphi_n\otimes\psi_n,\varphi\otimes\psi,\varphi_n\otimes\psi'_n,\varphi\otimes\psi'$ factorizing appropriately. If $q$ is a weak limit point of $\sequ nq$, $q_n\in Z^{s_n,t_n}_{\eta_n,\eta'_n}$ we know from lower weak semicontinuity of the norm that $\norm q\le1$. Further,  we get from $\norm {q_n}\le1$ 
\begin{displaymath}
  \absq{\scpro{\varphi\otimes\psi}{q\varphi\otimes\psi'}}=\lim_{n\to\infty}   \absq{\scpro{\varphi_n\otimes\psi_n}{q_n\varphi_n\otimes\psi'_n}}= 1
\end{displaymath}
what shows $\norm{q\varphi\otimes\psi}=\norm{\varphi\otimes\psi}$ and $q_n\varphi\otimes\psi\tlimitsto{\norm\p}{n\to\infty}q\varphi\otimes\psi$. On the other hand,
\begin{displaymath}
 P^{s,t}_\eta q=\wlim_{n\to\infty}P^{s_n,t_n}_{\eta_n}q_n =\wlim_{n\to\infty}q_n=q
\end{displaymath}
and, similarly, $q P^{s,t}_{\eta'}=q$ show that $q=z\unit\otimes\Pr{\psi',\psi}$ for some $z\in\NT$.    Since $Z^{s,t}_{\eta,\eta'}$ is a continuous  image of $\NT$ we obtain  $Z^{s_n,t_n}_{\eta_n,\eta'_n}\limitsto{}{n\to\infty} Z^{s,t}_{\eta,\eta'}$ in $\FG_{B_1(\B(\E_1))}$ where $B_1(\B(\E_1))$ is equipped with the weak topology.  By \cite[Theorem 6.6.7]{A:Bee93}  there is a measurable  mapping  $(s,t,\eta,\eta')\mapsto q^{s,t}_{\eta,\eta'}\in Z^{s,t}_{\eta,\eta'}$. 

$\NT$ identified with $[0,1)=\NR/\NZ$ acts in a natural manner on $\P(\E_1)\times\P(\E_1)$ through $(s,\eta,\eta')\mapsto (\eta\circ\sigma_{-s},\eta'\circ\sigma_{-s})$. Since $\gr t\sigma$  is $\sigma$-weakly continuous and periodic this action  is well-defined on $\NT$ and continuous. \cite[Theorem 6.6.7]{A:Bee93} shows  that there is a measurable mapping $\map T{\P(\E_1)\times\P(\E_1)}{\P(\E_1)\times\P(\E_1)}$ such that $T(\eta\circ\sigma_{-s},\eta'\circ\sigma_{-s})=T(\eta,\eta')=T(T(\eta,\eta'))$. Due to compactness of $\NT$ and continuity of $\sigma$  there is another map $\map S{\P(\E_1)\times\P(\E_1)}\NT$ for which $T(\eta,\eta')=(\eta\circ\sigma_{-S(\eta,\eta')},\eta'\circ\sigma_{-S(\eta,\eta')})$ for all $\eta,\eta'\in\P(\E_1)$. \nocite{A:Gli61}These maps are used to define
\begin{displaymath}
 \bar{q}^{s,t}_{\eta,\eta'}=\sigma_{S(\eta,\eta')}\klam{q^{s-S(\eta,\eta'),t-S(\eta,\eta')}_{T(\eta,\eta')}}
\end{displaymath}
such that $\bar{q}^{s+r,t+r}_{\eta\circ\sigma_{-r},\eta'\circ\sigma_{-r}}=\sigma_{r}(\bar{q}^{s,t}_{\eta,\eta'})$. Clearly, $\bar q$ remains measurable.

Further, $\P(\E_1)$ is a Standard Borel space such that there is a measurable total ordering $\prec$ on it. We set
\begin{displaymath}
  \bar{\bar q}^{s,t}_{\eta,\eta'}=\left\{
    \begin{array}[c]{c>{$~}l<$}
\bar{q}^{s,t}_{\eta,\eta'}&if $\eta\circ\sigma_{-S(\eta,\eta')}\prec\eta'\circ\sigma_{-S(\eta,\eta')}$\\
P^{s,t}_\eta& if $\eta=\eta'$\\
 \klam{\bar{q}^{s,t}_{\eta',\eta}}^*&otherwise   
\end{array}
\right.
\end{displaymath}
such that $\bar{\bar q}$ is a measurable section with  
\begin{displaymath}
  \bar{\bar q}^{s+r,t+r}_{\eta\circ\sigma_{-r},\eta'\circ\sigma_{-r}}=\sigma_{r}(\bar{\bar q}^{s,t}_{\eta,\eta'})\dmf{r\in\NR}
\end{displaymath}
and 
\begin{displaymath}
  \bar{\bar q}^{s,t}_{\eta,\eta'}=\klam{\bar{\bar q}^{s,t}_{\eta',\eta}}^*.
\end{displaymath}

Fix again the family   $\sgi t\eta$ used above and set $\eta_{s,t}=\eta_{t-s}\circ \sigma_{-s}$ as well as $\eta_{s,t',t}=\eta_{s,t'}\otimes\eta_{t',t}\otimes\eta_{t,s}$. For $\eta\in\P_{[s,t]}$ we  define 
\begin{displaymath}
  r^{s,t}_{\eta}=\bar{\bar q}^{s,t}_{\eta,\eta_{s,t}}
\end{displaymath}
if $s=\min\set{s':\eta\circ\sigma_{-S'(\eta)}\in\P_{[s',t]}\text{~for some~}t\ne s'}=\min\set{s':\exists t: (s,t)\in\tilde F(\eta\circ\sigma_{-S'(\eta)})}$ where $S'$ is defined as above but from the action of $\sigma$ on $\P(\E_1)$. Observe that Lemma \ref{lem:measurability F(eta)}, measurability of $S'$ and $\gr t\sigma$, and measurability of the map $\FG_\NT\ni Z\mapsto \min Z$ imply that $r$ is again measurable.  Clearly, we can extend this operator family in compliance with 
\begin{displaymath}
  r^{s,t'}_{\eta}r^{t',t}_{\eta}=r^{s,t}_{\eta}\bar{\bar q}^{s,t}_{\eta_{s,t},\eta_{s,t',t}}\dmf{t'\in(s,t)\in I_\NT,\eta\in\P_{[s,t']}\cap\P_{[t',t]}}
\end{displaymath}
and retaining  measurability of $r$.  It is easy to deduce that ${r}^{s+r,t+r}_{\eta\circ\sigma_{r}}=\sigma_{r}({r}^{s,t}_{\eta})$, too.

From this we can derive that $Q^{t}_{\eta,\eta'}={ r}^{0,t}_\eta\klam{{ r}^{0,t}_{\eta'}}^*$ fulfils  all of the relations \ref{mess0}--\ref{mess4} like follows. \ref{mess1} is clear. By construction, we have ${r}^{s,t}_{\eta}({r}^{s,t}_{\eta})^*=P^{s,t}_\eta$ and $({r}^{s,t}_{\eta})^*{r}^{s,t}_{\eta}=P^{s,t}_{\eta_{s,t}}$ such that $P^{s,t}_\eta{r}^{s,t}_{\eta}={r}^{s,t}_{\eta}={r}^{s,t}_{\eta}P^{s,t}_{\eta_{s,t}}$. This shows
\begin{eqnarray*}
  Q^{t}_{\eta,\eta'}Q^{t}_{\eta',\eta''}&=&{  r}^{0,t}_\eta\klam{{ r}^{0,t}_{\eta'}}^*{  r}^{0,t}_{\eta'}\klam{{ r}^{0,t}_{\eta''}}^*={  r}^{0,t}_\eta P^{t}_{\eta_{t}}\klam{{ r}^{0,t}_{\eta''}}^*={  r}^{0,t}_\eta \klam{{ r}^{0,t}_{\eta''}}^*=Q^{t}_{\eta,\eta''},
\end{eqnarray*}
i.e.\ \ref{mess3}. In a similar way, $Q^{t}_{\eta',\eta}=\klam{Q^{t}_{\eta,\eta'}}^*$ and $Q^{t}_{\eta,\eta}=P^t_\eta$ are immediate such that \ref{mess3} implies \ref{mess0}. To prove \ref{mess4}, fix $\eta,\eta'\in\P_{s,s+t}$. Then
\begin{eqnarray*}
  Q^s_{\eta,\eta'}\sigma_s(Q^t_{\eta\circ\sigma_{-s},\eta'\circ\sigma_{-s}})&=&{ r}^{0,s}_{\eta}\klam{{ r}^{0,s}_{\eta'}}^*\sigma_s\klam{{ r}^{0,t}_{\eta\circ\sigma_{-s}}\klam{{ r}^{0,t}_{\eta'\circ\sigma_{-s}}}^*}\\
&=&{ r}^{0,s}_{\eta}\klam{{ r}^{0,s}_{\eta'}}^*{ r}^{s,s+t}_{\eta}\klam{{ r}^{s,s+t}_{\eta'}}^*\\
&=&{ r}^{0,s}_{\eta}{ r}^{s,s+t}_{\eta}\klam{{ r}^{0,s}_{\eta'}{ r}^{s,s+t}_{\eta'}}^*\\
&=&{r}^{0,s+t}_{\eta}\bar{\bar q}^{0,s+t}_{\eta_{0,s+t},\eta_{0,s,s+t}}\klam{{r}^{0,s+t}_{\eta'}\bar{\bar q}^{0,s+t}_{\eta_{0,s+t},\eta_{0,s,s+t}}}^*\\
&=&{r}^{0,s+t}_{\eta}\bar{\bar q}^{0,s+t}_{\eta_{0,s+t},\eta_{0,s,s+t}}\klam{\bar{\bar q}^{0,s+t}_{\eta_{0,s+t},\eta_{0,s,s+t}}}^*\klam{{r}^{0,s+t}_{\eta'}}^*\\
&=&{r}^{0,s+t}_{\eta}P^{s+t}_{\eta_{0,s+t}}\klam{{r}^{0,s+t}_{\eta'}}^*\\
&=&{r}^{0,s+t}_{\eta}\klam{{r}^{0,s+t}_{\eta'}}^*=Q^{s+t}_{\eta,\eta'}.
\end{eqnarray*}
This completes the proof.
\end{proof}
\begin{corollary}
\label{cor:onlyonemeasurablestructure}
  Let $\E=(\sg t\E,(V_{s,t})_{s,t\in\NRp})$ be an algebraic continuous tensor product system. Then  all  measurable structures on $\sg t\E$, making the family $(V_{s,t})_{s,t\in\NRp}$ measurable, give rise to mutual isomorphic (measurable) product systems. 

Further, two (measurable) product systems are isomorphic iff they are isomorphic as algebraic product systems. 

More generally, two measurable product systems $\E$  and $\tilde \E$ are isomorphic iff there is a unitary $\map\theta{\E_1}{\tilde\E_1}$ such that $\theta^*(\A_{0,t})\theta=\tilde\A_{0,t}$ for a dense set of  $t\in(0,1)$.
\end{corollary}
\begin{proof}
  In the above proof we saw that  $\gr t\tau$ induces many measurable structures  turning $\sg t\E$ into a product system, depending on the choice of the vectors $\sgi t\psi$ and of $\map f{[0,1]}\NT$. If $\sg t\E$ is a (measurable) product system, we can even choose these vectors measurably and it is easy to see that $\sgi t{\theta'}$ becomes measurable by this. Clearly, this implies that $m$ is measurable and by   \cite[Corollary of Proposition 2.3]{OP:Arv89}  the corresponding $f$ can be chosen to be measurable. This shows that $\sgi t\theta$ is measurable and the product systems $\E$ and $\tilde\E$ are isomorphic.  The first assertion follows immediately.  

The second assertion follows from the first one and the fact that  any algebraic isomorphism $\map \theta\E{\E'}$ of product systems transfers a measurable structure from $\E$ to a measurable structure on $\E'$ compatible with the multiplication of $\E'$. 

For the proof of the third part observe that normality of $\theta^*(\cdot)\theta$ as well as continuity of the map $t\mapsto\A_{0,t}$ implies that  $\theta^*(\A_{s,t})\theta=\tilde\A_{s,t}$ for all $s,t\in[0,1]$. From this equality it follows that the product systems $\E'$ and $\tilde \E'$ constructed in the proof of Theorem \ref{th:intrinsicmeasurability} are the same. This shows that $\E$  and $\tilde \E$ are algebraically isomorphic and the preceding results complete the proof. 
\end{proof}
\begin{remark}
  \label{rem:onedimensional product systems}
From these results we see that for onedimensional product systems, all algebraic product structures  and all measurable structures are algebraically isomorphic.
\end{remark}

\subsection{Product Systems of $W^*$-Algebras}
\label{sec:psW*}

We want to follow the above developed lines a bit further and  summarize our view on  continuous tensor product systems of Hilbert spaces through continuous tensor product systems of $W^*$-algebras. Recall that a $W^*$-algebra is the abstract version of a von Neumann algebra, characterised as a $C^*$-algebra with unit which is  the dual of a Banach space. Normal states on a $W^*$-algebra are elements of this predual and isomorphisms of $W^*$-algebras map normal states into normal ones. The tensor product of $W^*$-algebras $\B,\B'$ is defined as the spatial tensor product (we can represent  both $W^*$-algebras faithfully on separate Hilbert spaces) \cite[section 1.22]{Sak71}. This abstract notation is useful for algebras like $L^\infty(\M)$ for which there is no canonical Hilbert space they act on.
\begin{definition}
\label{def:ps of vnas}
 We call a  family $(\B_{s,t})_{(s,t)\in I_{0,\infty}}$ of  $W^*$-algebras  \emph{continuous tensor product system of $W^*$-algebras} if 
 \begin{condition}{$W^*$}
     \item \label{W*1}For all $(s,t)\in I_{0,\infty}$, $r\in\NRp$ the $W^*$-algebra $\B_{s,t}$ is  isomorphic to $\B_{s+r,t+r}$ under $\beta^{s,t}_r$. Further    $\beta^{s+r,t+r}_{r'}\circ \beta^{s,t}_r=\beta^{s,t}_{r+r'}$ for all $(s,t)\in I_{0,\infty}$, $r,r'\in\NRp$.
\item \label{W*2}For all $(r,s),(s,t)\in I_{0,\infty}$ there is an isomorphism  $\map{\gamma_{r,s,t}}{\B_{r,s}\otimes \B_{s,t}}{\B_{r,t}}$  such that
  \begin{displaymath}
    \gamma_{s,s',t}(\Id_{\B_{s,s'}}\otimes\gamma_{s',t',t})=\gamma_{s,t',t}(\gamma_{s,s',t'}\otimes\Id_{\B_{t',t}})\dmf{(s,t)\in I_{0,\infty},(s',t')\in I_{s,t}}
  \end{displaymath}
and 
\begin{displaymath}
 \gamma_{s+r,s'+r,t+r}\circ(\beta^{s,s'}_r\otimes\beta^{s',t}_r)=\beta^{s,t}_r\circ \gamma_{s,s',t}\dmf{(s,s'),(s',t)\in I_{0,\infty},r\in\NRp}.
\end{displaymath}
 \end{condition} 
\end{definition}
 \begin{definition}
Two states $\eta,\eta'$ on a $W^*$-algebra $\A$ are \emph{equivalent} if the GNS-re\-pre\-sen\-ta\-tions $\pi_\eta$ and $\pi_{\eta'}$ are unitarily equivalent, i.e.\ there is a unitary $\map {U_{\eta,\eta'}}{H_\eta}{H_\eta'}$ with
\begin{equation}
\label{eq:interwiningGNS}
  U_{\eta,\eta'}\pi_\eta(a)=\pi_{\eta'}(a) U_{\eta,\eta'}\dmf{a\in\A}.
\end{equation}
A \emph{type of normal states} on  $\A$ is an equivalence class of normal states under this relation. 
\end{definition}
\begin{definition}
\label{def:stationary factorizing type of states}
Let $(\B_{s,t})_{(s,t)\in I_{0,\infty}}$ be a product system of $W^*$-algebras. A \emph{stationary factorizing type of states} is a family $(\S_{s,t})_{(s,t)\in I_{0,\infty}}$ of  sets of   normal states on $\B_{s,t}$ with the following properties.  
\begin{condition}{S}
\item \label{h0} All $\eta\in\S_{s,t}$ belong to the same type of states on $\B_{s,t}$.
\item \label{h3}For all $s,r\in\NRp$ the relation  $\eta\in\S_{0,s}$ is true iff $\eta\circ\beta^{0,s}_r\in\S_{r,r+s}$. 
\item \label{h2}For all $s,t\in\NRp$, $\eta_s\in\S_{0,s}$,   $\eta_t\in\S_{0,t}$ the state   $(\eta_s\otimes\eta_t\circ(\beta^{0,t}_s)^{-1})\circ\gamma_{0,s,s+t}^{-1}$ belongs to $\S_{0,s+t}$. 
\item \label{h1} For all $t\in\NRp$, $\eta\in\S_{0,t}$ the GNS Hilbert space $H^t_\eta$ is  separable. 
\end{condition}
\end{definition}
\begin{remark}
  The simplest situation is that  $\S_{s,t}$ is a complete  type of states on $\B_{s,t}$.   The given,  slightly more general definition is useful, e.g., in  Example \ref{ex:Tsirelson type III}.

Further, the separability condition  \ref{h1} is not  necessary from general reasons.  We add it here for compactness of notation.
\end{remark}
In the sequel, we use the short-hand notation $\eta_s\otimes\eta_t$ for the state $(\eta_s\otimes\eta_t\circ(\beta^{0,t}_s)^{-1})\circ\gamma_{0,s,s+t}^{-1}\in\S_{0,s+t}$ and likewise $a_s\otimes b_t$ for the operator $\gamma_{0,s,s+t}(a_s\otimes\beta^{0,t}_s(b_t))$.
\begin{lemma}
\label{lem:associativity tensor}
  These $\otimes$ operations are associative, i.e.
  \begin{eqnarray*}
    a_r\otimes (b_s\otimes c_t)&=&    (a_r\otimes b_s)\otimes c_t\dmf{r,s,t\ge0,a_r\in\B_{0,r},b_s\in\B_{0,s},c_t\in\B_{0,t}}\\
    \eta_r\otimes (\eta_s\otimes \eta_t)&=&    (\eta_r\otimes \eta_s)\otimes \eta_t\dmf{r,s,t\ge0,\eta_r\in\S_{0,r},\eta_s\in\S_{0,s},\eta_t\in\S_{0,t}}
  \end{eqnarray*}
\end{lemma}
\begin{proof}
  The first relation  is due to \ref{W*2} and the semigroup property of $\beta$ given by \ref{W*1}. The second relation follows from the first one by  applying the left hand side to the left hand one and the right hand side to the right hand one as well as noting that the elements $a_r\otimes(b_s\otimes c_t)$ are $\sigma$-weakly total in $\B_{0,r+s+t}$.\end{proof}
Define for all $s,t\in\NRp$, $\eta_s\in\S_{0,s}$, $\eta_t\in\S_{0,t}$ unitaries  $\map{V_{s,t}^{\eta_s,\eta_t}}{H_{\eta_s}\otimes H_{\eta_t}}{H_{\eta_s\otimes\eta_t}}$ by extension of
\begin{displaymath}
  V_{s,t}^{\eta_s,\eta_t}[a_s]_{\eta_s}\otimes[b_t]_{\eta_t}=[a_s\otimes b_t]_{\eta_s\otimes\eta_t}\dmf{a_s\in\B_{0,s},b_t\in\B_{0,t}}.
\end{displaymath}
\begin{proposition}
  \label{prop:tensorproductW*=ps}
Let $(\B_{s,t})_{(s,t)\in I_{0,\infty}}$ be a product systems of $W^*$-algebras,  $(\S_{s,t})_{(s,t)\in I_{0,\infty}}$ a stationary factorizing type of states on it and $\U=(\U^t)_{t\in\NRp}$,  $\U^t=(U^t_{\eta,\eta'})_{t\in\NRp,\eta,\eta'\in\S_{0,t}}$ families of unitaries     $\map{U^t_{\eta,\eta'}}{H^t_\eta}{H^t_{\eta'}}$ fulfilling \ref{eq:productUetaeta'}, \ref{eq:interwiningGNS} and 
  \begin{equation}
\label{eq:factorization Utetaeta'}
    U^{s+t}_{\eta_s\otimes\eta_t,\eta'_s\otimes\eta'_t}V_{s,t}^{\eta_s,\eta_t}=  V_{s,t}^{\eta'_s,\eta'_t}U^s_{\eta_s,\eta'_s}\otimes U^t_{\eta_t,\eta'_t}\dmf{s,t\in\NRp,\eta_r,\eta'_r\in\S_{0,r},r=s,t}.
  \end{equation}

Then the family  $\E^{\S,\U}=\sg t\E$ with  the  Hilbert spaces \begin{math}
  \E_t=H(\B_{0,t},\S_{0,t},\U^t)
\end{math}   equipped with the multiplication encoded by  $\sg{s,t}V$, 
\begin{displaymath}
\left(  V_{s,t}(\psi^s)\otimes(\tilde \psi^t)\right)_{\eta_s\otimes\eta_t}= V_{s,t}^{\eta_s,\eta_t}\psi^s_{\eta_s}\otimes\tilde \psi^t_{\eta_t}\dmf{s,t\in\NRp},
\end{displaymath}
is an algebraic  product system. 
\end{proposition}
\begin{proof}
By Lemma \ref{lem:defH(S)}, $\E_t$ is well-defined due to the properties of $\klam{U^t_{\eta,\eta'}}_{\eta,\eta'\in\S_{0,t}}$.  By \ref{h1}, $\E_t$ is separable. $V_{s,t}$ is well-defined because of \ref{eq:factorization Utetaeta'} as can be seen from
\begin{eqnarray*}
  \left( V_{s,t}(\psi^s)\otimes(\tilde \psi^t)\right)_{\eta'_s\otimes\eta'_t}&=& V_{s,t}^{\eta'_s,\eta'_t}\psi^s_{\eta'_s}\otimes\tilde \psi^t_{\eta'_t}\\
&=& V_{s,t}^{\eta_s,\eta_t}U^s_{\eta_s,\eta'_s}\psi^s_{\eta_s}\otimes U^t_{\eta_t,\eta'_t}\tilde \psi^t_{\eta_t}\\
&=&U^{s+t}_{\eta_s\otimes\eta_t,\eta'_s\otimes\eta'_t}V_{s,t}^{\eta_s,\eta_t}\psi^s_{\eta_s}\otimes\tilde\psi^t_{\eta_t}.
\end{eqnarray*}
Associativity of the resulting product is proved if we can show that 
\begin{displaymath}
  V_{r,s+t}^{\eta_r,\eta_s\otimes\eta_t}\circ(\unit_{H_{\eta_r}}\otimes V_{s,t}^{\eta_s,\eta_t})= V_{r+s,t}^{\eta_r\otimes\eta_s,\eta_t}\circ( V_{r,s}^{\eta_r,\eta_s}\otimes\unit_{H_{\eta_t}})
\end{displaymath}
for $r,s,t\ge0$, $\eta_p\in\S_{0,p}$, $p=r,s,t$. For this sake, select $a\in\B_{0,r}$, $b\in\B_{0,s}$, $c\in\B_{0,t}$ and derive using Lemma \ref{lem:associativity tensor}
\begin{eqnarray*}
\lefteqn{ V_{r,s+t}^{\eta_r,\eta_s\otimes\eta_t}\circ(\unit_{H_{\eta_r}}\otimes V_{s,t}^{\eta_s,\eta_t} ) [a]_{\eta_r}\otimes[b]_{\eta_s}\otimes[c]_{\eta_t}}&=& V_{r,s+t}^{\eta_r,\eta_s\otimes\eta_t}[a]_{\eta_r}\otimes(V_{s,t}^{\eta_s,\eta_t}[b]_{\eta_s}\otimes[c]_{\eta_t})
= V_{r,s+t}^{\eta_r,\eta_s\otimes\eta_t}[a]_{\eta_r}\otimes([b\otimes c)]_{\eta_s\otimes\eta_t})\\
&=& [a\otimes(b\otimes c)]_{\eta_r\otimes(\eta_s\otimes\eta_t)})= [(a\otimes b)\otimes c]_{(\eta_r\otimes\eta_s)\otimes\eta_t}\\
&=& V_{r+s,t}^{\eta_r\otimes\eta_s,\eta_t}\circ( V_{r,s}^{\eta_r,\eta_s}\otimes\unit_{H_{\eta_t}}) [a]_{\eta_r}\otimes[b]_{\eta_s}\otimes[c]_{\eta_t}.
\end{eqnarray*}
Thus    $\E^{\S,\U}$ is an  algebraic product system.
\end{proof}

It would be interesting to know under what circumstance $\E$ has a consistent  measurable structure. We have the following result on a necessary condition. 
\begin{proposition}
\label{prop:measurability W^*}
 There is a one parameter  automorphism group $\gr t{\smash{\hat\beta}}$ of $\B_{0,1}$, determined  by 
\begin{equation}
\label{eq:defbetaTt}
  \hat\beta_t(a_{1-t}\otimes b_{t})= b_{t}\otimes a_{1-t}\dmf{a_{1-t}\in\B_{0,1-t}, b_t\in\B_{0,t}}.
\end{equation}
If, for some family $\U$, the product system $\E$ of Proposition \ref{prop:tensorproductW*=ps} has a consistent measurable structure and $\pi_\eta$ is injective for some and thus all $\eta\in\S_{0,1}$ then $\gr t{\smash{\hat\beta}}$ is $w-*$-measurable.
\end{proposition} 
\begin{proof}
   $\hat\beta_t$ is normal since (omitting several  normal shifts $\beta^r_{s,t}$) $\hat\beta_t=\gamma_{0,t,1-t}\mathcal{F}_t\gamma_{0,1-t,t}^{-1}$, where  $\mathcal{F}_t$ is the unitarily implemented, thus normal, flip from $\mathcal{B}_{0,t}\otimes \mathcal{B}_{t,1}$ to  $\mathcal{B}_{0,1-t}\otimes \mathcal{B}_{1-t,1}$. It is a short computation similar
  to the proof of Proposition \ref{prop:tau_tiscontinuous} that $\hat
  \beta$ is a one parameter group of automorphisms.
  
  On the other side, if the GNS representations $\pi_\eta$,
  $\eta\in\S_{0,1}$, are faithful, measurability of $\gr
  t{\smash{\hat\beta}}$ is  necessary since it is the restriction
  of $\gr t\sigma$.
\end{proof}
\begin{remark}
  We conjecture that the condition is also sufficient but at the moment we have no proof for  this. But, for the special cases of type $\mathrm{I}$ factors $\B_{s,t}$ and abelian  $W^*$-algebras  Theorem \ref{th:intrinsicmeasurability} and Section \ref{sec:Polish Product System} respectively solve these problems. 

Similarly, we do not know anything about existence and uniqueness of $\E$ in its dependence on $\B$ and $\S$.
\end{remark}
\begin{example}
  A first example of this result was used in the proof of Theorem \ref{th:intrinsicmeasurability}. There the choice was $\B_{s,t}=\A_{s,t}$, $\S_{0,t}=\P_t$ with the tedious choice of $U^t_{\eta,\eta'}$ vis.\ $Q^t_{\eta,\eta'}$ to obtain  both  measurability and compatibility. The derived product system was the original one. Measurability of $\gr t\sigma$ yielded a consistent measurable structure on $\E$.   
\end{example}
\begin{example}
\label{ex:ps from faithful states}
  As another application, let $\E$ be a product system and   use as  $\S_{s,t}$ the sets of faithful normal states on $\A_{s,t}$. The unitaries  $U_{\eta,\eta'}$ are  defined by $U_{\eta,\eta'}[a]_\eta=[a{\varrho'}^{1/2}\varrho^{-1/2}]$  on a suitable  set of operators $a$ such that $[a]_\eta$ is dense in $H_\eta$. The above corollary gives us a product system $\sg t{\tilde\E}$, $\tilde\E_t=(H_\eta)_{\eta\in\S_{0,t}}$. By Example \ref{ex:GNSnormalstateB(H)} we know that $H_{\eta}\cong\E_t\otimes\E^*_t$. Thus, the product system $\tilde \E$ is just $\E\otimes \E^*$, where $\E^*=(\sg t{\E^*},\sg{s,t}V)$. Observe that $\E_s^*\otimes\E_t^*=(\E_s\otimes\E_t)^*$ by definition of the tensor product and $V_{s,t}$ is unitary on $\E_s^*\otimes\E_t^*$ too.

 Whether $\tilde \E\cong\E\otimes \E$ is also valid, as suggested by the discussion in Example \ref{ex:GNSnormalstateB(H)}, is not clear. This depends essentially on the question whether $\E\cong\E^*$, i.e.\ on the construction of antiunitary conjugations  $\sg tC$ on $\sg t\E$ such that $C_{s+t}=C_s\otimes C_t$.

Similarly, we could consider for a general product system of $W^*$-algebras the set $\S_{s,t}$  of faithful normal states on $\A_{s,t}$ and base the operators $U^t_{\eta,\eta'}$ on \textsc{Connes'} cocycle \cite{OP:Con73}. 
\end{example}
The derivation of product systems from product systems of  $L^\infty$ spaces is the subject of  section \ref{sec:randommeasuresandincrements}.  As preparation, we compute the flip $\gr t\tau$ for the product system  $\E^{\S,\U}$.
\begin{lemma}
  Suppose $\A$ is a $C^*$-algebra, $\eta$ a state on $\A$ and $\beta$ an automorphism of $\A$. Then there is a unique unitary  $\map{u^{\beta,\eta}}{H_\eta}{H_{\eta\circ\beta^{-1}}}$ defined by
\begin{displaymath}
  u^{\beta,\eta}[a]_\eta=[\beta(a)]_{\eta\circ\beta^{-1}}\dmf{a\in\A}.
\end{displaymath}
\end{lemma}
\begin{proof}
Clearly,
\begin{eqnarray*}
  \scpro{u^{\beta,\eta}[a]_\eta}{u^{\beta,\eta}[b]_\eta}_{\eta\circ\beta^{-1}}&=&\scpro{[\beta(a)]_{\eta\circ\beta^{-1}}}{[\beta(b)]_{\eta\circ\beta^{-1}}}_{\eta\circ\beta^{-1}}\\
&=&{\eta\circ\beta^{-1}}(\beta(a)^*\beta(b))\\
&=&\eta(\beta^{-1}(\beta(a^*b)))=\eta(a^*b)= \scpro{[a]_\eta}{[b]_\eta}_\eta.
\end{eqnarray*}
Since the range of $u^{\beta,\eta}$ is total, the operator is unitary.  
\end{proof}
\begin{lemma}
  Suppose $\S_{0,1}$ is a maximal type of states. Then  for all $t\in[0,1]$ and  $\eta\in \S_{0,1}$ also $\eta\circ\hat\beta_t\in\S_{0,1}$.
\end{lemma}
\begin{proof}
  We know that there are $\eta_t\in\S_{0,t}$, $\eta'_{1-t}\in\S_{0,1-t}$. Further, by the preceding lemma, we know that $\eta\sim\eta'$ implies $\eta\circ\hat\beta_t\sim\eta'\circ\hat\beta_t$. We conclude from that
  \begin{displaymath}
   \eta\sim \eta'_{1-t}\otimes\eta_t\sim \eta_t\otimes \eta'_{1-t}=(\eta'_{1-t}\otimes\eta_t)\circ\hat\beta_t\sim \eta\circ\hat\beta_t
  \end{displaymath}
and the proof is complete.
\end{proof}
\begin{lemma}
\label{lem:computation tau_t}
  Suppose additionally to the assumptions in  Proposition \ref{prop:tensorproductW*=ps} that 
  \begin{displaymath}
    U_{\eta\circ\hat\beta_t,\eta'\circ\hat\beta_t}u^{\hat\beta_{-t},\eta}=u^{\hat\beta_{-t},\eta'}U_{\eta,\eta'}.
  \end{displaymath}
Then
\begin{displaymath}
(\tau_t\psi)_\eta=U_{\eta\circ\hat\beta_{-t},\eta}u^{\hat\beta_{-t},\eta}\psi_\eta
\end{displaymath}
\end{lemma}
\begin{proof}
  Suppose $\psi\in\E_1$ fulfils for some $b_t\in\B_{0,t}$,
  $a_{1-t}\in\B_{0,1-t}$, $\eta_t\in\S_{0,t}$,
  $\eta'_{1-t}\in\S_{0,1-t}$ the relation $\psi_{\eta'_{1-t}\otimes
    \eta_t}=[a_{1-t}\otimes b_t]_{\eta'_{1-t}\otimes \eta_t}$. Then,
  denoting $\F$ the flip between $H_{\eta_t}$ and $H_{\eta'_{1-t}}$,
  we derive
\begin{eqnarray*}
  (\tau_t\psi)_{\eta_t\otimes \eta'_{1-t}}&=&V_{t,1-t}^{\eta_t,\eta_{1-t}}\F (V_{1-t,t}^{\eta'_{1-t},\eta_t})^*\psi_{\eta'_{1-t}\otimes\eta_t}\\
&=&V_{t,1-t}^{\eta_t,\eta'_{1-t}}\F (V_{1-t,t}^{\eta'_{1-t},\eta_t})^*[a_{1-t}\otimes b_t]_{\eta'_{1-t}\otimes\eta_t}\\
&=&V_{t,1-t}^{\eta_t,\eta'_{1-t}}\F [a_{1-t}]_{\eta'_{1-t}}\otimes[b_t]_ {\eta_t}\\
&=& V_{t,1-t}^{\eta_t,\eta_{1-t}}[b_t]_ {\eta_t}\otimes [a_{1-t}]_{\eta'_{1-t}}\\
&=&[b_t\otimes a_{1-t}]_{\eta_t\otimes \eta'_{1-t}}\\
&=&[\hat\beta_t( a_{1-t}\otimes b_t)]_{\eta_t\otimes \eta'_{1-t}}.
\end{eqnarray*}
We conclude for $\psi_{\eta'_{1-t}\otimes
  \eta_t}=[a]_{\eta'_{1-t}\otimes \eta_t}$ that
$\tau_t\psi_{\eta_t\otimes
  \eta'_{1-t}}=[\hat\beta_t(a)]_{\eta_t\otimes \eta'_{1-t}}$ or
$\tau_t\psi_{\eta_t\otimes
  \eta'_{1-t}}=u^{\hat\beta_t,\eta'_{1-t}\otimes\eta_t}\psi_{\eta'_{1-t}\otimes\eta_t}$.
From this the assertion follows easily. 
\end{proof}
\begin{remark}
  \label{rem:productsystemW*-algebras}
We could derive the  tensor product structure  even if the  algebras $\B_{r,s}$ and $\B_{s,t}$  commute, generate  $\B_{r,t}$, and there is a state $\eta$ such that $\eta(b_{r,s}b_{s,t})=\eta(b_{r,s})\eta(b_{s,t})$. In fact, such factorisations imply essentially that  $\B_{r,t}=\B_{r,s}\otimes \B_{s,t}$ \cite[Exercise 1(b) to section IV.5]{OP:Tak79}. An application of that result could be similar results for $C^*$-algebras, where the tensor product need not be uniquely determined. Use of $C^*$-algebras  would make it  more easy to define states, but we refrain from complicating things  here without need.

The above introduced structure could equally be defined with indeces in $[0,1]$ and using Proposition \ref{prop:psbyE_1}. This would allow us to consider all $W^*$-algebras as von Neumann subalgebras of a fixed $\B(\H)$. We could simplify  notations a little bit by considering  the algebras $\B_{0,t}$ only, but we think the above presentation shows a more clear division between shift (vis.\ $\beta$) and tensor product (vis.\ $\gamma$). In the proof of Theorem \ref{th:intrinsicmeasurability} we found already this separation useful.

Continuous tensor product systems of $W^*$-algebras should  be a basic structure for the construction of (generalised) quantum Markov processes and quantum Markov random fields, see \cite{Q:AF01a,Q:LS01,Q:Lie01b}.
\end{remark}

\subsection{Product systems and Unitary Evolutions}
\label{sec:ps2timeshifts}

Now we want to show how to use  product systems of $W^*$-algebras and the corresponding Hilbert spaces $H(\A,\S,\U)$ to construct for all product systems an $E_0$-semigroup it belongs to according to  equation \ref{eq:def ps from E_0(Arveson)}. Before we state the main result, we need some preparations.

Again, we use the sets $\Q_{s,t}$ as defined in Lemma \ref{lem:measurable P_t}. 
\begin{lemma}
  If $\E\ne\sgp t\NC$ then $\P(\E_1)\supsetneq\bigcup_{t\in(0,1)}\P_t$.
\end{lemma}
\begin{proof}
  We assume  that $\P(\E_1)=\bigcup_{t\in(0,1)}\P_t$. This would imply that  
  \begin{displaymath}
    \P(\E_1)=\bigcup_{s,s'\in\ND, 0<s<s'<1  
}\Q_{s,s'}. 
  \end{displaymath}
Obviously, each  $\Q_{s,s'}$ is closed and  we conclude from  Baires theorem of categories that one $\Q_{s,s'}$ contains an open ball. Denote its  centre by $\eta$, which should correspond to a vector $\psi$. Then we know 
\begin{displaymath}
  \scpro{\tilde\psi}{ab\tilde\psi}\scpro{\tilde\psi}{\tilde\psi}=\scpro{\tilde\psi}{a\tilde\psi}\scpro{\tilde\psi}{b\tilde\psi}\dmf{a\in\A_{0,s},b\in\A_{s',1}}
\end{displaymath}
for all $\tilde\psi$ in a ball around $\psi$. Inserting $\tilde\psi=\psi+\varepsilon\psi'$ with an arbitrary vector $\psi'$ and $\varepsilon>0$ small enough into this equation we arrive at
\begin{displaymath}
  \scpro{\psi+\varepsilon\psi'}{ab(\psi+\varepsilon\psi')}\scpro{\psi+\varepsilon\psi'}{\psi+\varepsilon\psi'}=\scpro{\psi+\varepsilon\psi'}{a(\psi+\varepsilon\psi')}\scpro{\psi+\varepsilon\psi'}{b(\psi+\varepsilon\psi')}.
\end{displaymath}
Comparing the terms of first order in $\varepsilon$ gives us
\begin{eqnarray*}
\lefteqn{(\scpro{\psi'}{ab\psi}+\scpro\psi{ab\psi'}) \scpro\psi\psi+\scpro\psi{ab\psi}(\scpro{\psi'}{\psi}+\scpro\psi{\psi'})}&=&(\scpro{\psi'}{a\psi}+\scpro\psi{a\psi'}) \scpro\psi{b\psi}+\scpro\psi{a\psi}(\scpro{\psi'}{b\psi}+\scpro\psi{b\psi'}).
\end{eqnarray*}
 Using this identity for $\psi'$ too ($\psi'$ was arbitrary) we find
\begin{eqnarray*}
\scpro{\psi'}{ab\psi}\scpro\psi\psi+\scpro\psi{ab\psi}\scpro{\psi'}{\psi}&=&\scpro{\psi'}{a\psi} \scpro\psi{b\psi}+\scpro\psi{a\psi}\scpro{\psi'}{b\psi}\emf{a\in\A_{0,s},b\in\A_{s',1}}
\end{eqnarray*}
what leads to
\begin{displaymath}
ab\psi\scpro\psi\psi+\scpro\psi{ab\psi}\psi= \scpro\psi{b\psi}a\psi+\scpro\psi{a\psi}b\psi\dmf{a\in\A_{0,s},b\in\A_{s',1}}.
\end{displaymath}
Multiplication of this  equation by another $b'\in\A_{s',1}$ and using the same formula for $ab'b\psi$ gives us with $\norm\psi=1$, $\eta(a)=\scpro{\psi}{a\psi}$,
\begin{eqnarray*}
  ab'b\psi+\eta(ab)b'\psi&=& \eta(b)ab'\psi+\eta(a)b'b\psi\\
\lefteqn{  -\eta(ab'b)\psi+\eta(b'b)a\psi +\eta(a)b'b\psi+\eta(ab)b'\psi}&=&- \eta(b)\eta(ab')\psi+\eta(b)\eta(b')a\psi+\eta(b)\eta(a)b'\psi+\eta(a)b'b\psi\\
\lefteqn{  -\eta(a)\eta(b'b)\psi+\eta(b'b)a\psi +\eta(a)\eta(b)b'\psi}&=&- \eta(b)\eta(a)\eta(b')\psi+\eta(b)\eta(b')a\psi+\eta(b)\eta(a)b'\psi\\
 \eta(b'b)(a\psi -\eta(a)\psi) &=& \eta(b')\eta(b)(a\psi -\eta(a)\psi)
\end{eqnarray*}
Since $\A_{s',1}$ is noncommutative, there are certainly $b,b'\in\A_{s',1}$ with $ \eta(b'b)\ne \eta(b')\eta(b)$ what implies that $a\psi =\eta(a)\psi$ for all $a\in\A_{0,s}$. Consequently, $\eta(a^*a)=\eta(a^*)\eta(a)$  which is impossible on the same grounds. This contradiction completes the proof.
\end{proof}

Originally, the following  was a combined result of \cite[Theorem 3.4]{Arv89}, \cite[Theorem 3.4]{OP:PR89} and \cite[Corollary 5.17]{OP:Arv90}. We aim here at a more explicit construction in terms of a Hilbert space $H(\A,\S,\U)$. 
Type I factors $\A\subseteq\B(\H)$ are von Neumann algebras isomorphic to some $\B(\H')$, i.e.\ there is a unitary $\map u{H}{H'\otimes\H''}$ for another Hilbert space $\H''$ such that $u\A u^*=\B(\H')\otimes\unit_{\H''}$. 
\begin{theorem}
  \label{th:E_0-semigroup}
  For all product systems $\E$ there is a separable Hilbert space $\H$ and a strongly continuous one-parameter group $\gr t\gamma$ of unitaries  together with two commuting type $\mathrm{I}$ factors $\A_{0)}$ and $\A_{[0}$ with $\B(\H)=\A_{0)}\vee\A_{[0}$ such that 
\begin{displaymath}
     \gamma_s^*\mathcal{A}_{[0}\gamma_s\subseteq\A_{[0}\dmf{s\ge0}.
\end{displaymath}
Further,
\begin{displaymath}
  \E'_t=\set{u\in\A_{[0}:\gamma_tb\gamma_t^*u=ub\forall b\in\A_{[0}}\dmf{t\ge0}
\end{displaymath}
together with the multiplication $u_s\otimes u_t=u_su_t$ is a product system  isomorphic to $\E$.
\end{theorem}
\begin{proof}
Again, we will construct $\H$ as some $H(\A,\S,\U)$. First, we define the corresponding $C^*$-algebra $\A$ as a so-called quasilocal algebra.

  We consider the $W^*$-algebras $\pa st\A$ given by  representations on  various Hilbert spaces $\E_r\otimes\E_{r'}$, $r,r'\in\NN$, as follows
  \begin{displaymath}
    \A_{s,t}=
    \left\{
      \begin{array}[c]{ll}
\underbrace{\unit_{\E_{s+r}}\otimes\B(\E_{t-s})\otimes\unit_{\E_{-t}}}_{\subseteq\B(\E_r)}\otimes\unit_{\E_{r'}}&\text{~if~}-r\le s<t<0\\
\underbrace{\unit_{\E_{s+r}}\otimes\B(\E_{-s})}_{\subseteq\B(\E_r)}\otimes\underbrace{\B(\E_t)\otimes\unit_{\E_{r'-t}}}_{\subseteq\B(\E_{r'})}&\text{~if~}-r\le s<0<t\le r'\\
\unit_{\E_r} \otimes\underbrace{\unit_{\E_{s}}\otimes\B(\E_{t-s})\otimes\unit_{\E_{r'-t}}}_{\subseteq\B(\E_{r'})} & \text{~if~} 0\le s<t\le r'   \end{array}
      \right.
  \end{displaymath}
All these embeddings to $\B(\E_r\otimes\E_{r'})$ are consistent with the natural embeddings 
\begin{displaymath}
  \B(\E_{r_1})\otimes\B(\E_{r'_1})\hookrightarrow \unit_{\E_{r_2-r_1}}\otimes\B(\E_{r_1})\otimes\B(\E_{r'_1})\otimes\unit_{\E_{r'_2-r'_1}}\subseteq\B(\E_{r_2})\otimes\B(\E_{r'_2})\dmf{a\in\B(\E_1)}
\end{displaymath}
for $r_1\le r_2$, $r_1'\le r_2'$. Thus there is a unique $C^*$-algebra $\A=\overline{\bigcup_{s<t}\A_{s,t}}$. This $C^*$-algebra can also be thought of as the (unique) infinite tensor product $\bigotimes_{p\in\NZ}\B(\E_1)$. The latter structure gives us injections $\sequz pj$, $\map{j_p}{\B(\E_1)}\A$ such that the images of $j_p,j_{p'}$, which are $\A_{p,p+1}$ and $\A_{p',p'+1}$, commute for $p\ne p'$, and $\bigcup_{p\in\NZ}j_p(\B(\E_1))$ generates $\A$. Further, there is a unique discrete automorphism group $\sequz n\Sigma$ determined  via
\begin{displaymath}
  \Sigma_nj_p(a)=j_{p+n}(a)\dmf{n,p\in\NZ}
\end{displaymath}
and the requirement that $\map{\Sigma_n}{\A_{m,m'}}{\A_{m+n,m'+n}}$ should be $\sigma$-strongly continuous. 

Next we want to extend $\sequz n\Sigma$ to an automorphism group $\gr t\Sigma$. Define
\begin{displaymath}
  \Sigma_tj_p(a)=    \left\{
      \begin{array}[c]{ll}
j_{p+\lfloor t\rfloor}(\sigma_t(a))&a\in\A_{0,1-t+\lfloor t\rfloor}\\
j_{p+\lfloor t\rfloor+1}(\sigma_t(a))&a\in\A_{1-t+\lfloor t\rfloor,1}
  \end{array}
      \right.,
\end{displaymath}
where $\gr t\sigma$ is the shift group defined in \ref{eq:defshiftautomorphisms}. Again, the additional requirement that $\map{\Sigma_t}{\A_{r,s}}{\A_{r+t,s+t}}$ should be $\sigma$-strongly continuous fixes $\Sigma_t$.  One obtains immediately $\Sigma_t^{-1}=\Sigma_{-t}$. Further, we obtain for  $0\le s,t$, $s+t\le 1$
\begin{displaymath}
  \Sigma_s( \Sigma_t(j_p(a)))=    \left\{
      \begin{array}[c]{ll}
j_{p}(\sigma_{s+t}(a))&a\in\A_{0,1-t-s}\\
j_{p+1}(\sigma_{s+t}(a))&a\in\A_{1-t-s,1}
  \end{array}
      \right.
 =\Sigma_{s+t}(j_p(a)).
\end{displaymath}
In the case $s+t\ge1$ the calculations are similar. This shows in general $\Sigma_s\circ\Sigma_t=\Sigma_{s+t}$.

If $\eta\in\P(\E_1)$ is a pure state there is a natural (locally normal) state $\eta^\infty$ on $\A$ given by  extension of
\begin{displaymath}
  \eta^\infty(a_{-n,m})=\eta^{\otimes (n+ m)}(a_{-n,m})\dmf{n,m\in\NN,a_{n,m}\in\A_{-n,m}}
\end{displaymath}
with regard to $\A_{-n,m}\cong\B(\E_n)\otimes\B(\E_m)\cong\B(\E_1)^{\otimes n+m}$.
The GNS representation space $H_\eta$ of $\eta^\infty$ is naturally isomorphic to an infinite tensor product Hilbert space. Choose $\psi\in\E_1$ with $\eta=\scpro\psi{\p\psi}$ and define a linear space
\begin{displaymath}
  \H_0={\mathrm{lh}\set{\bigotimes_{p\in\NZ}h_p:h_p=\psi\text{~for all but a finite number of~}p\in\NZ}}
\end{displaymath}
equipped with  the unique inner product such that 
\begin{displaymath}
  \scpro{\smash{\bigotimes_{p\in\NZ}}h_p}{\smash{\bigotimes_{p\in\NZ}}h'_p}=\prod_{p\in\NZ}\scpro{h_p}{h'_p}.
\end{displaymath}
This inner product is well-defined and $\H_0$ is a pre-Hilbert space. Its completion is denoted $\bigotimes_{p\in\NZ}^\psi\E_1=\overline{\H_0}$. Obviously, $\bigotimes_{p\in\NZ}^\psi\E_1$ is separable. The isomorphy between  $H_\eta$ and $\bigotimes_{p\in\NZ}^\psi\E_1$ is based on the fact that $\bigotimes_{p\in\NZ}^\psi\E_1\cong\bigotimes_{p\in\NZ,p<n}^\psi\E_1\otimes\E_{m-n}\otimes\bigotimes_{p\in\NZ,p>m}^\psi\E_1 $ naturally for all $n<m\in\NZ$. With this identification the maps $[a_{n,m}]_{\eta^\infty}\mapsto\bigotimes_{p\in\NZ,p<n}^\psi\psi\otimes(a_{n,m}\psi^{\otimes m-n})\bigotimes_{p\in\NZ,p>m}^\psi\psi$ form an inductive sequence which extends to $\A$. 

Next we want to establish operators $\map{U_{\eta,\tilde\eta}}{H_\eta}{H_{\tilde\eta}}$ which obey \ref{eq:productUetaeta'}. For that reason choose automorphisms  $\sigma_{\eta,\tilde\eta}$ on $\B(\E_1)$ such that $\tilde\eta\circ\sigma_{\eta,\tilde\eta}=\eta$. Further properties need to be fulfilled, but again we collect all of them at the end of the proof, see \ref{e00}--\ref{e03}  and show how to fulfil them.   Now the maps $\sigma_{\eta,\tilde\eta}^{\otimes{m-n}}$ on $\A_{n,m}$ possess a unique extension to an  automorphism $\Sigma_{\eta,\tilde\eta}$ of $\A$. Clearly, ${\tilde\eta}^\infty\circ\Sigma_{\eta,\tilde\eta}=\eta^\infty$. Setting $U_{\eta,\tilde\eta}[a]_{\eta^\infty}=[\Sigma_{\eta,\tilde\eta}(a)]_{{\tilde\eta}^\infty}$ we obtain
\begin{eqnarray*}
  \scpro{U_{\eta,\tilde\eta}[a]_{\eta^\infty}}{U_{\eta,\tilde\eta}[b]_{\eta^\infty}}_{\tilde\eta}&=&\scpro{[\Sigma_{\eta,\tilde\eta}(a)]_{\eta^\infty}}{[\Sigma_{\eta,\tilde\eta}(b)]_{\eta^\infty}}_{\tilde\eta}={\tilde\eta}^\infty(\Sigma_{\eta,\tilde\eta}(a^*)\Sigma_{\eta,\tilde\eta}(b))\\
&=&{\tilde\eta}^\infty(\Sigma_{\eta,\tilde\eta}(a^*b))=\eta^\infty(a^*b)=\scpro{[a]_{\eta^\infty}}{[b]_{\eta^\infty}}_\eta.
\end{eqnarray*}
Since $\Sigma_{\eta,\tilde\eta}$ is an automorphism the image of $U_{\eta,\tilde\eta}$ is full and  $U_{\eta,\tilde\eta}$ is unitary. Further, the relation  $\sigma_{\tilde\eta,\tilde\eta'}\circ\sigma_{\eta,\tilde\eta}=\sigma_{\eta,\tilde\eta'}$  would imply \ref{eq:productUetaeta'}. To get the Hilbert space of interest, we restrict our attention to the set $\P^\infty=\set{\eta^\infty:\eta\in\P'}$ of states on $\A$ where
 \begin{displaymath}
   \P'=\set{\eta\in\P(\E_1):\eta\notin\P^2\text{~and, if~}\eta\in\P_{1/2}\text{~then~}\eta=\eta\circ\sigma_{1/2}}.
 \end{displaymath}
In the following, we use also the notation  $\P'_t=\P_t\cap\P'$. Now  $\H=H(\A,\P^\infty,\U)$ is well-defined by Lemma \ref{lem:defH(S)}.

Now we want to  transfer the shift given by $\gr t\Sigma$ to $\H=H(\A,\P^\infty,\U)$. We define for $\eta\in\P_{1-t}'$
\begin{displaymath}
  (\gamma_t([a_{\tilde\eta}]_{\tilde\eta^\infty})_{\tilde\eta\in\P'})_{\eta^\infty\circ\Sigma_{-t}}=[\Sigma_t(a_\eta)]_{\eta^\infty\circ\Sigma_{-t}} \dmf{t\in\NR}.
\end{displaymath}
Suppose now for all $\eta,\tilde\eta\in\P'_{1-t}$ that $\sigma_t\circ\sigma_{\eta,\tilde\eta}=\sigma_{\eta\circ\sigma_{-t},\tilde\eta\circ\sigma_{t}}\circ\sigma_t$ and $\sigma_{\eta,\tilde\eta}\cong\sigma^{1-t}_{\eta,\tilde\eta}\otimes\sigma^t_{\eta,\tilde\eta}$  according to  $\B(\E_1)\cong\B(\E_{1-t})\otimes\B(\E_t)$. Then we obtain for $t\in[0,1]$, $p\in\NZ$ and $a\in\A_{1-t}$
\begin{eqnarray*}
  \Sigma_{\eta\circ\sigma_{-t},\tilde\eta\circ\sigma_{-t}}\circ\Sigma_t(j_p(a))&=&\Sigma_{\eta\circ\sigma_{-t},\tilde\eta\circ\sigma_{-t}}(j_p(\sigma_t(a)))\\
&=&j_p(\sigma_{\eta\circ\sigma_{-t},\tilde\eta\circ\sigma_{-t}}(\sigma_t(a)))=j_p(\sigma_t(\sigma_{\eta,\tilde\eta}(a)))
\end{eqnarray*}
and for  $a_t\in\A_{1-t,1}$
\begin{eqnarray*}
  \Sigma_{\eta\circ\sigma_{-t},\tilde\eta\circ\sigma_{-t}}\circ\Sigma_t(j_p(a))&=&\Sigma_{\eta\circ\sigma_{-t},\tilde\eta\circ\sigma_{-t}}(j_{p+1}(\sigma_t(a)))\\
&=&j_{p+1}(\sigma_{\eta\circ\sigma_{-t},\tilde\eta\circ\sigma_{-t}}(\sigma_t(a)))=j_{p+1}(\sigma_t(\sigma_{\eta,\tilde\eta}(a))).
\end{eqnarray*}
Local $\sigma$-strong continuity of $\gr t\Sigma$ and the fact that $\Sigma_n$ commutes with both $\Sigma_{\eta\circ\sigma_{-t},\tilde\eta\circ\sigma_{-t}}$ and $\Sigma_t$, show
\begin{displaymath}
  \Sigma_{\eta\circ\sigma_{-t},\tilde\eta\circ\sigma_{-t}}\circ\Sigma_t=\Sigma_t\circ\Sigma_{\eta,\tilde\eta}\dmf{t\in\NR,\eta,\tilde\eta\in\P'_{1-t}}
\end{displaymath}
Consequently, $U_{\eta\circ\sigma_{-t},\tilde\eta\circ\sigma_{-t}}[\Sigma_t(a)]_{\eta^\infty\circ\sigma_{-t}}=[\Sigma_{\eta,\tilde\eta}\circ\Sigma_t(a)]_{\eta^\infty\circ\sigma_{-t}}=[\Sigma_t\circ\Sigma_{\eta,\tilde\eta}(a)]_{\eta^\infty\circ\sigma_{-t}}$ and we have shown that the definition of $\gr t\gamma$ does not depend on the choice of $\eta$. Further,
\begin{eqnarray*}
\lefteqn{\scpro {\gamma_t([a_\eta]_{\eta^\infty})_{\eta\in\P'}}{\gamma_t([b_\eta]_{\eta^\infty})_{\eta\in\P'}}_{H( \A,\P^\infty,\U)}}&=&\scpro{[\Sigma_t(a_\eta)]_{(\eta\circ\sigma_{-t})^\infty}}{[\Sigma_t(b_\eta)]_{(\eta\circ\sigma_{-t})^\infty}}\\
&=&\eta^\infty\circ\Sigma_{-t}(\Sigma_t(a^*_\eta)\Sigma_t(b_\eta))\\
&=&\eta^\infty(a^*_\eta b_\eta)=\scpro {([a_\eta]_{\eta^\infty})_{\eta\in\P'}}{[b_\eta]_{\eta^\infty})_{\eta\in\P'}}_{H(\A,\P,\U)}
\end{eqnarray*}
shows that $\gamma_t$ is in fact well-defined and unitary. $\gamma_t^*=\gamma_{-t}$ and $\gamma_{s+t}=\gamma_s\gamma_t$ are immediate from the corresponding properties of $\gr t\Sigma$ and $\gr t\sigma$.

Choosing $\sigma_{\eta,\tilde\eta}$ measurably depending on $\eta,\tilde\eta$, we find that
\begin{eqnarray*}
t\mapsto \scpro {([a_\eta]_{\eta^\infty})_{\eta\in\P'}}{\gamma_t([b_\eta]_{\eta^\infty})_{\eta\in\P'}}_{H(\A,\P^\infty,\U)}&=&\scpro{[a_{\eta\circ\sigma_t}]_{(\eta\circ\sigma_t)^\infty}}{[\Sigma_t(b_\eta)]_{(\eta\circ\sigma_{-t})^\infty}}\\
&=&\eta^\infty\circ\Sigma_{-t}(\sigma_{\eta,\eta\circ\sigma_t}(a^*_\eta)\Sigma_t(b_\eta))\\
&=&\eta(\sigma_{\eta,\eta\circ\sigma_t}(a^*_\eta) b_\eta)
\end{eqnarray*}
is measurable due to $\sigma$-strong continuity of $\gr t\sigma$. Thus $\gr t\gamma$ is weakly measurable and due to the fact that $\bigotimes_{p\in\NZ}^\psi\E_1$ and thus $H(\A,\P^\infty,\U)$ is separable, $\gr t\gamma$ is continuous. 

Of course, $\A=\A^-\vee\A^+$ where $\A^-=\overline{\bigcup_{t<0}\A_{t,0}}$ and $\A^+=\overline{\bigcup_{t<0}\A_{0,t}}$ and both subalgebras are commuting. Let $H_\eta^+$ being the GNS Hilbert space of the restriction of $\eta$ to  $\A^+$ and define $H^-_\eta$ similarly. Now consider the map $\map{V^0_\eta}{H_\eta}{H_\eta^-\otimes H_{\eta}^+}$,
\begin{displaymath}
  V^0_\eta[a^-a^+]_{\eta^\infty}=[a^-]_{\eta^\infty}\otimes [a^+]_{\eta^\infty}\dmf{a^+\in\A^+,a^-\in\A^-}.
\end{displaymath}
From
\begin{eqnarray*}
  \scpro{V^0_\eta[a^-a^+]_{\eta^\infty}}{V^0_\eta[b^-b^+]_{\eta^\infty}}&=&\scpro{[a^-]_{\eta^\infty}\otimes [a^+]_{\eta^\infty}}{[b^-]_{\eta^\infty}\otimes [b^+]_{\eta^\infty}}\\&=&\eta^\infty((a^-)^*b^-)\eta^\infty((a^+)^*b^+)\\
&=&\eta^\infty((a^-)^*b^-(a^+)^*b^+)\\&=&\eta^\infty((a^-a^+)^*b^-b^+)\\
&=&\scpro{[a^-a^+]_{\eta^\infty}}{[b^-b^+]_{\eta^\infty}}
\end{eqnarray*}
we derive that $V^0_\eta$ is  well-defined and unitary. Further, $\Sigma_{\eta,\tilde\eta}(\A^\pm)=\A^\pm$ shows that $U_{\eta,\tilde\eta}V^0_\eta=V^0_{\tilde\eta}(U_{\eta,\tilde\eta}\otimes U_{\eta,\tilde\eta})$ and $U_{\eta,\tilde\eta}\H^\pm_\eta=\H^\pm_{\tilde\eta}$. The latter gives us Hilbert spaces $\H^\pm$ such that  there is  a unitary $\map{V^0}{H(\A,\P,\U)}{\H^-\otimes \H^+}$. Further, set $\A_{[0}=V^0(\unit_{\H^-}\otimes\B(\H^+))(V^0)^*$ and see $\gamma_t\A_{[0}\gamma_t^*\subset \A_{[0}$. Thus  $\gr t\gamma$ reduces to an $E_0$-semigroup $\sg t{\alpha^\gamma}$ on this  type $\mathrm{I}$ factor.

Now we want to establish that $\E'$ is isomorphic to $\E$. By Proposition \ref{prop:psbyE_1}, it is enough to establish isomorphy on $\sgi t{\E'}$ and $\sgi t\E$. To understand $\gr t\gamma$ better, fix some $\eta=\scpro\psi{\p\psi}\in\P'$ such that $\H\cong \bigotimes_{p\in\NZ}^\psi\E_1$ as remarked above. If $0<t\le1$  and  $\eta_{1-t}\in\P'_{1-t}$, consider $a=\prod_{p\in\NZ}\sigma_{\eta_{1-t},\eta}(a_{0,{1-t}}^pa_{{1-t},1}^p)$ for operators $a_{0,{1-t}}^p\in\A_{p,p+{1-t}}$, $a_{{1-t},1}^p\in\A_{p+{1-t},p+1}$ with only finitely many operators  different from $\unit$. By the GNS construction, we find
\begin{displaymath}
  [a]_{\eta^\infty}\cong\bigotimes_{p\in\NZ}^\psi \underbrace{x_{0,{1-t}}^p\otimes x_{{1-t},1}^p}_{\in\E_{1-t}\otimes\E_t\cong\E_1}
\end{displaymath}
where $x_{0,{1-t}}^p\cong[a^p_{0,{1-t}}]_{\eta_{1-t}}$. Now it is easy to see that 
\begin{eqnarray*}
  \Sigma_{\eta\circ\sigma_t,\eta}(\Sigma_t(a))&=&\prod_{p\in\NZ}\sigma_{\eta_{1-t}\circ\sigma_t,\eta}(\sigma_t(a_{0,{1-t}}^pa_{{1-t},1}^{p-1}))\\
&=&\prod_{p\in\NZ}\sigma_{\eta_{1-t},\eta}(\sigma_{\eta_{1-t}\circ\sigma_t,\eta_{1-t}}(\sigma_t(a_{0,{1-t}}^pa_{{1-t},1}^{p-1})))
\end{eqnarray*}
such that
\begin{equation}
\label{eq:squeezed shift}
  \gamma_t [a]_{\eta^\infty}\cong\bigotimes_{p\in\NZ}^\psi \underbrace{U_tx_{{1-t},1}^{p-1}\otimes U_{1-t}x_{0,{1-t}}^p}_{\in\E_t\otimes\E_{1-t}\cong\E_1}
\end{equation}
for some unitaries $U_t\in\B(\E_t)$, $U_{1-t}\in\B_{E_{1-t}}$ such that $\sigma_{\eta_{1-t}\circ\sigma_t,\eta_{1-t}}$ corresponds to $U_t\otimes U_{1-t}$. This shows, regarding $\H^+\cong \bigotimes_{p\in\NN}^\psi\E_1\cong\E_t\otimes\E_{1-t} \otimes\bigotimes_{p\in\NN\setminus\set0}^\psi\E_1$ that $\alpha^\gamma_t(\B(\H^+))=\unit_{\E_t}\otimes\B(\E_{1-t})\otimes\B(\bigotimes_{p\in\NN\setminus\set0}^\psi\E_1)$. Now it  is immediate to see that  $\E'_t$ and $\E_t$ are  isomorphic also under preservation of tensor products, at least for $0\le s,t \le1$, $s+t\le 1$ since setting $\H_1=\bigotimes_{p\in\NN\setminus\set0}^\psi\E_1$
\begin{displaymath}
  \alpha^\gamma_{s+t}(\B(\H^+))=\unit_{\E_{s+t}}\otimes\B(\E_{1-s-t})\otimes\B(\H_1)=\unit_{\E_s}\otimes\underbrace{\unit_{\E_t}\otimes\B(\H_1)}_{=\alpha^\gamma_s(\B(\H^+))}
\end{displaymath}
This shows that $\E$ and $\E'$ are algebraically isomorphic. Using Corollary \ref{cor:onlyonemeasurablestructure} we conclude they are isomorphic as product system.

All the above  derivations are true provided we find a family  $(\sigma_{\eta,\tilde\eta})_{\eta,\tilde\eta\in\P}$ of  $\sigma$-strongly continuous automorphisms of $\B(\E_1)$ such that
\begin{enumerate}
\item \label{e00}$\tilde\eta\circ\sigma_{\eta,\tilde\eta}=\eta$ for all $\eta,\tilde\eta\in\P'$.
\item \label{e01}$\sigma_{\tilde\eta,\tilde\eta'}\circ\sigma_{\eta,\tilde\eta}=\sigma_{\eta,\tilde\eta'}$ for all $\eta,\tilde\eta,\tilde\eta'\in\P'$.
\item \label{e02}$\sigma_{1-t}\circ\sigma_{\eta,\tilde\eta}=\sigma_{\eta\circ\sigma_t,\tilde\eta\circ\sigma_t}\circ\sigma_{1-t}$ for all $t\in[0,1]$, $\eta,\tilde\eta\in\P'_t$.
\item \label{e04} If $t\in[0,1]$ and $\eta,\tilde\eta\in\P'_t$ then, regarding  $\E_1\cong\E_t\otimes \E_{1-t}$,  $\sigma_{\eta,\tilde\eta}$ factorizes as $\sigma_{\eta,\tilde\eta}\cong\sigma^t_{\eta,\tilde\eta}\otimes\sigma^{1-t}_{\eta,\tilde\eta}$ for two $\sigma$-strongly continuous automorphisms $\sigma^t_{\eta,\tilde\eta}$ on $\B(\E_t)$ and $\sigma^{1-t}_{\eta,\tilde\eta}$ on $\B(\E_{1-t})$.
\item\label{e03} The map $(\eta,\tilde\eta)\mapsto\sigma_{\eta,\tilde\eta}$ is pointwise $\sigma$-strongly measurable.
\end{enumerate}
We construct now  automorphisms in the form  $\sigma_{\eta,\tilde\eta}=V_{\eta,\tilde\eta}\p V_{\eta,\tilde\eta}^*$ for unitaries $V_{\eta,\tilde\eta}$. Before we proceed to the definition of $V_{\eta,\tilde\eta}$, we need additional ingredients.

Suppose we are given a  complete orthonormal system $\sequ n{e}$ in some $\E_t$.  We derive from a pure state  $\eta$ a vector $f_0^\eta$ with $\eta(\Pr{f_0^\eta})=1$. Then  the Gram-Schmidt orthogonalisation procedure for the sequence $\sequ n{f^\eta}$ with the chosen $f_0^\eta$ and $f_{n+1}^\eta=e_n$, $n\in\NN$,  yields   a new  complete orthonormal system   $\sequ n{e^\eta}$. Obviously, if $\sequ n{e}$ and $f_0^\eta$ are chosen measurably then $\sequ n{e^\eta}$ is measurable too. We will assume this is true by choosing a measurable section of complete orthonormal systems through $\sgi t\E$ in advance.

Further, $\P(\E_{1/2})$ is a  Standard Borel spaces. Thus there exists a measurable total ordering  $\succ_\P$ on it, which we fix. Moreover, let $\map l\NN{\NN\times\NN}$ be a fixed one-to-one map with $l(0)=(0,0)$ and denote by $l_1,l_2$ its components.

For every state $\eta\in\P'$ we define now an orthonormal system $\sequ n{\bar e^\eta}$ like follows. First we consider  $\eta=\eta^1_t\otimes\eta^2_{1-t}\in\P_t$. If  $t<1/2$, we define $\bar e^\eta_n=e^{\eta^1_t}_{l_1(n)}\otimes e^{\eta^2_{1-t}}_{l_2(n)}$ whereas for $t>1/2$ we set $\bar e^\eta_n=e^{\eta^1_t}_{l_2(n)}\otimes e^{\eta^2_{1-t}}_{l_1(n)}$. In the case $t=1/2$, where $\eta^1_{1/2}=\eta^2_{1/2}$, we define  $\bar e^\eta_n=e^{\eta^1_{1/2}}_{l_1(n)}\otimes e^{\eta^2_{1/2}}_{l_2(n)}$.  Lastly, we set $\bar e_n^\eta=e_n^\eta$ for all $\eta\in\P(\E_1)\setminus\P^1$. Then our goal is accomplished by setting 
\begin{displaymath}
  V_{\eta,\tilde\eta} \bar e^\eta_n=\bar e^{\tilde\eta}_n\dmf{n\in\NN}.
\end{displaymath}
Namely, \ref{e00} follows from $l(0,0)=0$ and the properties of $\sequ n{e^\eta}$. \ref{e01} is easy to verify. The relations \ref{e02} and  \ref{e04} follow from $V_{\eta,\tilde\eta}e^{\eta^1_t}_{n_1}\otimes e^{\eta^2_t}_{n_2}=e^{{\tilde\eta}^1_t}_{n_1}\otimes e^{{\tilde\eta}^2_t}_{n_2}$. Finally, measurability follows from  measurability of $F(\eta)$ and measurability of $\eta\mapsto\bar e^\eta_n$.  This completes the proof.
\end{proof}
\begin{remark}
  A look at \ref{eq:squeezed shift} shows  that the shift on $\H$ is not just a usual shift on an infinite product Hilbert space. This is due to the fact that $\eta$ is not a $\gr t\Sigma$ invariant state on $\A$. May be, with this idea one is able to do the construction without reference to a whole bundle $H_\eta$ of Hilbert spaces just  on a single Hilbert space $\bigotimes^\psi_{p\in\NZ}\E_1$. Nevertheless, it should not be easy to fulfil the compatibility relations on  $U_t,U_{1-t}$ needed in equation \ref{eq:squeezed shift} for $\gr t\gamma$ to be a unitary group.   Here the flexible structure of $H(\A,\P^\infty,\U)$ helps a lot.
\end{remark}
\begin{remark}
\label{rem:factorization sets as invariant}
  The sets $F(\eta)$ have an interesting  structure. In fact, it is easy to see that $\D_\E=\set{F(\eta):\eta\in\P(\E_1)}$ is an  invariant of $\E$. E.g., it is easy to work out that in the type $I_n$ case $\D_\E=\set{(0,1)}$ for $n=0$ and  $\D_\E=\FG_{(0,1)}$ for $n\ge1$ what implies the same formulae for  type $\mathrm{II}_{n}$ product systems since the latter contain the former as subsystems.  From the above work it seems that at least in the type $\mathrm{III}$ case $\D_\E$ contains a lot of information.  We leave here just the indication that this structure might be interesting for further study. 
\end{remark}

\subsection{Additional Results on Measurability}
\label{sec:measurability}

At the end of this section we provide some results referee to above. The following existence of measurable direct integral representations   was needed in section \ref{sec:direct integrals}. For a Standard Borel space $Y$, let $L^\infty(Y)$ denote  the $C^*$-algebra of bounded Borel measurable functions on $Y$ normed by the supremum norm. A representation $\pi$ of $L^\infty(Y)$ is called \emph{normal} if $\pi(\sup_{n\in\NN} f_n)=\sup_{n\in\NN}\pi(f_n)$ for all bounded  increasing sequences $\sequ nf\subset L^\infty(Y)$.
\begin{lemma}
  \label{lem:measurabilitydirectintegral}
  Suppose $X,Y$ are  Standard Borel spaces and $(H_x)_{x\in X}$ is a measurable family of separable Hilbert spaces carrying a family  $(\pi_x)_{x\in X}$ of normal representations of $L^\infty(Y)$ such that $x\mapsto \pi_x(a)$ is measurable for all $a\in L^\infty(Y)$. Then there is a measurable family $(H^y_x)_{x\in X,y\in Y}$ of Hilbert spaces and a measurable family $(\mu_x)_{x\in X}$ of probability measures on $Y$ such that  $H_x\cong\int^\oplus_Y\mu_x(\d y)H^y_x$.
\end{lemma}
\begin{proof}
  We describe first an algorithm to derive a direct integral decomposition  $\H\cong\int^\oplus_Y\mu(\d y)H^y$ from  a representation $\pi$ of $L^\infty(Y)$ on $\H$ without reference to Zorn's lemma like \cite{Par92,Sak71,BR87}. So let us  fix a faithful normal state $\eta$ on  $\B(\H)$ and a sequence $\sequ n{\psi^0}\subset\set{\psi\in\H:\norm\psi=1}$ which is dense in $\set{\psi\in\H:\norm\psi=1}$. Then we set $k=0$, $H_0=\H$ and iterate
{\newenumi{(\arabic{enumi})}
\begin{enumerate}
  \item Choose the first $n_k\in\NN$ such that 
    \begin{displaymath}
      \eta(\Pr{\psi^k_{n_k}})\ge 1/2\max\set{\eta(\Pr\psi):\psi\in H_k,\norm\psi=1}.
    \end{displaymath}
  \item Compute $H_{k+1}$ as orthogonal complement of $\pi(L^2(\M))\psi^k_{n_k}$.
  \item  Define a new sequence  $\sequ n{\psi^{k+1}}\subset\set{\psi\in H_{k+1}:\norm\psi=1}$ by projecting $\psi^k_{n_k}$ onto $H_{k+1}$, normalising it if the result is a non zero vector or discard this vector otherwise.
  \item Increment $k$ and go to step (1).
  \end{enumerate}
} Depending on whether for some $k\in\NN$ it happens $H_k=\set0$ or not, there result  finitely or infinitely many iterations. For the sake of simplicity,  let us  assume an infinite iteration. From \cite[Proposition 7.3]{Par92} we know that $\pi\restriction_{H_k\ominus H_{k+1}}$ is isomorphic to the canonical representation $\pi_k$ of $L^\infty(\M)$ on $L^2(\mu_k)$, $\mu_k$ being defined through $\int f\d\mu_k=\scpro{\psi^k_{n_k}}{\pi(f)\psi^k_{n_k}}$, $f\in L^\infty(\M)$. 

We want to prove now that $\bigoplus_{k\in\NN}(H_k\ominus H_{k+1})=\H$, or $\Pr{H_k}\limitsto s{k\to\infty}0$.  Define the Hilbert subspace $H_\infty=\bigcap_{k\in\NN}H_k$ by $\Pr{H_\infty}=\slim_{k\to\infty}\Pr{H_k}$, which exists by monotony of $\sequ kH$ and suppose $H_\infty\ne\set0$. Then we obtain from faithfulness of $\eta$ some $\psi\in H_\infty$, $\norm \psi=1$ with $\eta(\Pr\psi)=c>0$. Since $H_\infty\subseteq H_k$ for all $k\in\NN$, we obtain $\eta(\Pr{\psi^k_{n_k}})\ge c/2$ for all $k$. Now, the sequence $\sequp k{\psi^k_{n_k}}$ is orthonormal and $\sum_{k\in\NN}\Pr{\psi^k_{n_k}}\le\bigvee_{k\in\NN}(\unit-\Pr{H_k})=\unit-\Pr{H_\infty}\le\unit$. Thus
\begin{displaymath}
1=\eta(\unit)\ge\eta(\sum_{k\in\NN}\Pr{\psi^k_{n_k}})\ge\sum_{k\in\NN}\eta(\Pr{\psi^k_{n_k}})\ge\sum_{k\in\NN}c/2=\infty
\end{displaymath}
which is a contradiction. Thus $H_\infty=\set0$.

As a consequence, $\pi\restriction_\H$ is isomorphic to the direct sum $\bigoplus_{k\in\NN}\pi\restriction_{H_k\ominus H_{k+1}}\cong \bigoplus_{k\in\NN}\pi_k$. $\pi_k$ depends only on the measure $\mu_k$. Taking $\mu=\sum_{k\in\NN}2^{-k}\mu_k$, we get  $\bigoplus_{k\in\NN}L^2(\mu_k)\cong L^2(\mu')$, with the $\sigma$-finite measure $\mu'=\sum_{k\in\NN}2^{-k}\int \mu_k(\d y)\delta_{y,k}$ being defined on $Y\times\NN$.  Under these unitary equivalences, $\pi$ transforms into the canonical representation of $L^\infty(Y)$ on $L^2(\mu')$. Disintegrating $\mu'$ with respect to the first component in $Y\times \NN$ (which has distribution $\mu$) gives a kernel $p$. Like in section \ref{sec:direct integrals} we derive $H_y=L^2(p(y,\cdot))$.

At the end, we want to  use this algorithm to obtain measurable direct integral decompositions. So fix a measurable family $(H_x)_{x\in X}$ of Hilbert spaces. The measurable structure allows us to choose a measurable section $(\eta_x)_{x\in X}$ of faithful normal states on $\B(H_x)$. Then step (1)--(4) above lead to  Hilbert subspaces $H^x_k$, $k=1,\dots$. Since there exists a countable generating algebra inside $L^\infty(Y)$ and the representations $\pi^x$ are normal, the corresponding  measures $\mu^x_k$ depend measurably on $x$, i.e.\ $x\mapsto \int f\d\mu^x_k$ is measurable for all $f\in L^\infty(Y)$. Now disintegration of measures on Standard Borel spaces is a measurable operation, see e.g.\ \cite[Lemma 6.11]{Tsi00b}. This shows that the family $(H_x^y)_{x\in X,y\in Y}$ is measurable and  $\mu^x$ is measurable too. The proof is complete.
\end{proof}
In the beginning we mentioned that our definition of measurability in product systems differs from  that one used in literature  previously. Here  connect both approaches.
\begin{lemma}
  \label{lem:both measurability definitions}
Let $\H$ be a Standard Borel space with a measurable projection $\map p{\H}{(0,\infty)}$ such that $p^{-1}(\set t)=\H_t$, $t>0$ are Hilbert spaces and $\H$ is Borel isomorphic to $(0,\infty)\times l^2(\NN)$ under a  fibrewise unitary map. Suppose further that there is a measurable associative  multiplication on $\H$ such that all restrictions to $\H_s\times\H_t$ are bilinear and correspond to a unitary   $\map{W_{s,t}}{\H_s\otimes\H_t}{\H_{s+t}}$. Set $\H_0=\NC$, define  $W_{0,t},W_{t,0}$ for all $t\ge0$ as trivial operations and equip $\sg t\H$ with the measurable structure generated by the preimages of the measurable sections through $(0,\infty)\times l^2(\NN)$.  Then $(\sg t\H,\sg{s,t}W)$  is a product system.

Conversely, if $\E$ is a product system with $\dim\E_t=\infty$, $t>0$, then there  is a Standard Borel structure on  $\H=\bigcup_{t>0}\E_t$ such that  the fibre map $\map p\H{(0,\infty)}$, $p|\E_t\equiv t$, is measurable  and the multiplication determined  by  $x_sy_t=V_{s,t}x_s\otimes y_t$, $s,t>0$, is  measurable and associative.  Moreover, the above procedure leads again to the  product system $\E$. 
\end{lemma}
\begin{proof}
The properties of the  measurable structure on $\sg t\H$ are  provided in \cite[Proposition 1.15]{Arv89}. Thus it is sufficient to show \ref{eq:associativity Ust} and measurability of $\sg{s,t}W$. Suppose $x_s\in\H_s$, $y_t\in\H_s$, $z_r\in\H_r$. Then by associativity of the multiplication
\begin{eqnarray*}
  W_{s,t+r}(\unit_{\H_s}\otimes W_{t,r}) x_s\otimes y_t\otimes z_r&=& W_{s,t+r}x_s\otimes y_t z_r =x_s(y_t z_r)\\
&=&(x_sy_t) z_r= W_{s+t,r}x_s y_t\otimes z_r\\
&=&W_{s+t,r}(W_{s,t}\otimes\unit_{\H_r}) x_s\otimes y_t\otimes z_r.
\end{eqnarray*}
 Using  the remarks before Lemma 6.18 in \cite{Arv89}, multiplication is measurable iff for measurable sections $\sg tx, \sg ty,\sg tz$ $(s,t)\mapsto x_sy_t$ is measurable, i.e.\ $(s,t)\mapsto \scpro{x_sy_t}{z_{s+t}}$ is measurable. This is exactly the same as measurability of $\sg{s,t}W$.

 For the proof of the converse direction, let $\E$ be a product system  with $\dim\E_t=\infty$, $t>0$. By the Gram-Schmidt orthogonalisation procedure we obtain measurable sections $\sgs t{h^n}$, $n\in\NN$ through $\sgs t\E$ such that for all $t>0$ $(h^n_t)_{n\in\NN}$ is a complete orthonormal system in $\H_t$. Sending $h_t\in\H_t$ to $(t,\sequ n{\scpro{h_t}{e^n_t}})$ we get the desired bundle isomorphism of $\H=\bigcup_{t>0}\E_t$ to $(0,\infty)\times l^2(\NN)$ which we can make an Borel isomorphism by imposing the  preimage  Borel structure on $\H$. Multiplication is introduced by $x_sy_t=V_{s,t}x_s\otimes y_t$, $s,t>0$. A similar calculation as above shows that \ref{eq:associativity Ust} makes this multiplication  associative. Measurability was already shown above. 

It is clear that constructing a product system out of $\H$ we arrive again at $\E$. This completes the proof.
\end{proof}
\section{Construction of Product Systems from General Measure Types}
\label{sec:randommeasuresandincrements}
We want to shed some more light on the construction of  product systems by $\E_t=L^2(\M_t)$ with $\M_t$ being a measure type on a suitable space $X_t$. Our key  observation  is that $L^2(\M)$ is just a bundle of GNS Hilbert spaces $(H_\mu)_{\mu\in\M}$. So all we need to look at are conditions for using Proposition \ref{prop:tensorproductW*=ps} and conditions for measurability like mentioned in Proposition \ref{prop:measurability W^*} and given in Theorem \ref{th:intrinsicmeasurability}. These conditions are build in the first subsection in  such a way to be applicable in the following subsections, dealing with examples for this construction.  
\subsection{General Results}
\label{sec:Polish Product System}
Before we present the general theorem, we need some preparatory lemmata which are more or less folklore. We add them since the necessary conditions derived in these lemmata show that the construction  is the most general one we can hope for in the present context.  

First we want to know  which  spaces of square integrable functions are separable.
\begin{lemma}
  \label{lem:L^2separable=Polish}
  For a probability  measure $\mu$  on some measurable space $(X,\XG)$ the space  $L^2(\mu)$ is separable iff  $L^\infty(X,\mu)$ is isomorphic as $W^*$-algebra to some $L^\infty(X',\mu')$ for  a Borel probability measure  $\mu'$ on a Standard Borel  space $X'$.  
\end{lemma}
\begin{proof}
 If $L^2(\mu)$ is separable  choose sets $\sequ nY\subset\XG$ such that $\sequp n{\chfc{Y_n}}$ is total (in $L^2(\mu)$) in the set $\set{\chfc Y:Y\in\XG}$. Then $\map iX{\set{0,1}^\NN=X'}$, $i(x)=(\chfc{Y_n}(x))_{n\in\NN}$ is a map from $X$ to the Polish space $X'$, which is clearly measurable. Set $\mu'=\mu\circ i^{-1}$ and define $\map j{L^\infty(E,\mu')}{L^\infty(X,\mu)}$, $j(f)=f\circ i$.  Assume $f=0$ $\mu'$-a.s.\ Then 
  \begin{displaymath}
    \int \mu(\d x)\abs{j(f)(x)}=\int\mu'(\d x')\abs{f(x')}=0,
  \end{displaymath}
 i.e.\ $j(f)=0$ $\mu$-a.s.\ and  $j$ is well-defined. Moreover, $j$ is a normal homomorphism. Similar we  prove that it is injective. Assume it is not onto. Then $j(L^\infty(X',\mu'))$ is not dense  in $L^2(\mu)$. Since    $j(L^\infty(X',\mu'))$ contains the total set $\sequp n{\chfc{Y_n}}$, this  is a contradiction.

Since any Standard Borel space contains a countable algebra of sets generating the Borel sets, the reverse direction is obvious. This completes the proof.
\end{proof}
Further, the operator algebraic notion of equivalence and the measure theoretic one coincide.
\begin{proposition}
  \label{prop:both equivalences}
Let  $\mu,\mu'$ be two probability measures on a measurable space $(X,\XG)$. Then they are equivalent as measures iff they are equivalent as states on $L^\infty(\M)$ for any measure type $\M$, $\mu_1,\mu_2\ll\M$.
\end{proposition}
\begin{proof}
  Since the unitary $U_{\mu,\mu'}$ is multiplication by a measurable function (see equation \ref{eq:defUmumu'}), the \emph{only if} direction is immediate. For the proof of the \emph{if} direction, assume the states induced by $\mu,\mu'$ are equivalent. The corresponding intertwining operator $U$ defines a vector $\psi=    U[\unit]_\mu\in L^2(\mu')$. Then 
  \begin{displaymath}
    \scpro\psi{\pi_{\mu'}(f)\psi}=\scpro{U[\unit]_\mu}{\pi_{\mu'}(f)U[\unit]_\mu}=\scpro{[\unit]_\mu}{U^*\pi_{\mu'}(f)U[\unit]_\mu}=\scpro{[\unit]_\mu}{\pi_{\mu}(f)[\unit]_\mu}=\int f\d\mu
  \end{displaymath}
shows that the normal state induced by $\psi$ on $\pi_{\mu'}(L^\infty(X,\XG))$ is $\mu$. Suppose $Y\in\XG$ fulfils  $\mu(Y)=0$. This implies $[\chfc Y]_\mu=0$ and we derive all  $f\in L^\infty(X,\XG)$
\begin{displaymath}
  \pi_{\mu}(\chfc Y)[f]_\mu=\pi_{\mu}(\chfc Y)\pi_{\mu}(f)[\unit]_\mu=\pi_{\mu}(f)\pi_{\mu}(\chfc Y)[\unit]_\mu=\pi_{\mu}(f)[\chfc Y]_\mu=0
\end{displaymath}
such that $\pi_{\mu}(\chfc Y)=0$. We obtain $\pi_{\mu'}(\chfc Y)=U\pi_{\mu}(\chfc Y) U^*=0$ what shows that
  \begin{displaymath}
  \mu'(Y)= \scpro{[\unit]_{\mu'}}{\pi_{\mu'}(\chfc Y)[\unit]_{\mu'}}=0
  \end{displaymath}
and $\mu'\ll\mu$. $\mu\ll\mu'$ follows by symmetry.
\end{proof}
Next we want to derive conditions on $\mu$ ensuring that $L^\infty(\mu)$ factorizes as   $L^\infty(\mu)\cong L^\infty(\mu_1)\otimes L^\infty(\mu_2)(\cong L^\infty(\mu_1\otimes\mu_2))$. Before, we state a more general lemma.
\begin{lemma}
  \label{lem:L^2tensorproduct=measureproduct}
Let  $X,X'$ be  Standard Borel spaces and  $\M,\M'$ be measure types on them. Then   $L^\infty(X,\M)\cong L ^\infty(X',\M')$ as $W^*$-algebras iff there is a measurable   map $\map iXX'$ with
  \begin{displaymath}
    \M'=\M\circ i^{-1}
  \end{displaymath}
such that for all $f\in L^\infty(\M)$ there is some $g\in L^\infty(\M')$ with
\begin{displaymath}
  g(i(x))=f(x)\dmf{\M-\text{a.a.~}x\in X}.
\end{displaymath}
Especially, the assertion holds if $i$ is a bijection.
\end{lemma}
\begin{proof}
 Observe that  $g\mapsto g\circ i$ defines a normal unital homomorphism from $L^\infty(X',\M')$ to $L ^\infty(X,\M)$. Our task is to prove it is injective and onto. The former was already proved. The second conditions is exactly the requirement that this homomorphism has full image. This shows the \emph{if} direction.

The \emph{only if} statement is trivial in the case that $X$ (and thus $X'$) are countable, since we can identify the points in $X$ with the minimal projections in $L^\infty(X,\M)$ (upto sets of measure zero). If $X$ and $X'$ are uncountable,  we find  a  Borel  isomorphism $i_1$ from $\set{0,1}^\NN$ to $X'$. This isomorphism gives us  a sequence of sets $\sequ n{Y'}\subset\XG'$ being the images of the sets $\set{\sequ kx:x_n=1}$. Under an isomorphism $\map j{L^\infty(X,\M)}{ L ^\infty(X,\M)}$, $j^{-1}(\chfc {Y'_n})$ is a  projection in $L^\infty(\M)$. Thus we may choose $Y_n\in\XG$ such that $\chfc{Y_n}=j^{-1}(\chfc {Y'_n})$. Naturally, we define $\map{i_2}{X}{\set{0,1}^\NN}$ by $i_2(x)=(\chfc{Y_n}(x))_{n\in\NN}$.  Choose $i=i_1\circ i_2$, then normality of $j$ and the monotone class theorem   establish the assertion. 
\end{proof}
\begin{corollary}
 Let $\M,\M_1,\M_2$ be measure types on   Standard Borel spaces $X,X_1,X_2$ respectively. Then $L^\infty(X,\M)\cong L ^\infty(X_1\times X_2,\M_1\otimes\M_2)$ as $W^*$-algebras iff there is a measurable   map $\map i{X_1\times X_2}X$ with
  \begin{displaymath}
    \M=(\M_1\otimes\M_2) \circ i^{-1}
  \end{displaymath}
such that for all $f\in L^\infty(\M_1\otimes\M_2)$ there is some $g\in L^\infty(\M)$ with
\begin{displaymath}
  g(i(x_1,x_2))=f(x_1,x_2)\dmf{\M_1-\text{a.a.~}x_1\in X_1,\M_2-\text{a.a.~}x_2\in X_2}.\proofend
\end{displaymath}
\end{corollary}
In the same way, one can prove
\begin{lemma}
\label{lem:measurable automorphism Linfty}
  $\sgi t\beta$ is a $w-*$-measurable family of $\ast$-automorphisms of $L^\infty(\M)$ iff there is a measurable function $\map h{[0,1]\times X}{X}$ with
  \begin{displaymath}
    \beta_t(f)(x)=f(h(t,x))\dmf{t\in[0,1],\M-a.a.~x\in X}.\proofend
  \end{displaymath}
\end{lemma}
\begin{definition}
  \label{def:generalstationaryfactorizingmeasuretype}
  Let  $(X_{s,t})_{(s,t)\in I_{0,\infty}}$ be a family of Standard Borel  spaces with the following properties:
  \begin{condition}{X}
  \item \label{x1} For all $(s,t)\in I_{0,\infty}$, $r\in\NRp$ the space  $X_{s,t}$ is isomorphic to $X_{s+r,t+r}$ under a measurable  map $\mathfrak{s}^{s,t}_r$ with $\mathfrak{s}^{s,t}_{r+r'}=\mathfrak{s}^{s+r,t+r}_{r'}\circ\mathfrak{s}^{s,t}_r$,
\item \label{x2}For $(r,s),(s,t)\in I_{0,\infty}$ there is
 a  measurable   map $\map{\mathfrak{i}_{r,s,t}}{X_{r,s}\times X_{s,t}}{X_{r,t}}$ with 
  \begin{displaymath}
  \mathfrak{i}_{r,s',t}\circ ( \mathfrak{i}_{r,s,s'}\times\Id_{X_{s',t}})= \mathfrak{i}_{r,s,s'}\circ(\Id_{X_{r,s}}\times \mathfrak{i}_{s,s',t}) \dmf{(r,s),(s,s'),(s',t)\in I_{0,\infty}}
  \end{displaymath}
and
\begin{displaymath}
 \mathfrak{i}_{s+r,s'+r,t+r}\circ (\mathfrak{s}^{s,s'}_r\times \mathfrak{s}^{s',t}_r)=\mathfrak{s}^{s,t}_r\circ \mathfrak{i}_{s,s',t}\dmf{(s,s'),(s',t)\in I_{0,\infty},r\in\NRp}.
\end{displaymath}
  \end{condition}
Then a \emph{stationary factorizing measure type on $(X_{s,t})_{(s,t)\in I_{0,\infty}}$} is a family $\MG=(\M_{s,t})_{(s,t)\in I_{0,\infty}}$  of measure types, each $\M_{s,t}$ on $X_{s,t}$, with 
\begin{equation}
\label{eq:measureequivalenceundershift}
  \M_{s+r,t+r}=\M_{s,t}\circ (\mathfrak{s}^{s,t}_r)^{-1}\dmf{(s,t)\in I_{0,\infty},r\in\NRp}
\end{equation}
and 
\begin{equation}
\label{eq:measurefactorizationunderembedding}
  \M_{r,t}= \M_{r,s}\otimes \M_{s,t}\circ \mathfrak{i}_{r,s,t}^{-1}\dmf{(r,s),(s,t)\in I_{0,\infty}}.
\end{equation}
\end{definition}
\begin{remark}
  Relation \ref{eq:measurefactorizationunderembedding} is weaker than  expected from a factorisation property. We need additionally a kind of non-overlap condition like provided in \ref{eq:factorization f under i} below. Nevertheless, the weaker construction is also of interest. E.g., one could look for  Levy-Khintschine type characterisations  of convolution semigroups of measure types on $\NR$ corresponding to  families $\sg t\varphi$ of probability measures on $\NR$ fulfilling
  \begin{displaymath}
    \varphi_{s+t}\sim\varphi_s\ast\varphi_t.
  \end{displaymath}
\end{remark}
Now we can  present the main result of this section.
\begin{theorem}
  \label{th:generalproductsystemfrom stationaryfactorizingmeasuretype}
  Let $(\M_{s,t})_{(s,t)\in I_{0,\infty}}$ be a stationary factorizing measure type on $(X_{s,t})_{(s,t)\in I_{0,\infty}}$. Suppose that  for all $(r,s),(s,t)\in I_{0,\infty}$ and $f_{r,s}\in L^\infty(\M_{r,s})$, $f_{s,t}\in L^\infty(\M_{s,t})$ there is a function $f_{r,t}\in L^\infty(\M_{r,t})$ such that
\begin{equation}
\label{eq:factorization f under i}
  f_{r,t}\circ\mathfrak{i}_{r,s,t}=f_{r,s}\otimes f_{s,t} \dmf{\M_{r,s}\otimes\M_{s,t}-\mbox{a.s.}}.
\end{equation}

 Then $\E^\MG=\sg t\E$ with $\E_t=L^2(\M_{0,t})$, $t\in\NRp$, and  the multiplication 
\begin{displaymath}
  (V_{s,t}(\psi^s)\otimes (\psi^t))_{\mu\otimes(\mu'\circ\mathfrak{s}^{-1}_s)}=(\psi^s_\mu\otimes\psi^t_{\mu'})\circ \mathfrak{i}_{0,s,s+t}\dmf{\M_{0,s+t}-a.s.}
\end{displaymath}
defines an algebraic  product system.

 $\E^\MG$ is a product system iff   additionally to the above assumptions  there is  a measurable function $\map h{[0,1]\times X_{0,1}}{X_{0,1}}$ fulfilling 
\begin{eqnarray*}
    \mathfrak{i}_{0,1-t,1}(x_{1-t},\mathfrak{s}^{0,t}_{1-t}(x_t))=h(t,\mathfrak{i}_{0,1-t,1}(x_t,\mathfrak{s}^{0,1-t}_{t}(x_{1-t})))\emf[.]{t\in[0,1],\M_{0,p}-\mbox{a.a.}~x_p\in X_{0,p},p=t,1-t}
\end{eqnarray*}
\end{theorem}
Following \cite{Tsi00b}, we equip for any measure type $\M$ the set $L^0(\M)$ of all $\M$-equivalence classes of measurable maps with the topology of convergence in measure. It is  metrized, e.g.,  by
\begin{displaymath}
  d_\mu(f,f')=\int \frac{\abs{f-f'}}{1+\abs{f-f'}}\d\mu
\end{displaymath}
for some $\mu\in \M$. Observe that the topology does not depend on the choice of $\mu$. Further, this introduces a concept of measurability which is equivalent for Standard Borel spaces to existence of a jointly measurable version and, for $L²$-valued maps, to measurability as Hilbert space valued map.   
\begin{proof}
 We want to use Proposition \ref{prop:tensorproductW*=ps} and set $\B_{s,t}=L^\infty(X_{s,t},\M_{s,t})$. $\M_{s,t}$ defines a family $\S_{s,t}$ of (faithful) normal states on $\B_{s,t}$. Further, the isomorphisms  $\beta_t^{r,s}$ and $\gamma_{r,s,t}$ are  implemented by $\mathfrak{s}^{r,s}_t$ and  $\mathfrak{i}_{r,s,t}$ respectively. The first assertion follows by  application of Proposition \ref{prop:tensorproductW*=ps},  the above lemmata as well as the fact that
 \begin{displaymath}
   \frac{\d \mu_1\otimes\mu'_1}{\d \mu_2\otimes\mu'_2}(x,y)=   \frac{\d \mu_1}{\d \mu_2}(x)\frac{\d\mu'_1}{\d\mu'_2}(y)
 \end{displaymath}
what implies \ref{eq:factorization Utetaeta'}.

This second result relies additionally on Proposition \ref{prop:measurability W^*}, Theorem \ref{th:intrinsicmeasurability} and Lemma \ref{lem:measurable automorphism Linfty}. The first and the last  establish necessity of the existence of $h$. For sufficiency, observe that the Lemma guarantees that  the one parameter automorphism  group $\gr t{\smash{\hat\beta}}$ defined in Proposition \ref{prop:measurability W^*} is weakly continuous, i.e.\ $t\mapsto \int\d\mu\hat\beta_t(f)$ is continuous for all $f\in L^\infty(\M_{0,1})$, $\mu\in\M_{0,1}$. Since the Radon-Nikodym derivative $(\mu,\mu')\mapsto \frac{\d\mu}{\d\mu'}$ is measurable \cite[Lemma 6.11]{Tsi00b} as well as multiplication in $L^0(\M)$, we derive that $t\mapsto \sqrt{\frac{\d\mu\circ\hat\beta_t}{\d\mu}}\hat\beta_t(f)$ is measurable. On the other hand, we know that from Lemma \ref{lem:computation tau_t} that 
\begin{math}
  u^{\hat\beta_t}f= \hat\beta_t(f)
\end{math}
what results in  
\begin{displaymath}
  \sqrt{\frac{\d\mu\circ\hat\beta_t}{\d\mu}}\hat\beta_t(f)=
\tau_tf
\end{displaymath}
and  shows that there is a consistent measurable structure on $\E^\MG$. 
\end{proof}
The description of a consistent measurable structure needs a little  more work. We introduce  $\sigma$-fields $\F^T_{s,t}\subseteq\mathfrak{X}_{0,T}$, $T\in\NRp$, $(s,t)\in I_{0,T}$ in the following way.
$\F^T_{0,t}$ consists of all sets $Y\subseteq X_{0,T}$ for which there is a Borel set $Y'\subset{X}_{0,t}$ with
\begin{displaymath}
  \chfc Y\circ\mathfrak{i}_{0,t,T}=\chfc{Y'}\otimes \chfc{X_{t,T}}
\end{displaymath}
$\M_{0,t}\otimes\M_{t,T}$-a.s. Similarly, $Y\in\F^T_{t,T}$ fulfils 
\begin{displaymath}
  \chfc Y\circ\mathfrak{i}_{0,t,T}=\chfc{X_{0,t}}\otimes \chfc{Y'}
\end{displaymath}
and we define $\F^T_{s,t}=\F^T_{s,T}\cap\F^T_{0,t}$.

We have canonical embeddings $L^\infty(\M_{0,t})\cong L^\infty(X_{0,T},\F^T_{0,t},\M_{0,t})\subseteq L^\infty(\M_{0,t})$. With these embeddings (which are consistent by associativity \ref{x2} of the family $\mathfrak{i}$), there exist \emph{consistent} families $\sg t\mu$ of probability measures $\mu_t\in\M_{0,t}$ in the sense that
\begin{displaymath}
  \mu_t\restriction L^\infty(\M_{0,s})=\mu_s\dmf{0<s<t}.
\end{displaymath}
\begin{corollary}
  \label{cor:generalproductsystemfrom stationaryfactorizingmeasuretype,ms}
 The measurable structure determined  by all sections $\sg t\psi$  such that for all $T>0$ and there exists a measurable function $\map f{\NRp\times X_{0,T}}\NC$ with 
\begin{displaymath}
  f(t,\mathfrak{i}_{0,t,T}(x_{t},x_{T-t}))=(\psi_t)_{\mu_t}(x_t)\dmf{0< t<T,\M_{0,t}-\text{a.a.}~x_t,\M_{t,T}-\text{a.a.}~x_{T-t}}
\end{displaymath}
for a consistent family $\sg t\mu$ of measures on $\sgp t{X_{0,t}}$ turns  $\E^\MG=\sg t\E$ into a product system.
\end{corollary}
\begin{proof}
Clearly, for all $\psi\in L²(\M_{0,t})$ there is some measurable section such that $\psi_t=\psi$. It remains to prove that the multiplication on $\E^\MG$ is measurable.  From Proposition \ref{prop:psbyE_1} we derive that it is enough to restrict to $\sgii stV$. We obtain from consistency of $\sg t\mu$ for three measurable sections $\psi¹,\psi²,\psi³$ corresponding to measurable functions $f¹,f²,f³$ on $[0,1]\times X_{0,1}$ that
\begin{displaymath}
  \scpro{\psi¹_{s+t}}{V_{s,t}\psi²_s\otimes \psi³_t}=\int \overline{f^1(s+t,x)}\sqrt{\frac{\d\mu_{s+t}\otimes\mu_{1-s-t}}{\d\mu_s\otimes\mu_t\otimes\mu_{1-s-t}}(x)}f²(s,x)f³(t,h(s,x))\mu_1(\d x).
\end{displaymath}
Similar to  the proof  that the algebra valued map  $(s,t)\mapsto\A_{s,t}$ is  continuous (see Proposition  \ref{prop:singlepoint}), we derive  that the family $(\F¹_{s,t})_{(s,t\in I_{0,1})}$ is measurable as defined in \cite{Tsi00b}. Further, this result and  measurability of $\grp t{\smash{\hat\beta_t}}$ imply for consistent families $\sg t\mu$ that the maps $(s,t)\mapsto\mu_{s+t}\otimes\mu_{1-s-t}$ and $(s,t)\mapsto\mu_s\otimes \mu_t\otimes\mu_{1-s-t}$ are measurable.    Measurability of Radon-Nikodym derivatives \cite[Lemma 6.11]{Tsi00b} and  measurability of integrals  complete  the proof.
\end{proof}
\begin{remark}
\label{rem:restricted ps of Standard Borel spaces}
 By Proposition \ref{prop:psbyE_1} and an analogous construction  for a continuous product system of Standard Borel measure spaces, it would be enough to have a family $(X_{s,t})_{(s,t)\in I_{0,1}}$ of Standard Borel spaces with the corresponding  structures restricted to $[0,1]$.
\end{remark}
Standard Borel spaces are either finite, isomorphic to $\NN$ or Borel isomorphic to $\set{0,1}^\NN\cong[0,1]$. We want to show now that the second is impossible within the structure of stationary factorizing measure types on Standard Borel spaces.
\begin{proposition}
  Suppose that $(\M_{s,t})_{(s,t)\in I_{0,\infty}}$ is a stationary factorizing measure type on Standard Borel spaces $(X_{s,t})_{(s,t)\in I_{0,\infty}}$ which are all countable and fit into the above theorem. Then $\M_{s,t}=\set{\delta_{x_{s,t}}}$ for some family $(x_{s,t})_{(s,t)\in I_{0,\infty}}$.
\end{proposition}
\begin{proof}
  Since $L^\infty(\M_{s,t})$ is generated by its minimal projections, there are no nontrivial one parameter automorphism groups  of it. Consequently, the automorphism group  $\gr t{\smash{\hat\beta}}$  defined in Proposition  \ref{prop:measurability W^*} is trivial. Thus, ${\hat\beta}_{1/2}$ maps $L^\infty(\M_{0,1},\F^1_{0,1/2})$ both into itself and $L^\infty(\M_{0,1},\F^1_{1/2,1})$. I.e., both algebras coincide. On the other side, there is a $\mu\in\M_{0,1}$ such that $\F^1_{0,1/2}$ and $\F^1_{1/2,1}$ are independent. Consequently, these $\sigma$-fields are $\M_{0,1}$-a.s.\ trivial and the same is true of $\F^1_{0,1}$. Therefore, $\M_{0,1}$ and thus each of the  $\M_{s,t}$ consists of a single Dirac measure. 
\end{proof}
\begin{corollary}
  If one of the measure types $(\M_{s,t})_{(s,t)\in I_{0,\infty}}$ consists of discrete measures, the resulting product system is isomorphic to $\sgp t\NC$.\proofend
\end{corollary}
\begin{corollary}
  Let $(\M_{s,t})_{(s,t)\in I_{0,\infty}}$ be  a stationary factorizing measure type on $(X_{s,t})_{(s,t)\in I_{0,\infty}}$ which  fits into the above theorem. Then all the measure types $\M_{s,t}$ have at most one atom.
\end{corollary}
\begin{proof}
  We could restrict $\M_{s,t}$ to the space of its atoms $X^0_{s,t}$. Provided $X^0_{s,t}\ne\emptyset$,  the proposition applies to the pair $X^0,\MG^0$, $\M_{s,t}^0=\set{\mu:\mu\sim \M_{s,t}\restriction X^0_{s,t}}$. Since $\M_{s,t}$ is positive on its atoms, this shows that $X^0_{s,t}$ is a singleton.  
\end{proof}
\begin{remark}
  If one (and thus all) $\M_{s,t}$ has an atom, say $x_{s,t}$, it is easy to see that 
\begin{displaymath}
    (u_t)_\mu=\mu(\set{x_{0,t}})^{-1/2}\chfc{\set{x_{0,t}}}
  \end{displaymath}
defines a unit of $\E^\MG$. In this sense, product systems from random sets are typical (as long as $L$ is compact). 

On the other side, $\M_{s,t}$ need not have an atom for units in $\E^\MG$ to exists, see  Example \ref{ex:ps from Brownian increments}. 
\end{remark}
\begin{remark}
  It would simplify matters considerably if all product systems could be derived as $\E^\MG$. It is not clear whether  this is possible but we conjecture it is not for the following reason. 

Namely, $\E^\MG_t$ carries a natural conjugation (i.e.\ antiunitary) mapping $f\in L^2(\mu)$ into its complex conjugate $\overline{f}$ which extends to $L^2(\M_{0,t})$ since $U_{\mu,\mu'}$ is multiplication by a real (in fact positive) function. Consequently, if $\E\cong\E^\MG$ then the conjugate product system $\E^*$ (see Example \ref{ex:ps from faithful states}) is isomorphic to $\E$. Therefore, any product system with   $\E\not\cong \E^*$ (which are likely to exist in view of the now many counterexamples in product systems) would provide an example for a product system which is not of  of the form $\E^\MG$. It could even not derived from the set of faithful states on a product system of $W^*$-algebras (see end of Example \ref{ex:ps from faithful states}) since these carry the intrinsic modular conjugations.   

Another equivalent, but seemingly useless, formulation is the following: There exists a product system of maximal abelian $W^*$-algebras inside $\sgp t{\B(\E_t)}$. 
\end{remark}
\subsection{Product Systems from Random Sets}
\label{sec:psfromRACS}

This section is dedicated to  the proof of Proposition \ref{prop:statfactmeastype=ps}. Thereby, we  have to  analyse  the events  $\set{Z:(\set t\times L)\cap Z\ne\emptyset}$ for $t\in[0,1]$. First, we  consider  singletons $L$ here.
\begin{lemma}
  \label{lem:threedifferentpossibilities}
  For a quasistationary quasifactorizing measure $\mu$ on $\FG_{[0,1]}$ there are only two possibilities
\begin{enumerate}
\item \label{all0}$\mu(\set{Z:t\in Z})=0$ for all $t\in[0,1]$, or
\item \label{all1}$\mu(\set{Z:t\in Z})=1$ for all $t\in[0,1]$.
\end{enumerate}
The latter happens for $\mu=\delta_{[0,1]}$ only.
\end{lemma}
\begin{proof}
From quasistationarity it is clear that  it suffices to prove that 
\begin{enumerate}
\setcounter{enumi}2\item \label{allin(0,1)}$0<\mu(\set{Z:t\in Z})<1$ for all $t\in[0,1]$
\end{enumerate}
cannot happen, so let us  assume that \ref{allin(0,1)} is true.  This means that $\sgi tP$, $P_t=\chfc{\set{Z:t\in Z}}(\cdot)$ is a set of nontrivial projections in $L^2(\mu)$. 

Separability of  $L^2(\mu)$ shows that there is a countable set  $Q\subset[0,1]$ such that $(P_q)_{q\in Q}$ is weakly dense in $\sgi tP$. Assume we have a sequence $\sequ nq\subset Q$ with $P_{q_n}\limitsto{w}{n\to\infty}P_t$ for some $t$, this is equivalent to $\chfc{\set{Z:q_n\in Z}}\limitsto{\mu}{n\to\infty}\chfc{\set{Z:t\in Z}}$. In order to show $\lim_{n\to\infty}q_n=t$ assume $\lim_{n\to\infty}q_n=t'\ne t$. We introduce $\sigma$-fields $\F_t^-$ and $\F_t^+$ as the $\mu$-completion of the $\sigma$-fields generated by the random variables $Z\mapsto Z\cap[0,t)$ and $Z\mapsto Z\cap(t,1]$ respectively. From  $(Z_1\cup Z_2)\cap[0,t)=Z_1\cap[0,t)$, $(Z_1\cup Z_2)\cap(t,1]=Z_2\cap(t,1]$, and   equation \ref{eq:quasifactorizingmeasure}  we derive that $\F_t^-$ and $\F_t^+$ are independent under $\mu_{0,t}\ast\mu_{t,1}$. Thus $t'<t$ implies $\chfc{\set{Z:t\in Z}}=\mu-\lim_{n\to\infty}\chfc{\set{Z:q_n\in Z}}$ is $\F_{(t+t')/2}^-$-measurable. On the other side, it is $\F_{(t+t')/2}^+$-measurable and therefore $\mu_{0,(t+t')/2}\ast\mu_{(t+t')/2,1}$-independent of itself. This happens only if it is almost surely constant contradicting nontriviality of $P_t$. The same arguments show that for each $t\in(0,1)$ there cannot be two sequences $\sequ nq,\sequ n{q'}\subset Q$, $q_n\uparrow_{n\to\infty} t$, $q'_n\downarrow_{n\to\infty} t$ with $\mu-\lim_{n\to\infty}\chfc{\set{Z:q_n\in Z}}=\chfc{\set{Z:t\in Z}}=\mu-\lim_{n\to\infty}\chfc{\set{Z:q'_n\in Z}}$. Due to density of  $(P_q)_{q\in Q}$ in $\sg tP$ this amounts to say for each $t\in[0,1]$ either $\inf\set{d_w(P_s,P_t):s<t}>0$ or $\inf\set{d_w(P_s,P_t):s>t}>0$ if $d_w$ is some metric for the weak topology on the unit ball  $B_1(\B(L^2(\mu)))$. Quasistationarity of $\mu$ implies that either of these relations is true  for all $t\in[0,1]$ simultaneously. Without loss of generality, assume $\inf\set{d_w(P_s,P_t):s<t}=0$ for all $t$, i.e.\ $\chfc{\set{Z:t\in Z}}$ is $\F_t^-$-measurable.

So, fix $t\in(0,1)$ and $f\in L^\infty(\mu_{0,t})$ and denote $1_t$ a version of $\chfc{\set{Z:t\in Z}}$ which is measurable with respect to the $\sigma$-field generated by $Z\cap[0,t)$. Consider $Y=\set{Z:Z\cap(t,t+1/n)\ne\emptyset\forall n\in\NN}\subset\set{t\in Z}$. Since $Y\in\F_t^+$ is independent of $\set{Z:t\in Z}\in\F_t^-$, we conclude that $Y$ is a null set. Therefore, for all $t\in[0,1)$ the random set  $Z\cap(t,1]$ is closed almost surely. 

 The Baire theorem on categories tells us that there is at least one $\varepsilon>0$ such that the closure of $\set{t:\mu(\set{Z:t\in Z})>\varepsilon}$ contains an interval, around $t_0$ say. Thus we find  a sequence $\sequ nt$, $t_n\downarrow t_0$, $t_n>t_0$ with $\mu(\set{Z:t_n\in Z})\ge\varepsilon$. Consequently,
\begin{eqnarray*}
 \mu(\set{Z:t_n\in Z \text{~for infinitely many ~}n\in\NN})&=&\mu(\bigcap_{n\in\NN}\bigcup_{k\ge n}\set{Z:t_k\in Z})\\
&=&\lim_{n\to\infty}\mu(\bigcup_{k\ge n}\set{Z:t_k\in Z})\\
&\ge& \limsup_{n\to\infty}\mu(\set{Z:t_n\in Z})\ge\varepsilon.
\end{eqnarray*}
But the set $\set{Z:t_n\in Z \text{~for infinitely many ~}n\in\NN}$ is contained in $Y$, which is a null set. This  contradiction shows that  \ref{allin(0,1)} is not possible. 

If \ref{all1} is valid, $\mu$-a.s.\ $q\in Z$ for all $q\in\mathbb{Q}\cap[0,1]$. Since $Z$ is closed, $Z=[0,1]$ almost surely and $\mu=\delta_{[0,1]}$ results.
\end{proof}
\begin{remark}
  In fact, such a non-atomicity  result is valid at a much more abstract level. All what is needed is that we can almost surely recover $x_{t'},x_{t''}$ from $\mathfrak{i}_{0,t'',1}(\mathfrak{i}_{0,t',t''}(x_{t'},y),x_{t''})$ for arbitrary $t'<t''$.  
\end{remark}
\begin{lemma}
  Let $\mu$ be a quasistationary quasifactorizing measure on $\FG_{[0,1]\times L}$. Then there is a maximal closed set $L'\subseteq L$ such that any measurable $Y\subseteq L$ with $\mu(\set{Z:(\set t\times Y)\cap Z\ne\emptyset})>0$ for some $t\in[0,1]$ fulfils $Y\cap L'\ne\emptyset$. For this $L'$ almost surely $[0,1]\times L'\subseteq Z$. 
\end{lemma}
\begin{proof}
  We define 
  \begin{displaymath}
    L'=\set{l\in L:\mbox{$\mu(\set{Z:(t,l)\in Z})>0$ for some $t\in[0,1]$}}.
  \end{displaymath}
Like in the above lemma, it is possible  to restrict ourselves to one $t$, say $t=1/2$ and $\mu(\set{Z:(t,l)\in Z})=1$. Then it is easy to see that $L'$ is closed. 

It remains to prove that  $\mu(\set{Z:(\set t\times L\setminus L')\cap Z\ne\emptyset})=0$ for all $t\in[0,1]$. Since $L\setminus L'$ is again locally compact, let us assume that $L'$ is empty and $\mu(\set{Z:(\set t\times L)\cap Z\ne\emptyset})>0$ as well as $\mu(\set{Z:(t,l)\in Z})=0$ for all $t\in[0,1]$, $l\in L$. By $\sigma$-compactness of $L$, we derive that there is a compact $K\subset L$ with $\mu(\set{Z:(\set t\times K)\cap Z\ne\emptyset})>0$ for all $t\in[0,1]$. Since $\pi(Z\cap([0,1]\times K))$ is again a quasistationary quasifactorizing random closed set for the projection $\map\pi{[0,1]\times L}{[0,1]}$, we derive from the above lemma that $\mu(\set{Z:(\set t\times K)\cap Z\ne\emptyset})=1$ for all $t\in[0,1]$. Dissecting $K$ we find a decreasing sequence $\sequ nK$, $K_n\subset K$, with $\mathrm{diam} K_n\le 2^{-n}$ and $\mu(\set{Z:(\set t\times K_n)\cap Z\ne\emptyset})=1$ for all $t\in[0,1]$, $n\in\NN$. Since $\bigcap_{n\in\NN}K_n=\set l$ for some $l\in L$ by compactness of $K$ it results $\mu(\set{Z:(\set t\times\set l)\cap Z\ne\emptyset})=1$ which contradicts  $\mu(\set{Z:(t,l)\in Z})=0\forall l\in L$. This completes the proof.
\end{proof}
\begin{proof}[ of Proposition \ref{prop:statfactmeastype=ps}]
\label{page:proof prop:statfactmeastype=ps}
We need  to   establish the conditions for application of Theorem \ref{th:generalproductsystemfrom stationaryfactorizingmeasuretype}. 

So we can assume the analogue of \ref{all0}. Define the maps $\map {\mathfrak{i}_{r,s,t}}{\FG_{[r,s]\times L}\times\FG_{[s,t]\times L}}{\FG_{[r,t]\times L}}$  by $\mathfrak{i}_{r,s,t}(Z_1,Z_2)=Z_1\cup Z_2$ and $\mathfrak{s}_r(Z)=Z+r$. We need to show that for all $f\in L^\infty(\M_{r,s})$ there is some $g\in L^\infty(\M_{r,t})$ such that 
\begin{displaymath}
  g(Z_{r,s}\cup Z_{s,t})=f(Z_{r,s})\dmf{\M_{r,s}-\text{a.a.~}Z_{r,s},\M_{s,t}-\text{a.a.~}Z_{s,t}}
\end{displaymath}
and a similar condition for all $f\in L^\infty(\M_{s,t})$ since this implies the general condition. The proof for the latter is the same, so take $f\in L^\infty(\M_{r,s})$ and set $g(Z)=f(Z\cap([r,s]\times L))$. We know for $\M_{r,s}\otimes\M_{s,t}$--a.a.\ $(Z_1,Z_2)$ that $(\set s\times L)\cap Z_1=\set s\times L'=(\set s\times L)\cap Z_2$ where $L'$ is taken from the above lemma. Therefore, $(Z_1\cup Z_2)\cap([r,s]\times L)=Z_1$  $\M_{r,s}\otimes\M_{s,t}$--a.s.\ and this part is proven. 

For proving measurability, set $\map h{[0,1]\times\FG_{[0,1]\times L}}\FG_{[0,1]\times L}$
\begin{displaymath}
  h(t,Z)=Z+t\dmf{t\in[0,1],Z\in\FG_{[0,1]}}
\end{displaymath}
where the RHS was defined at the beginning of section \ref{sec:theoryrs2ps}, on page  \pageref{sec:theoryrs2ps}. Clearly,   $h$ is measurable and   Theorem  \ref{th:generalproductsystemfrom stationaryfactorizingmeasuretype} completes the proof.
\end{proof}
\subsection{Product Systems from Random Measures }
\label{sec:random measure}

 Consider the space $\MG_\NRp$, the set of Radon measures on $\NRp$. Under  the vague  topology, $\MG_\NRp$ is a Polish space \cite{Kal86}. Similar results are valid  for the spaces $\MG_{[s,t]}$ of finite measures on $[s,t]$. Further, a natural choice of $\map{\mathfrak{i}_{r,s,t}}{\MG_{[r,s]}\times\MG_{[s,t]}}{\MG_{[r,t]}}$ is
 \begin{displaymath}
   \mathfrak{i}_{r,s,t}(\varphi,\varphi')=\varphi+\varphi'.
 \end{displaymath}
Further, the shift  is implemented as $\mathfrak{s}_t\varphi(B)=\varphi(B-t)$. Having fixed this, we call a family $\M_{s,t}$ of measure types on $\MG_{s,t}$ \emph{stationary factorizing measure type on $\MG_\NRp$} if it fulfils \ref{eq:measureequivalenceundershift} and \ref{eq:measurefactorizationunderembedding}.  Similar to Proposition \ref{prop:statfactmeastype=ps} one can prove
\begin{proposition}
\label{prop:statfactmeastypemeasure=ps}
 Every stationary factorizing measure type $\M$ on $\MG_\NRp$ gives rise to a product system $(\E^\M_t)_{t\in[0,1]}$.\proofend
\end{proposition}
\begin{remark}
  Observe that the map $\varphi\mapsto\supp\varphi\in\FG_\NRp$ is measurable. Thus  every stationary factorizing measure type $\M$ on $\MG_\NRp$ yields a stationary factorizing measure type $\M$ on $\FG_\NRp$. Despite the cases discussed in section \ref{sec:rs2ps} there is only one additional one: It may happen that $\supp\varphi_{s,t}=[s,t]$ almost surely. This may lead to nontrivial product systems. E.g.\ if $\varphi$ is an (infinitely divisible) Gamma Random Measure \cite[Exercise 6.1.2]{DVJ88} the resulting product system is of type $\mathrm{I}_\infty$. 
\end{remark}
\begin{example}
Conversely, some quasistationary  quasifactorizing measures  on $\FG_\NRp$ can be canonically related to distributions of random measures. 
 
 E.g., a set $\set{t_1,\dots,t_n}\in\FG^f_{[0,t]}$ corresponds in a one-to-one fashion to the measure $\varphi=\delta_{t_1}+\cdots+ \delta_{t_n}$. Thus we could define  the Poisson process  $\Pi_\ell$ as well as distribution of  a random measure.

  In the same spirit, the zero set of a Brownian motion (see Example \ref{ex:zerosBt1}) corresponds almost surely one-to-one to the restriction of a certain  Hausdorff measure to this set \cite{C:Per81}.  
\end{example}
\subsection{Product Systems from Random Increment Processes}
\label{sec:random increments}
We start with an example.
\begin{example}
\label{ex:ps from Brownian increments}
  For  a Brownian motion $\sg tB$ we define $\E_t=L^2(\mathscr{L}((B_s)_{0\le s\le t}))$. Independence of the increments of Brownian motion and the Markov property show $\E_{s+t}\cong{\E_s\otimes\E_t}$ by 
  \begin{displaymath}
    (B_s)_{0\le s\le t+t'})\leftrightarrow((B_s)_{0\le s\le t},(B_{s+t}-B_t)_{0\le s\le t'})).
  \end{displaymath}
The well-known Wiener-It\^o-Segal chaos decomposition states that $\E_t$ is isomorphic to $\Gamma(L^2([0,t],\ell))$. The isomorphism is given by multiple It\^o-Integrals: 
  \begin{displaymath}
    \sequ nf\mapsto\sum_{n\in\NN}\int^\infty_0\int^{t_1}_0\cdots\int  ^{t_{n-1}}_0f(t_1,\dots,t_n)\d B_{t_n}\d B_{t_2}\d B_{t_1}.
  \end{displaymath}
(the series on the right consist of iterated ordinary It\^o-Integrals). Clearly, this isomorphism respects the independence of increments, what  shows that $\E\cong\Gammai(\NC)$. On the other side, we could generate $\E_t$ in an  equivalent manner from the distribution of the Brownian  increments $(B_s-B_r)_{0\le r<s\le t}$. 
\end{example}
This examples leads us to consider  increment processes. Define $\IG_\NRp$ to be the set of all maps $\map \xi{I_{0,\infty}}\NR$ with
\begin{equation}
\label{eq:increments}
  \xi(r,t)=\xi(r,s)+\xi(s,t)\dmf{(r,s),(s,t)\in I_{0,\infty}}  
\end{equation}
and the property that the map $(s,t)\mapsto \xi(s,t)$ is (jointly) c\`adl\`ag, i.e.\ it is  right continuous and has limits from the left everywhere. Similarly, we have also $\IG_{s,t}$, $(s,t\in I_{0,\infty})$, equipped with 
\begin{displaymath}
  \mathfrak{i}_{r,s,t}(\xi,\xi')(s',t')=\left\{
    \begin{array}[c]{c>{~$}l<$}
\xi(s',t')&if $t'\le s$\\
\xi(s',s)+\xi'(s,t')&if $s'< s < t'$\\
\xi'(s',t')&if $s'\ge s$
 \end{array}
\right.
\end{displaymath}
and the obvious shift. Again, there is the canonical notion of a  \emph{stationary factorizing measure type  on $\IG_\NRp$}. 

The set of c\`adl\`ag functions is Polish \cite[Theorem 12.2]{Bil68}, thus one can prove similar to Proposition \ref{prop:statfactmeastype=ps}
\begin{proposition}
\label{prop:statfactmeastypeincrements=ps}
 For every stationary factorizing measure type $\M$ on $\IG_\NRp$ 
the bundle $\E^\M=\sg t{\E^\M}$, $\E^\M_t=L^2(\M_{0,t})$ is a product system.\proofend
\end{proposition}
\begin{remark}
  Clearly, every increment process $\xi$ corresponds to an ordinary process $\sg t{\xi'}$ with c\`adl\`ag paths, $\xi'_t=\xi(0,t)$ via $\xi_{s,t}=\xi_t'-\xi_s'$. Conversely, any process $\sg t{\xi'}$ defines an increment process $\xi$ in this manner, although  we loose the information about  $\xi'_0$ this way. Identifying a Radon measure $\varphi$ with its distribution function $\xi'_t=\varphi([0,t))$, we see that the set of product systems based on  increment processes includes  in fact the set of product systems induced by random measures.     
\end{remark}
\begin{remark}
  An interesting generalisation are  multiplicative increments, where in  \ref{eq:increments} the $+$ is replaced by $\cdot$:
\begin{displaymath}
  \xi(r,t)=\xi(r,s)\xi(s,t)\dmf{(r,s),(s,t)\in I_{0,\infty}}  
\end{displaymath}
where $\xi$ takes values in some Polish semigroup. We just mention some examples.

If $\xi$ is $[0,\infty)$-valued, this reveals equation \ref{eq:arst}. But, $\xi(s,t)$ need  not to be a bounded random variable. If $\xi$ is positive, we can take logarithms and come  back to \ref{eq:increments}. On the other side, $\omega(s,t)=\chfc{\set{(0,\infty)}}(\xi(s,t))$ fulfils \ref{eq:prodomega}. Assuming sufficient continuity of $\xi$  this relates such  multiplicative increment processes to  a pair of a random set  and an additive increment process (which is possibly $\infty$).   
 
 Another example would be $\NT$-valued increment processes which  correspond to adapted shift-equivariant unitaries (i.e.\ the analogue of \ref{eq:Pstshift} is fulfilled). Such unitaries are in one-to-one correspondence with automorphisms in the same manner like adapted  shift-equivariant projections correspond to product subsystems, so this structure could be important to study automorphisms of product systems.

If the range of $\xi$ is not a group, we cannot recover an ordinary process $\sg t {\xi'}$ from its increments. The simplest example is $\set{0,1}$ with multiplication. Then increment processes correspond to bisets introduced in Section \ref{sec:nonsep}, see Proposition \ref{prop:charbisets} below, under continuity assumptions to  closed sets. 

Summarisingly, increment processes deserve further explorations, may be with modification of the continuity properties as the next example and the next section suggest. 
\end{remark}
  Which types of product systems  could arise in these constructions? Obviously, both constructions allow type $\mathrm{I}$. But we recover also the type $\mathrm{II}$ examples from section \ref{sec:perfect} connected with Hausdorff measures, identifying $Z$ with $\mathscr{H}^h(Z\cap\cdot)$, where $\map h\NRp\NRp$ is the scaling function or the increment process $\xi(s,t)=\mathscr{H}^h(Z\cap(s,t])$.  As  \cite{Tsi03,Tsi00b} show, we can also  derive    type $\mathrm{III}$ product systems from  increment processes.  
\begin{example}
\label{ex:Tsirelson type III}
 \textsc{Tsirelson}  considered in \cite{Tsi00b} generalised Gaussian processes $W$ with covariance
 \begin{displaymath}
   \NE W(f)W(g)=\int \ell(\d t)\ell(\d s)f(s)\kappa(t-s)g(t)=\scpro fg_\kappa
 \end{displaymath}
where  $\kappa$ is a positive definite function with $\int\ell(\d t)|\kappa(t)|<\infty$ and $\kappa(t)=\frac1{t\abs{\ln^{\alpha} t}}$ for small positive $t$ and some $\alpha>1$. 

We can create a Gaussian increment process $(\xi_{s,t})_{(s,t\in I_{0,\infty})}$ by
\begin{displaymath}
  \xi(s,t)=W(\chfc{[s,t]})\dmf{(s,t)\in I_{0,\infty}}
\end{displaymath}
and totality  of  indicator functions in the reproducing kernel Hilbert space $H_\kappa$ of $\kappa$ (i.e.\ $H_\kappa$ is the completion of $L^2(\NRp)$ in $\norm\cdot_\kappa$) shows that $\xi$ contains the same information as $W$. 

Clearly, $\xi$ is stationary and we see that the respective laws of the restrictions $\xi_{s,t}\in\IG_{s,t}$ of $\xi$ to $I_{s,t}$ form a stationary factorizing measure type. A problem is whether $\xi$ has a  c\`adl\`ag modification since  standard conditions like \cite[Theorem 15.7]{Bil68} are  not directly applicable.  At least, $(s,s',t,t')\mapsto \NE \xi(s,t)\xi(s',t')$ is continuous since $(s,t)\mapsto\chfc{[s,t]}\in H_\kappa$ is. Thus $\xi$ is continuous in probability and  $\xi$ is determined by countably many values $\xi(s,t)$. This shows that there is a Standard Borel structure supporting this probability measure and by Lemma \ref{lem:L^2separable=Polish} $\L²(\mathscr{L}(\xi))$ is separable. Since the Standard Borel structure enters the proof of Theorem \ref{th:generalproductsystemfrom stationaryfactorizingmeasuretype} only through that lemma, $L²(\M)$ is a product system.  

Observe that in the present example it is enough to consider instead of the  whole measure type only the Gaussian measures from this measure type. The latter set is easily shown to fulfil \ref{h0}--\ref{h1} in Definition \ref{def:stationary factorizing type of states}.
\end{example}

\section{Beyond Separability:  Random Bisets}
\label{sec:nonsep}
 In this section we want to discuss briefly, what problems would appear if we dropped the separability condition on the fibres of a product system. If  $\E_1$ is not any more  separable, the notions of weak and strong measurability are not any more identical and  Proposition \ref{prop:singlepoint} and Proposition  \ref{prop:tau_tiscontinuous} may break down. Further, this may happen for algebraic product systems. Therefore, we want  to understand what we really have to add to the relations 
$$
     \mathrm{P}_{r,t} =\mathrm{P}_{r,s}  \mathrm{P}_{s,t}\dmf{(r,s),(s,t)\in I_{0,1}}\eqno\let\upshape\relax\ref{eq:Prst}
$$
 to get a version of Theorem \ref{th:RACS}, i.e.\ random \emph{closed} sets,  in more general situations too. Thus  we think it is instructive to analyse these  relations  on their own. 
In Theorem \ref{th:RACS} we analysed projection families which were additionally adapted what implied continuity in Proposition \ref{prop:singlepoint}. We are interested, to what extent the result stays valid without the product structure.  From a projection family  with \ref{eq:Prst} we build new projections as 
\begin{eqnarray}
\label{eq:defPplus}
  \mathrm{P}_s^+&=&\bigvee_{\varepsilon>0}\mathrm{P}_{s,s+\varepsilon}\dmf{0\le s<1}\\
\label{eq:defPminus}
  \mathrm{P}_s^-&=&\bigvee_{\varepsilon>0}\mathrm{P}_{s-\varepsilon,s}\dmf{0< s\le1}.
\end{eqnarray}
\begin{proposition}
\label{prop:RACS on Hilbert space}
    Let $\H$ be a Hilbert space and $(\mathrm{P}_{s,t})_{(s,t)\in I_{0,1}}$  be  projections in $\B(\H)$ which fulfil \ref{eq:Prst} and $\mathrm{P}_{0,1}\ne0$. Then for all normal states $\eta$ on $\B(\H)$ there exists a unique probability measure $\mu_\eta$ on $\FG_{[0,1]}$ with$$
    \mu_\eta(\set{Z:Z\cap[s_i,t_i]=\emptyset, i=1,\dots,k})=\eta(\mathrm{P}_{s_1,t_1}\cdots \mathrm{P}_{s_k,t_k})\dmf{(s_i,t_i)\in I_{0,1}}
\eqno\let\upshape\relax\ref{eq:mueta}
$$
iff $\mathrm{P}_s^+=\mathrm{P}_s^-$ for all $s\in(0,1)$.
 
Then, if $\eta$ is faithful,  the correspondence 
\begin{displaymath}
  \chfc{\set{Z:Z\cap[s,t]=\emptyset}}\mapsto
\mathrm{P}_{s,t}\dmf{(s,t)\in I_{0,1}}
\end{displaymath}
extends to an injective  normal representation $J_{\mathrm{P}}$  of $L^\infty(\mu_\eta)$ on $\B(\H)$ with  image $\set{\mathrm{P}_{s,t}:{(s,t)\in I_{0,1}}}''$. Moreover,  $\mu_{\eta'}\ll\mu_\eta$ for all  normal states $\eta'$ on $\B(\H)$.
\end{proposition}
\begin{proof}
  The proof of the \emph{if} direction parallels essentially that of Theorem \ref{th:RACS}. The only thing we need as substitute of Proposition \ref{prop:singlepoint} is the relation $\slim_{n\to\infty}\mathrm{P}^\circ_{(s-1/n,t+1/n)}=\mathrm{P}_{s,t}$. Since we have for $\varepsilon>0$ that
  \begin{displaymath}
\mathrm{P}_{s-\varepsilon,t+\varepsilon}\le\mathrm{P}^\circ_{(s-\varepsilon,t+\varepsilon)}\le\mathrm{P}_{s-\varepsilon/2,t+\varepsilon/2}\le\mathrm{P}_{s,t}
\end{displaymath}
we need only to show $\slim_{n\to\infty}\mathrm{P}_{s-1/n,t+1/n}=\mathrm{P}_{s,t}$. From definition, \ref{eq:Prst}, assumption,  and $\sigma$-strong continuity of multiplication we derive 
\begin{eqnarray*}
\slim_{n\to\infty}\mathrm{P}_{s-1/n,t+1/n}&=&\slim_{n\to\infty}\mathrm{P}_{s-1/n,s}\mathrm{P}_{s,t}\mathrm{P}_{t,t+1/n}\\
&=&\mathrm{P}_s^-\mathrm{P}_{s,t}\mathrm{P}_t^+\\
&=&\mathrm{P}_s^+\mathrm{P}_{s,t}\mathrm{P}_t^-\\
&=&\slim_{n\to\infty}\mathrm{P}_{s,s+1/n}\mathrm{P}_{s,t}\mathrm{P}_{t-1/n,t}\\
&=&\slim_{n\to\infty}\mathrm{P}_{s,s+1/n}\mathrm{P}_{s,s+1/n}\mathrm{P}_{s+1/n,t-1/n}\mathrm{P}_{t-1/n,t}\mathrm{P}_{t-1/n,t}\\
&=&\slim_{n\to\infty}\mathrm{P}_{s,s+1/n}\mathrm{P}_{s+1/n,t-1/n}\mathrm{P}_{t-1/n,t}\\
&=&\mathrm{P}_{s,t}.
\end{eqnarray*}

For the proof of the \emph{only if} statement observe that validity of \ref{eq:mueta} for one faithful $\eta$ implies existence of $J_{\mathrm{P}}$. Then for any $s\in(0,1)$ we derive from  closedness of $Z\in\FG_{[0,1]}$ that $  \set{Z:s\notin Z}=\set{Z:Z\cap[s,s+1/n]=\emptyset}$ and 
\begin{equation}
  \label{eq:sets determined from right}
  \chfc{\set{Z:s\notin Z}}=\lim_{\varepsilon\downarrow0}\chfc{\set{Z:Z\cap[s,s+\varepsilon]=\emptyset}}.
\end{equation}
Normality of $J_{\mathrm{P}}$ shows $J_{\mathrm{P}}(\chfc{\set{Z:s\notin Z}})=\mathrm{P}_s^+$. Similar arguments work for $\mathrm{P}_s^-$ such that
\begin{displaymath}
   \mathrm{P}_s^- =J_{\mathrm{P}}(\chfc{\set{Z:s\notin Z}})=\mathrm{P}_s^+
\end{displaymath}
proves the \emph{only if} statement.
\end{proof}
Since  the interpretation by (random) closed sets forces the condition  $\mathrm{P}_s^+=\mathrm{P}_s^-$ we are left with the problem of finding a similar interpretation of the relations \ref{eq:Prst} only. 
\begin{definition}
  \label{def:fulldc*}
  Let $\C_{0,1}$ be the $C^*$-algebra generated by $\unit$ and commuting projections $(\mathrm{P}_{s,t})_{(s,t)\in I_{0,1}}$
fulfilling the relations   \ref{eq:Prst}.
\end{definition}
\begin{remark}
  We refrain from showing that $\C_{0,1}$ is well-defined. This   is operator algebraic folklore since projections are always contractive.
\end{remark}
The essential problem with closed \emph{sets} can be seen from  \ref{eq:sets determined from right}: The fact $s\notin Z$ poses a restriction on the behaviour of $Z$ on the right of $s$ as well as  on its left. On the other hand, the relations \ref{eq:Prst} separate easily into those in the interval $[0,s]$ and those in the interval $[s,1]$. Consequently, we should allow that $s$ belongs to $Z$ ``from the right'' but not ``from the left''. This induces two sets  $L_Z,R_Z$ collecting the points \emph{not} in $Z$ ``from the left (right)''.  Such constructs we call \emph{bisets}. Of course, any   set $Z\subseteq[0,1]$ induces  a pair $(L_Z,R_Z)$ of subsets of $[0,1]$ such that 
\begin{eqnarray*}
  s\in L_Z&\iff&Z\cap[s,t]=\emptyset \text{~for some $t>s$~}\\
  t\in R_Z&\iff&Z\cap[s,t]=\emptyset \text{~for some $s<t$~}.
\end{eqnarray*}
What is important here is the fact that we dealt with \emph{closed} sets and we should deal with something like \emph{closed} bisets. Therefore, $L_Z,R_Z$ should be open in some sense. Since we are concerned with the analysis of \ref{eq:Prst} only, we drop the word closed in the sequel.
\begin{definition}
\label{def:biset}
  A \emph{biset} is a pair $(L,R)$ of subsets of $[0,1]$ with the following properties:
\begin{condition}{B}
\item\label{b1} $s\in L$ implies $[s,s+\varepsilon)\subset L$ for some $\varepsilon>0$.
\item\label{b2} $t\in R$ implies $(t-\varepsilon,t]\subset R$ for some $\varepsilon>0$.
\item\label{b3} $[s,t)\subset L$ implies $(s,t)\subset R$.
\item\label{b4} $(s,t]\subset R$ implies $(s,t)\subset L$.
\end{condition} 
The set of all bisets in $[0,1]$ is denoted $\BG_{[0,1]}$.
\end{definition}
We list several equivalent ways to describe bisets.
\begin{proposition}
\label{prop:charbisets}
The following objects  stay in 1-1 correspondence.
\begin{enumerate}
\item Functions $\map\omega{I_{0,1}}{\set{0,1}}$ such that 
  \begin{equation}\label{eq:prodomega}
    \omega(r,s)\omega(s,t)=\omega(r,t)\dmf{(r,s),(s,t)\in I_{0,1}},
  \end{equation}
\item bisets $(L,R)$, and
\item transitive relations $M\subset[0,1]^2$ with the additional property that $(r,t)\in M$ implies both $r<t$
and $(r',t')\in M$ for all $(r',t')\in I_{r,t}$. 
\end{enumerate}
Namely, 
\begin{equation}
\label{eq:defomegaLR}
  \omega_{(L,R)}(s,t)=\left\{
    \begin{array}[c]{cl}
1&\text{~if~$[s,t)\subset L$ and $(s,t]\subset R$}\\
0&\text{~otherwise}
    \end{array}
\right.
\end{equation}
and $(s,t)\in M_{(L,R)}$ iff $[s,t)\subset L$ and $(s,t]\subset R$.
\end{proposition}
\begin{proof}
 Any $\omega$ determines two sets $L_\omega=\set{s\in[0,1):\exists t\in(0,1]:\omega(s,t)=1}$, $R_\omega=\set{t\in(0,1]:\exists s\in[0,1):\omega(s,t)=1}$. If $s\in L_\omega$, there is some $t>s$ such that $\omega(s,t)=1$. Thus this is true  for all $t=s+\varepsilon$, where $\varepsilon$ small is chosen suitably. This shows that $L_\omega$ fulfils \ref{b1}, \ref{b2}--\ref{b4} are proven similarly.

Fix a biset $(L,R)$. We prove now that  $\omega=\omega_{(L,R)}$ fulfils $(L,R)=(L_\omega,R_\omega)$. Suppose $s\in L$. Then there is  some $\varepsilon>0$ such that $[s,s+\varepsilon)\subset L$ and, due to \ref{b3}, $(s,s+\varepsilon)\subset R$, i.e.\ $[s,s+\varepsilon/\pi)\subset L$ and $(s,s+\varepsilon/\pi]\subset R$. Conversely, $s\in L_\omega$ shows
$[s,s+\varepsilon)\subset L$, i.e.\ $s\in L$. This proves $L=L_\omega$ and similarly, $R=R_\omega$. 

The  proof for the  relation $M_{L,R}$ is a consequence of the one-to-one correspondence between $M$ and $\omega$, $(s,t)\in M$ iff $\omega(s,t)=1$.
\end{proof}
\begin{remark}
For general  $Z\subseteq[0,1]$ the pair $(L_Z,R_Z)$ introduced above  need not be a biset since \ref{b1},\ref{b2} may not be fulfilled.  On the other hand, the function $\omega_Z$ defined by 
\begin{displaymath}
  \omega_Z(s,t)=\chfc{\set{Z':Z'\cap[s,t]=\emptyset}}(Z)\dmf{(s,t)\in I_{0,1}}
\end{displaymath}
fulfils \ref{eq:prodomega} for all $Z$. This means that we cannot distinguish $Z$ from a biset by this hitting function. This can be easily seen for $Z=\NQ\cap[0,1]$.

 At least,  $\FG_{[0,1]}$  can be mapped to $\BG_{[0,1]}$ via  $L_Z=[0,1)\setminus Z$, $R_Z=(0,1]\setminus Z$. Conversely, a biset $(L,R)$ corresponds to  a closed set in this way iff $L\cap(0,1)=R\cap(0,1)$.  
\end{remark}
Recall that a function $f$ on a metric space $X$ is upper semicontinuous if for all $x\in X$ $\liminf_{x'\to x}f(x')\ge f(x)$.  
\begin{lemma}
\label{lem:biset upper semicontinuous closed}
If $\map\omega{I_{0,1}}{\set{0,1}}$ fulfils \ref{eq:prodomega} then there is  a closed set  with $\omega=\omega_Z$ iff $\omega$ is upper semicontinuous. 
\end{lemma}
\begin{proof}
  Suppose $Z$ is closed. Then $\omega_Z(s,t)=1$ iff $Z\cap[s,t]=\emptyset$. This shows  $Z\cap[s-\varepsilon,t+\varepsilon]=\emptyset$ for some $\varepsilon>0$ such that $\omega_Z(s',t')=1$ if $|s-s'|<\varepsilon$, $|t-t'|<\varepsilon$, what means upper continuity.

On the other side, assume $\omega$ fulfils \ref{eq:prodomega} and is upper semicontinuous. This shows that both $L$ and $R$ are open such that they coincide on $(0,1)$ by \ref{b3} and \ref{b4}. Setting $Z=[0,1]\setminus(L\cup R)$ we obtain $\omega=\omega_Z$. This completes the proof.
\end{proof}
We equip $\BG_{[0,1]}$ with the coarsest topology making $(L,R)\mapsto\omega_{(L,R)}(s,t)$ continuous for all $(s,t)\in I_{0,1}$. 
\begin{proposition}
  \label{propPstu=bisets}
$\C_{0,1}$ is isomorphic to $C(\BG_{[0,1]})$.
\end{proposition}
\begin{proof}
 By the  well-known Gelfand representation theorem for abelian $C^*$-algebras \cite[Theorem 2.1.11A]{BR87} it follows that $\C_{0,1}$ is isomorphic to  the continuous functions on the set of characters on $\C_{0,1}$. Characters are  multiplicative linear functionals $\chi$ on $\C_{0,1}$, i.e.\ $\chi(ab)=\chi(a)\chi(b)$ for all $a,b\in \C_{0,1}$. Since $\C_{0,1}$  is generated by the projections $\mathrm{P}_{s,t}$, $(s,t)\in I_{0,1}$,  we have to determine all possible maps $(s,t)\mapsto\chi(\mathrm{P}_{s,t})$. Clearly, since $\mathrm{P}_{s,t}$ is a projection, $\chi(\mathrm{P}_{s,t})\in\set{0,1}$. Thus any  character of $\C_{0,1}$ is   determined by a  map $\map\omega{\set{(s,t)\in I_{0,1}}}{\set{0,1}}$ with \ref{eq:prodomega} and each such map defines a character of $\C_{0,1}$. Proposition \ref{prop:charbisets} completes the proof.
\end{proof}
\begin{corollary}
  \label{cor:Baire}
  Suppose $\H$ is a Hilbert space and $(\mathrm{P}_{s,t})_{(s,t)\in I_{0,1}}$  are  projections on $\B(\H)$ fulfilling  \ref{eq:Prst}. Then  there exists for all (normal) states $\eta$ on $\B(\H)$ a unique Baire  measure $\mu_\eta$ on $\BG_{[0,1]}$ with 
  \begin{equation}
    \label{eq:muetaBaire}
    \mu_\eta(\set{(L,R):\omega_{(L,R)}(s_i,t_i)=1})=\eta(\mathrm{P}_{s_1,t_1}\cdots \mathrm{P}_{s_k,t_k})\dmf{(s_i,t_i)\in I_{0,1}}
  \end{equation}
\end{corollary}
\begin{proof}
Since $\mathrm{lh}\set{\mathrm{P}_{s_1,t_1}\cdots \mathrm{P}_{s_k,t_k}:(s_i,t_i)\in I_{0,1}}$ is dense in $\C_{0,1}$, there is a unique Baire measure $\mu_\eta^b$ on $\BG_{[0,1]}$ with \ref{eq:muetaBaire}.  This completes the proof.
\end{proof}
\begin{remark}
  From the conditions \ref{b1}, \ref{b2} one obtains that $(0,1)\setminus(L\cap R)\in\FG_{(0,1)}$  and $L\setminus R$ and $R\setminus L$ are both countable. Since both $\FG_{(0,1)}$ and $\NR^\NN\cup\bigcup_{n\in\NN}\NR^n$ are Standard Borel spaces, there is a Standard Borel structure on $\BG_{[0,1]}$ too. If  this Standard Borel structure would be generated by the maps $(L,R)\mapsto\omega_{(L,R)}(s,t)$, $0\le s<t\le1$, we could extend the corollary to Borel measures. 
\end{remark}
\begin{remark}
\label{rem:nonseparable no faithful}
By the above note, a biset corresponds to a closed set $\cmpl{(L\cap R)}$. Consequently, we can associate with every (even nonseparable) product system  distributions of  random closed  sets, coming from the projections $\tilde{\mathrm{P}}_{s,t}=\mathrm{P}_s^-\mathrm{P}_{s,t}\mathrm{P}_t^+$. This corresponds to \cite[Lemma 2.9]{Tsi03}, where only times $t$ with $\mathrm{P}_s^\circ=\unit$ were considered in equation \ref{eq:mueta}. 

If $\E_1$ is nonseparable, weak and strong measurability do not coincide any more. Thus, depending on the choice of the measurability concept, $\gr t\tau$ may be only weakly continuous and Proposition \ref{prop:singlepoint} would fail.    Thus, it is not clear whether $ \mathrm{P}_{s,t}^{\Us,\circ}=  \mathrm{P}_{s,t}^{\Us}$ and we loose the direct connection between the random set and the space of units. Nevertheless, these projections encode an invariant structure in (nonseparable) product systems.  

But, the main obstacle for fully generalising Theorem \ref{th:RACS} to  nonseparable product systems is  that  there is no faithful normal state on $\B(\H)$ if $\H$ is not separable. Thus, there is no analogue of the statement concerned with  faithful normal states and  we cannot hope for a generalisation of \ref{eq:equZs+t}.  
\end{remark}
We want to close this section by showing that $\BG_{[0,1]}$ has besides its  topological  structure also   an order structure similar to $\FG_{[0,1]}$. Define, as natural
\begin{displaymath}
  (L,R)\preceq_\BG (L',R')\iff L\supseteq L'\text{~and~}R\supseteq R'\iff\omega_{(L,R)}\le\omega_{(L',R')}\iff M_{(L,R)}\subseteq M_{(L',R')}.
\end{displaymath}
\begin{proposition}
  \label{propbisetlattice}
  $(\BG_{[0,1]},\preceq_\BG)$ is a complete lattice, i.e.\ each family $((L_i,R_i))_{i\in I}$ has a unique least upper  and greatest lower bound.
\end{proposition}
\begin{proof}
We apply again \cite[Theorem I.6]{S:Bir84}, so we need to prove only the existence of a greatest lower bound.
  Identify $(L,R)$ with the function $\omega_{(L,R)}$ defined in \ref{eq:defomegaLR}. Then $\bigwedge_{i\in I}(L_i,R_i)$
should correspond to a function less than $\min_{i\in I}\omega_{L_i,R_i}=\prod_{i\in I}\omega_{L_i,R_i}$.  Since the latter   function fulfils \ref{eq:prodomega}, it corresponds to the  greatest lower bound of $((L_i,R_i))_{i\in I}$.
\end{proof}
\begin{remark}
  For the least upper bound, identify $(L,R)$ with the relation $M_{(L,R)}$. Then the transitive hull $M$ of 
$\bigcup_{i\in I}M_{(L_i,R_i)}$  is that  relation  which corresponds to  the least upper bound  of $((L_i,R_i))_{i\in I}$.
\end{remark}

\section{An Algebraic Invariant of Product Systems}
\label{sec:algebraic invariant}
In this last section we analyse an invariant related to ideas used in  \cite{Tsi00} for the special product systems described in Example \ref{ex:Tsirelson type III}. Especially, we want to show that this invariant yields no additional information for the classification of type $\mathrm{II}$ product systems.
 
An \emph{elementary set} $F\subset[0,1]$ has the form $F=\bigcup_{i=1}^n[s_i,t_i]$, let denote $\FG^e_{[0,1]}$  the set of all elementary sets in $[0,1]$. Consider any product system $\E=\sg t\E$. For each elementary set  $F=\bigcup_{i=1}^n[s_i,t_i]$ we can define the von Neumann subalgebra $\A_F=\bigvee_{i=1}^n\A_{s_i,t_i}\subseteq\B(\E_1)$, where $\A_{s,t}$ was given in \ref{eq:defAst}. Based upon this,  we  define for every $F\in\FG_{[0,1]}$ the von Neumann  subalgebra
\begin{displaymath}
  \A_F=\bigwedge_{F'\in\FG^e_{[0,1]},F'\supseteq F}\A_{F'}\subseteq\B(\E_1),
\end{displaymath}
which is consistent with the previous definition for  $F\in\FG^e_{[0,1]}$. 

In \cite{Tsi00b} there was considered for a sequence $\sequ nF\subset\FG^e_{[0,1]}$ norms of  differences  like $\norm{\eta_n^{\psi_1}-\eta_n^{\psi_2}}$ where $\eta^\psi_n$ is the restriction of the normal state $\eta^\psi=\scpro\psi{\cdot\psi}$ to $\A_{F_n}$. Thereby
\begin{displaymath}
  \norm{\eta_n-\eta'_n}=\max\set{
\abs{\eta(a)-\eta'(a)}:a\in\A_{F_n},\norm a\le 1}.
\end{displaymath}
In section \ref{sec:2measurability} we defined already a topology on the set $\N(\H)$ of all von Neumann subalgebras of $\B(\H)$.
\begin{proposition}
  \label{prop:convergence norms=convergence algebras}
  Suppose $\sequ nF$ converges to  $F\in\FG_{[0,1]}$ from above, i.e\ $F_{n+1}\subseteq F_n$ for all $n\in\NN$ and $\bigcap_{n\in\NN}F_n=F$. Then $\sequp n{\A_{F_n}}$ converges to $\A_F$. 

Moreover, the following statements are equivalent
\begin{enumerate}
\item \label{tsiinv1}$\norm{\eta_n-\eta'_n}\limitsto{}{n\to\infty}0$ for all pure normal states $\eta,\eta'$ on $\B(\E_1)$.
\item \label{tsiinv2}$\norm{\eta_n-\eta'_n}\limitsto{}{n\to\infty}0$ for all normal states $\eta,\eta'$ on $\B(\E_1)$.
\item \label{tsiinv3}$\A_F=\NC\unit$.
\end{enumerate}

Further, for disjoint $F,F'\in\FG_{[0,1]}$ the algebras $\A_F$ and $\A_{F'}$ commute and $\A_{F\cup F'}=\A_F\vee\A_{F'}$. For general  $F,F'\in\FG_{[0,1]}$ the relation  $\A_{F\cap F'}=\A_{F}\wedge \A_{F'}$ is true.
\end{proposition}
\begin{proof}
Define for $\varepsilon>0$ and $F\in\FG_{[0,1]}$ another set $F_\varepsilon=[0,1]\cap\bigcup_{t\in F}[t-\varepsilon,t+\varepsilon]\in\FG^e_{[0,1]}$. The first observation is that $\A_F=\bigwedge_{\varepsilon>0}\A_{F_\varepsilon}$. This is proven by  $\A_F'=\bigvee_{s,t:(s,t)\cap F=\emptyset}\A_{s,t}$, which follows from the same relation for elementary $F$. It is easy to see then that 
\begin{displaymath}
  (\bigwedge_{\varepsilon>0}\A_{(F)_\varepsilon})'=\bigvee_{s,t:(s,t)\cap(F)_\varepsilon=\emptyset\forall \varepsilon>0}\A_{s,t}.
\end{displaymath}
 But, $(s,t)\mapsto\A_{s,t}$ is continuous, thus 
 \begin{displaymath}
   \bigvee_{s,t:(s,t)\cap(F)_\varepsilon=\emptyset\forall \varepsilon>0}\A_{s,t}=\bigvee_{s,t:(s,t)\cap F=\emptyset}\A_{s,t}.
 \end{displaymath}
If $\sequ nF$ converges  to $F$ from above, then it  converges in the Hausdorff metric on $\FG_{[0,1]}$ \cite[Corollary 3 of Theorem 1-2-2]{C:Mat75}. Thus, for each $\varepsilon>0$ and large  $n\in\NN$, $F_n\subset (F)_\varepsilon$.  Consequently, $ \A_{F_n}\subseteq\A_{(F)_\varepsilon}$ for large $n$. This implies $\bigcap_{n\in\NN}\A_{F_n}\subseteq\A_F$. On the other hand, $\A_{F_n}\supseteq\A_F$ and the first statement is proven.

Clearly, \ref{tsiinv2} implies \ref{tsiinv1}.  Conversely, convexity of the norm shows that \ref{tsiinv1} implies \ref{tsiinv2}.

Suppose $\A_F\ne\NC\unit$. Then there exist two normal states $\eta,\eta'$ on $\B(\E_1)$ and $a\in\A_F$, $\norm a\le 1$ such that 
$\eta(a)\ne\eta'(a)$. From $\A_{F_n}\limitsto{}{n\to\infty}\A_F$ we derive $\sequ na$, $a_n\in\A_{F_n}$, $\norm{a_n}\le1$ with $\lim_{n\to\infty}a_n=a$. Continuity of $\eta-\eta'$ implies $\lim_{n\to\infty}\abs{\eta(a_n)-\eta'(a_n)}=\abs{\eta(a)-\eta'(a)}>0$. Thus \ref{tsiinv2} implies \ref{tsiinv3}.  
 
On the other hand, $\eta-\eta'$ is uniformly weakly continuous on the unit ball of $\B(\E_1)$ and zero on $\NC\unit$. Thus $\A_{F_n}\limitsto{}{n\to\infty}{\NC\unit}$  implies for all $\varepsilon>0$ for large enough $n$ that  $\abs{\eta(a_n)-\eta'(a_n)}\le\varepsilon$ for all $a_n\in\A_{F_n}$, $\norm{a_n}\le1$. 

If $F,F'$ are disjoint there exist elementary sets $F_1\supset F$, $F'_1\supset F'$ with $F_1\cap F_1'=\emptyset$. This shows that $\A_{F\cap F'}=\A_F\vee\A_{F'}$ and that $\A_F$ and $\A_{F'}$ commute since $\A_{F_1}$ and $\A_{F_1'}$  do so. 

For general $F,F'\in\FG_{[0,1]}$ we have  $\A_F=\bigwedge_{\varepsilon>0}\A_{F_\varepsilon}$ and $\A_{F'}=\bigwedge_{\varepsilon>0}\A_{F'_\varepsilon}$. Thus we derive from monotony of $\varepsilon\mapsto\A_{F_\varepsilon}$ that
\begin{displaymath}
  \A_{F}\wedge \A_{F'}=\bigwedge_{\varepsilon>0}\A_{F_\varepsilon}\wedge\bigwedge_{\varepsilon'>0}\A_{F'_{\varepsilon'}}=\bigwedge_{\varepsilon>0}\A_{F_\varepsilon}\wedge\A_{F'_\varepsilon}=\bigwedge_{\varepsilon>0}\A_{F_\varepsilon\cap F'_\varepsilon}=\bigwedge_{\varepsilon>0}\A_{(F\cap F')_\varepsilon}=\A_{F\cap F'}.
\end{displaymath}
This completes the proof.
\end{proof}
\begin{remark}
  The  construction of  the map $F\mapsto \A_F$ reminds of  the construction of a  Hausdorff measure, see Definition \ref{eq:defHausdorff}.  We can copy that definition by saying $\B(F)=\bigvee_{\varepsilon>0}\B_\varepsilon(F)$ with 
  \begin{displaymath}
     \B_\varepsilon(F)= \bigwedge\set{\bigvee_{i\in\NN}\A_{B_i}:\text{$\sequ iB$ are balls with  $d(B_i)\le\varepsilon$ and $\bigcup_{i\in\NN}B_i\supseteq F$}},
  \end{displaymath}
but there are some observations.  Defining $\A^\circ_G=\bigvee_{s,t:(s,t)\subseteq G}\A_{s,t}$ for open $G\subseteq[0,1]$, the fact that we can cover any open set by countably many arbitrarily small balls implies that   $  \B_\varepsilon(F)$ is independent from $\varepsilon>0$ and
\begin{displaymath}
  \B(F)=\bigwedge_{G\subseteq[0,1]:G\text{~open~},G\supseteq F}\A^\circ_G\dmf{F\text{~Borel set in~}[0,1]}.
\end{displaymath}
For compact $F$, it is easy to see that an open covering set can be chosen as interior of an elementary set and we get $\A_F=\B(F)$. Using the more general $\B(F)$ instead of $\A_F$ is more problematic because of  worse convergence properties. It is even not clear, whether $F\cap F'=\emptyset$ implies that $\B(F)$ and $\B(F')$ commute. For our purpose, it is enough to work with the map $F\mapsto\A_F$.
\end{remark}
\begin{corollary}
  \label{cor:Hausdorff-like algebra measure}
 Let $\E$ be a product system. Then   the map $\map{\kappa_\E}{\FG_{[0,1]}}{\set{0,1}}$ defined through
  \begin{displaymath}
    \kappa_\E(F)=\left\{
    \begin{array}[c]{c>$l<$}
0&~if~$\A_F=\NC\unit$\\
1&~otherwise
    \end{array}\right.
  \end{displaymath}
 is an invariant of the product system $\E$.\proofend
\end{corollary}
\begin{example}
  \label{ex:AFforFock}
  Suppose $\E=\Gammai(\NC^d)$. Then $\A_F=\B(\Gamma(L^2(F\times\set{1,\dots,d},\ell\restriction_F\otimes\#)))$. Especially, 
\begin{displaymath}
    \kappa_\E(F)=\left\{
    \begin{array}[c]{c>$l<$}
0&~if~$\ell(F)=0$\\
1&~otherwise
    \end{array}\right.
  \end{displaymath}
This is proven like follows. Define $\B_E=\B(\Gamma(L^2(E,\ell,\NC^d)))\otimes\unit_{\Gamma(L^2(E,\ell,\NC^d))}\subset\B(\E_1)$ for any Borel set $E$. Then $\B_E'=\B_{\cmpl E}$ and $\B_E=\A_E$ for $E\in\FG^e_{[0,1]}$. Observe that  $\B_E$ is generated in the $\sigma$-weak topology  by the  set  $\set{\W(f):f=f\chfc E~\ell\text{--a.s.}}$ of \emph{Weyl operators}  $\W(f)$, defined by 
\begin{displaymath}
  \W(f)\psi_h=\e^{-\frac12\norm f^2-\scpro 
hf}\psi_{h+f}\dmf{h\in L^2([0,1],\ell)}.
\end{displaymath}
If $\sequ nE  $ is increasing, we get $\bigvee_{n\in\NN} \B_{E_n}=\B_{\bigcup_{n\in\NN}E_n}$ since  $f\chfc{E_n}\limitsto{}{n\to\infty}f\chfc{\bigcup_{m\in\NN}E_m}$ and $f\mapsto\W(f)$ is strongly continuous \cite[section 20]{Par92}. Taking commutants respectively complements, we get  $\bigwedge_{n\in\NN} \B_{E_n}=\B_{\bigcap_{n\in\NN}E_n}$ for decreasing $\sequ nE$. This proves $\A_F=\B_F$ for all closed $F$. 
\end{example}
\begin{example}
  \label{ex:kappaFforRACS}
Let  $\M$ be a stationary factorizing measure type on $\FG_{[0,1]}$. We want  to compute $\kappa(F)$ for the product system $\E^\M$,  using relation  \ref{tsiinv1} from Proposition \ref{prop:convergence norms=convergence algebras}.
  \begin{lemma}
    In this situation, 
    \begin{displaymath}
      \kappa_{\E^\M}(F)=\left\{
    \begin{array}[c]{c>$l<$}
0&~if~$Z\cap F=\emptyset$~$\M$-a.s.\\
1&~otherwise
    \end{array}\right.\dmf{F\in\FG_{[0,1]}}
    \end{displaymath}
  \end{lemma}
  \begin{proof}
Take the unit $u$, $u_\mu=\mu(\set\emptyset)^{-1/2}\chfc{\set\emptyset}(\cdot)$ as considered in Corollary \ref{cor:psrstypeII}, the corresponding normal state $\eta^u$ and any other normal state $\eta$ on $\B(\E^\M_1)$. We want to look for sets $\sequ nF\subset\FG^e_{[0,1]}$, $F\in\FG_{[0,1]}$ with $F_n\downarrow F$ and $\norm{\eta_n-\eta^u_n}\limitsto{}{n\to\infty}0$. Since $\eta^u$ is the pure state with $\eta^u(\mathrm{P}^u_{0,1})=1$ the latter happens if and only if $\eta(\mathrm{P}^u_{F_n})\limitsto{}{n\to\infty}1$, defining $\mathrm{P}^u_{\bigcup_{i=1}^n[s_i,t_i]}=\prod_{i=1}^n\mathrm{P}^u_{s_i,t_i}$. By construction,
\begin{displaymath}
1=\lim_{n\to\infty} \eta(\mathrm{P}^u_{F_n})=\lim_{n\to\infty} \mu_\eta(\set{Z:Z\cap F_n=\emptyset})= \mu_\eta(\set{Z:Z\cap F=\emptyset}).
\end{displaymath}
 This shows that $\kappa_{\E^\M}(F)=0$ if and only if $\mu^{\E^\M,u}_\eta(\set{Z:Z\cap F=\emptyset})=1$ for all normal states $\eta$ on $\B(L^2(\M_{0,1}))$. But we know from Corollary \ref{cor:psrstypeII} that the latter is equivalent to $Z\cap F=\emptyset$~$\M$-a.s.\ which completes the proof.
  \end{proof}
Since we know that $\Pi_\ell(\set{Z:Z\cap F=\emptyset})=\e^{-\ell(F)}$, this is a direct generalisation of the result in the preceding example.
\end{example}
In exactly the same way we find 
\begin{corollary}
   For all type $\mathrm{II}$ product systems $\E$, all units $u$ of $\E$ the map $\kappa_\E$ computes as follows: 
    \begin{displaymath}
      \kappa_\E(F)=\left\{
    \begin{array}[c]{c>$l<$}
0&~if~$Z\cap F=\emptyset$~$\M^{\E,u}$-a.s.\\
1&~otherwise
    \end{array}\right.\dmf{F\in\FG_{[0,1]}}.\proofend
    \end{displaymath}
\end{corollary}
\begin{remark}
\label{rem:zerocapacitiesofM^u=M^v}
Thus, $\kappa_\E$ depends on the null sets of the capacity of any  $\mu\in\M^{\E,u}$ only. Therefore, it cannot characterize type $\mathrm{II}$ product systems fully.  

On the other hand, regarding the problem whether $\M^{\E,u}=\M^{\E,v}$ for all units $u,v$ of a type \textrm{II} product system or not, we see that at least the capacities $F\mapsto\mu(\set{Z:Z\cap F=\emptyset})$  have the same null sets for $\mu\in\M^{\E,u}$ and $\mu\in\M^{\E,v}$.  In spite of the fact that the capacities characterize $\mu$ by the Choquet theorem, this is a strong indication that  actually $\M^{\E,u}=\M^{\E,v}$ for all units $u,v$. 
\end{remark}
\begin{remark}\label{rem:dependence degree}
  We want to quantify how dense  $F\in\mathscr{F}_{[0,1]}$ with $\kappa_\E(F)=0$ could be and set
  \begin{displaymath}
    \mathrm{deg}_\E=\sup\set{\dim_{\mathscr{H}}F:\kappa_\E(F)=0}.
  \end{displaymath}
From  the above corollary and Note \ref{rem:Kahane}, we derive that in the case of  a stationary factorizing measure type on $\FG_{[0,1]}$ this degree is just $1$ minus the essential supremum of the Hausdorff dimension of a realisation. Further, the larger the degree, the less information about a realisation is needed to reconstruct it. I.e., the more dependence is inside the measure (type), see the announcement in the introduction on page \pageref{page:dependence measure}. 

Note that to differentiate the family of product systems constructed  in \cite{Tsi00b}, we should use  scaling functions like $h_\alpha(u)=\frac u{\abs{\ln u}^{\alpha-1}}$  for the Hausdorff measures instead of $h_\alpha(u)=u^\alpha$. It would be interesting to know for this family whether there is a family of type $\mathrm{II}$ product system $\E$ with the same $\kappa_\E$'s. 
\end{remark}
\begin{remark}
  We could go this way further. First, we could combine $\kappa$ with the invariant of marked product subsystems and consider the invariant   $\set{(\F,\M^\F,\kappa_\F):\F\text{~is subsystem of ~}\E}$ carrying again a natural order structure. But, we conjecture 
     \begin{displaymath}
      \kappa_\F(F)=\left\{
    \begin{array}[c]{c>$l<$}
0&~if~$Z\cap F=\emptyset$~$\M^{\F}$-a.s.\\
1&~otherwise
    \end{array}\right.\dmf{F\in\FG_{[0,1]}}
    \end{displaymath}
such that we do not  get anything new compared to  the invariant from Theorem   \ref{th:S(E)islattice}. 

Another extension could be  provided by  \textsc{Connes} classification of hyperfinite $W^*$-factors \cite{OP:Con73}.  So instead of just  asking whether $\A_F=\NC$ or not we could also use the classification of $\A_F$. First there may be the possibility that    $\A_F$ is not  a factor, since the intersection resp.\ the union of a monotone sequence of type $\mathrm{I}$ factors is only a factor, if a certain state fulfils a mixing condition \cite[Theorem 2.6.10]{BR87}.   Second there may be the possibility that the intersection is a factor, but not type $\mathrm{I}$. Yet, we have no explicite form of the algebras $\A_F$ besides for  type $\mathrm{I}$   product systems.
\end{remark}
\begin{remark}
\label{rem:vnaproductsubsystems=lattice}  Another extension, much more in the spirit of the present  paper,  would be to look at continuous tensor  product subsystems of $W^*$-algebras $(\B_{s,t})_{(s,t)\in I_{0,\infty}}$ in the product system $(\A_{s,t})_{(s,t)\in I_{0,\infty}}$.  We have as prominent examples
  \begin{enumerate}
  \item $\B_{s,t}=\A_{s,t}$: the algebra  itself,
\item $\B_{s,t}=L^\infty(\M_{s,t})$ acting in $\E^\M$,
\item the abelian $W^*$-subalgebras $\B_{s,t}=\set{\mathrm{P}^u_{s',t'}:(s',t')\in I_{s,t}}''$ and $\B_{s,t}=\set{\mathrm{P}^\Us_{s',t'}:(s',t')\in I_{s,t}}''$ if the projections $\mathrm{P}^u_{s',t'}$ are defined by canonical extension of \ref{eq:defPust} and \ref{eq:defPUst}, and
\item  $\B_{s,t}=\set{\mathrm{P}^u_{s',t'}:(s',t')\in I_{s,t}}'$ and $\B_{s,t}=\set{\mathrm{P}^\Us_{s',t'}:(s',t')\in I_{s,t}}'$ connected with  the direct integral representations in section \ref{sec:direct integrals}.
  \end{enumerate}
Of course, this list is not complete. Consider,  as example, in the case $\E=\Gammai(\NC^2)$ the families 
\begin{displaymath}
  \B^{\alpha}_{s,t}=\set{\W(\alpha f\oplus \beta\overline f):f=f\chfc{[s,t]}\text{$\ell$-a.e.~}}''\dmf{(s,t)\in I_{0,1}},
\end{displaymath}
where  $\alpha,\beta>0$ are fixed under the restriction  $\alpha^2-\beta^2=1$ and $\overline f$ denotes the complex conjugate of $f$. These families consist of   type III factors which cannot happen in  any of the above listed possibilities.

 The definitions of $\A_F$ and $\kappa(F)$ extend easily to von Neumann subalgebras, e.g.
\begin{displaymath}
  \B_F=\bigwedge_{F'\in\FG^e_{[0,1]},F'\supseteq F}\B_{F'}\subseteq\B(\E_1),
\end{displaymath}
with $\B_{\bigcup_{i=1}^n[s_i,t_i]}=\bigvee_{i=1}^n\B_{s_i,t_i}$. This way we could enrich the lattice of   product subsystems of $W^*$-algebras   of the  product system   $(\A_{s,t})_{(s,t)\in I_{0,\infty}}$ of $W^*$-algebras like in Theorem  \ref{th:spect}. Since we have  already collected enough invariants of product systems, we postpone this discussion to future  work.
\end{remark}

\section{Conclusions and Outlook}
\label{sec:conclusio}

In this paper, we provided a strong connection between continuous product systems of Hilbert spaces (of type $\mathrm{I}$ or $\mathrm{II}$) and the distribution of quasistationary quasifactorizing random sets on $[0,1]$. Summarisingly, we showed that the ideas developed in \cite{Tsi03,Tsi00,Tsi00b,Tsi99a,Tsi99b,C:TV98} to construct new examples of product systems lead to many new insights into the general structure of product systems.  These result are based essentially on the following new ingredients in the theory of product systems:
\begin{itemize}
\item reduction of the structure of product systems to a single Hilbert space $\E_1$ in Proposition \ref{prop:psbyE_1}
\item determination of  a \emph{continuous} unitary shift group on $\E_1$ in Proposition \ref{prop:tau_tiscontinuous}
\item association of commuting projections in $\E_1$ with a measure type of random sets, see Theorem \ref{th:RACS}.   
\item associating a commutative algebra with the related  direct integral representation, see Theorem \ref{th:directintegralps} 
\end{itemize}

To overview the productivity of the related ideas, we want to mention here all the  invariants of a product system $\E$ used or derived in the present work:
\def\labelenumi{\arabic{enumi}.~}
\begin{enumerate}
\item The measure type $\M^{\E,\Us}$ provided in Theorem \ref{th:ps2randomset}.
\item The lattice $\mathscr{S}(\E)$ together with the map $m_\E$ declared in Theorem \ref{th:spect}.
\item The direct integral representation of $\E$  related to the representation $J_{\mathrm{P}^\Us}$ of the abelian $W^*$-algebra $L^\infty(\M^{\E,\Us})$, see  Theorem \ref{th:directintegralps}.
\item The spectral decomposition of the flip group $\gr t\tau$, see Remark \ref{rem:spectral theory of flip as invariant}.
\item The leading spectral type of the representation $J_{\mathrm{P}^\Us}$ discussed in Note \ref{rem:HahnHellingerinvariant}.
\item The set of factorisation sets derived in the proof of Theorem \ref{th:E_0-semigroup}, see Note \ref{rem:factorization sets as invariant}.
\item The map $\kappa_\E$ from Corollary  \ref{cor:Hausdorff-like algebra measure}.
\item The lattice of product subsystems of von Neumann algebras  as discussed in Note \ref{rem:vnaproductsubsystems=lattice}.
\end{enumerate}
Although we were able, besides recording this list, to settle questions like the intrinsic measurable structure of product systems and the relation between type $\mathrm{II}$ and type $\mathrm{III}$, we are far from understanding the general structure of product systems. Therefore, we want to close this work with a  summary of  some questions we would like to answer in future.
\begin{enumerate}
\item At the moment, most important seems to us  to answer the question whether  the automorphisms of an arbitrary  product system act  transitively on the normalized units. We have  strong indications (see Proposition \ref{prop:unitandunitalprojections} and Note  \ref{rem:zerocapacitiesofM^u=M^v}) that the answer is affirmative. This  would have far-reaching consequences: 
  \begin{enumerate}
\item Like in Theorem  \ref{th:ps2randomset}  it would result that $\M^{\E,u}=\M^{\E,v}$ for any two units, see Note \ref{rem:ME,u as invariant}. This would  give us additional knowledge about the map $m_\E$ on the  lattice $\mathscr{S}(\E)$, provided in Theorem  \ref{th:spect}  (Note  \ref{rem:refinementS(E)}). Further, Corollary \ref{cor:S(E)forE^M}could be extended to all stationary factorizing measure types.
\item It would be the finishing touch on the characterisation of weak dilations of $\sigma$-strongly and uniformly continuous semigroups of unital completely positive maps on $\B(\H)$ for any separable Hilbert space \cite{OP:Bha01}. 
\item Any stationary factorizing measure type  $\M$ on $\FG_\NT$ would be an invariant of the  product system $\E^\M$ (Note \ref{rem:measuretypeasitsowninvariant}).
\item The direct integral representation given in Theorem \ref{th:directintegralps} used for $J_{\mathrm{P}^u}(L^\infty(\M^{\E,u}))$ would yield a much better structure theory, see Note \ref{rem:structuretheorybydirectintegral}.
  \end{enumerate}
\item A related topic is to determine the structure of automorphisms of type $\mathrm{II}$ product systems and examine the consequences on the structure of automorphisms of type $\mathrm{III}$ product systems completing the work of \textsc{Arveson} \cite[section 8]{Arv89}.
\item We think it should be quite instructive to determine all subsystems of the product system derived from the measure type of the sets of zeros of a Brownian motion, see Example \ref{ex:zerosBt}, as well as the automorphism group of this product system, see Example \ref{ex:automorphism for Bt0}. Similarly, the product systems of stationary factorizing measure type of random sets on $[0,1]\times L$ with noncompact $L$ deserves further study.
\item What structures from stochastic calculus have counterparts for quasistationary quasifactorizing measures? What type of  quantum stochastic calculus exists for  product systems or $E_0$-semigroups which are not type $\mathrm{I}$?
\item Are there product systems not of type $\mathrm{I}_{0}$ or $\mathrm{I}_{1}$ without any nontrivial product subsystems? Certainly, they should be type $\mathrm{II}_0$ or $\mathrm{III}$. 
\item One shall study  the structure  of  stationary  factorizing measure types on other types of Standard Borel  spaces  and of the corresponding  product systems determined by Theorem \ref{th:generalproductsystemfrom stationaryfactorizingmeasuretype}. 
\item Theorem \ref{th:directintegralps} indicates that factorizing Hilbert spaces over more complex structures than intervals of $\NRp$ like $\set{(t,Z):t\ge0, Z\in\FG_{[0,t]}}$ are important. Are there examples for $\NRp\times\NRp$ despite the more or less trivial ones derived from Poisson processes in the plane? 
\item To what extent do the results of the present paper  have counterparts for   product systems of  Hilbert modules? What about product systems of Hilbert spaces over the reals or quaternions?   
\item Is the conjugate product system $\E^*$ introduced in Example \ref{ex:ps from faithful states} isomorphic to $\E$? More specifically, are there product systems which could not be derived from stationary factorizing measure types on general Standard Borel spaces?
\item Is there a structure theory of   product systems of  $W^*$-algebras? 
\end{enumerate}

\paragraph{Acknowledgement} This work  profited very much from discussion with \textsc{B.V.~Rajarama Bhat} and  \textsc{P.~M\oe rters}. Without their help it would not have reached this state. \textsc{B.V.~Rao} pointed out a wrong statement in a former version of this paper.  Further, we want to thank \textsc{L.~Accardi}, \textsc{V.~Belavkin}, \textsc{St.~Bareto}, \textsc{R.~Dalang}, \textsc{F.~Fidaleo}, \textsc{U.~Franz}, \textsc{H.~Führ}, \textsc{M.~Lindsay}, \textsc{B.~Kümmerer}, \textsc{A.~Martin}, \textsc{M.~Skeide}, \textsc{B.~Tsirelson}, \textsc{H.~von Weizsäcker}, \textsc{G.~Winkler},  \textsc{O.~Wittich} and \textsc{J.~Zacharias} for interesting and stimulating  discussions on this topic  and  bibliographic information.  

This work was supported  by a DAAD-DST Project based Personnel exchange Programme  and INTAS grant N 99-00545.
\bibliographystyle{myplain}
\bibliography{qm-art,qm-lnm,quanten,abuch,kbuch,klass,punktpr,pbuch,qp8,qp9,malliavin,ag,qbuch,oper,stat,anal,here}
\clearpage
\end{document}